\documentclass[10pt,twoside]{article}

\usepackage{amssymb,amsmath,makeidx}
\usepackage[T1]{fontenc}
\usepackage{latexsym}
\usepackage{exscale}
\usepackage{pb-diagram}
\usepackage[dvips]{graphicx}
\usepackage{xr}

\usepackage{xcolor,pgf,tikz}

\newcommand{\al}{\alpha}
\newcommand{\alvec}{{\vec{\al}}} 
\newcommand{\alvecbe}{{\vec{\al}^\frown\be}} 
\newcommand{\alvecga}{{\vec{\al}^\frown\ga}} 
\newcommand{\alvecpr}{{\vec{\al}'}}
\newcommand{\alvecpl}{{\vec{\al}^+}}
\newcommand{\alvecst}{{\vec{\al}^\star}}
\newcommand{\alti}{\tilde{\alpha}}

\newcommand{\altipr}{\tilde{\alpha}^\prime}
\newcommand{\alpr}{{\alpha^\prime}}
\newcommand{\alnod}{{\al_0}}

\newcommand{\ale}{{\al_1}}
\newcommand{\alevec}{{\alvec_1}}

\newcommand{\ali}{{\al_i}}

\newcommand{\alivec}{{\alvec_i}}
\newcommand{\alvecrestri}{{\alvec{\scriptscriptstyle{\restriction {i}}}}}

\newcommand{\alie}{{\al_{i+1}}}

\newcommand{\alimin}{{\al_{i-1}}}

\newcommand{\aln}{{\al_n}}
\newcommand{\alne}{{\al_{n+1}}}
\newcommand{\alnvec}{{\alvec_n}}

\newcommand{\alnmin}{{\al_{n-1}}}

\newcommand{\alnminvec}{{\alvec_{n-1}}}

\newcommand{\alm}{{\al_m}}
\newcommand{\alme}{{\al_{m+1}}}
\newcommand{\almn}{{\al_{m+n}}}

\newcommand{\albar}{\bar{\al}}
\newcommand{\alnbar}{\overline{\aln}}
\newcommand{\alstar}{{\al^\star}}
\newcommand{\alstarvec}{{\alvec^\star}}

\newcommand{\alplus}{{\al^+}}
\newcommand{\alhat}{{\widehat{\al}}}
\newcommand{\be}{\beta}

\newcommand{\beti}{\tilde{\be}}
\newcommand{\bevec}{{\vec{\be}}}

\newcommand{\bevecpr}{{\vec{\be}'}}

\newcommand{\bebar}{\bar{\be}}

\newcommand{\bepr}{{\beta^\prime}}

\newcommand{\ga}{\gamma}
\newcommand{\gahat}{{\widehat{\ga}}}
\newcommand{\gapr}{{\ga^\prime}}
\newcommand{\gati}{\tilde{\gamma}}

\newcommand{\gavec}{{\vec{\ga}}}
\newcommand{\gavecpr}{{\vec{\ga}'}}

\newcommand{\gastar}{{\ga^\star}}
\newcommand{\Ga}{\Gamma}
\newcommand{\de}{\delta}
\newcommand{\depr}{{\de^\prime}}
\newcommand{\deti}{\tilde{\delta}}

\newcommand{\devec}{{\vec{\de}}}
\newcommand{\tildedevec}{{\vec{\tilde{\de}}}}

\newcommand{\De}{\Delta}

\newcommand{\eps}{\varepsilon}

\newcommand{\epsn}{\varepsilon_0}
\newcommand{\epsom}{\varepsilon_{\Omega+1}}

\newcommand{\epsal}{\varepsilon_{\al+1}}

\newcommand{\iovec}{{\vec{\iota}}}
\newcommand{\la}{\lambda}
\newcommand{\lapr}{{\la^\prime}}
\newcommand{\ka}{\kappa}
\newcommand{\om}{\omega}
\newcommand{\omq}[1]{{\bar{\omega}^{#1}}}
\newcommand{\Om}{\Omega}
\newcommand{\Omme}{\Om_{m+1}}
\newcommand{\Ommz}{\Om_{m+2}}
\newcommand{\Omkz}{\Om_{k+2}}
\newcommand{\Omke}{\Om_{k+1}}
\newcommand{\Omiz}{\Om_{i+2}}
\newcommand{\Omie}{\Om_{i+1}}

\newcommand{\rhoi}{\rho_i}

\newcommand{\nupr}{{\nu^\prime}}

\newcommand{\muti}{{\tilde{\mu}}}
\newcommand{\sivec}{{\vec{\si}}}

\newcommand{\varsivec}{{\vec{\varsigma}}}
\newcommand{\varsivecpr}{{\vec{\varsigma}^\prime}}
\newcommand{\varsivecnod}{{\vec{\varsigma}_0}}
\newcommand{\varsivecti}{{{\vec{\tilde{\varsigma}}}}}
\newcommand{\varsivecstar}{{\vec{\varsigma^\star}}}
\newcommand{\varsivecde}{{\vec{\varsigma_\de}}}
\newcommand{\sibar}{{\bar{\si}}}
\newcommand{\sipr}{\si^\prime}

\newcommand{\siti}{\tilde{\si}}

\newcommand{\sirstar}{{\si_r^\star}}
\newcommand{\siestar}{{\si_e^\star}}
\newcommand{\tauti}{\tilde{\tau}}
\newcommand{\tauvec}{{\vec{\tau}}}
\newcommand{\tauvecrestrk}{{\tauvec{\scriptscriptstyle{\restriction {k}}}}}

\newcommand{\tauevec}{\vec{\tau}_1}
\newcommand{\taunvec}{\vec{\tau}_n}
\newcommand{\taupr}{{\tau^\prime}}
\newcommand{\taubar}{\bar{\tau}}

\newcommand{\tautipr}{{\tauti^\prime}}

\newcommand{\tauipr}{{\tau^\prime_i}}
\newcommand{\taukpr}{{\tau^\prime_k}}
\newcommand{\taunpr}{{\tau^\prime_n}}

\newcommand{\taustar}{{\tau^\star}}

\newcommand{\tauistar}{{\tau_i^\star}}
\newcommand{\taui}{{\tau_i}}

\newcommand{\taukstar}{{\tau_k^\star}}
\newcommand{\taunstar}{{\tau_n^\star}}

\newcommand{\taunminstar}{{\tau_{n-1}^\star}}

\newcommand{\si}{\sigma}
\newcommand{\tht}{\vartheta}

\newcommand{\xivec}{{\vec{\xi}}}

\newcommand{\xivecpr}{{\vec{\xi}'}}
\newcommand{\ze}{\zeta}
\newcommand{\zepr}{{\zeta^\prime}}
\newcommand{\zevec}{{\vec{\ze}}}
\newcommand{\zevecy}{{\vec{\ze}}^y}
\newcommand{\zevecpr}{{\vec{\ze}'}}
\newcommand{\etapr}{{\eta^\prime}}

\newcommand{\etavec}{{\vec{\eta}}}

\newcommand{\etaij}{{\eta_{i,j}}}

\newcommand{\etauv}{{\eta_{u,v}}}

\renewcommand{\phi}{\varphi}

\newcommand{\N}{{\mathbb N}}
\newcommand{\Hz}{{\mathbb P}}
\newcommand{\Hznod}{{\mathbb P}_0}
\newcommand{\Lz}{{\mathbb L}}
\newcommand{\Ez}{{\mathbb E}}
\newcommand{\Ezone}{{\mathbb E}_1}
\newcommand{\On}{{\mathrm{Ord}}}
\newcommand{\card}{{\mathrm{Card}}}

\newcommand{\CNF}{{\mathrm{\scriptscriptstyle{CNF}}}}
\newcommand{\ANF}{{\mathrm{\scriptscriptstyle{ANF}}}}

\newcommand{\NF}{{\mathrm{\scriptscriptstyle{NF}}}}

\newcommand{\Lim}{\mathrm{Lim}}
\newcommand{\Limnod}{\mathrm{Lim}_0}
\newcommand{\Image}{\mathrm{Im}}

\newcommand{\finsub}{\subseteq_\mathrm{fin}}
\newcommand{\logend}{{\mathrm{logend}}}
\newcommand{\sumend}{{\mathrm{end}}}

\newcommand{\leo}{\le_1}

\newcommand{\lo}{<_1}

\newcommand{\klex}{<_\mathrm{\scriptscriptstyle{lex}}}
\newcommand{\kglex}{\le_\mathrm{\scriptscriptstyle{lex}}}
\newcommand{\glex}{>_\mathrm{\scriptscriptstyle{lex}}}

\newcommand{\lepw}{\le_\mathrm{\scriptscriptstyle{pw}}}

\newcommand{\phib}{\bar{\phi}}

\newcommand{\C}{{\operatorname{C}}}

\newcommand{\Cnm}{{\C^n_m}}

\newcommand{\Cnem}{{\C^{n+1}_m}}
\newcommand{\Cnmalbe}{\Cnm(\al,\be)}
\newcommand{\Cnemalbe}{\Cnem(\al,\be)}

\newcommand{\thtnm}{\tht^n_m}

\newcommand{\thtnem}{\tht^{n+1}_m}

\newcommand{\thtnnmin}{\tht^n_{n-1}}

\newcommand{\thtnk}{\tht^n_k}

\newcommand{\thtnek}{\tht^{n+1}_k}

\newcommand{\thtnmal}{\thtnm(\al)}
\newcommand{\thtnemal}{\thtnem(\al)}

\newcommand{\thtm}{\tht_m}

\newcommand{\thtn}{\tht_n}
\newcommand{\thtnod}{\tht_0}
\newcommand{\thte}{\tht_1}
\newcommand{\thti}{\tht_i}

\newcommand{\thtk}{\tht_k}
\newcommand{\thtke}{\tht_{k+1}}
\newcommand{\thtmal}{\tht_m(\al)}

\newcommand{\thtt}{\tht^\tau}
\newcommand{\thtti}{\tht^\taui}

\newcommand{\thtal}{\tht^\al}
\newcommand{\thtali}{\tht^\ali}
\newcommand{\thtbe}{\tht^\be}
\newcommand{\thts}{\tht^\si}

\newcommand{\T}{{\operatorname{T}}}

\newcommand{\Tm}{{\operatorname{T}_m}}

\newcommand{\Tnm}{{\operatorname{T}^n_m}}
\newcommand{\Tnn}{{\operatorname{T}^n_n}}

\newcommand{\Tnem}{{\operatorname{T}^{n+1}_m}}

\newcommand{\Tnme}{{\operatorname{T}^n_{m+1}}}

\newcommand{\Tt}{{\operatorname{T}^\tau}}

\newcommand{\Ttvec}{{\operatorname{T}^\tauvec}}

\newcommand{\laTalvec}{{\la\operatorname{-T}^\alvec}}
\newcommand{\Tsivec}{{\operatorname{T}^\sivec}}

\newcommand{\ltvec}{{\operatorname{l}^\tauvec}}

\newcommand{\ltvecal}{{\operatorname{l}^{\tauvec^\frown\al}}}

\newcommand{\Ttrestral}{{\operatorname{T}^\tau_{\restriction_\al}}}

\newcommand{\Ttrestrga}{{\operatorname{T}^\tau_{\restriction_\ga}}}
\newcommand{\Ttrestralplus}{{\operatorname{T}^\tau_{\restriction_{\alplus}}}}
\newcommand{\Talrestralplus}{{\operatorname{T}^\al_{\restriction_{\alplus}}}}

\newcommand{\Tal}{{\operatorname{T}^\al}}
\newcommand{\Talmn}{{\operatorname{T}^{\al_{m+n}}}}
\newcommand{\Tali}{{\operatorname{T}^\ali}}

\newcommand{\Tga}{{\operatorname{T}^\ga}}

\newcommand{\Ts}{{\operatorname{T}^\si}}
\newcommand{\Tts}{{\operatorname{T}^\tau}{\mbox{\raisebox{0.38ex}{$\scriptstyle [\si]$}}}}
\newcommand{\Tsr}{{\operatorname{T}^\si}{\mbox{\raisebox{0.38ex}{$\scriptstyle [\rho]$}}}}
\newcommand{\Talga}{{\operatorname{T}^\al}{\mbox{\raisebox{0.38ex}{$\scriptstyle [\ga]$}}}}

\newcommand{\Pm}{{\operatorname{P}_m}}
\newcommand{\Ptm}{{\operatorname{P}^\tau_m}}

\newcommand{\Ptnod}{{\operatorname{P}^\tau_0}}
\newcommand{\Pme}{{\operatorname{P}_{m+1}}}
\newcommand{\Pnm}{{\operatorname{P}^n_m}}

\newcommand{\Pt}{{\operatorname{P}^\tau}}

\newcommand{\starm}{{\ast_m}}

\newcommand{\start}{{\ast^\tau}}

\newcommand{\startm}{{\ast^\tau_m}}

\newcommand{\starnod}{{\ast_0}}

\newcommand{\starnm}{{\ast^n_m}}

\newcommand{\thetam}{{\theta_m}}
\newcommand{\thetame}{{\theta_{m+1}}}
\newcommand{\thetake}{{\theta_{k+1}}}
\newcommand{\thetanm}{\theta^n_m}
\newcommand{\thetann}{\theta^n_n}
\newcommand{\thetanne}{\theta^n_{n+1}}
\newcommand{\thetannmin}{\theta^n_{n-1}}
\newcommand{\thetanme}{\theta^n_{m+1}}
\newcommand{\thetank}{\theta^n_k}
\newcommand{\thetanke}{\theta^n_{k+1}}

\newcommand{\thetanem}{\theta^{n+1}_m}

\newcommand{\Thetanm}{\Theta^n_m}
\newcommand{\Thetann}{\Theta^n_n}

\newcommand{\Thetannmin}{\Theta^n_{n-1}}
\newcommand{\Thetanme}{\Theta^n_{m+1}}

\newcommand{\htarg}[1]{\operatorname{ht}_{#1}}

\newcommand{\htt}{{\operatorname{ht}_\tau}}

\newcommand{\htal}{{\operatorname{ht}_\al}}

\newcommand{\subtm}{{\operatorname{Sub}^\tau_m}}

\newcommand{\subtme}{{\operatorname{Sub}^\tau_{m+1}}}

\newcommand{\subtnod}{{\operatorname{Sub}^\tau_0}}

\newcommand{\tal}{{\operatorname{t}^\tau_\al}}

\newcommand{\ttau}{{\operatorname{t}^\al_\tau}}

\newcommand{\tlocaln}{(\tau=\al_0,\ldots,\al_n=\al)}

\newcommand{\lh}{{\operatorname{lh}}}

\newcommand{\kvnod}{\kappa^{{\scriptscriptstyle{\mathbf{()}}}}}

\newcommand{\kval}{\kappa^\alvec}

\newcommand{\kvalbe}{\kappa^\alvec_\be}

\newcommand{\kvalga}{\kappa^\alvec_\ga}

\newcommand{\kvalde}{\kappa^\alvec_\de}

\newcommand{\kvga}{\kappa^\gavec}

\newcommand{\laal}{{\la_\al}}

\newcommand{\laaln}{{\la_{\aln}}}
\newcommand{\laalnbe}{{\la^\aln_\be}}
\newcommand{\laalnminaln}{{\la^\alnmin_\aln}}

\newcommand{\iotal}{\iota_{\tau,\al}}

\newcommand{\iotga}{\iota_{\tau,\ga}}

\newcommand{\iotalinv}{{\iota_{\tau,\al}^{-1}}}

\newcommand{\zetal}{{\ze^\tau_\al}}

\newcommand{\latau}{{\la^\tau}}
\newcommand{\lasi}{{\la^\si}}
\newcommand{\latal}{{\la^\tau_\al}}

\newcommand{\latga}{{\la^\tau_\ga}}

\newcommand{\latbe}{{\la^\tau_\be}}

\newcommand{\pist}{\pi_{\si,\tau}}

\newcommand{\pigaal}{\pi_{\ga,\al}}
\newcommand{\pirhosi}{\pi_{\rho,\si}}
\newcommand{\pirhosiinv}{\pi^{-1}_{\rho,\si}}

\newcommand{\pizesiinv}{\pi^{-1}_{\ze,\si}}

\newcommand{\pigaalinv}{\pi^{-1}_{\ga,\al}}

\newcommand{\laspistal}{{\la^\si_{\pist(\al)}}}

\newcommand{\laalbe}{{\la^\al_\be}}

\newcommand{\laalbetal}{{\la^\al_{\be^\tal}}}

\newcommand{\Xij}{X_{i,j}}

\newcommand{\Xuv}{X_{u,v}}

\newcommand{\Xpr}{X^\prime}

\newcommand{\Zij}{Z_{i,j}}

\newcommand{\Zuv}{Z_{u,v}}

\newcommand{\Zpr}{Z^\prime}

\newcommand{\tildev}{\tilde{V}}

\newcommand{\tildey}{\tilde{Y}}
\newcommand{\tildeypr}{\tilde{Y}^\prime}
\newcommand{\tildeypl}{\tilde{Y}^+}

\newcommand{\yti}{\tilde{y}}

\newcommand{\tildez}{\tilde{Z}}
\newcommand{\tildezij}{\tilde{Z}_{i,j}}
\newcommand{\tildezijpr}{\tilde{Z}^\prime_{i,j}}
\newcommand{\dbltildezij}{\tilde{\tilde{Z}}_{i,j}}

\newcommand{\tildezuv}{\tilde{Z}_{u,v}}
\newcommand{\tildezpr}{\tilde{Z}^\prime}
\newcommand{\tildezstar}{\tilde{Z}^\star}
\newcommand{\tildezprnod}{\tilde{Z}^\prime_0}

\newcommand{\mc}{{\operatorname{mc}}}

\newcommand{\lf}{{\operatorname{lf}}}

\newcommand{\Par}{\operatorname{Par}}
\newcommand{\Part}{\operatorname{Par}^\tau}
\newcommand{\Parsi}{\operatorname{Par}^\si}

\newcommand{\PA}{{\operatorname{PA}}}
\newcommand{\BH}{{\operatorname{BH}}}
\newcommand{\kpom}{{\operatorname{KP}\!\om}}
\newcommand{\kplnod}{{\operatorname{KP}\!\ell_0}}
\newcommand{\kplr}{{\operatorname{KP}\!\ell^r}}
\newcommand{\pioneonecanod}{\Pi^1_1{\operatorname{-CA}_0}}
\newcommand{\pionetwocanod}{\Pi^1_2{\operatorname{-CA}_0}}
\newcommand{\idlessomega}{{\operatorname{ID}_{<\om}}}
\newcommand{\idn}{{\operatorname{ID}_n}}
\newcommand{\R}{{\cal R}}

\newcommand{\Rone}{{{\cal R}_1}}
\newcommand{\Rtwo}{{{\cal R}_2}}
\newcommand{\Rthree}{{{\cal R}_3}}
\newcommand{\Ctwo}{{{\cal C}_2}}

\newcommand{\Rn}{{{\cal R}_n}}
\newcommand{\Rom}{{{\cal R}_\om}}

\newcommand{\Ronepl}{{{\cal R}_1^+}}
\newcommand{\Rtwopl}{{{\cal R}_2^+}}

\newcommand{\Core}{{\operatorname{Core}}}

\newcommand{\bardot}{\bar{\cdot}}

\newcommand{\qed}{\mbox{ }\hfill $\Box$\vspace{2ex}}
\newcommand{\imp}{\Rightarrow}
\newcommand{\aeq}{\Leftrightarrow}
\newcommand{\andsp}{\:\&\:}
\newcommand{\veesp}{\:\vee\:}
\newcommand{\sub}{\subseteq}

\newcommand{\set}[2]{\{ #1 \:|\: #2\}}

\newcommand{\singleton}[1]{\{ #1 \}}

\newlength{\hilflh}
\newcommand{\hilfminus}[1]{
  \settowidth{\hilflh}{$#1-$}\mbox{$#1-\hspace{-0.5\hilflh}
  \makebox[0pt]{\raisebox{0.24\hilflh}{$#1\cdot$}}\hspace{0.5\hilflh}$}}
\newcommand{\minusp}{\mathbin{\mathchoice {\hilfminus{\displaystyle}}
  {\hilfminus{\textstyle}}{\hilfminus{\scriptstyle}}
  {\hilfminus{\scriptscriptstyle}}}}

\newtheorem{theo}{Theorem}[section]
\newtheorem{cor}[theo]{Corollary}
\newtheorem{lem}[theo]{Lemma}
\newtheorem{defi}[theo]{Definition}
\newtheorem{prop}[theo]{Proposition}
\newtheorem{conv}[theo]{Convention}
\newtheorem{rmk}[theo]{Remark}
\newtheorem{claim}[theo]{Claim}

\newcommand{\oneinf}{1^\infty}
\newcommand{\siinf}{\si^\infty}
\newcommand{\tauinf}{\tau^\infty}
\newcommand{\alinf}{\al^\infty}

\newcommand{\chicheck}{\check{\chi}}
\newcommand{\chis}{\chi^\si}
\newcommand{\chit}{\chi^\tau}
\newcommand{\chitcheck}{\check{\chi}^\tau}
\newcommand{\chial}{\chi^\al}

\newcommand{\chialncheck}{\check{\chi}^\aln}

\newcommand{\mutal}{\mu^\tau_\al}

\newcommand{\mualn}{\mu_\aln}
\newcommand{\mus}{\mu^\si}
\newcommand{\mut}{\mu^\tau}
\newcommand{\mutbe}{\mu^\tau_\be}
\newcommand{\mutali}{\mu^\tau_\ali}

\newcommand{\mutaln}{\mu^\tau_\aln}

\newcommand{\muga}{\mu_\ga}
\newcommand{\mualbe}{\mu^\al_\be}

\newcommand{\mualbetal}{\mu^\al_{\be^\tal}}

\newcommand{\muspistal}{{\mu^\si_{\pist(\al)}}}
\newcommand{\rhosi}{\varrho^\si}
\newcommand{\rhot}{\varrho^\tau}
\newcommand{\rhotal}{\varrho^\tau_\al}
\newcommand{\rhoalmutal}{\varrho^\al_{\mutal}}
\newcommand{\rhoalbe}{{\varrho^\al_\be}}

\newcommand{\rhoalnga}{{\varrho^\aln_\ga}}

\newcommand{\MNF}{{\mathrm{\scriptscriptstyle{MNF}}}}
\newcommand{\mNF}{{\mathrm{\scriptscriptstyle{NF}}}}
\newcommand{\Mz}{{\mathbb M}}

\newcommand{\trs}{{\mathrm{ts}}}

\newcommand{\trst}{{\mathrm{ts}^\tau}}
\newcommand{\trsal}{{\mathrm{ts}^\al}}

\newcommand{\trsaln}{{\mathrm{ts}^\aln}}
\newcommand{\trsrelal}{{\mathrm{ts}[\al]}}

\newcommand{\cspr}{{\mathrm{cs}^\prime}}

\newcommand{\rs}{{\mathrm{rs}}}
\newcommand{\ers}{{\mathrm{ers}}}
\newcommand{\rsij}{{\mathrm{rs}_{i,j}}}
\newcommand{\ersij}{{\mathrm{ers}_{i,j}}}
\newcommand{\rsistar}{{\mathrm{rs}_{i^\star}}}
\newcommand{\rsarg}[1]{{\mathrm{rs}_{#1}}}

\newcommand{\tc}{{\mathrm{tc}}}
\newcommand{\TC}{{\mathrm{TC}}}
\newcommand{\laTC}{{\la\operatorname{-TC}}}
\newcommand{\latTC}{{(\la,t)\operatorname{-TC}}}

\newcommand{\RS}{\mathrm{RS}}
\newcommand{\laRS}{\la\operatorname{-RS}}

\newcommand{\xiRS}{\xi\operatorname{-RS}}
\newcommand{\RSt}{\mathrm{RS}^\tau}
\newcommand{\letwo}{\le_2}
\newcommand{\lethree}{\le_3}
\newcommand{\ktwo}{<_2}
\newcommand{\kthree}{<_3}
\newcommand{\lei}{\le_i}

\newcommand{\lSeq}{\mathrm{lSeq}}

\newcommand{\dom}{\mathrm{dom}}

\newcommand{\domkval}{{\mathrm{dom}({\kval})}}
\newcommand{\domkvga}{{\mathrm{dom}({\kvga})}}
\newcommand{\domnuval}{{\mathrm{dom}({\nuval})}}
\newcommand{\domnuvga}{{\mathrm{dom}({\nuvga})}}

\newcommand{\dpf}{\mathrm{dp}}

\newcommand{\dpval}{{\mathrm{dp}_\alvec}}

\newcommand{\dpvga}{{\mathrm{dp}_\gavec}}

\newcommand{\nual}{\nu^\al}

\newcommand{\nuval}{\nu^\alvec}

\newcommand{\nuvalbe}{\nu^\alvec_\be}
\newcommand{\nuvalga}{\nu^\alvec_\ga}

\newcommand{\nuvga}{\nu^\gavec}

\newcommand{\alcp}[1]{{\al_{#1}}}
\newcommand{\becp}[1]{{\be_{#1}}}

\newcommand{\decp}[1]{{\de_{#1}}}
\newcommand{\sicp}[1]{{\si_{#1}}}
\newcommand{\sicppr}[1]{{\si^\prime_{#1}}}
\newcommand{\taucp}[1]{{\tau_{#1}}}
\newcommand{\taucppr}[1]{{\tau^\prime_{#1}}}

\newcommand{\zecp}[1]{{\ze_{#1}}}
\newcommand{\zecpy}[1]{{\ze^y_{#1}}}

\newcommand{\rcp}[1]{{r_{#1}}}

\newcommand{\alticp}[1]{{\alti_{#1}}}
\newcommand{\tauticp}[1]{{\tauti_{#1}}}
\newcommand{\siticp}[1]{{\siti_{#1}}}

\newcommand{\rhoargs}[2]{\varrho^{#1}_{#2}}
\newcommand{\alvecrestrarg}[1]{{\alvec{\scriptscriptstyle{\restriction {#1}}}}}
\newcommand{\bevecrestrarg}[1]{{\bevec{\scriptscriptstyle{\restriction {#1}}}}}

\newcommand{\devecrestrarg}[1]{{\devec{\scriptscriptstyle{\restriction {#1}}}}}

\newcommand{\tpr}{{t^\prime}}

\newcommand{\beucp}[1]{\be^{#1}}

\newcommand{\ordvalue}{\mathrm{o}}
\newcommand{\ov}{\mathrm{o}}

\newcommand{\ordvalalvec}{\mathrm{o}^\alvec}
\newcommand{\ordcp}[1]{{\mathrm{o}_{#1}}}

\newcommand{\ktc}{<_\mathrm{TC}}
\newcommand{\letc}{\le_\mathrm{TC}}

\newcommand{\cml}{\operatorname{cml}}
\newcommand{\gbo}{\operatorname{gbo}}
\newcommand{\predec}{\operatorname{pred}}
\newcommand{\predecs}{\operatorname{Pred}}
\newcommand{\succs}{\operatorname{Succ}}

\newcommand{\maxeta}{\eta_{\operatorname{max}}}

\newcommand{\TS}{\operatorname{TS}}
\newcommand{\laTS}{\la\operatorname{-TS}}

\newcommand{\TSe}{\operatorname{TS}^1}
\newcommand{\TSt}{\operatorname{TS}^\tau}
\newcommand{\TSalvec}{\operatorname{TS}^\alvec}
\newcommand{\TSal}{\operatorname{TS}^\al}
\newcommand{\TSaln}{\operatorname{TS}^\aln}

\newcommand{\mts}{\operatorname{mts}}
\newcommand{\mtsal}{\operatorname{mts}^\al}
\newcommand{\mtsale}{\operatorname{mts}^\ale}

\newcommand{\mtsga}{\operatorname{mts}^\ga}

\newcommand{\hop}{\operatorname{h}}
\newcommand{\homega}{\operatorname{h}_\om}
\newcommand{\hbe}{\operatorname{h}_\be}
\newcommand{\hga}{\operatorname{h}_\ga}

\newcommand{\maxmucovtau}{\operatorname{max-cov}^\tau}
\newcommand{\minmucoval}{\operatorname{min-cov}^\al}

\newcommand{\sk}{\operatorname{sk}}
\newcommand{\skbe}{\operatorname{sk}_\be}

\newcommand{\ec}{\operatorname{ec}}
\newcommand{\me}{\operatorname{me}}
\newcommand{\mepl}{\operatorname{me}^+}

\newcommand{\refcp}[1]{\operatorname{ref}_{#1}}

\def\vec#1{\mathchoice{\mbox{\boldmath$\displaystyle#1$}}
{\mbox{\boldmath$\textstyle#1$}}
{\mbox{\boldmath$\scriptstyle#1$}}
{\mbox{\boldmath$\scriptscriptstyle#1$}}}

\newcommand{\bs}{{\mathrm{bs}}}
\newcommand{\bspr}{\mathrm{bs}^\prime}

\newcommand{\ups}{\upsilon}
\newcommand{\upsseg}{{\upsilon\mathrm{seg}}}

\newcommand{\Romannumeral}[1]{\uppercase\expandafter{\romannumeral #1\relax}}

\begin{document}

\title{Pure $\Sigma_2$-Elementarity beyond the Core
\footnote{Preprint of article in press:
G. Wilken, Pure $\Sigma_2$-elementarity beyond the core, Annals of Pure and Applied Logic 172 (2021), https://doi.org/10.1016/j.apal.2021.103001
 }}

\author{Gunnar Wilken\\
Structural Cellular Biology Unit\\
Okinawa Institute of Science and Technology\\
1919-1 Tancha, Onna-son, 904-0495 Okinawa, Japan\\
{\tt wilken@oist.jp}
}

\maketitle

\begin{abstract}
\noindent We display the entire structure $\Rtwo$ coding $\Sigma_1$- and $\Sigma_2$-elementarity on the ordinals.
This will enable the analysis of pure $\Sigma_3$-elementary substructures.
\end{abstract}

\section{Introduction}
Intriguing theorems that demonstratively serve as examples for mathematical incompleteness have been a subject of great interest ever since
G\"odel established his incompleteness theorems \cite{G31}, showing that Hilbert's programme (see \cite{Hilbert}, and for a recent
account \cite{Z06}) could not be executed as originally expected.
Some of the most appealing such theorems exemplifying incompleteness are the Paris-Harring\-ton theorem \cite{PH}, 
Goodstein sequences \cite{Goodstein, KirbyParis, WW13}, Kruskal's theorem \cite{Kr60, RW93}, 
its extension by Friedman \cite{S85}, the graph minor theorem by Robertson and Seymour, see \cite{FRS87} and, for a general reference on 
so-called concrete mathematical incompleteness, Friedman's book \cite{Fr11}. The term concrete incompleteness refers to natural mathematical theorems
independent of significantly strong fragments of $\mathrm{ZFC}$, Zermelo-Fraenkel set theory with the axiom of choice. 
Technique and theorems on phase transition from provability to unprovability that bear on methods and results from analytic number theory 
were developed by Weiermann \cite{We03,We09} in order to further understand phenomena of mathematical independence, with interesting
contributions from others, e.g.\ Lee  \cite{L14}. 
A natural generalization of Kruskal's theorem was shown by Carlson \cite{C16}, using elementary patterns of resemblance \cite{C01} as 
basic structures of nested trees. 
The embeddability of tree structures into one another plays a central role in the area of well-quasi orderings, cf.\ \cite{C16}, 
maximal order types of which can be measured using ordinal notation systems from 
proof theory, cf.\ \cite{Schm79} and the introduction to \cite{WW11}. 

Elementary patterns of resemblance (in short: patterns) of order $n$ are finite structures of orderings $(\le_i)_{i\le n}$ where 
$\le_0$ is a linear ordering and $\le_1,\ldots,\le_n$ are forests such that 1) $\le_{i+1}\subseteq\le_i$ and
2) $a\le_{i+1}b$ whenever $a\le_i b\le_i c$ and $a\le_{i+1}c$ for all $i<n$ and all $a,b,c$ in the universe of the pattern. 
Patterns that do not contain further non-trivial functions or relations are called pure patterns.  
Kruskal's theorem establishes that the collection of finite trees is well-quasi ordered with respect to inf-preserving embeddings
\cite{Kr60}. When it comes to patterns, the natural embeddings are coverings, i.e.\ $\le_0$-embeddings that maintain
the relations $\le_i$ for $i=1,\ldots,n$. Carlson's theorem \cite{C16} states that pure patterns of order $2$ are well-quasi ordered
with respect to coverings. This is provable in the extension of the base theory of reverse mathematics, $\mathrm{RCA}_0$, of recursive comprehension, 
by the uniform $\Pi^1_1$-reflection principle for $\kplnod$, a subsystem of set theory axiomatizing a mathematical
universe that is a limit of admissible sets, see \cite{Ba75,C16}. In \cite{W18} we give a proof, which is independent of \cite{C17}, 
of the fact that Carlson's theorem is unprovable in $\kplnod$, or equivalently, $\pioneonecanod$. 
This latter subsystem of second order number theory, where induction is restricted to range over 
sets and set comprehension is restricted to $\Pi^1_1$-formulae, plays a prominent role in reverse mathematics, see Simpson \cite{S09} 
for the relevance of these theories in mathematics. 
Pohlers \cite{P98} provides an extensive exposition of various subsystems of set theory and second order number theory,
equivalences and comparisons in strength, and their proof-theoretic ordinals. Note that in \cite{P98} $\kplnod$ goes by the name 
$\kplr$.

\subsection*{\boldmath Discovery of patterns, the structure $\Rone$, and Cantor normal form\unboldmath}
Patterns were discovered by Carlson \cite{C00} during model construction that verifies the consistency of epistemic arithmetic with 
the statement \emph{I know I am a Turing machine}, thereby proving a conjecture by Reinhardt, see Section 21.2 of \cite{W} for a short
summary.
The consistency proof uncovers, through a demand for $\Sigma_1$-elementary substructures via a well-known finite set criterion, see
Proposition 1.2 of \cite{C99}, here part 1 of Proposition \ref{letwocriterion}, 
a structure $\Rone=\left(\On;\le,\leo\right)$ of ordinals, 
where $\le$ is the standard linear ordering on the ordinals and the relation $\leo$ is defined recursively in $\be$ by
\[\al\leo\be:\aeq (\al;\le,\leo) \preceq_{\Sigma_1} (\be;\le,\leo).\]
We therefore have $\al\leo\be$ if and only if $(\al;\le,\leo)$ and $(\be;\le,\leo)$ satisfy the same $\Sigma_1$-sentences over the language 
$(\le,\leo)$ with parameters from $\al=\{\ga\mid\ga<\al\}$, where $\Sigma_1$-formulas are quantifier-free formulas preceded by finitely many existential quantifiers.
$\Rone$ was analyzed in \cite{C99} and shown to be recursive and periodic in multiples of the ordinal $\epsn$, the proof-theoretic ordinal
of Peano arithmetic ($\PA$).
Pure patterns of order $1$ comprise the finite isomorphism types of $\Rone$ and provide ordinal notations that denote their unique pointwise minimal coverings 
within $\Rone$, an observation Carlson elaborated in \cite{C01}. 
In $\Rone$ the relation $\le_1$ can still be described by a relatively simple recursion formula, see Proposition \ref{Roneformula} below, 
given in terms of Cantor normal form notation, as was carried out in \cite{C99}.

Ordinal notations in Cantor normal form, indicated by the notation $=_\CNF$, are built up from $0,+$, and $\om$-exponentiation, 
where $\om$ denotes the least infinite ordinal, see \cite{P09,Sch77} for reference. 
We write $\al=_\CNF\om^\ale+\ldots+\om^\aln$ if $\al$ satisfies the equation with $\ale\ge\ldots\ge\aln$.
If $\al$ is represented as a sum of weakly decreasing additive principal numbers (i.e.\ powers of $\om$) 
$\rho_1\ge\ldots\ge\rho_m$, we also write $\al=_\ANF\rho_1+\ldots+\rho_m$, where we allow $m=0$ to cover the case $\al=0$.
$\epsn$ is the least fixed point of $\om$-exponentiation, as in general the class function
$\al\mapsto\eps_\al$ enumerates the class $\Ez:=\{\xi\mid\xi=\om^\xi\}$ of epsilon numbers. 
Note that Cantor normal form notations with parameters from an initial segment, say $\tau+1$, of the ordinals provide notations for the 
segment of ordinals below the least epsilon number greater than $\tau$. 

Defining $\lh(\al)$, the length of $\al$, to be the greatest $\be$ such that $\al\leo\be$, if such $\be$ exists,
and $\lh(\al):=\infty$ otherwise, the recursion formula for $\Rone$ reads as follows:

\begin{prop}[\cite{C99}]\label{Roneformula} For $\al=_\CNF\om^\ale+\ldots+\om^\aln<\epsn$, where $n>0$ and $\aln=_\ANF\rho_1+\ldots+\rho_m$, 
we have \[\lh(\al)=\al+\lh(\rho_1)+\ldots+\lh(\rho_m).\]
\end{prop}

While, as shown in \cite{C99}, the structure $\Rone$ becomes periodic in multiples of $\epsn$, with the proper multiples of $\epsn$ characterizing those ordinals $\al$ which satisfy $\lh(\al)=\infty$, in $\Rtwo$ this isomorphic repetition of the interval
$[1,\epsn]$ with respect to additive translation holds only up to the ordinal $\epsn\cdot(\om+1)$, since the pointwise least $<_2$-pair is
$\epsn\cdot\om<_2\epsn\cdot(\om+1)$ and hence $\lh(\epsn\cdot(\om+1))=\epsn\cdot(\om+1)$, as was shown in \cite{CWc}.

\subsection*{G\"odel's program and patterns of embeddings}
G\"odel's program to solve incompleteness progressively by the introduction of large cardinal axioms, see \cite{F96,K09}, 
together with the heuristic correspondence between large cardinals and ordinal notations through reflection properties and embeddings, 
motivated Carlson to view patterns as a programmatic approach to ordinal notations, as we pointed out earlier in \cite{W18} and \cite{W}. 
According to his view, the concept of patterns 
of resemblance, in which binary relations code the property of elementary substructure, can be generalized to patterns 
of embeddings that involve codings of embeddings, and would ultimately tie in with inner model theory in form of an ultra-fine structure.   
It seems plausible that $\Rn$-patterns ($n<\om$) just suffice to analyze set-theoretic systems of $\Pi_m$-reflection ($m<\om$), 
see \cite{R94,PS}.
 
\subsection*{\boldmath Additive patters of order $1$ and the structure $\Ronepl$\unboldmath}
The extension of $\Rone$ to a relational structure $\Ronepl=\left(\On;0,+;\le,\leo\right)$, containing the graphs of $0$ and ordinal addition,  
gives rise to additive patterns of order $1$ and is the central object of study in \cite{C01,W06,W07a,W07b,W07c,CWa,A15}.
Such extensions are of interest when it comes to applications of patterns to ordinal analysis.
Following Carlson's recommendation in \cite{C01}, \cite{W06} gives a sense of the relation $\leo$ in $\Ronepl$, moreover, 
an easily accessible characterization of the Bachmann-Howard structure in terms of $\Sigma_1$-elementarity is given in full detail. 
The Bachmann-Howard ordinal characterizes the segment of ordinals below the 
proof-theoretic ordinal of Kripke-Platek set theory with infinity, $\kpom$, which axiomatizes the notion of (infinite) admissible set, 
cf.\ Barwise \cite{Ba75}, or equivalently $\operatorname{ID}_1$, the theory of non-iterated positive inductive definitions, cf.\ \cite{P09}. 
A benefit from such characterizations, as pointed out in the introduction to \cite{W11}, is an illustrative semantics for Skolem hull 
notations, cf.\ \cite{P91}, the type of ordinal notations that have been most-useful in ordinal analysis so far, by the description of 
ordinals given in hull-notation as least solutions to $\Sigma_1$-sentences over the language $(0,+,\le,\leo)$.
The pattern approach to ordinal notations, on the other hand, is explained in terms of ordinal arithmetic in the style of \cite{RW93}.

\cite{W07c} presents an ordinal assignment producing pointwise minimal instantiations of $\Ronepl$-patterns in the ordinals. 
In turn, pattern-characterizations are assigned to ordinals given in classical notation. In \cite{CWa} normal forms are considered, 
also from the viewpoint of \cite{C01}. The elementary recursive isomorphism between Skolem-hull notations 
and pattern notations for additive patterns of order $1$, established in \cite{W07a,W07b,W07c}, was of considerable interest to  
proof-theorists, as it explained the new concept of patterns in a language familiar to proof-theorists.
It was shown in \cite{W07b,W07c} that notations in terms of $\Ronepl$-patterns characterize the proof-theoretic ordinal of $\pioneonecanod$.
By simply incorporating (the graph of) ordinal addition into the pattern's language, patterns of order one considerably increase in strength.
A similar phenomenon for classical notations is observed in the opposite direction in \cite{SchS} when withdrawing basic functions like 
addition from a notation system, as doing so causes a collapse, in case of a system for $\pioneonecanod$ down to $\epsn$.
On the other hand, the analysis of additive patterns of order $1$ is easily modified when adding (the graphs of) other basic arithmetic 
functions for expressive convenience, as we pointed out in earlier work. 
Such a procedure does not increase the order type of the resulting notation system any further, cf.\ Section 9 of \cite{C01}. 

\subsection*{Skolem-hull based ordinal notations and arithmetic}
The analysis of $\Ronepl$ in \cite{W07b} and the assignment of minimal ordinal solutions to additive patterns of order $1$ in \cite{W07c}
and \cite{CWa} rely on a toolkit of ordinal 
arithmetic based on relativized Skolem hull-based notation systems $\Tt$ and was introduced in \cite{W07a}. 
The parameter $\tau$ is intended to be either $1$ or an arbitrary epsilon number, 
which we often denote as $\tau\in\Ezone=\Ez\cup\{1\}$.
This type of ordinal notation system builds upon contributions over many decades mainly by Bachmann, Aczel, Feferman, Bridge, 
Sch\"utte, Buchholz, Rathjen, and Weiermann (see \cite{B55,Br75,B75,B81,R90,RW93}, and the introduction of \cite{W07a}). 
Preparatory variants of systems $\Tt$ that can be extended to stronger ordinal notation systems and arithmetic were provided in 
\cite{WW11}, starting with work by Buchholz and Sch\"utte in \cite{BSch76}. 
The so far most powerful notation systems of this classical type were introduced by Rathjen, 
resulting in a notation system for $\pionetwocanod$ \cite{R05}, and have quite recently been elaborated more completely for an analysis 
of the provably recursive functions of reflection by Pohlers and Stegert \cite{PS}.

We now go into more detail regarding ordinal arithmetic sufficiently expressive to display pattern notations on the basis of order $1$ 
(additive) or order $2$ (pure).
Let $\Om^{>\tau}$ be the least regular cardinal greater than $\card(\tau)\cup\aleph_0$. We also simply write $\Om$ for $\Om^{>\tau}$
when the dependence of $\Om$ on $\tau\in\Ezone$ is easily understood from the context. For convenience we set $\T^0:=\T^1=\T$ 
for the original notation system that provides notations for the ordinals below the proof theoretic ordinal of theories 
such as $\pioneonecanod$, $\kplnod$, and $\idlessomega$, the latter being the theory of finitely times-iterated inductive definitions, 
cf.\ e.g.\ \cite{P98}. In the presence of constants for all ordinals less than $\tau$ and ordinal addition, 
stepwise collapsing, total, injective, unary functions $(\thti)_{i<\om}$, where 
$\thtnod=\thtt$ is relativized to $\tau$, give rise to unique terms for all ordinals in the notation system $\Tt$.
The relativized systems $\Tt$ we are referring to were carefully introduced in Section 3 of \cite{W07a}, including complete proofs, 
but for the reader's convenience we provide a review of ordinal arithmetic in terms of the systems $\Tt$ in the next section, 
including the formal definition of $\tht$-functions. 
However, for now, descriptions given here should suffice to understand the larger picture.   
  
The Veblen function \cite{V1908} is a binary function in which the first 
argument indicates a lower bound of the fixed point level. 
While $\varphi(0,\cdot)$ enumerates the additive principal numbers, i.e.\ the class of ordinals
greater than $0$ that are closed under ordinal addition, which we often denote as $\Hz$, $\varphi(1,\cdot)$ enumerates 
$\eps$-numbers starting with $\epsn$. In general, $\al\mapsto\varphi(\xi+1,\al)$ enumerates the fixed points of $\varphi(\xi,\cdot)$, 
and $\varphi(\la,\cdot)$ enumerates the common fixed points of all $\varphi(\xi,\cdot)$, $\xi<\la$, where $\la$ is a limit ordinal. 
The fixed-point free variant $\phib$ of the Veblen function omits fixed-points, meaning that $\phib(0,\cdot)$ enumerates those 
additive principal numbers that are not epsilon numbers, $\phib(1,\cdot)$ enumerates all epsilon numbers $\al$ such that $\al<\eps_\al$,
and so on. Therefore, the first argument of the function $\phib$ denotes the exact fixed-point level. 

Now, the fixed point-level of $\Ga$-numbers, i.e.\ ordinals $\al$ such that $\al=\varphi(\al,0)$, 
cannot be expressed unless we introduce a ternary Veblen function.
In the context of $\Tt$-systems this problem of finding names for increasing fixed-point levels is resolved by arithmetization 
through higher $\tht_j$-functions,
relying on the regularity properties of the $\Om_i$, in that they cannot be reached by enumeration functions from below. 

With ordinal addition being part of the set of functions over which Skolem hulling is performed, the function $\xi\mapsto\thtt(\xi)$, 
where $\xi<\Om$, enumerates all additive principal numbers in the interval $[\tau,\Om)$ that are not epsilon numbers greater than $\tau$, 
for the sake of uniqueness of notation. 
Fixed points of increasing level, in the sense explained above in the context of the Veblen function, are in the range of $\thtt$, 
however, as enumerations starting from multiples of $\Om$. For example, $\xi\mapsto\thtt(\Om+\xi)$, where $\xi<\Om$, enumerates
the ordinals $\{\eps_\al\mid\tau,\al<\eps_\al<\Om\}$.

At each multiple of $\Om$ in the domain of $\thtt$ the enumeration of 
ordinals below $\Om$ of next higher fixed point level begins, again in fixed point omitting manner, which means that $\xi<\thtt(\De+\xi)$
whenever $\De$ is such a multiple of $\Om$ and $\xi<\Om$. 
Letting $\Om_0:=\tau$, $\Om_1:=\Om^{>\tau}$, and generally $\Om_{i+1}$ to be the least (infinite) regular cardinal greater than $\Om_i$ 
for $i<\om$, this algorithm of denoting ordinals allows for a powerful mechanism extending the fixed-point free variant of the Veblen 
function.
Generalizing from $\thtt=\tht_0$, the function $\tht_i$ enumerates additive principal numbers starting with $\Om_i$ in fixed-point free manner: 
\begin{lem}[cf.\ 4.2 of \cite{W07a}]\label{tenumlem} Let $i<\om$.
\begin{enumerate}
\item[a)] $\tht_i(0)=\Om_i$.
\item[b)] For any $\al<\Om_{i+1}$ we have $\tht_i(1+\al)=\omq{\Om_i+\al}$,
\end{enumerate}
where $\xi\mapsto\omq{\xi}$ denotes the enumeration function of the additive principal numbers that are not epsilon numbers.
\end{lem}
Consequently, fixed-points of increasing level are enumerated by $\tht_i$ with arguments starting from multiples of $\Om_{i+1}$, 
which in turn are (additively) composed of values of the function $\tht_{i+1}$. Note that we have $\tht_i=\tht^{\Om_i}$ for $i<\om$,
where the right hand side $\tht$-function is a $\tht_0$-function that is relativized to $\Om_i$. 

\medskip

\noindent{\bf Examples.} Assuming $\tau=1$, we give a few instructive examples of values of $\tht$-functions: 
$\tht_0(0)=1=\om^0=\varphi(0,0)$, $\tht_0(1)=\om=\varphi(0,1)$, $\tht_0(\epsn)=\om^{\epsn+1}$, 
$\tht_0(\tht_1(0))=\tht_0(\Om)=\epsn=\varphi(1,0)$, 
$\tht_0(\Om+\epsn)=\eps_{\epsn}$, $\tht_0(\Om+\Om)=\varphi(2,0)$, $\tht_0(\Om^2)=\Ga_0$ where $\Om=\aleph_1$ and $\Om^2=\tht_1(\tht_1(0))$,
and $\tht_0(\tht_1(\tht_2(0)))=\tht_0(\eps_{\Om+1})$ denoting the Bachmann-Howard ordinal. 
The proof-theoretic ordinal of the theory of $n$-times iterated inductive definitions $\mathrm{ID}_n$, $0<n<\om$, is 
$|\mathrm{ID}_n|=\tht_0(\tht_1(\ldots(\tht_{n+1}(0))\ldots))$.
For $n<\om$ the ordinal $\tht_0(\tht_1(\ldots(\tht_{n+1}(0))\ldots))$ is equal to the least $\leo$-predecessor of the pointwise minimal
$\ktwo$-chain of $(n+2)$-many ordinals, as was first shown in \cite{CWc}.
Finally, $\sup\{\tht_0(\ldots\tht_n(0)\ldots)\mid n<\om\}$ is the proof-theoretic ordinal of $\mathrm{ID}_{<\om}$, $\kplnod$, and
$\pioneonecanod$, and is the limit of expressibility of ordinals in terms of pure patterns of order $2$, see \cite{W18}.

Taking the functions $\thti$ for granted for now, we can formally define $\Tt$ and state a rough lemma on the collapsing nature of the
functions $\thti$.
\begin{defi}[\boldmath Inductive Definition of $\Tt$, cf.\ 3.22 of \cite{W07a}\unboldmath]\label{Ttdefi}
\mbox{}
\begin{itemize}
\item $\tau\sub\Tt$
\item $\xi,\eta\in\Tt\imp\xi+\eta\in\Tt$
\item $\xi\in\Tt\cap\Omiz\imp\thti(\xi)\in\Tt$ for all $i<\om$.
\end{itemize}
\end{defi}

\begin{lem}[cf.\ 3.30 of \cite{W07a}] For every $i<\om$ the function
$\thti\restriction_{\Tt\cap\Omiz}$ is $1$-$1$ and has values in $\Hz\cap[\Om_i,\Omie)$. 
\end{lem}

Thus every ordinal in $\Tt$ can be identified with a unique term
built up from parameters below $\tau$ using $+$ and the functions $\thti$ $(i<\om)$, if we assume that additive compositions are given
in additive normal form. The intersection $\Tt\cap\Om$ turns out to be an initial segment of ordinals, which for $\tau=1$ is the proof-theoretic
ordinal of $\pioneonecanod$.
The essential reason for $\Tt\cap\Om$ to be an ordinal is that all proper components $\eta\in\Tt\cap\Om$ used to compose the unique notation 
of some $\xi\in\Tt\cap\Om$ satisfy $\eta<\xi$, see Lemma \ref{segmentlem} and Theorem \ref{tnmthm} (3.12 and 3.14 of \cite{W07a}). 
For a term of shape $\al:=\thtt(\De+\xi)\in\Tt$, where $\De$ is a multiple of $\Om$ and $\xi<\Om$, 
as was mentioned earlier, this means that besides $\xi<\al$ we also have $\eta<\al$ for every subterm $\eta$ of $\De$ such that $\eta<\Om$.

The choice of regular cardinals for the $\Om_i$ is just convenient for short proofs. One can choose recursively regular ordinals instead
(see \cite{R93}), at the cost of more involved proofs to demonstrate recursiveness (see \cite{Schl93}).
The usage of regular cardinals causes gaps in notation systems constructed in a similar fashion as the systems $\Tt$, where the first such
gap is $[\Tt\cap\Om,\Om)$. 
Patterns arise in a more dynamic manner from local reflection properties. 
This is one major reason that explains why, 
despite the elegance and brevity of their definition, their calculation is quite involved. Another, intrinsically connected, major reason
explaining the subtleties in pattern calculation is the necessity to calculate or analyze in terms of the connectivity components of the
relations $\le_i$. We will return to the discussion of connectivity components later.

As explained above, we developed an expressive ordinal arithmetic in \cite{W07a} on the basis of Skolem hull term systems for initial segments 
$\tauinf:=\Tt\cap\Om^{>\tau}$ of $\On$, that was further extended in \cite{CWa}, \cite{CWc}, and \cite{W18}. 
The following definition transfinitely iterates the closure under $\tau\mapsto\tauinf=\Tt\cap\Om^{>\tau}$ continuously through all of $\On$.
Note that we have chosen the Greek lowercase letter $\ups$ (upsilon) in order to avoid ambiguity of notation. 
\begin{defi}[modified 9.1 of \cite{W07a}]\label{upsilondefi}
Let $(\ups_\iota)_{\iota\in\On}$ be the sequence defined by
\begin{enumerate}
\item $\ups_0:=0$,
\item $\ups_{\xi+1}:=\ups_\xi^\infty=\T^{\ups_\xi}\cap\Om^{>\ups_\xi}$,
\item $\ups_\la:=\sup\{\ups_\iota\mid\iota<\la\}$ for $\la\in\Lim$, i.e.\ $\la$ a limit ordinal.
\end{enumerate}
\end{defi}
In Corollary 5.10 of \cite{W07b} we have shown that the maximal $<_1$-chain in $\Ronepl$ is $\Image(\ups)\setminus\{0\}$. 
The results from \cite{W07b}, \cite{W07c}, and from Sections 5 and 6 of \cite{CWa}
were applied in a case study for an ordinal analysis of $\operatorname{ID}_{<\om}$, \cite{Pae}, following Buchholz' style of 
operator-controlled derivations \cite{B92}. Proof-theoretic analysis based on patterns as originally intended by Carlson uses 
$\be$-logic \cite{Cb}. 
As we pointed out in \cite{W}, the intrinsic semantic content via the notion of elementary substructure directly on the ordinals in 
connection with their immediate
combinatorial characterization and the concise elegance of their definition on the one hand, and the remarkable intricacies in their
calculation, revealing mathematical depth not yet fully understood, on the other hand, create the impression that patterns contribute 
to the quest for natural well-orderings, cf.\ \cite{F}. 

\subsection*{\boldmath The structure $\Rtwo$ and its core, notes on $\Rtwopl$ and $\Rthree$\unboldmath}
Returning to the discussion of pure patterns of order $2$, let \[\Rtwo=\left(\On;\le,\leo,\letwo\right)\] be the structure of ordinals 
with standard linear ordering $\le$ and partial orderings $\le_1$ and $\le_2$, simultaneously defined by induction on $\be$ in  
\[\al\le_i\be:\aeq \left(\al;\le,\le_1,\le_2\right) \preceq_{\Sigma_i} \left(\be;\le,\le_1,\le_2\right)\]
where $\preceq_{\Sigma_i}$ is the usual notion of $\Sigma_i$-elementary substructure (without bounded quantification), as analyzed thoroughly
below the ordinal $\upsilon_1=\oneinf$ of $\kplnod$ in \cite{CWc,W17,W18}.

$\le_i$-relationships between ordinals considered as members of structures of patterns can be verified using finite-set criteria.
The following criterion can be extended successively to provide criteria for relations $\le_i$, $i<\om$, see e.g.\ Proposition 3.9 of \cite{W}.

\begin{prop}[cf.\ 7.4 of \cite{CWc}]\label{letwocriterion} Let $\al,\be$ be such that $\al<\be$ and $X,\tildey,Y$ be 
finite sets of ordinals such that $X,\tildey\subseteq\al$ and $Y\subseteq[\al,\be)$. Consider the following properties: 
\begin{enumerate}
\item $X<\tildey<\al$ and there exists an isomorphism $h:X\cup\tildey\stackrel{\cong}{\longrightarrow} X\cup Y$. 
\item For all finite $\tildey^+$, where $\tildey\sub\tildey^+\sub\al$, $h$ can be extended to an isomorphism $h^+$ such that
\[h^+: X\cup\tildey^+\stackrel{\cong}{\longrightarrow} X\cup Y^+\]
for a suitable superset $Y^+\subseteq\be$ of $Y$.
\end{enumerate}
If for all such $X$ and $Y$ there exists a set $\tildey$ that satisfies property 1, then we have $\al<_1\be$.
If $\tildey$ can be chosen so that additionally property 2 holds, then we even have $\al<_2\be$.  
\end{prop} 
{\bf Proof.} A proof is given in Section 7 of \cite{CWc} and in greater detail in \cite{W} (Propositions 3.1 and 3.6). 
\qed

Note that for any $\al\in\On$ the class $\{\be\mid\al\le_i\be\}$ is closed (under limits) for $i=1,2$ and a closed interval for $i=1$.
The criterion for $\le_1$ can be applied to see that in $\Rtwo$ we have $\al\leo\al+1$ if and only if $\al\in\Lim$ and for every 
$<_2$-predecessor $\be$ of $\al$ the ordinal $\al$ is a proper supremum of $<_2$-successors of $\be$, see Lemma \ref{loalpllem}.
Another basic observation is that whenever $\al<_2\be$, $\al$ must be the proper supremum of an infinite $<_1$-chain, i.e.\ the order type
of the set of $<_1$-predecessors of $\al$ must be a limit ordinal, see Lemma \ref{ktwoinflochainlem}. 
Another useful elementary observation proved in \cite{CWc} and in \cite{W} (Lemma 3.7) is the following

\begin{lem}[7.6 of \cite{CWc}]\label{letwoupwlem}
Suppose $\alpha<_2\beta$, $X\subseteq_\mathrm{fin}\alpha$, and $\emptyset\not=Y\subseteq_\mathrm{fin}[\alpha,\beta)$.
\begin{enumerate}
\item There exist cofinally many $\tilde{Y}\subseteq\beta$ such that $X\cup\tilde{Y}\cong X\cup Y$.
More generally, for any $Z\subseteq_\mathrm{fin}\alpha$ with $X<Z$, if 
$\al\models\forall x\exists\tildez\;(x<\tildez\wedge\mbox{``}X\cup Z\cong X\cup\tildez\mbox{''})$ then this also holds in $\be$.
\item Cofinally in $\al$, copies $\tildey\subseteq\al$ of $Y$ can be chosen that besides satisfying $X<\tildey$ and 
$X\cup\tildey\cong X\cup Y$ also 
maintain $\leo$-connections to $\be$: For any $y\in Y$ such that $y\lo\be$ the corresponding $\yti$ satisfies $\yti\lo\al$.\qed
\end{enumerate}
\end{lem}

$\Core(\Rtwo)$, the \emph{core} of $\Rtwo$, i.e.\ the union of pointwise minimal instantiations of all finite isomorphism types of $\Rtwo$, 
was analyzed in \cite{W18} and shown to coincide (in domain) with the initial segment of the ordinals below $\oneinf$. 
Note that any (finite) subset of $\Rtwo$ gives rise to a substructure of $\Rtwo$; hence, it represents a (finite) isomorphism type of $\Rtwo$.
It was shown in \cite{W18} that collections of pure patterns of order $2$ as defined above and of finite isomorphism types of $\Rtwo$
coincide, and that each pure pattern of order $2$ has a \emph{unique isominimal} representative $P\subset\Rtwo$,
where a finite substructure $Q\subset\Rtwo$ is isominimal if and only if $Q\lepw R$ for every $R\subset\Rtwo$ such that $Q\cong R$. 
$Q\lepw R$ means that $Q$ is pointwise less than or equal to $R$ with respect to increasing enumerations. The notions core and isominimality
where introduced by Carlson in \cite{C01}. 
In the case of patterns of order $1$ the increased strength resulting from basic arithmetic functions such as addition, is matched by pure 
patterns of order $2$, as we have shown in \cite{CWc} that any pure pattern of order $2$ has a covering below $\oneinf$, 
the least such ordinal. According to a conjecture by Carlson, this phenomenon of compensation holds generally for all orders.      

Despite the fact that the cores of $\Ronepl$ and $\Rtwo$ cover the same initial segment of the ordinals, their structures differ considerably. While, as pointed out in \cite{W}, the core of $\Ronepl$ shows a great deal of uniformity, reminding one of Girard's notion of dilator, 
cf.\ \cite{G81}, and giving rise to Montalb\'an's Question 27 in \cite{Mo11}, the core of $\Rtwo$ is a structure, the regularity of which is 
far less obvious, due to the absence of uniformity provided by ordinal addition.
In this context the ordinal $\oneinf$ is obtained as a collapse when weakening additive pattern notations of order $2$ on the basis of the 
structure $\Rtwopl=\left(\On;0,+;\le,\leo,\letwo\right)$ to pure patterns of order $2$ arising in $\Rtwo$. 
We claim, as in \cite{W}, that the segment of countable ordinals denoted by the Skolem-hull notation system derived from the first 
$\om$-many weakly inaccessible cardinals covers (the domain of) $\Core(\Rtwopl)$. 
Note that the notation system based on the first weakly inaccessible cardinal matches the proof-theoretic ordinal of the set theory 
$\operatorname{KPI}$, which axiomatizes an admissible universe that is a limit of admissible sets, and which is equivalent to the system 
$\De^1_2$-$\operatorname{CA}+\operatorname{BI}$ of second order number theory, first analyzed by J\"ager and Pohlers in \cite{JP82}, 
cf.\ also \cite{P98}. See Buchholz' seminal work \cite{B92} for an analysis of $\operatorname{KPI}$ via operator controlled derivations. 
The analysis of $\Rtwopl$ is a topic of future work and requires a generalization of ordinal arithmetical methods, the beginning of which
is outlined in \cite{WW11}. Note that \cite{C09} discusses patterns of order 2 in a general way over Ehrenfeucht-Mostowski structures; 
however, on the basis of modified relations $\leo$, $\letwo$, see Section 5 of \cite{C09}.

In this article, we generalize the approach taken in \cite{CWc,W17,W18}, which in turn naturally extends the arithmetical 
analysis of pure $\Sigma_1$-elementarity given by Carlson in \cite{C99}, to arithmetically characterize the relations $\leo$ and $\letwo$
in all of $\Rtwo$, not just its core.
Define \[I:=\{\iota\in\On\mid\iota>1\mbox{ and not of a form }\iota=\la+1\mbox{ where }\la\in\Lim\}\]
and let the expression $\iota\minusp1$ for $\iota\in\On$ denote $\iota_0$ if $\iota=\iota_0+1$ for some $\iota_0$, and simply $\iota$ otherwise.
The following theorem is an immediate corollary of Theorem \ref{maintheo} of the present article.
\begin{theo}[\boldmath Maximal $\ktwo$-chain in $\Rtwo$\unboldmath]\label{maxchaintheo} $\mbox{ }$
\begin{enumerate}
\item The sequence $(\ups_\iota)_{\iota>0}$
is a $\lo$-chain through the ordinals, supporting the maximal $\ktwo$-chain through the ordinals, which is $(\ups_\iota)_{\iota\in I}$.
Here maximality means that for any $\iota\in I$ (actually for any ordinal $\iota$), the enumeration of all $\ktwo$-predecessors of 
$\ups_\iota$ is given by $(\ups_\xi)_{\xi\in I\cap\iota}$. 

\item For $\iota\in\On\setminus I$ the ordinal $\ups_\iota$ is 
$\ups_{\iota\minusp1}$-$\leo$-minimal, i.e.\ there does not exist any $\al\in(\ups_{\iota\minusp1},\ups_\iota)$ such that $\al\lo\ups_\iota$, 
and does not have any $\ktwo$-successor. 

\item For every $\iota\in I$ the ordinal $\ups_\iota$ is a $\ktwo$-predecessor of every proper multiple of $\ups_{\iota\minusp1}$ greater
than $\upsilon_\iota$.
The ordinals of the form $\ups_\la$ where $\la\in\Lim$ comprise the set of suprema of $\ktwo$-chains of limit order type. \qed
\end{enumerate}
\end{theo}

We give a simple illustration of Theorem \ref{maxchaintheo} below, where $\la$ stands for any $\la\in\Lim$. 
All indicated ordinals are $\leo$-connected, blue edges except for the connection $\ups_\la\ktwo\ups_{\la+1}$ indicate the maximal 
$\ktwo$-chain of part 1 of the theorem. 
The \emph{least} $\ktwo$-successor of, for instance, $\ups_2$ is $\ups_2+\ups_1$, which in turn does not
have any $\ktwo$-successors itself. The least $\ktwo$-successor of $\ups_3$ is $\ups_3+\ups_2$, and so forth, 
while the least $\ktwo$-successor of $\ups_\la$ is $\ups_\la\cdot2$ and the least $\ktwo$-successor of $\ups_{\la+2}$ is 
$\ups_{\la+2}+\ups_{\la+1}$, etc.
Note that according to the theorem, $\ups_1$ and $\ups_{\la+1}$ do not have any $\ktwo$-successors.
The ordinal $\ups_1$ is $\leo$-minimal, and the greatest $\lo$-predecessor of $\ups_{\la+1}$ is $\ups_\la$.

\color{blue}
\begin{pgfpicture}{0cm}{0cm}{15cm}{2cm}

\pgfxyline(1,1.1)(1,0.9)

\pgfxyline(2,1)(4,1)

\pgfxyline(2,1.1)(2,0.9)
\pgfxyline(3,1.1)(3,0.9)
\pgfxyline(4,1.1)(4,0.9)
\pgfxyline(4,1)(4.8,1)

\pgfputat{\pgfxy(2.5,1.25)}{\pgfbox[center,base]{\color{red}${\scriptscriptstyle <_2}$}}
\pgfputat{\pgfxy(3.5,1.25)}{\pgfbox[center,base]{\color{red}${\scriptscriptstyle <_2}$}}
\pgfputat{\pgfxy(4.5,1.25)}{\pgfbox[center,base]{\color{red}${\scriptscriptstyle <_2}$}}

\pgfxyline(7.2,1)(8,1)
\pgfputat{\pgfxy(7.5,1.25)}{\pgfbox[center,base]{\color{red}${\scriptscriptstyle <_2}$}}

\pgfxyline(8,1.1)(8,0.9)
\pgfxyline(9,1.1)(9,0.9)

\pgfxyline(10,1)(12,1)
\pgfxyline(10,1.1)(10,0.9)
\pgfxyline(11,1.1)(11,0.9)
\pgfxyline(12,1.1)(12,0.9)

\pgfxyline(12,1)(12.8,1)
\pgfputat{\pgfxy(12.5,1.25)}{\pgfbox[center,base]{\color{red}${\scriptscriptstyle <_2}$}}

\pgfputat{\pgfxy(10.5,1.25)}{\pgfbox[center,base]{\color{red}${\scriptscriptstyle <_2}$}}
\pgfputat{\pgfxy(11.5,1.25)}{\pgfbox[center,base]{\color{red}${\scriptscriptstyle <_2}$}}

\pgfxycurve(8,1)(8.25,1.2)(8.75,1.2)(9,1)
\pgfstroke
{\color{red}
\pgfputat{\pgfxy(8.5,1.25)}{\pgfbox[center,base]{\color{red}${\scriptscriptstyle <_2}$}}}

\pgfxycurve(8,1)(8.5,1.8)(9.5,1.8)(10,1)
\pgfstroke
{\color{red}
\pgfputat{\pgfxy(9,1.8)}{\pgfbox[center,base]{\color{red}${\scriptscriptstyle <_2}$}}}

{\color{black}
\pgfputat{\pgfxy(1,0.5)}{\pgfbox[center,base]{$\ups_1$}}
\pgfputat{\pgfxy(2,0.5)}{\pgfbox[center,base]{$\ups_2$}}
\pgfputat{\pgfxy(3,0.5)}{\pgfbox[center,base]{$\ups_3$}}
\pgfputat{\pgfxy(4,0.5)}{\pgfbox[center,base]{$\ups_4$}}
\pgfputat{\pgfxy(6,0.5)}{\pgfbox[center,base]{$\cdots$}}
\pgfputat{\pgfxy(8,0.5)}{\pgfbox[center,base]{$\ups_{\la}$}}
\pgfputat{\pgfxy(9,0.5)}{\pgfbox[center,base]{$\ups_{\la+1}$}}
\pgfputat{\pgfxy(10,0.5)}{\pgfbox[center,base]{$\ups_{\la+2}$}}
\pgfputat{\pgfxy(11,0.5)}{\pgfbox[center,base]{$\ups_{\la+3}$}}
\pgfputat{\pgfxy(12,0.5)}{\pgfbox[center,base]{$\ups_{\la+4}$}}
\pgfputat{\pgfxy(13,0.5)}{\pgfbox[center,base]{$\cdots$}}
}
\end{pgfpicture}

\color{black}

Theorem \ref{maintheo} also shows that the gap $[\ups_1,\ups_2)$ contains a $\lo$-chain of order type $\ups_2$ that is $\lo$-connected 
to $\ups_2$, and the gap $[\ups_{\la+1},\ups_{\la+2})$ contains a $\lo$-chain of order type $\ups_{\la+2}$ that is $\lo$-connected to 
$\ups_{\la+2}$, as in general for any successor ordinal $\iota$ the interval $[\ups_\iota,\ups_{\iota+1})$ contains a $\lo$-chain of 
order type $\ups_{\iota+1}$ that is $\lo$-connected to $\ups_{\iota+1}$.

Theorem \ref{maxchaintheo} provides overview and general structure of $\Rtwo$, and follows from a detailed arithmetical 
characterization of the relations $\leo$ and $\letwo$ in all of $\Rtwo$,
shown here through a generalization of Theorem 7.9 and Corollary 7.13 of \cite{CWc}
from the initial segment $\oneinf=\ups_1$ to all of $\On$, by
Theorem \ref{maintheo}, which for $i=1,2$ explicitly describes the $\le_i$-predecessors, and Corollary \ref{maincor}, which describes 
the $\le_i$-successors of a given ordinal. 
As a byproduct, a flaw in the proof of Theorem 7.9 of \cite{CWc} is corrected here in the generalized version, see also the paragraph
preceding Theorem \ref{maintheo}. It is worth mentioning that Theorem \ref{maxchaintheo} corrects the claim made in the first paragraph of
Section 8 of \cite{CWc}, where the role of the ordinals $\{\ups_{\la+1}\mid\la\in\Lim\}$ was overlooked.
The exact description of $<_i$-predecessors and $\le_i$-successors for $i=1,2$ (in particular the greatest one whenever such exists) of an 
ordinal $\al$ relies on its tracking chain $\tc(\al)$ and the notion of maximal extension ($\me$) of tracking chains. These notions are 
carefully elaborated in Section \ref{tstcsec} in the generalized form needed for an analysis of the entire structure $\Rtwo$, and can be 
derived from the term decomposition of $\al$ in Skolem hull notation. Tracking chains and their (maximal) extensions essentially make 
the surrounding (nested) $\le_i$-connectivity components visible, in which $\al$ is located. 

Incomplete fragments of the general big picture of $\Rtwo$ described above, repeat in a cofinal manner of growing complexity throughout 
all of $\Rtwo$, with the union of pointwise minimal isomorphic copies of all finite patterns comprising $\Core(\Rtwo)$, 
the universe of which is $\ups_1$ as shown in \cite{W18}. The $\lei$-relationships claimed to hold in Theorem \ref{maintheo} are verified
using Proposition \ref{letwocriterion} by transfinite induction through the ordinals. Conversely, starting from the backbone $\Image(\ups)$
provided by Theorem \ref{maxchaintheo} and applying reflection properties given by Lemma \ref{letwoupwlem} and the converse of 
Proposition \ref{letwocriterion} in a transfinitely iterated manner, we may eventually arrive back at (an isomorphic copy of) $\Rtwo$,
cf.\ \cite{C01,C09} for an elaboration of such an approach. As a consequence of Theorem \ref{maintheo} the converse of Propostion 
\ref{letwocriterion} holds in $\Rtwo$.  

The results established here enable us to show in \cite{W} that on the initial segment $\ups_{\om^2+2}$ the structures $\Rtwo$
and $\Rthree$ agree, where $\Rthree=\left(\On;\le,\leo,\letwo,\lethree\right)$ and
\[\al\le_i\be:\aeq \left(\al;\le,\leo,\letwo,\lethree\right) \preceq_{\Sigma_i} \left(\be;\le,\leo,\letwo,\lethree\right)\]
simultaneously for $i=1,2,3$ and recursively in $\be$, while \[\ups_{\om^2}\kthree\ups_{\om^2+2}\] is the least occurrence of a 
$\kthree$-pair in $\Rthree$. 
A detailed arithmetical analysis of the structure $\Rthree$, using extended arithmetical means as for the analysis of $\Rtwopl$, 
is the subject of ongoing work that will also be based on the present article.

\subsection*{Organization of this article, stand-alone readability}
The present article generalizes \cite{CWc}, with several improvements and corrections, and provides the necessary preparation for the
result in \cite{W} (Section 21.4) and for future work on $\Rthree$. 
\cite{W} starts out from an earlier, slightly longer, version of this introduction and reviews basic insights around patterns along 
with an outline of how they were discovered. 
The interested reader not yet familiar with patterns will find Sections 21.2 and 21.3 of \cite{W} helpful before reading 
the present article in detail. However, previous knowledge of \cite{W} is not required for understanding this article, which is intended
to be readable as stand-alone text. 
Clearly, for proofs and more details the reader is ultimately referred to the cited previous work on the subject, however, we have decided
to renew several proofs in order to increase accessibility of the work.

Section 2 provides a review of Skolem-hull based ordinal arithmetic developed for the analysis of patterns.

Section 3 introduces the reader to the generalized concepts of tracking sequences and chains. These provide the necessary machinery to 
describe the nested structure of $\leo$- and $\letwo$-connectivity components of $\Rtwo$. Besides the required generalization as compared to
\cite{CWc}, the exposition contains considerably more explanations, examples, and motivations, as well as renewed and simplified proofs.
The conceptual and technical improvement and simplification of Section 4 of \cite{CWc} in particular was first published in \cite{W17},
used in \cite{W18}, and is reviewed and extended for application in the present article. 
  
Section 4 contains the complete description of $\Rtwo$. The proof of the main theorem, Theorem \ref{maintheo}, begins with an overview
called \emph{proof map} and contains simple examples for the various cases that need to be discussed.

\section{Preliminaries}
Here we give a review of ordinal notational and arithmetical tools developed in \cite{W07a}, Section 5 of \cite{CWa}, and \cite{CWc}. 
For a reference to basic ordinal arithmetic we recommend Pohlers' book on proof theory, \cite{P09}. 

Let $\Limnod$ denote the class $\{0\}\cup\Lim$ where $\Lim$ is the class of limit ordinals. 
Generally, for a set or class $X$ of ordinals we define $\Lim(X)$\index{$\Lim(X)$} to be the set of all $\al\in X$ that are proper 
suprema of subsets of $X$.
By $\Hz$ we denote the class of additive 
principal numbers, i.e.\ the image of the $\om$-exponentiation function $\{\om^\eta\mid\eta\in\On\}$,  
and we write $\Hznod:=\{0\}\cup\Hz$ for the class of all ordinals that are closed under ordinal addition.
By $\Lz$ we denote the class $\{\om^\eta\mid\eta\in\Lim\}=\Lim(\Hz)$. 
By $\Mz$ we denote the class of multiplicative principal numbers, i.e.\ nonzero ordinals closed under ordinal multiplication.
Note that $\Mz\subseteq\Hz$ since \[\Mz=\{1\}\cup\{\om^{\om^\eta}\mid\eta\in\On\}.\]
As in the introduction, let $\Ez=\{\eta\mid\eta=\om^\eta\}$ denote the class of all epsilon numbers, 
i.e.\ fixed points of $\om$-exponentiation, and $\Ezone:=\{1\}\cup\Ez$. 
By $\Ez^{>\de}$ we denote the class of epsilon numbers that are greater than $\de$. Similarly, $\Hz^{\le\eta}$ denotes the set of 
additive principal numbers that are less than or equal to $\eta$, etc.

Representations in normal form are sometimes explicitly marked, such as additive normal form 
($\ANF$, weakly decreasing additive principal summands) or multiplicative normal form of ordinals $\al\in\Hz$
($\MNF$, weakly decreasing multiplicative principal factors), 
or Cantor normal form ($\CNF$) itself,
i.e.\ $\al=_\CNF\om^{\al_1}+\ldots+\om^{\al_n}$ such that $\al_1\ge\ldots\ge\al_n$, $n\ge 0$. Here $n=0$ covers the case $\al=0$.  
For $\al=_\ANF\al_1+\ldots+\al_n$ we define $\sumend(\al):=\al_n$ and set $\sumend(0):=0$,
while $\mc(\al)$ denotes the greatest additive component of $\al$, i.e.\ $\ale$ if $n>0$ and $0$ otherwise.
For an ordinal $\al=_\ANF\ale+\ldots+\aln$ the notation $\al=_\NF\be+\ga$ is a shorthand for $\be=_\ANF\ale+\ldots+\alnmin$ and $\ga=\aln$
if $n>0$, and $\be=\ga=0$ otherwise, for completenss.
This notation implies that $\ga=\sumend(\al)$.
For $\al\in\Hz$ such that $\al=_\MNF\ale\cdot\ldots\cdot\aln$ the notation $\al=_\NF\be\cdot\ga$ indicates that 
$\be=_\MNF\ale\cdot\ldots\cdot\alnmin$ and $\ga=\aln$ if $n>1$, and $\be=1$, $\ga=\al$ otherwise, for completeness.
This has the effect that the last multiplicative principal factor in the multiplicative normal form of $\al$ (also written as $\lf(\al)$) 
is equal to $\ga$. 

We define the ordinal $\al\minusp\beta$ as usual, namely in case of $\al\le\be$ to be $0$, while in case of $\al>\be$ to be the least $\ga$
s.t.\  $\al=\ga+\be$ if such $\ga$ exists, and $\al$ otherwise.
If $\al\le\be$ then we write $-\al+\be$\index{$-\al+\be$} for the unique $\ga$ such that $\al+\ga=\be$.
By $(1/\ga)\cdot\al$ we denote the least ordinal $\de$ such that $\al=\ga\cdot\de$, whenever such an ordinal exists.
We write $\al\mid\be$\index{$\al\mid\be$} if $\be$ is a (possibly zero) multiple of $\al$, i.e. $\exists\xi\,(\be=\al\cdot\xi)$.
i.e.\ $\al_1$ if $\al>0$ and $0$ otherwise.
$\logend(\al)$\index{logend} is defined to be $0$ if $\al=0$ and $\al_n$ if $\al=_\CNF\om^{\al_1}+\ldots+\om^{\al_n}$.
For $\al\in\Hznod$ we also write $\log(\al)$ instead of $\logend(\al)$.

For a function $f$ and a subset $X$ of its domain we denote the image of $X$ under $f$ by $f[X]=\set{f(x)}{x\in X}$.
Inequalities like $X<Y$ or $\al<X$ where $X,Y$ are sets of ordinals mean the conjunction of all inequalities taking each element of the concerning sets.
For sets $X$ and $Y$ we denote the set $\set{x}{x\in X\andsp x\not\in Y}$ by $X\setminus Y$.
Intervals of ordinals are often written in the following way: $(\al,\be)=\{\ga\mid\al<\ga<\be\}$, $[\al,\be]=\{\ga\mid\al\le\ga\le\be\}$,
and mixed forms $(\al,\be]$ and $[\al,\be)$ are defined analogously. Clearly, we also simply have $\al=\{\ga\mid\ga<\al\}$.

Sequences of ordinals (also called ordinal vectors) are often written as $\alvec=(\ale,\ldots,\aln)$, and appending an ordinal $\al$ is 
written as $\al^\frown\alvec=(\al,\ale,\ldots,\aln)$, $\alvec^\frown\al=(\ale,\ldots,\aln,\al)$. 
Concatenation of $\alvec$ with $\bevec=(\be_1,\ldots,\be_m)$ is writen as $\alvec^\frown\bevec=(\ale,\ldots,\aln,\be_1,\ldots,\be_m)$.
Similar notation is used for sequences of ordinal vectors.

\subsection{Stepwise collapsing functions \boldmath $(\thti)_{i\in\om}$ \unboldmath}

\cite{W06} can be seen as an introduction to the kind of ordinal arithmetic which is reviewed here, however, with ordinal addition
\emph{and} $\om$-exponentiation as basic functions over which Skolem hulling is performed, and with only one collapsing function $\tht$
instead of a family $(\thti)_{i\in\om}$ of stepwise collapsing functions.
\cite{W06} therefore only covers notations for ordinals below the Bachmann-Howard ordinal, but provides an elementary and quickly
accessible treatise of the Bachmann-Howard structure both in terms of hull notations and additive patterns of order one. 
However, we do not assume knowledge of \cite{W06} for the understanding of this article. 
The main reference for this subsection is Section 3 of \cite{W07a}, which contains detailed proofs.
Here we give an outline of the construction of systems $\Tt$, which along the lines of the introduction reduces to the formal 
definition of the system $(\thti)_{i\in\Om}$ of collapsing functions. The construction can be seen as a straightforward direct limit
construction, since $\Tt$ and $(\thti)_{i\in\Om}$ are obtained from increasing initial segments, i.e.\ through a limit of end-extensions.
Once established, the usefulness of these relativized notations becomes apparent through algebraic exactness and the absence of normal
conditions that would have to be permanently verified as is the case when one works with non-injective collapsing functions.
The usage of functions $\thti$ for increasing $i$ will allow us to arithmetically characterize pure patterns of an increasing maximal number 
of nestings of $\letwo$, as indicated in the examples following Lemma \ref{tenumlem}.

Let us fix the general setting for the definition of relativized hull notations as in the introduction:
\emph{Let $\tau\in\Ezone$, set $\Om_0:= \tau$,
suppose that $\Om_1$ is an uncountable regular cardinal number
greater than $\tau$, and let $\Om_{i+1}$ for $i\in(0,\om)$ be the regular cardinal successor of $\Om_i$.}

\subsubsection{\boldmath The notation systems $\Tnm$\unboldmath}
The following definition taken from \cite{W07a} is fundamental for the construction of systems $\Tt$. It provides the necessary 
support framework for the definition of $\thtnm$-functions, which in turn yield the desired $\thti$-functions via successive end-extension
for $n\to\om$. We finally obtain $\Tt$ as a direct limit of systems $(\T^n_0)_{n<\om}$. Apart from their role as auxiliary systems,
for $\tau=1$
the system $\T^1_0$ provides a notation system suited for an analysis of Peano arithmetic ($\PA$), as $\T^1_0\cap\Om=\epsn$, 
and $\T^{n+1}_0$ provides a notation system for the theory $\idn$ of $n$-times iterated inductive definitions, for $n\in(0,\om)$.

\begin{defi}[3.1 of \cite{W07a}]
Let $n\in(0,\om)$. Descending from $m=n-1$ down
to $m=0$ we define sets of ordinals $\Cnmalbe$\index{$\Cnm$} where $\be<\Omme$ and ordinals $\thtnm(\al)$\index{$\thtnm$} by simultaneous recursion on $\al<\Om_{m+2}$.

\noindent For each $\be<\Omme$ the set $\Cnmalbe$ is defined inductively by
\begin{itemize}
\item $\Om_m\cup\be\sub\Cnmalbe$
\item $\xi,\eta\in\Cnmalbe\;\imp\;\xi+\eta\in\Cnmalbe$
\item $\xi\in\Cnmalbe\cap\Om_{k+2}\;\imp\;\thtnk(\xi)\in\Cnmalbe$ for $m<k<n$
\item $\xi\in\Cnmalbe\cap\al\;\imp\;\thtnm(\xi)\in\Cnmalbe$.
\end{itemize}
Having defined $\thtnm(\xi)$ for all $\xi<\al$ and $\Cnmalbe$ for every $\be<\Omme$ we set
\[\thtnmal:=\min(\set{\xi<\Om_{m+1}}{\Cnm(\al,\xi)\cap\Omme\sub\xi\wedge\al\in\Cnm(\al,\xi)}
                  \cup\singleton{\Omme}).\]
\end{defi}

Note that we have $\thtnm(0)=\Om_m$ for $m<n$.
The function $\thtnnmin$ is not a proper collapsing function simply because $\Om_n\not\in\Cnmalbe$ which will become
clear in the sequel.
The next two lemmas follow immediately from the above definition.
\begin{lem}[3.2 of \cite{W07a}] Let $\al, \al_1, \al_2, \ga<\Ommz$ and $\be, \be_1, \be_2, \de<\Omme$.
\begin{enumerate}
\item[a)] If $\de\sub\Cnmalbe$ then $\Cnm(\al,\de)\sub\Cnmalbe$.
\item[b)] For $\al_1\le\al_2$ and $\be_1\le\be_2$ we have $\Cnm(\al_1,\be_1)\sub\Cnm(\al_2,\be_2)$.
\item[c)] We have $\Cnmalbe=\bigcup_{\ga<\al}\Cnm(\ga,\be)$ for $\al\in\Lim$ and similarly we have
$\Cnmalbe=\bigcup_{\de<\be}\Cnm(\al,\de)$ for $\be\in\Lim$.
\item[d)] $\card(\Cnmalbe)<\Omme$.
\end{enumerate}
\end{lem}
In the following we write $\Cnm(\Ommz,\be)$ for $\bigcup_{\al<\Ommz}\Cnmalbe$ and
$\Cnm(\al,\Omme)$ for $\bigcup_{\be<\Omme}\Cnmalbe$.

\begin{lem}[3.3 of \cite{W07a}]\label{thtseglem} Let $\al<\Ommz$ and $\be<\Omme$.
\begin{enumerate}
\item[a)] $\thtnmal=\Cnm(\al,\thtnmal)\cap\Omme$.
\item[b)] $\thtnmal\in\Hz\cap[\Om_m,\Omme]$.
\item[c)] Let $\xi=_\ANF\xi_1+\ldots+\xi_l$. Then $\xi\in\Cnmalbe$ iff $\xi_1,\ldots,\xi_l\in\Cnmalbe$.
\end{enumerate}
\end{lem}

As in Lemma \ref{tenumlem} additive principal numbers which are not epsilon numbers can be characterized as follows.
Recall that the function $\omq{\cdot}$ enumerates the non-epsilon additive principal numbers.
\begin{lem}[3.5 of \cite{W07a}]\label{enumlem} Let $m<n$. For every $\al<\Omme$ we have
\[\thtnm(1+\al)=\omq{\Om_m+\al}.\]
\end{lem}

\begin{defi}[3.6 of \cite{W07a}] For $m<n$ we define
\[\Tnm:=\Cnm(\Ommz,0)\] and set $\Tnn:=\Om_n$ for convenience.\index{$\Tnm$}
\end{defi}

\begin{lem}[3.7 of \cite{W07a}] For $m<n$ the set $\Tnm$ is inductively characterized as follows:
\begin{itemize}
\item $\Om_m\sub\Tnm$
\item $\xi,\eta\in\Tnm\;\imp\;\xi+\eta\in\Tnm$
\item $\xi\in\Tnm\cap\Omkz\;\imp\;\thtnk(\xi)\in\Tnm$ for $m\le k<n$.
\end{itemize}
\end{lem}

The next lemma the most important claim of which is that the $\tht$-functions are collapsing functions depends
on the regularity of the cardinals $\Om_n$ where $0<n<\om$. 
In \cite{W06} we showed that the $\tht$-function defined there is a total collapsing function on the segment
$\epsom$, the least epsilon number greater than $\Om$. 
The analogy here is that $\thtnm$ is a total collapsing function on the set $\Cnm(0,\Omme)\cap\Ommz$
(for $m<n-1$). The latter set is actually the largest segment of ordinals having notations in $\Tnme$.

\begin{lem}[3.8 of \cite{W07a}]\label{collapslem} Let $m<n$. For all $\al\in\Cnm(0,\Omme)\cap\Ommz$ we have
\[\thtnmal<\Omme\;\mbox{ and }\;\thtnmal\not\in\Cnm(\al,\thtnmal).\]
\end{lem}

\begin{cor}[3.9 of \cite{W07a}]\label{idsetslem} For $m<n$ we have
\[\Tnm\sub\Tnme=\Cnm(0,\Omme)=\Cnm(\Ommz,\Omme).\]
\end{cor}

In order to compare $\thtnm$-terms we need to detect the additive principal parts of ordinals in $\Tnme$.
This is done by the following definition. The symbol $\finsub$ indicates a finite subset.

\begin{defi}[3.10 of \cite{W07a}] Let $m<n$. By recursion on the definition of $\Tnme$ we define $\Pnm(\xi)\finsub\Omme$\index{$\Pnm$} for
every $\xi\in\Tnme$.
\begin{itemize}
\item $\Pnm(\xi):=\singleton{\xi_1,\ldots,\xi_r}$, if $\xi=_\ANF\xi_1+\ldots+\xi_r<\Omme$
\item $\Pnm(\xi):=\Pnm(\xi_1)\cup\Pnm(\xi_2)$, if $\xi_1,\xi_2\in\Tnme$ and $\xi=_\NF\xi_1+\xi_2>\Omme$
\item $\Pnm(\xi):=\Pnm(\eta)$, if $\xi=\thtnk(\eta)$, $\eta\in\Tnme\cap\Omkz$, $m<k<n$.
\end{itemize}
We define $\xi^\starnm:=\max((\Pnm(\xi)\setminus\Om_m)\cup\singleton{0})$\index{$\starnm$} for $\xi\in\Tnme$.
\end{defi}

According to the next lemma, by this definition $\Pnm(\xi)$ is uniquely determined for every
$\xi\in\Tnme$. 
The lemma provides a criterion for the comparison of ordinals within $\Tnm$
which is elementary recursive in $\Om_m$.

\begin{lem}[3.11 of \cite{W07a}]\label{propslem}
For $n>0$ and $m\in\singleton{0,\ldots,n-1}$ we have
\begin{enumerate}
\item[a)] $\Pnm$ is well defined.
\item[b)] Let $\al<\Ommz$ and $\be<\Omme$. For every $\xi\in\Tnme$ we have
\[\xi\in\Cnmalbe \aeq \Pnm(\xi)\sub\Cnmalbe.\]
\item[c)] $\al^\starnm<\thtnmal$ for all $\al\in\Tnme\cap\Ommz$.
\item[d)] The restriction of $\thtnm$ to $\Tnme\cap\Ommz$ is $1$-$1$.
We have
\[\thtnmal<\thtnm(\ga) \aeq \left(\al<\ga \andsp \al^\starnm<\thtnm(\ga)\right) \vee \thtnmal\le\ga^\starnm\]
for all $\al,\ga\in\Tnme\cap\Ommz$.
\end{enumerate}
\end{lem}

As mentioned in the introduction, below $\Omme$ the Skolem hulls $\Cnm(\al,\be)$ form initial segments of ordinals.
This property is essential in their role to provide the basis for an ordinal notation system.

\begin{lem}[3.12 of \cite{W07a}]\label{segmentlem}
For $m<n$ we have
\begin{enumerate}
\item[a)] $\Cnm(\al,\al^\starnm+1)=\Cnm(\al,\thtnmal)$ for all $\al\in\Tnme\cap\Ommz$.
\item[b)] $\Cnmalbe\cap\Omme\in\On$ for all $\al<\Ommz$ and all $\be<\Omme$.
\end{enumerate}
\end{lem}

For $m<n$ and $k\in\singleton{m,\ldots,n-1}$ the ordinal
$\thetank$ defined below is the supremum of all ordinals in the segment $[\Om_k,\Omke)$ that have
a notation within $\Tnm$. Moreover, according to the theorem below, $\thetanm$ is the maximal segment of ordinals having a notation in $\Tnm$.

\begin{defi}[3.13 of \cite{W07a}] \label{thetadefi} Let $n>0$.
For $k<\om$ we define \index{$\Thetanm$}\index{$\thetanm$}
\[\Thetannmin(k):={\thtnnmin}^{(k)}(0)\;\mbox{ and }\;\thetannmin:=\sup_{k<\om}\,\Thetannmin(k)\]
where ${\thtnnmin}^{(k)}(0)$ denotes the $k$-fold application of $\thtnnmin$ to $0$.
Descending from $m=n-2$ down to $m=0$ we define
\[\Thetanm(k):=\thtnm(\Thetanme(k))\;\mbox{ and }\;\thetanm:=\sup_{k<\om}\,\Thetanm(k).\]
For convenience of notation we set $\Thetann(k):=\Thetannmin(k)$ for $k<\om$, $\thetann:=\thetannmin$ 
(since the $\thtnnmin$-function is not a collapsing function),
and $\thetanne:=0$.
\end{defi}

\begin{theo}[3.14 of \cite{W07a}]\label{tnmthm}
For $m<n$ we have \[\Tnm\cap\Omme=\thetanm.\]
\end{theo}

\begin{cor}[3.15 of \cite{W07a}]\label{thtsubscor}
Let $m<n$, $\al<\Ommz$, and $\be<\Omme$. Then all terms of a shape $\thtnk(\eta)$ in $\Cnmalbe$ satisfy $\eta<\thetanke$.
\end{cor}

\subsubsection{\boldmath The notation system $\Tt$\unboldmath}

The following lemma is central to see that we can obtain $\Tt$ as a direct limit of the $\T^n_0$.

\begin{lem}[3.20 of \cite{W07a}]\label{endextlem} For all $m\le n$ we have $\thetanm<\thetanem$, and for $m<n$
\[\Cnmalbe\cap\thetanme=\Cnemalbe\cap\thetanme \quad\&\quad \thtnmal=\thtnemal\]
for all $\al<\thetanme$ and all $\be<\Omme$.
\end{lem}

\begin{cor}[3.21 of \cite{W07a}]
Let $m<n$. $\Tnm\sub\Tnem$ and for $k$ such that $m\le k<n$ the functions $\thtnk$ and $\thtnek$\index{$\thetam$}
agree on $\Tnm\cap\Omkz$.
\end{cor}

\begin{defi}[\boldmath $\thtm$ and $\Tm$\unboldmath, 3.22 of \cite{W07a}, cf.\ Def.\ \ref{Ttdefi}]\label{Tmdef} For $m<\om$ we set
\[\thetam:=\sup_{n>m}\thetanm,\] 
and we define a function $\thtm:\thetame\to\Omme$\index{$\thtm$} by \[\thtmal:=\thtnmal\] for
$\al<\thetame$ where $n>m$ is large enough to satisfy $\al<\thetanme$.
$\Tm$\index{$\Tm$} is defined inductively as follows:
\begin{itemize}
\item $\Om_m\sub\Tm$
\item $\xi,\eta\in\Tm\;\imp\;\xi+\eta\in\Tm$
\item $\xi\in\Tm\cap\Omkz\;\imp\;\thtk(\xi)\in\Tm$ for $k\ge m$.
\end{itemize}
\end{defi}

It is immediate from the previous lemma that the functions $\thtm$ are well defined.
The well-definedness of $\Tm$, which means that $\sup(\Tm\cap\Omkz)\le\thetake$ where $m\le k$, follows from the next theorem.
This theorem establishes the systems of relativized ordinal notations based on Skolem hull operators we aim for.

\begin{theo}[3.23 of \cite{W07a}] For $m<\om$ we have $\Tm=\bigcup_{n>m}\Tnm$ and
\[\Tm\cap\Omme=\thetam=\sup_{n\ge m}\thtm(\cdots(\thtn(0))\cdots).\]
\end{theo}

\begin{cor}[3.24 of \cite{W07a}] For every $k\ge m$ we have $\sup(\Tm\cap\Omkz)=\thetake$.
For $m<n$ we have \[\Tnm\cap\Omme=\thetanm=\thtm(\cdots(\thtn(0))\cdots).\]
\end{cor}

Notice that by the end-extension property shown in the previous lemma and theorem it follows that
each $\Tm$ again gives rise to a notation system with parameters from $\Om_m$ that provides a unique term
for every ordinal which is element of some $\Tnm$ where $n>m$ (refining the clause for ordinal addition
with a normal form condition as mentioned in the introduction).
The comparison of $\tht$-terms
in $\Tm$ can be done within a sufficiently large fragment $\Tnm$ where $n>m$.
The notation system $\Tm$ as well as the criterion for the comparison of its elements
are now easily seen to be elementary recursive in $\Om_m$.
From now on we will only need to consider the notation system $\T_0$.

\begin{conv}[cf.\ 3.25 of \cite{W07a}] In our setting the ordinal notations are relativized to $\tau$.
Later on we will indicate this explicitly in writing $\thtt$\index{$\thtt$} and $\Tt$\index{$\Tt$} instead of $\tht_0$ and $\T_0$.
As already defined in the introduction, we have \[\tauinf:=\Tt\cap\Om=\theta_0.\]
\end{conv}

The notion of term height in the following sense is often useful in inductive proofs. It recovers the least fragment in which to
find the given notation (modulo the trivial embedding of a system $\Tnm$ into $\Tm$).
\begin{defi}[3.26 of \cite{W07a}]\label{htdef} We define a function $\htt:\Tt\to\om$\index{$\htt$} as follows:
\[\htt(\al):=\left\{\begin{array}{ll}
   m+1 & {\small \mbox{ if } m=\max\set{k}{\mbox{there is a subterm of $\al$ of shape $\thtk(\eta)$}}}\\
   0   & {\small \mbox{ if such $m$ does not exist.}}
   \end{array}\right.\]
\end{defi}

\begin{lem}[3.27 of \cite{W07a}]\label{htlem}
For $\al<\Tt\cap\Om_1$
\[\htt(\al)=\min\set{n}{\al<\thtnod(\cdots(\thtn(0))\cdots)}.\]
$\htt$ is weakly increasing successively as the indices of occurring $\tht$-functions increase.
\end{lem}

The following notion of subterm is crucial for the comparison of $\tht$-terms. Subterms of lower cardinality become parameters when
comparing $\tht$-terms of higher cardinality, a natural property if we take into account that $\tht$-functions are collapsing functions.
\begin{defi}[3.28 of \cite{W07a}]\label{subtermdefi}
We define sets of subterms $\subtm(\al)$\index{$\subtm$} for $m<\om$ and notations $\al$ in $\Tt$ by recursion on
the build up of $\Tt$:
\begin{itemize}
\item $\subtm(\al):=\singleton\al$ for parameters $\al<\tau$
\item $\subtm(\al):=\singleton{\al}\cup\subtm(\xi)\cup\subtm(\eta)$
for $\xi,\eta\in\Tt$ s.t. $\al=_\NF\xi+\eta>\tau$
\item $\subtm(\al):=\left\{\begin{array}{cl}
                                  \singleton\al\cup\subtm(\xi) & \mbox{ if } \;k\ge m\\
                                  \singleton\al & \mbox{ if } \;k<m
                                \end{array}\right.$\\
for $\al=\thtk(\xi)$ where $\xi\in\Tt\cap\Omkz$.
\end{itemize}
We define the additive principal parts of level $m$ of $\al\in\Tt$ by\index{$\Pm$}\index{$\starm$}
\[\Pm(\al):=\subtm(\al)\cap\Hz\cap[\Om_m,\Om_{m+1}) \mbox{ and } \al^\starm:=\max\left(\Pm(\al)\cup\singleton0\right).\]
The set of parameters $<\tau$ used in the unique term denoting some $\al\in\Tt$ is denoted by\index{$\Part$}
\[\Part(\al):=\subtnod(\al)\cap\tau.\]
\end{defi}

\begin{conv}[3.29 of \cite{W07a}]
In order to make the setting of relativization explicit we write $\Ptm$\index{$\Ptm$} for $\Pm$ and ${}^\startm$ for ${}^\starm$.
We write $\Pt$\index{$\Pt$} for $\Ptnod$, and instead of ${}^\starnod$ we will also write ${}^\start$.\index{$\start$}
\end{conv}

\begin{rmk}[Remarks following 3.29 of \cite{W07a}]\mbox{ }
\begin{enumerate}
\item $\subtm(\al)$ consists of the subterms of $\al$ where parameters below $\tau$ as well as subterms of a shape
$\thtk(\eta)$ with $k<m$ are considered atomic.
\item $\Pm$ consists of the $\thtm$-subterms of $\al$ which are not in the
scope of a $\thtk$-function with $k<m$.
\item By lemma \ref{propslem} part c) and the end extension properties shown above it follows that the notion ${}^\starm$ is consistent with the notion ${}^\starnm$ where $m<n$ on the common
domain. It also follows that
\[\al^\starm=\max\left(\subtme(\al)\cap\Hz\cap[\Om_m,\Om_{m+1})\cup\singleton0\right).\]
\item The notion $\Pm$ takes more subterms into consideration than $\Pnm$ since $\subtm$ also decomposes
$\thtm$-subterms. However, $\Pme$ is consistent with $\Pnm$ on the common domain.
\item In order to clarify the definition of $\Part$ consider the following examples: $\operatorname{Par}^{\epsn}(\om+1)=\singleton{\om+1}$ and
$\operatorname{Par}^{\epsn}(\epsn+\om+1)=\singleton{1,\om}$.
\end{enumerate}
\end{rmk}

The following lemma concerning $\tht$-terms within $\Tt$ and their comparison will be
used frequently without further mention. 

\begin{lem}[3.30 of \cite{W07a}]\label{thttcomplem}
For $m<\om$ the function
$\thtm\restriction_{\Tt\cap\Ommz}$ is $1$-$1$ and has values in $\Hz\cap[\Om_m,\Omme)$.
Let $\al,\ga\in\Tt\cap\Ommz$. Then $\al^\starm<\thtm(\al)$ and
\[\thtm(\al)<\thtm(\ga) \;\aeq\; \left(\al<\ga\wedge\al^\starm<\thtm(\ga)\right)\vee\thtm(\al)\le\ga^\starm.\]
\end{lem}

\begin{rmk}\label{extrarmk} Note that in particular we have
\begin{enumerate}
\item $\Pt(\al)$ is the set of all subterms of $\al$ of a form $\thtt(\xi)$ for some $\xi$.
\item $\al^\start=\max(\Pt(\al)\cup\singleton{0})$.
\item $\al=\al^\start$ whenever $\al=\thtt(\xi)$ for some $\xi$.
\item $\Part(\al)$ is the set of parameters $<\tau$ used in the unique term denoting $\al$.
\end{enumerate}
\end{rmk}

Note that the ordinal defined by a $\tht$-term, say $\tht(\xi)$, is characterized as the least 
$\theta>\xi^\star$ that is closed under parameters and basic functions such that $\tht(\ze)<\theta$ for all $\ze<\xi$ satisfying 
$\ze^\star<\theta$, cf.\ Lemma 4.10 of \cite{W06}.

\subsection{Localization}\label{localizationsec}

The notion of \emph{localization} introduced in Section 4 of \cite{W07a} and refined in \cite{CWa} is to be understood in terms of 
closure properties or fixed point levels in the sense discussed in the introduction. As an example, given an ordinal $\al\not\in\Ez$ 
we ask for the greatest epsilon number $\be<\al$, if such exists. The ordinal $\be$ is of higher fixed point level and has stronger 
closure properties than $\al$, as $\be$ is closed under $\om$-exponentiation. The refinement also considers degrees of limit point thinning,
i.e.\ (transfinitely) iterated applications of the operation $\Lim(\cdot)$.

We now formally fix a notation already used in the introduction, and define $\alplus$ for given $\al=\thtt(\De+\eta)$, which is the 
least ordinal greater than $\al$ that is of the same fixed point level as $\al$, as becomes clear shortly. 
All ordinals mentioned are assumed to be represented in $\Tt$.

\begin{conv}[cf.\ 4.1 of \cite{W07a}]\label{thtargconv} $\al=\thtt(\De+\eta)$ automatically means that $\Om_1\mid\De$ and $\eta<\Om_1$. 
In a situation where some $\al=\thtt(\De+\eta)$ is fixed, we will write $\alplus$ for $\thtt(\De+\eta+1)$. 
\end{conv} 

We will apply this notation frequently and use Greek capital letters to indicate that part of the argument which is
a (possibly zero) multiple of $\Om_1$.
Note that for $\De$, to be a proper multiple of $\Om$ means to be a sum of $\tht_1$-terms. $\eta$ can be additively composed of $\thtt$-terms
and parameters below $\tau$.
Recall Lemma \ref{tenumlem} from the introduction, which characterizes additive principal numbers that are not epsilon numbers in the 
interval $(\tau,\tauinf)$. Epsilon numbers in the interval $(\tau,\tauinf)$ are characterized as follows.

\begin{lem}[4.3 of \cite{W07a}] For $\al=\thtt(\De+\eta)\in\Tt$ we have $\al\in\Ez^{>\tau}$ if and only if $\De>0$.
\end{lem}

The (informal) notion of fixed point level is justified by the following lemma. We sometimes also call $\De$ the fixed point level of an
ordinal $\thtt(\De+\eta)$.

\begin{lem}[4.4 of \cite{W07a}]\label{fixlevcarlem} 
For $\al=\thtt(\De+\eta)\in\Tt$ we have $\eta=\sup_{\si<\eta}\thtt(\De+\si)$ if and only if $\eta=\thtt(\Ga+\rho)$ with $\Ga>\De$ 
and $\eta>\De^\start$.
\end{lem}

An immediate consequence of Lemma \ref{thttcomplem} shows that for any $\al$ of the form $\thtt(\De+\eta)$ the interval $(\al,\alplus)$ 
does not contain ordinals of fixed point level greater than or equal to $\De$. 

\begin{lem}[4.5 of \cite{W07a}]
For $\al=\thtt(\De+\eta)\in\Tt$ and $\thtt(\Ga+\rho)\in(\al,\alplus)$ we have $\Ga<\De$.
\end{lem}

Recalling part 3 of Remark \ref{extrarmk}, we are now prepared to make sense of the notion of localization.

\begin{defi}[4.6 of \cite{W07a}]\label{localizationdef}
Let $\al=\thtt(\De+\eta)\in\Tt$.
We define a finite sequence of ordinals as follows:
Set $\al_0:=\tau$. Suppose $\al_n$ is already defined and $\al_n<\al$. Let
$\al_{n+1}:=\thtt(\xi)\in\Pt(\al)\setminus(\al_n+1)$ where $\xi$ is maximal.
This yields a finite sequence $\tau=\al_0<\ldots<\al_n=\al$ for some $n<\om$ which we call the
\boldmath{\bf $\tau$-localization}\unboldmath\ of $\al$.\index{localization!$\tau$-localization of $\al$}  
\end{defi} 

For example, let $\al=\thtnod(\thte(\tht_2(0)))$ be the Bachmann-Howard ordinal,
$\be$ be the least $\Ga$-number greater than $\al$, i.e\ $\be=\Ga_{\al+1}=\thtal(\Om^2)=\thtnod(\thte(\thte(0))+\al)$, 
and $\ga=\eps_{\be+1}=\thtbe(\thte(0))=\thtnod(\thte(0)+\be)$ be the least epsilon number greater than $\be$.
Then the $1$-localization of $\ga$ is $(1,\al,\be,\ga)$, the $\al$-localization of $\ga$ is $(\al,\be,\ga)$, 
whereas the $\be$-localization of $\ga$ is just $(\be,\ga)$ and the $\ga$-localization of $\ga$ is simply $(\ga)$,
the trivial localization.
 
Note that $\al_1,\ldots,\al_n$, in case of $n\ge 2$, forms by definition a sequence of $\thtt$-terms of strictly decreasing arguments.
The following lemma summarizes properties of localization starting from this observation.

\begin{lem}[4.7, 4.8, and 4.9 of \cite{W07a}]\label{localipic}
Let $\al=\thtt(\De+\eta)\in\Tt$, $\al>\tau$ and $\tlocaln$ be the $\tau$-localization of $\al$ where $\al_i=\thtt(\De_i+\eta_i)$ for $i=1,\ldots,n$. Then 
\begin{itemize} 
\item[a)] For $i<n$ and any $\be=\thtt(\Ga+\rho)\in(\al_i,\al_{i+1})$ we have $\Ga+\rho<\De_{i+1}+\eta_{i+1}$.
\item[b)] $(\De_i)_{1\le i\le n}$ forms a strictly descending sequence of multiples of $\Om_1$.
\item[c)] For $i<n$ the sequence $(\al_0,\ldots,\al_i)$ is the $\tau$-localization of $\al_i$.
\end{itemize} 
The guiding picture for localizations is \[\tau<\al_1<\ldots<\al_n=\al<\alplus=\al_n^+<\ldots<\al_1^+.\]
\end{lem}

The notion of \emph{fine localization} is obtained by iterated application of the operator $\al\mapsto\albar$
from \cite{W07a}, defined there in the proof of Lemma 8.2, and extended in Section 5 of \cite{CWa} as follows.
While the predecessor in the $\tau$-localization of an ordinal $\al\in(\tau,\tauinf)$ is the greatest ordinal below $\al$ 
that has a (strictly) greater fixed point level than $\al$, if such ordinal exists, and $\tau=\alnod$ otherwise, the predecessor
$\albar$ of $\al$ in the $\tau$-fine-localization of $\al$ is the greatest ordinal below $\al$ that is of the same fixed point level and of a 
degree of limit point thinning greater than or equal to that of $\al$, if such ordinal exists, and the predecessor of $\al$ in its 
$\tau$-localization of $\al$ otherwise. 
According to Corollary 6.3 of \cite{CWa} a pattern notation (of least cardinality) for an ordinal $\al\in\Core(\Ronepl)$, i.e.\ $\al<\oneinf$, 
is obtained by closure of the set $\{0,\al\}$ under additive decomposition, the function $\lh$ for $\Ronepl$, and the operator $\bardot$.
Closure under additive decomposition means that for any ordinal $\be=_\ANF\be_1+\ldots+\be_n$ in the set, also the ordinals 
$\be_i, \be_1+\ldots+\be_i$ for $i=1,\ldots,n$ are in the set.

\begin{defi}[5.1 of \cite{CWa}]\label{barop}
Let $\al=\thtt(\De+\eta)\in\Tt$, $\al>\tau$, and $\tau=\al_0,\ldots,\al_n=\al$
be its $\tau$-localization.
$\albar\in[\tau,\al)$ is defined as follows.
\begin{itemize}
\item Suppose $\eta$ is of the form $\etapr+\eta_0$ where $\eta_0\in\Hz$ and either $\etapr=0$ or
$\eta=_\NF\etapr+\eta_0$. Then if either $\eta_0=1$ or $\etapr<\sup_{\si<\etapr}\thtt(\De+\si)$
we set $\albar:=\thtt(\De+\etapr)$.
\item In all remaining cases we let $\albar:=\al_{n-1}$.
\end{itemize}
\end{defi}

In case of $\al\in\Lz$, i.e.\ for $\al$ a limit of additive principal numbers, the above definition is consistent with the definition 
given in \cite{W07a}. For $\al\in\Hz\setminus\Lz$, that is, for $\al$ of a shape $\al=\om^{\alpr+1}$ the above definition yields 
$\albar=\om^\alpr$.
For further clarification, consider the following examples, in which we assume that $\thtnod=\thtt$ for $\tau=1$: 
\begin{enumerate}
\item $\overline{\epsn\cdot\om}=\epsn$, where $\epsn\cdot\om=\thtnod(\thtnod(\thte(0)))$, 
\item $\overline{\eps_{\om+\om}}=\eps_\om$, where $\eps_{\om+\om}=\thtnod(\thte(0)+\om+\om)$ and $\om=\thtnod(\thtnod(0))$, 
\item $\overline{\eps_{\om^2}}=1$, where $\eps_{\om^2}=\thtnod(\thte(0)+\om^2)$ and $\om^2=\thtnod(\thtnod(0)+\thtnod(0))$, 
\item $\overline{\eps_{\Ga_0+1}}=\Ga_0=\overline{\eps_{\Ga_0+\om}}$, where 
  $\eps_{\Ga_0+1[+\om]}=\thtnod(\thte(0)+\thtnod(\thte(\thte(0)))[+\om])$, 
\item $\overline{\eps_{\Ga_0+2}}=\eps_{\Ga_0+1}$, where $\eps_{\Ga_0+2}=\thtnod(\thte(0)+\thtnod(\thte(\thte(0)))+\thtnod(0))$,
\item $\overline{\Ga_{\om^\om+\om^2}}=\Ga_{\om^\om}$, where 
  $\Ga_{\om^\om[+\om^2]}=\thtnod(\thte(\thte(0))+\thtnod(\thtnod(\thtnod(0)))[+\om^2])$.
\end{enumerate}
See also Lemma \ref{simplebarlem} for a more general statement about ordinals $\al\in\Hz\setminus\Mz$.

As pointed out in \cite{CWa}, the iterated application of $\bardot$ to $\al$ leads to $\al_{n-1}$ after finitely many steps,
from there to $\al_{n-2}$ and finally to $\tau$. This follows from the lemmas concerning $\tau$-localization cited below.

\begin{lem}[5.2 of \cite{CWa}]\label{loclexordlem} Let $\al,\be\in(\tau,\tauinf)\cap\Hz$ with $\al<\be$. For their $\tau$-localizations 
$\tau=\al_0,\ldots,\al_n=\al$ and $\tau=\be_0,\ldots,\be_m=\be$ we have
\[\alvec :=(\al_1,\ldots,\al_n)\klex(\be_1,\ldots,\be_m)=:\bevec .\]
\end{lem}

\begin{lem}[5.3 of \cite{CWa}]\label{subtermloclem}
Let $\al=\thtt(\xi)\in\Tt$ with $\tau$-localization $\alvec :=(\al_0,\ldots,\al_n)$
and $\be\in\Pt(\xi)$.
If there is $\al_i\in\Pt(\be)$ where $0\le i\le n$ then $(\al_0,\ldots,\al_i)$ is an initial
sequence of the $\tau$-localization of any $\ga\in[\be,\al]$.
\end{lem}

\begin{lem}[5.4 of \cite{CWa}]\label{albarloclem}
Let $\al\in(\tau,\tauinf)\cap\Hz$ with $\tau$-localization $(\al_0,\ldots,\al_n)$.
Then the $\tau$-locali\-zation of $\albar$ is $(\al_0,\ldots,\al_{n-1}=\albar)$ if $\albar\in\Pt(\al)$
or $(\al_0,\ldots,\al_{n-1},\albar)$ otherwise. In the latter case for $\al=\thtt(\De+\eta)$ and
$\albar=\thtt(\Ga+\rho)$ we have
\[(\De+\eta)^\start=(\Ga+\rho)^\start.\]
\end{lem}

\begin{defi}[5.5 of \cite{CWa}]
The \boldmath{\bf $\tau$-fine-localization}\unboldmath\   of $\al\in(\tau,\tauinf)\cap\Hz$ is defined to be either
$(\tau,\al)$ if $\albar=\tau$ or the $\tau$-fine-localization of $\albar$ concatenated with $\al$.
\end{defi}

\begin{lem}[5.6 of \cite{CWa}]\label{finelocinitialseqlem}
The $\tau$-localization of $\al$ is a subsequence of the $\tau$-fine-localization of $\al$.
\end{lem}

\begin{lem}[5.9 of \cite{CWa}]\label{klexfineloclem} 
Let $\al,\be\in[\tau,\tauinf)\cap\Hz$ with $\tau$-fine-localizations $\alvec=(\al_0,\ldots,\al_n)$ and $\bevec=(\be_0,\ldots,\be_m)$ be given.
Then $\al<\be$ if and only if $\alvec\klex\bevec$.
\end{lem}

In order to obtain a formal notion of limit point thinning in the context of $\Tt$, we first characterize the function $\log$
introduced earlier in terms of $\Tt$. The indicated correction corrects for an unintended flaw in the original formulation. 
\begin{lem}[corrected 4.10 of \cite{W07a}]\label{logendcharlem}
For $\al=\thtt(\De+\eta)\in\Tt$ we have
\[\log(\al)=\left\{\begin{array}{ll}
                        \al & {\small \mbox{ if }\De>0}\\
                        \eta+1 & {\small \mbox{ if } \De=0 \mbox{ and } \eta=\de+n \mbox{ such that } \de\in\Ez^{>\tau}, n<\om}\\
                        (-1+\tau)+\eta & {\small \mbox{ otherwise.}}
                      \end{array}\right.\]           
\end{lem}

We now define an operator $\ze$ that given a setting of relativization $\tau$ and an ordinal $\al$ in the image of the enumeration function 
$\eta\mapsto\thtt(\De+\eta)$, $\eta<\tauinf$, outputs the \emph{degree of limit point thinning of $\al$}, written as $\zetal$.
\begin{defi}[4.11 of \cite{W07a}]\label{zetaldefi}
For $\al=\thtt(\De+\eta)\in\Tt$ we define 
\[\zetal:=\left\{\begin{array}{ll}
                    \logend(\eta) & {\small \mbox{ if } \eta<\sup_{\si<\eta}\thtt(\De+\si)}\\
                    0 & {\small \mbox{ otherwise.}}
                 \end{array}\right.\] 
\end{defi}

The connection between the operators $\bardot$ and $\ze$ cannot be stated yet. We first need to introduce the notions of 
\emph{base transformation}, \emph{cofinality operator}, and \emph{translation}.

\subsection{Base transformation, cofinality operators, and translation}\label{arithmeticsubsec}

Here we still cite results from \cite{W07a,CWa}, but we summarize as in Subsection 2.3 of \cite{W07c}.

\emph{Base transformation}, as introduced and treated in detail in Section 5 of \cite{W07a}, is the crucial notion which allows us to express
essential uniformity properties in both the development of a strong ordinal arithmetic and the formulation of characteristic properties of elementary substructure on the ordinals. It enables a precise comparison of ordinals modulo their relativizations. 

\begin{defi}[\boldmath $\Tts$, $\pist$\unboldmath]
Let $\tau\in\Ez$ and $\si\in\singleton{1}\cup\Ez\cap\tau$.
\begin{itemize}
\item $\Tts:=\set{\al\in\Tt}{\Part(\al)\sub\si}$. 
\item The base transformation $\pist:\Tts\to\Ts$ maps $\al\in\Tts$ to the ordinal one obtains from the term
representation of $\al$ in $\Tt$ by substituting every occurring function $\thtt$ by $\thts$, i.e.\ $\pist(\al):=\al[\thtt/\thts]$. 
\end{itemize}
\end{defi}

\begin{lem}[5.3, 5.5, and 5.6 of \cite{W07a}] 
$\pist$ is a $(<,+)$-isomorphism, and base transformation commutes with the process of
localization as well as with the degree of limit point thinning (i.e.\ $\pist(\zetal)=\ze^\si_{\pist(\al)}$).
For $\al\in\Tts$ we have $\Part(\al)=\Parsi(\pist(\al))$.
\end{lem}

The \emph{cofinality operators $\iotal$ and $\la^\tau$}, introduced in Section 7 of \cite{W07a}, allow for the exact classification of the
\emph{cofinality properties} of additive principal numbers in $\Tt\cap\Om_1$ which are also crucial in the analysis of $\leo$ and $\letwo$. 
It will become clear in lemma \ref{cofinlem} what is meant by cofinality properties of ordinals.

\begin{defi}[\boldmath$\Ttrestral$, $\iotal$, $\latal$\unboldmath]\label{lataldefi}
Let $\al=\thtt(\De+\eta)\in\Tt$. 
\begin{itemize}
\item The restriction of $\Tt$ to $\al$ is defined by
$\Ttrestral:=\set{\be\in\Tt}{\be^\start<\al}$.
\item If $\De>0$ we define $\iotal:\Ttrestral\to\Tal$ by recursion on the definition of $\Tt$ 
by transforming every $\Tt$-term for an ordinal $\xi<\al$ into the parameter $\xi$, defining $\iotal(\xi)$ for 
any non-principal $\xi>\al$ homomorphically, and substituting every $\thtke$-function by $\thtk$ (where $\tht_0=\thtal$).
\footnote{Note that $\Ttrestral$ does not contain any $\thtt$-terms which are greater than or equal to $\al$.}
\item If $\al>\tau$ we define
\[\latal:=\left\{\begin{array}{ll}
      \iotal(\De)+\zetal & 
        \mbox{{\small if $\al\in\Ez$}} \\
      \zetal & \mbox{{\small otherwise.}}
      \end{array}\right.\]
\end{itemize}
\end{defi}

As can be seen from the above definition, the operator $\la$ measures not only the fixed point level,
but also the degree of limit point thinning. Easy examples are $\la_{\epsn}=\epsn$, $\la_{\eps_1}=\eps_1$, while
$\la_{\eps_\om}=\eps_\om+1$ and $\la_{\eps_{\om^2}}=\eps_{\om^2}+2$, where we have omitted the superscript $\tau=1$ and used
the enumeration function $\al\mapsto\eps_\al$ of the epsilon numbers.

For completeness we summarize basic properties of the $\iota$-operator that also justify proofs by induction on $\htt$ and $\htal$,
respectively, as given in Definition \ref{htdef}, and provide an important upper bound.

\begin{lem}[7.2 and 7.3 of \cite{W07a}]\label{iotalem} 
Let $\al=\thtt(\De+\eta)\in\Tt$ where $\De>0$.
\begin{enumerate}
\item[a)] $\iotal$ is a $(<,+)$-isomorphism of $\Ttrestral$ and $\Tal$, and we have $\iotal(\De)<\al^+$.
\item[b)] Suppose $\xi\in\Ttrestral$, $\xi>\al$. Then $\htal(\iotal(\xi))<\htal(\xi)$. 
\item[c)] $\xi\in\Tal\;\imp\;\htt(\iotalinv(\xi))\le\max\singleton{\htt(\al), \htal(\xi)+1}$.
\end{enumerate}
\end{lem}

The interplay between base transformation and the operator $\iotal$ is described as follows. 

\begin{lem}[7.8 of \cite{W07a}]\label{iopiinterplaylem}
Let $\ga,\al\in\Tt\cap\Om_1$ be epsilon numbers such that $\tau<\ga<\al$. We have
\[\iotga=\pigaal\circ\iotal\restriction_\Ttrestrga,\] i.e.\ the following diagram is commutative:
\[\begin{diagram}
\node{\Ttrestrga}\arrow[2]{e,t}{\iotga}\arrow{se,b}{\iotal}\node[2]{\Tga}\\
\node{}\node{\Talga}\arrow{ne,b}{\pigaal}
\end{diagram}\]
In particular, the image of $\iotal\restriction_\Ttrestrga$ is $\Talga$.
\end{lem}

The notion of term translation was introduced in Section 6 of \cite{W07a}. 
For given epsilon number $\al=\thtt(\De+\eta)\in\Tt$ (where $\De>0$) it is possible to define (relativized elementary recursive) 
\emph{partial translation procedures} ${}^\tal:\Tt\to\Tal$ and ${}^\ttau:\Tal\to\Tt$ between the systems
$\Tt$ and $\Tal$, see Definition 6.2 of \cite{W07a}. According to Lemma 6.3 of \cite{W07a},
${}^\tal$ and ${}^\ttau$ are correct on $\Ttrestralplus$ and $\Talrestralplus$, respectively. 
These two restrictions contain the same ordinals, and $\alplus=\thtt(\De+\eta+1)=\thtal(\De)$.
The details are quite technical and can be taken for granted here as we only use translation implicitly in a straightforward manner.

The translation procedure ${}^\ttau$ allows us to consider the ordinal $\latal$ as $\Tt$-term since $\iotal(\De)<\al^+$. 
We will omit translation superscripts since their application can easily be understood from the context.

\begin{lem}[7.6, 7.7, and 7.10 of \cite{W07a}]\label{latallem} Suppose $\al=\thtt(\De+\eta)\in\Tt$, $\al>\tau$.
\begin{enumerate}
\item[a)] We have $\latal=0$ if and only if $\al\not\in\Lz$. 
\item[b)] $\latal<\al^+$. 
\item[c)] $\htal(\latal)<\htt(\al)$ in case of $\De>0$.
\item[d)] If $\De>0$, i.e.\ $\al\in\Ez$, and $\be=\thtt(\Ga+\rho)\in(\al,\alplus)$. Then we have
\[\latbe=\laalbetal=\laalbe.\]
\item[e)] If $\tau\in\Ez$ and $\si\in\Ez\cap\tau$ such that $\al\in\Tts$ then 
\[\latal\in\Tts \quad\mbox{ and }\quad \pist(\latal)=\laspistal,\]
i.e.\ the following diagram is commutative:
\[\begin{diagram}
\node{\Tts\cap(\tau,\tauinf)\cap\Hz}\arrow[2]{e,t}{\latau}\arrow{s,l}{\pist}\node[2]{\Tts}\arrow{s,r}{\pist}\\
\node{\Ts\cap(\si,\siinf)\cap\Hz}\arrow[2]{e,b}{\lasi}\node[2]{\Ts}
\end{diagram}\]
\end{enumerate}
\end{lem}

The characterization of the (intuitive notion of) cofinality properties of an ordinal can now be stated precisely.
The next two lemmas settle that the cofinality properties of an additive principal number $\al\in\Tt$
are exactly classified by $\latal$.
First of all, note that immediately by definition $\latal=0$ if and only if $\al$
is not a limit of additive principal numbers.

\begin{lem}[8.1 of \cite{W07a}]\label{cofinlem}
Let $\al=\thtt(\De+\eta)\in\Tt$, $\al>\tau$, be given. Then for every $\la<\latal$ we have
\[\al=\left\{\begin{array}{ll}
      \sup\set{\ga\in(\tau,\al)\cap\Ez}{\pigaalinv(\latga)\ge\la}&
        \mbox{{\small if $\al\in\Lim(\Ez)$}} \\[2ex]
      \sup\set{\ga\in(\tau,\al)\cap\Hz}{\latga\ge\la} & \mbox{{\small otherwise.}}
      \end{array}\right.\]
\end{lem}

Fine-localization commutes with base transformation and in the sense of the following lemma
with translation. Part 1 of the following lemma, which slightly generalizes Lemma 8.2 of \cite{W07a},
gives a characterization of $\albar$ in terms of cofinality properties of $\al$.
The lemma also justifies that the operator $\bardot$ does not carry a superscript $\tau$: 
For $\be\in(\al,\alinf)\subseteq(\tau,\tauinf)$, $\al\in\Tt\cap\Ez^{<\Om}$, the ordinal $\bebar$ is invariant of the term 
representation of $\be$, whether given in $\Tt$ or in $\Tal$. 

\begin{lem}[5.7 of \cite{CWa}] \label{finelocbasicpropslem} Let $\al=\thtt(\De+\eta)\in\Tt$, $\al>\tau$.
\begin{enumerate}
\item We have
\[\albar=\left\{\begin{array}{ll}
\max\left(\set{\ga\in(\tau,\al)\cap\Ez}{\pigaalinv(\latga)\ge\latal}\cup\singleton{\tau}\right) & \mbox{ if } \al\in\Ez\\
\max\left(\set{\ga\in(\tau,\al)\cap\Hz}{\latga\ge\latal}\cup\singleton{\tau}\right) & \mbox{ otherwise.}
\end{array}\right.\]
\item The operator $\bardot$ commutes with base transformation, i.e.\ the following diagram
is commutative:

\[\begin{diagram}
\node{\Tts\cap(\tau,\tauinf)\cap\Hz}\arrow[3]{e,t}
{\bardot \scriptscriptstyle{\:in\:}\Tt}\arrow{s,l}{\pist}\node[3]{\Tts\cap[\tau,\tauinf)\cap\Hz}\arrow{s,r}{\pist}\\
\node{\Ts\cap(\si,\siinf)\cap\Hz}\arrow[3]{e,b}
{\bardot \scriptscriptstyle{\:in\:}\Ts}\node[3]{\Ts\cap[\si,\siinf)\cap\Hz}
\end{diagram}\]

That is to say, for $\al\in\Tts\cap(\tau,\tauinf)\cap\Hz$ we have
\[\albar\in\Tts\quad\mbox{and}\quad\overline{\pist(\al)}=\pist(\albar).\]

\item For $\al\in\Ez$ the operator $\bardot$ commutes with the translation mappings $\tal$ and $\ttau$:

\[\begin{diagram}
\node{\Tt\cap(\al,\al^+)\cap\Hz}\arrow[2]{e,t}{\tal}\arrow{s,r}{\:\bardot \scriptscriptstyle{\:in\:}\Tt}
\node[2]{\Tal\cap(\al,\al^+)\cap\Hz}\arrow[2]{e,t}{\ttau}\arrow{s,r}{\:\bardot \scriptscriptstyle{\:in\:}\Tal}
\node[2]{\Tt\cap(\al,\al^+)\cap\Hz}\arrow{s,r}{\:\bardot \scriptscriptstyle{\:in\:}\Tt}\\
\node{\Tt\cap[\al,\al^+)\cap\Hz}\arrow[2]{e,b}{\tal}
\node[2]{\Tal\cap[\al,\al^+)\cap\Hz}\arrow[2]{e,b}{\ttau}
\node[2]{\Tt\cap[\al,\al^+)\cap\Hz}
\end{diagram}\]

\item If $\al\in\Ez$ with $\tau$-fine-localization
$(\al_0,\ldots,\al_n)$ and $\be\in(\al,\al^+)$ with $\al$-fine-localization
$(\be_0,\ldots,\be_m)$, then
\[\tau=\al_0,\ldots,\al_n=\al=\be_0^\ttau,\ldots,\be_m^\ttau=\be^\ttau\]
is the $\tau$-fine-localization of $\be$.
\end{enumerate}
\end{lem}

For more details regarding term representations with respect to (fine) localization, closure and cofinality properties, and 
invariance with respect to base transformation and translation, see \cite{W07a,CWa}. In particular Lemma 5.8 of both \cite{W07a}
and \cite{CWa} is informative, and Lemma 5.10 of \cite{CWa} characterizes $\om$-exponentiation in $\Tt\cap(\tau,\tauinf)$, 
also with respect to $\tau$-fine-localization and base transformation. Here we only need a simplified version of this latter lemma.

\begin{lem}[cf.\ 5.10 of \cite{CWa}]\label{simplebarlem}
Let $\ze<\al=\om^\ze\in(\tau,\tauinf)$. We have
\[\al=\thtt\left(1+(-\tau+\ze)\minusp \operatorname{e}_\ze\right)\]
where $\operatorname{e}_\ze:=\left\{\begin{array}{ll}
1&\mbox{ if }\ze=\eps+n\mbox{ for some }\eps\in\Ez\mbox{ and }n<\om\\
0&\mbox{ otherwise.}\end{array}\right.$

\medskip
\noindent For any $\si\in\Ez\cap\tau$ such that $\ze\in\Tts$ we have 
\[\pist(\al)=\thts\left(1+(-\si+\pist(\ze))\minusp \operatorname{e}_\ze\right)=\om^{\pist(\ze)}.\]

\medskip
\noindent Suppose $\ze=_\CNF\om^{\ze_1}+\ldots+\om^{\ze_k}$ with $k>1$.
We have \[\zetal=\ze_k \quad\mbox{  and }\quad\albar=\om^\zepr\quad\mbox{ where }\quad\zepr=\om^{\ze_1}+\ldots+\om^{\ze_{k-1}}.\]
\end{lem} 

This shows that for $\al\in\Hz\setminus\Mz$ the operator $\bardot$ cancels the last factor in the multiplicative normal form of $\al$, 
and the degree of limit point thinning of $\al$ is indeed $\ze_k$, which for $\al\in\Hz\setminus\Lz$ is $0$.
In other words, for $\al=_\NF\eta\cdot\xi\in\Hz-\Mz$ we have $\albar=\eta$ and $\zetal=\log(\log(\xi))=\logend(\log(\al))$ 
for any $\tau\in\Ezone$ such that $\al\in(\tau,\tauinf)$.

\section{Generalized Tracking Sequences and Chains}\label{tstcsec}
In this section we establish connections with earlier work on the structure $\Rtwo$, in particular \cite{CWc}, \cite{W17}, and \cite{W18}, 
and generalize notions originally developed there to be applicable to the entire ordinal structure $\Rtwo$. 
In order to make this article more accessible we review quite extensively the relevant notions and results from \cite{CWc} with 
improvements from the book chapter \cite{W17} and to some extent from \cite{W18}. 

The notions of \emph{tracking sequences and chains}, which were introduced in \cite{CWc}, are motivated easily, despite the relatively 
involved technical preparations and definitions required to formulate them. The approach is natural when analyzing well-partial orderings
and was already chosen by Carlson to calculate the structure $\Rone$ in \cite{C99}: one considers enumeration functions of connectivity
components. In the case of $\Rone$, there is only one such function, called $\ka$, the domain of which is $\epsn+1$, since the 
$\epsn$-th component of $\Rone$ was shown in \cite{C99} to be $[\epsn,\infty)$, i.e.\ $\ka_{\epsn}=\epsn<_1\al$ for all $\al>\epsn$, 
which is sometimes written as $\epsn<_1\infty$. 
Within connectivity components one can consider components relative to the root of the component. In $\Rone$ this is 
again very easy: $[\epsn+1,\epsn\cdot2)\cong [0,\epsn)$ as in general $\On\cong [\al+1,\infty)$ for all $\al$.
Considering the nesting of connectivity components now allows us to locate any ordinal in the segment $\epsn\cdot\om$ by a 
finite sequence of $\ka$-indices. This means that we can \emph{track down} an ordinal in terms of nested connectivity components. 
The ordinal $\epsn\cdot\om$ is the least supremum of an infinite $\lo$-chain, and it is therefore not surprising that the
least $\ktwo$-pair of $\Rtwo$ is $\epsn\cdot\om\ktwo\epsn\cdot(\om+1)$, see \cite{CWc} and Theorem 3.8 of \cite{W} for a direct proof. 
Note that therefore, $\Rone$ and $\Rtwo$ agree on the initial segment $\epsn\cdot(\om+1)$.

In $\Ronepl$, as compared to $\Rone$, the situation is dynamized, and the role of nonzero multiples of $\epsn$ is taken by 
$\Image(\upsilon)\setminus\{0\}$ with $\upsilon$ as in Definition \ref{upsilondefi},
as was shown in \cite{W07b}. There, the enumeration function $\ka$ is therefore relativized to $\ka^\al$,
so that $\ka^\al$ enumerates the $\al$-$\leo$-minimal ordinals, ordinals greater than or equal to $\al$ that do not have 
$\lo$-predecessors strictly greater than $\al$. 
The index of the maximum $\al$-$\leo$-minimal component that is $\leo$-connected to $\al$ is defined to be $\la_\al$, 
and thereby the domain of $\ka^\al$ is $\la_\al+1$. As was shown in \cite{W07b}, if $\al\in\Tt\cap\Hz\cap(\tau,\tauinf)$ 
for some $\tau\in\Ezone$,
then $\laal=\latal$, i.e.\ the index $\laal$ of the largest $\al$-relativized $\leo$-connectivity component that is $\leo$-connected to
$\al$ is equal to the cofinality of $\al$ given by the cofinality operator defined in \ref{lataldefi}, namely $\latal$.
In short, $\al\leo\ka^\al_\latal$, and the problem of finding the largest $\be$ such that $\al\leo\be$, called $\lh(\al)$, 
the length of $\al$, is reduced to the calculation of $\lh(\ka^\al_\latal)$.  
For a summary of the results of \cite{W07b} see Subsection 2.4 of \cite{W07c}.

In $\Rtwo$ any additive principal number $\al<\upsilon_\om$ is the maximum of a maximal finite chain 
\begin{equation}\label{principalchain}
\al_0\lo\al_1\ktwo\cdots\ktwo\al_n=\al,
\end{equation} 
for suitable $n<\om$, such that $\al_{i-1}$ is the greatest $\ktwo$-predecessor of $\al_i$ for $i=2,\ldots,n$, and $\al_0$ is
$\leo$-minimal. Note that if $n>0$, $\al_0$ is therefore the least $\lo$-predecessor of $\al_1$, and 
according to Lemma \ref{ktwoinflochainlem} we indeed have $\al_0\not\ktwo\ale$. 
Note further that for every $i<n$, the ordinal $\alie$ is $\ali$-$\letwo$-minimal, i.e.\ $\be\le\ali$ for any $\be$ such that $\be\ktwo\alie$.
The existence of such sequences was shown in \cite{CWc} for $\al\le\oneinf=\ups_1$.
The \emph{tracking sequence} for such $\al$ characterizes the sequence $\al_0,\ldots,\aln$ in terms of nested connectivity components by
providing their indices.
In order to track down additively decomposable ordinals also, the process is iterated, leading to a \emph{tracking chain}: a finite
sequence of tracking sequences. For additively decomposable ordinals the iterated descent along greatest 
$\ktwo$-predecessors does not in general suffice to locate the ordinal in terms of nested $\le_i$-components, $i=1,2$.
Analyzing $\R$-structures in terms of enumerations of connectivity components lays bare their internal structure and reveals their regularity. 
An evaluation function $\ov$, for \emph{ordinal}, is defined in order to recover the actual ordinal from its description by indices. 

\subsection{\boldmath Ordinal operators for $\Rtwo$\unboldmath}\label{ordopsec}

While the $\la$-operator, which was first introduced in \cite{W07a}, could be motivated independently as a cofinality operator, 
as we did in the previous section, the $\mu$-operator, first introduced in \cite{CWc}, is specifically tailored for index calculation in 
$\Rtwo$. 
The ($\tau$-relativized) ordinal operator $\al\mapsto\latal$ returning the index of the largest (relative) $\leo$-connectivity 
component in the context of $\Ronepl$ is recovered in $\Rtwo$ in an analogue way, see part a) of Lemma \ref{rhomulamestimlem}; 
however, the nesting of $\letwo$-connectivity components 
within a $\leo$-component has to be considered as well. Each ordinal $\al$ such that $\al\ktwo\be$ for some $\be$ is the supremum of an 
infinite $\lo$-chain, along which $\letwo$-connectivity components can occur:

\begin{lem}[7.5 of \cite{CWc}]\label{ktwoinflochainlem}
If $\alpha<_2\beta$ then $\alpha$ is the $\mathrm{sup}$ of an infinite $<_1$-chain.
\end{lem}
{\bf Proof.} For any $\rho<\alpha$ we have $\beta\models\exists x\:\forall y>x\:(\rho<x<_1 y)$.
Hence the same holds true in $\alpha$. We obtain $\rho_1<_1\rho_2<_1\rho_3<_1\ldots<_1\alpha$.
\qed
 
The ($\tau$-relativized) ordinal operator $\al\mapsto\mutal$ and an enumeration function $\xi\mapsto\nu^\al_\xi$ where $\xi\le\mutal$
are defined so as to keep track of such infinite $\lo$-chains accommodating $\letwo$-connectivity components. 
$\mutal$ characterizes the order type of the $\lo$-chain leading to 
the greatest newly arising $\ktwo$-component with the refinement that this $\lo$-chain is subject to a notion of relativized 
$\letwo$-minimality. 
This means that $\nu^\al$ does not enumerate ordinals within the newly arising $\ktwo$-components apart from their roots. 
To clarify this latter statement, call such root $\be$, which (apart from a possible additve offset to place it into a surrounding 
component) is equal to $\nu^\al_\xi$ for some $\xi$ and its greatest $\letwo$-successor 
$\ga$. Then the enumeration function $\nu^\al$ omits all ordinals in the interval $(\be,\ga]$ and in case of $\xi<\mutal$ continues with 
the least ordinal $\de$ that is $\leo$-connected to the greatest newly arising $\letwo$-component, which is indexed by $\mutal$.

It should further be noted that any non-trivial $\letwo$-component is itself supremum of an infinite $\lo$-chain, which has the consequence that  
the function $\nu$ also enumerates those ordinals on such $\lo$-(sub-)chains that do not have any $\ktwo$-successor themselves but lead to 
the next non-trivial $\ktwo$-component. This entails that the image of $\mu^\tau$ consists of additive principal numbers.

Easy examples, again in the setting $\tau=1$ omitted as superscript, are $\mu_{\epsn}=\om$, which is the order type of the $\lo$-chain
leading to the ordinal $\epsn\cdot\om$, which is the least ordinal that has a $\ktwo$-successor.
While we still have $\mu_{\eps_\om}=\om$ (but $\lh(\eps_\om)=\eps_\om\cdot(\om+1)+1$), 
the index $\mu_{\phi(2,0)}=\om^2$ leads to the least ordinal that has two $\ktwo$-successors,
namely $\phi(2,0)\cdot\om^2\ktwo\phi(2,0)\cdot(\om^2+i)$, $i=1,2$.
The index $\mu_{\phi(\om,0)}=\om^\om$ governs the chain 
\[\phi(\om,0)\lo\ldots\lo\phi(\om,0)\cdot\om^\om\ktwo\phi(\om,0)\cdot(\om^\om+\om)\lo\phi(\om,0)\cdot(\om^\om+\om)+1=\lh(\phi(\om,0)),\]
where along the infinite $\lo$-chain of order type $\om^\om$ we find ordinals $\al$ such that $\al\ktwo\al+\phi(\om,0)\cdot n$ for every
$n<\om$, nested in a cofinally increasing manner. This provides an example for the following elementary observation, a complete proof of
which is also given in \cite{W} (Lemma 3.3): 

\begin{lem}[7.2 of \cite{CWc}]\label{loalpllem}
\begin{enumerate} 
\item In $\Rone$ we have (see \cite{C99}) \[\al\leo\al+1\quad\aeq\quad\al\in\Lim.\]
\item In $\Rtwo$ we have 
      \[\al\leo\al+1\quad\aeq\quad\al\in\Lim\;\andsp\;\forall\be(\be\ktwo\al\imp\al=\sup\set{\ga<\al}{\be\letwo\ga}).\]
\end{enumerate}
\end{lem}

The ordinal $\al:=\phi(\epsn,0)$ gives rise to an infinite $\lo$-chain leading to its largest $\letwo$-component at 
$\be:=\nual_{\mu_\al}=\al\cdot\epsn$, where $\mu_\al=\epsn$ and superscript $\tau=1$ is suppressed. This $\letwo$-component is
\[\be\ktwo\be+\be=:\ga\lo\ga+\epsn\cdot\om=:\de\ktwo\de+\epsn\] and contains an infinite $\lo$-chain below $\de$. 
The $\letwo$-component $\de\ktwo\de+\epsn$ is not new, as the interval $(\ga,\de+\epsn]$ is isomorphic to the initial segment 
$\epsn\cdot(\om+1)+1$.

Values of the $\mu$-operator look very canonical so far, but a subtlety arises for the first time when we consider the prominent ordinal
$\Ga_0=\min\{\al\mid\al=\phi(\al,0)\}=\thtnod(\thte(\thte(0)))$. While $\phi(\om,0)\cdot\om^\om$ was the least example for an ordinal $\al$
such that there exist ordinals $\be$ and $\ga$ satisfying $\al\ktwo\be\lo\ga$, the ordinal $\Ga_0$ is the least $\lo$-predecessor of the least
ordinal $\al$ such there exist $\be,\ga$ such that $\al\ktwo\be\lo\ga$ \emph{and} $\al\ktwo\ga$. 
We find $\al=\Ga_0^2\cdot\om$, $\be=\Ga_0^2\cdot(\om+1)$, and $\ga=\Ga_0^2\cdot(\om+1)+\Ga_0$.
The \emph{index} $\mu_{\Ga_0}$ leading to $\al$ is $\Ga_0\cdot\om$. $\al=\nu_{\mu_{\Ga_0}}$ is the supremum of the first infinite chain in 
$\Rtwo$ of a form \[\al_1\ktwo\be_1\lo\al_2\ktwo\be_2\lo\al_3\ktwo\be_3\lo\ldots\] 
of alternating $\lo$- and $\ktwo$-connections. Here $\al_1=\nu_{\Ga_0}=\Ga_0^2$ is the root of a $\letwo$-component that contains an element
($\be_1$) apart from the root that is $\lo$-connected to the greater $\letwo$-components  rooted in the ordinals $\al_2,\al_3,\ldots,\al$. 
The ordinals  $\al_1,\al_2,\al_3,\ldots$ are the least witnesses of such a phenomenon in $\Rtwo$. It turns out that there is a simple criterion 
for indices of the $\nu$-function of whether this phenomenon occurs or not. If we call the set of ordinals $\{\ga\mid\Ga_0\leo\ga\leo\al\}$
the \emph{main line} (here for the $\leo$-connectivity component rooted in $\Ga_0$), then we can say that the $\letwo$-components rooting 
in the $\al_i$ \emph{fall back onto the main line}, namely at the ordinals $\be_i$. 
This way of expressing this characteristic phenomenon in general in $\Rtwo$ has been used sometimes 
in earlier work. \cite{CWc} starts with the so-called \emph{indicator function} $\chi$, indicating the occurrence of the just described phenomenon. For the reader's convenience we review definition and key properties of $\chi$.

\begin{defi}[3.1 of \cite{CWc}]\label{indicatorchi}
For $\tau\in\Ez$ the {\bf\boldmath  indicator function\unboldmath} $\chit:\Tt\to\singleton{0,1}$\index{$\xit$@$\chit$} is defined by
\begin{itemize}
\item $\chit(\xi):=0$ for parameters $\xi<\tau$
\item $\chit(\tau):=1$
\item $\chit(\eta+\xi):=\chit(\xi)$ if $\eta+\xi>\tau$ is in normal form
\item Let $i<\om$ and $\xi=\De+\eta\in\dom(\thti)$ where $\eta<\Omega_{i+1}\mid\De$ with $\xi>0$ in case of $i=0$.
      \begin{itemize}
      \item $\chit(\thti(\xi)):=\chit(\De)$ if $\eta=\sup_{\si<\eta}\thti(\De+\si)$ or $\logend(\eta)=0$
      \item $\chit(\thti(\xi)):=\chit(\xi)$ otherwise.
      \end{itemize}
\end{itemize}
Let $\chitcheck:\Tt\to\singleton{0,1}$ be the dual indicator function, i.e.\ $\chitcheck:=1-\chit$.\index{$\xitc$@$\chitcheck$}
\end{defi}

\begin{lem}[3.2 of \cite{CWc}]\label{chibasetrafolem}
Let $\si,\tau\in\Ez$, $\si<\tau$, and
$\al\in\Tts$. Then $\chis(\pist(\al))=\chit(\al)$,
i.e.\ the following diagram is commutative:
\[\begin{diagram}
\node{\Tts}\arrow[2]{e,t}{\chit}\arrow{se,b}{\pist}\node[2]{\singleton{0,1}}\\
\node{}\node{\Ts}\arrow{ne,b}{\chis}
\end{diagram}\]
The analogue statement holds for $\chitcheck$.
\end{lem}

\begin{lem}[3.3 of \cite{CWc}]\label{chiinvlem}
Let $\tau\in\Ez$ and $\al=\thtt(\De+\eta)>\tau$. 
\begin{enumerate}
\item[a)] $\chit(\al)$ is equal to each of the following: $\chit(\be+\al)$ for all $\be<\tauinf$, $\chit(\logend(\al))$, $\chit(\om^\al)$, $\chit(\be\cdot\al)$ for all $\be\in(0,\tauinf)$,
$\chit\left((1/\be)\cdot\al\right)$ for all $\be\in\Hz^{<\al}$, and $\chit(\latal)$.
\item[b)] If $\al\in\Ez$ then for all $\xi\in\Ttrestral$ such that $\chial(\iotal(\xi))=1$ we have $\chit(\xi)=0$. 
\end{enumerate}
\end{lem}

In order to complete our list of instructive examples of values of the ordinal operator $\mu$, consider the Bachmann-Howard ordinal
$\al:=\thtnod(\thte(\tht_2(0)))$, which is the least $\lo$-predecessor of the least $\ktwo$-chain of the form $\be\ktwo\ga\ktwo\de$ in 
$\Rtwo$ (in fact, any $\Rn$ for $n>1$). $\al$ has fixed point level $\epsom=\thte(\tht_2(0))$, and setting $\tau:=1$ we have 
\[\mutal=\iotal(\epsom)=\epsal,\quad \be=\nu^\al_{\epsal}=\epsal, \quad\ga=\epsal\cdot\om,
  \quad\mbox{ and }\quad\de=\epsal\cdot(\om+1).\]  

With this preparation we can now review the precise definition of $\mu$:

\begin{defi}[3.4 of \cite{CWc}]\label{muindex}
Let $\tau\in\Ezone$ and $\al\in(\tau,\tauinf)\cap\Ez$, say
$\al=\thtt(\De+\eta)$ where $\De=\Omega_1\cdot(\la+k)$ such that $\la\in\Limnod$ and $k<\om$.
We define
\[\mutal:=\om^{\iotal(\la)+\chial(\iotal(\la))+k}.\]\index{$\mutal$}
\end{defi}

The next lemma justifies inductive proofs along $\htt$. The more refined estimation is useful when
dealing with localizations. The subsequent algebraic lemmas concerning the notions of translation and base transformation 
will later be used without explicit mention.

\begin{lem}[3.5 of \cite{CWc}]\label{muestimlem} 
$\htal(\mutal)\le\htal(\latal)<\htt(\al)$ and $\mutal,(\mutal)^+<\alplus$.
\end{lem}

\begin{lem}[3.6 of \cite{CWc}]\label{mutranslem} 
Let $\al=\thtt(\De+\eta)\in\Ez^{>\tau}$. For every $\be=\thtt(\Ga+\rho)\in(\al,\alplus)\cap\Ez$ we have
\[\mutbe=\mualbetal=\mualbe.\]
\end{lem}

\begin{lem}[3.7 of \cite{CWc}]\label{mubasetrafolem}
Let $\si,\tau\in\Ez$, $\si<\tau$, and
$\al=\thtt(\De+\eta)\in\Tts\cap(\tau,\tauinf)\cap\Ez$. Then 
\[\mutal\in\Tts\quad\mbox{ and }\quad\pist(\mutal)=\muspistal,\]
i.e.\ the following diagram is commutative:
\[\begin{diagram}
\node{\Tts\cap(\tau,\tauinf)\cap \Ez}\arrow[2]{e,t}{\mut}\arrow{s,l}{\pist}\node[2]{\Tts}\arrow{s,r}{\pist}\\
\node{\Ts\cap(\si,\siinf)\cap\Ez}\arrow[2]{e,b}{\mus}\node[2]{\Ts}
\end{diagram}\]
\end{lem}

\begin{lem}[3.8 of \cite{CWc}]\label{munoncofinlem}
Let $\tau\in\Ezone$ and $\al\in\Ez\cap(\tau,\tauinf)$ and $\ga\in\Ez\cap(\albar,\al)$. Then we have 
\[\pigaalinv(\mu^\tau_\ga)\le\mutal.\]
\end{lem}

With $\tau$ and $\al$ as in the above definition of $\mutal$, the ordinal operator $\rhoargs{\al}{\xi}$, 
where $\xi\le\mutal$, denotes the index of that
$\nu^\al_\xi$-relativized $\leo$-connectivity component which contains the largest $\letwo$-successor of $\nu^\al_\xi$, 
where $\nu^\al_\xi$ is the $\xi$-th such newly arising $\letwo$-component in the $\al$-th component, assuming
that there is no surrounding component causing an additive offset. This latter phenomenon is taken care of by the notion of 
tracking chain, i.e.\ nesting of tracking sequences, and will be discussed later. 

\begin{defi}[3.9 of\cite{CWc}]\label{rhoindex}
Let $\tau\in\Ez$ and $\al<\tauinf$ where $\logend(\al)=\la+k$ such that $\la\in\Limnod$ and $k<\om$.
We define
\[\rhotal:=\tau\cdot(\la+k\minusp\chit(\la)).\]\index{$\rhotal$}
\end{defi}

\begin{lem}[3.10 of\cite{CWc}] $\rhotal\le\tau\cdot\logend(\al)$ and $\htt(\rhotal)\le\max\singleton{1,\htt(\al)}$.
\end{lem}

\begin{lem}[3.11 of\cite{CWc}]\label{rhobasetrafolem} 
Let $\si,\tau\in\Ez$, $\si<\tau$, and  $\al\in\Tts\cap\tauinf$. Then we have 
\[\rhotal\in\Tts\quad\mbox{ and }\quad\pist(\rhotal)=\rhoargs{\si}{\pist(\al)},\]
i.e.\ the following diagram is commutative:
\[\begin{diagram}
\node{\Tts\cap\tauinf}\arrow[2]{e,t}{\rhot}\arrow{s,l}{\pist}\node[2]{\Tts}\arrow{s,r}{\pist}\\
\node{\Ts\cap\siinf}\arrow[2]{e,b}{\rhosi}\node[2]{\Ts}
\end{diagram}\]
\end{lem}

The lemma below shows the interrelations between the operators from the previous section (\cite{W07a}) and the new ones (\cite{CWc}).
Note in particular part a), where the new index operators $\mu$ and $\varrho$ fall into place with $\iota$.

\begin{lem}[3.12 of\cite{CWc}, corrected in \cite{W17}]\label{rhomulamestimlem}
Let $\tau\in\Ezone$ and $\al=\thtt(\De+\eta)\in(\tau,\tauinf)\cap\Ez$. Then we have
\begin{enumerate}
\item[a)] $\iotal(\De)=\rhoalmutal$ and hence $\latal=\rhoalmutal+\zetal$.
\item[b)] $\rhoalbe\le\latal$ for every $\be\le\mutal$. For $\be<\mutal$ such that 
\footnote{This condition is missing in \cite{CWc}. However, that inequality was only applied under this condition, cf.\ Def.\ 5.1 and L.\ 5.7 of \cite{CWc}.}  
$\chial(\be)=0$ we even have $\rhoalbe+\al\le\latal$.
\item[c)] If $\mutal<\al$ we have $\al\le\latal<\al^2$, while otherwise 
\[\max\left(\Ez^{\le\mutal}\right)=\max\left(\Ez^{\le\latal}\right).\]
\item[d)] If $\latal\in\Ez^{>\al}$, we have $\mutal=\latal\cdot\om$ in case of $\chial(\latal)=1$, and $\mutal=\latal$ otherwise.
\end{enumerate}
\end{lem}

Note that whenever $\mutal\in\Ez^{>\al}$, we have $\mutal=\iotal(\De)=\mc(\latal)$, where $\mc$ denotes the largest additive component
as introduced in the previous section.
An easy example already mentioned earlier might be helpful to illustrate part a): The ordinal $\eps_\om$ is the $\eps_\om$-th $\leo$-minimal
ordinal in $\Rtwo$, $\ka_{\eps_\om}=\eps_\om$,
$\mu_{\eps_\om}=\om$, $\nu^{\eps_\om}_\om=\eps_\om\cdot\om$, 
$\varrho^{\eps_\om}_\om=\eps_\om$, $\ze_{\eps_\om}=1$, and $\la_{\eps_\om}=\eps_\om+1$. 

\subsubsection*{Extending the domain of ordinal operators}

Note that for $\tau\in\Ezone$ we have $\tauinf\in\Image(\ups)$ and $(\tau,\tauinf)\cap\Image(\ups)=\emptyset$. 
Recall Definitions \ref{lataldefi}, \ref{muindex}, \ref{rhoindex}, and \ref{indicatorchi} of the ordinal operators $\la$, $\mu$, 
$\varrho$, and $\chi$, respectively. 
The (partial) extension of their domain to $\Image(\ups)$ is motivated by Theorem \ref{maxchaintheo}.
\begin{defi}\label{lamuextdefi}

For $\tau\in\Ezone$ we extend the definitions of the ordinal operators $\la$, $\mu$, $\varrho$, and $\chi$ as follows. 
\[\mu^\tau_{\tauinf}:=(\tauinf)^\infty=:\la^\tau_{\tauinf}.\]
\[\varrho^\tau_{\tauinf}:=\tauinf\quad\mbox{ and }\quad\chit(\tauinf):=0.\]
For $\la\in\Lim$ we set
\[\la_{\ups_\la}:=\ups_{\la+1}.\]
\end{defi}

We thus obtain \[\mu_{\ups_{\la+k}}=\ups_{\la+k+1}=\la_{\ups_{\la+k}}\] for $\la\in\Limnod$ and $k\in(0,\om)$. 
Note that the expressions $\la_0$ and $\mu_{\ups_\la}$ for $\la\in\Limnod$ remain undefined.

\subsection{\boldmath Tracking sequences relative to limit $\ups$-segments\unboldmath}
\noindent
As mentioned earlier, by definition the $\ka$-function enumerates $\leo$-minimal ordinals in $\Rone$ as well as in $\Rtwo$, 
and the $\nu$-functions essentially enumerate $\letwo$-components relative to the component they occur in.
Precise definitions will be given later.
Returning to the chain (\ref{principalchain}) of ordinals below any additive principal number (less than $\ups_\om$), $\al_0$ is
therefore an element in the image of $\ka$, while $\ali$ for $i=1,\ldots,n$ is an element of the image of that $\nu$-function
which relates to the context given by $\alnod,\ldots,\alimin$. This context dependence is indicated by a superscript
$\nu^\alvec$ where the vector $\alvec$ stands for the sequence of \emph{indices} of the ordinals $\alnod,\ldots,\alimin$. 
For arbitrary ordinals we will also need context dependent enumeration functions $\ka^\alvec$ of relativized $\leo$-components,
so that $\ka^{()}$ becomes another notation for simply $\ka$.

So, given an ordinal $\al$ in $\Rtwo$, we want to calculate the indices of the nested $\le_i$-components it occurs in. For
additive principal $\al$ this will be a sequence, called a \emph{tracking sequence}, or $\trs(\al)$, and for arbitrary $\al$ it 
will be a \emph{tracking chain}, or $\tc(\al)$, a sequence of tracking sequences, where the first element of each tracking sequence
is a $\ka$-index (relativized from the second tracking sequence in the chain on), and the later indices in each sequence are $\nu$-indices.

It turns out that tracking sequences can easily be characterized using the $\mu$-operator, and we denote the independently defined
set (or class) of sequences by $\TS$, in relativized form by $\TSt$.  
We cite the following definition of $\TSt$ from \cite{CWc}, thereby correcting a flaw
in the original formulation that caused a deviation from the intended meaning. 

\begin{defi}[corrected 4.2 of \cite{CWc}]\label{TSdefi}
Let $\tau\in\Ezone$. A nonempty sequence $(\ale,\ldots,\aln)$ of ordinals below $\tauinf$ is called a 
{\bf  \boldmath$\tau$-tracking sequence\unboldmath}\index{$\tau$-tracking sequence} if 
\begin{enumerate}
\item $(\ale,\ldots,\alnmin)$ is either empty or $\ale,\ldots,\alnmin\in\Ez$ and is such that $\tau<\ale<\ldots<\alnmin$.
\item $\aln\in\Hz$ and is such that $\aln\ge\tau$ if $n=1$ and $\aln>1$ if $n>1$.
\item $\alie\le\mutali$ for every $i\in\singleton{1,\ldots,n-1}$.
\end{enumerate} 
By $\TSt$\index{$\tst$@$\TSt$} we denote the set of all  $\tau$-tracking sequences. For convenience we set $\TS^0:=\TSe$.
\end{defi}

According to Lemma \ref{muestimlem} the length of a tracking sequence is bounded in terms of the largest index of $\tht$-functions in the term representation of the first element of the sequence.

\begin{defi}[cf.\ 4.3 of \cite{CWc}]\label{RSdefi} The set of sequences obtained from $\TSt$ by erasing the last entry in each sequence, 
is denoted by $\RSt$ and called the set of 
{\bf \boldmath$\tau$-reference sequences\unboldmath}. 
\end{defi}

Note that ($\tau$-)reference sequences are, apart from the empty sequence, sequences $(\ale,\ldots,\aln)$ of strictly increasing epsilon 
numbers greater than $\tau$ that are \emph{$\mu$-covered}, that is $\al_{i+1}\le\mu^{\al_{i-1}}_{\ali}$ for $i=1,\ldots,n-1$,
setting $\alnod:=\tau$. Note that $\mu^{\al_{i-1}}_{\ali}=\mu^\tau_\ali$ via translation, see Lemma \ref{mutranslem}.
We therefore sometimes omit the superscript when a suitable relativization parameter can be understood from the context.
$\mu$-coverings were first explicitly discussed in Subsection 3.1 of \cite{W17}; we will return to this notion shortly.
Reference sequences $\alvec$ characterize the contexts of relativization needed for the definition of $\ka^\alvec$ and $\nu^\alvec$,
where $\nu$ requires $\alvec$ not to be the empty sequence.

We now modify the definitions of $\TS$, which was defined to be just $\TSe$ in ealier work, and $\RS$, which used to be defined as $\RS^1$, 
to apply to all of $\On$. Since $\ups_\om$ is the supremum of the first infinite $\ktwo$-chain
in $\Rtwo$, an additional parameter $\la\in\Limnod$ comes into the picture as we now need to relate to the interval $[\ups_\la,\ups_{\la+\om})$
that we want to consider tracking sequences in. The initial sequence $(\ups_{\la+1},\ldots,\ups_{\la+m})$ for tracking
sequences leading into the interval $[\ups_{\la+m},\ups_{\la+m+1})$, which might appear redundant at first sight, is needed if $m>0$, 
since $\ups_{\la+1}$ is the $\ka$-index (relative to $\ups_\la$ if $\la\in\Lim$) that specifies the $\leo$-component 
(relative to $\ups_\la$ if $\la\in\Lim$) in which nested $\letwo$-components are specified by the tracking sequence, 
cf.\ Theorem \ref{maxchaintheo}. 

\begin{defi}[\boldmath$\laRS$ and $\laTS$\unboldmath]\label{laRSlaTSdefi}
For $\la\in\Limnod$ we define
\begin{enumerate}
\item $\alvec\in\laRS$ if and only if $\alvec$ is of a form\footnote{For the sake of notational simplicity, we deviate from the usual 
convention that a vector $\alvec$ has components $\ale,\ldots,\aln$ where $n$ is the length of $\alvec$, whenever an explicit redefinition 
is given.} 
$(\ups_{\la+1},\ldots,\ups_{\la+m})^\frown(\ale,\ldots,\aln)$ where $m<\om$ 
and either $n=0$ or  $n>0$ and $\ale,\ldots,\aln\in\Ez$ such that $\ale<\ldots<\aln$, where
$\ale\in(\ups_{\la+m},\ups_{\la+m+1})$, and $\alie\le\mu_{\ali}$ for $i=1,\ldots,n-1$.
\item $\alvec^\frown\be\in\laTS$ if and only if either $\alvec=()$ and $\be\in\Hz\cap[\ups_\la,\ups_{\la+1}]$ or 
$()\not=\alvec=(\al_1,\ldots,\aln)\in\laRS$ and $\be\in\Hz\cap(1,\mu_{\aln}]$.
\end{enumerate}
We define the extended class $\RS$ (respectively $\TS$) as the union of all $\laRS$ ($\laTS$, respectively).
\end{defi}

Note that the above definition accords with the already defined $\RS$ and $\TS$: 
$\RS$ comprises the sequences in $0$-$\RS$ with elements below $\ups_1=\oneinf$ and $\TS$ comprises the sequences below 
$\ups_1$ in $0$-$\TS$. For $\la\in\Lim$ we have $(\ups_\la)\in\laTS\setminus\laRS$.
While the parameter $\la$ in $\laRS$ and $\laTS$ explicitly mentions the interval $[\ups_\la,\ups_{\la+\om})$ in which the 
sequences reside, $\la$ can always be recovered from the first element in the sequence. 
The empty sequence $()\in\laRS$ clearly does not cause any problem.
 	  
Tracking sequences for additive principal numbers were first introduced in Definition 3.13 of \cite{CWc}. Here we first state the assignment
of the proper tracking sequence of a multiplicative principal number, which is based on the notion of localization reviewed in Subsection 
\ref{localizationsec}.

\begin{defi}[cf.\ 3.5 of \cite{W17}]\label{trsofmzdefi}
Let $\tau\in\Ezone$ and $\al\in\Mz\cap(\tau,\tauinf)$ with $\tau$-localization $\tau=\al_0,\ldots,\al_n=\al$.
The {\bf \boldmath tracking sequence of $\al$\unboldmath} above $\tau$\index{tracking sequence}, 
$\trst(\al)$\index{$\trstmz$}, is defined as follows.
If there exists the largest index $i\in\{1,\ldots,n-1\}$ such that $\al\le\mutali$, then 
\[\trst(\al):=\trst(\ali)^\frown(\al),\]
otherwise $\trst(\al):=(\al)$. We extend the definition of $\trst$ to \[\trst(\tauinf):=(\tauinf).\]
\end{defi}

We will apply $\trst$ to ordinals in $\Mz$ in particular when defining the evaluation function $\ov$. 
The extension of the assignment $\trst$ to all additive principal numbers in $(\tau,\tauinf)$ is somewhat more involved.

\begin{defi}[3.15 of \cite{W17}]\label{trsofhzdefi}
Let $\tau\in\Ezone$ and $\al\in[\tau,\tauinf)\cap\Hz$.
The tracking sequence of $\al$ above $\tau$\index{tracking sequence}, $\trst(\al)$\index{$\trsthz$}, is defined as in Definition \ref{trsofmzdefi} 
if $\al\in\Mz^{>\tau}$, and otherwise recursively in the multiplicative decomposition of $\al$ as follows.
\begin{enumerate}
\item If $\al\le\tau^\om$ then $\trst(\al):=(\al)$.
\item Otherwise. Then $\albar\in[\tau,\al)$ and $\al=_\mNF\albar\cdot\be$ for some $\be\in\Mz^{>1}$.
      Let $\trst(\albar)=(\ale,\ldots,\aln)$ and set $\al_0:=\tau$.\footnote{As verified in parts 1 and 2 of 
      Lemma \ref{trsbasicpropslem} we have $\be\le\al_n$.}
      \begin{enumerate}
      \item[2.1.] If $\aln\in\Ez^{>\al_{n-1}}$ and $\be\le\mutaln$ then $\trst(\al):=(\ale,\ldots,\aln,\be)$.
      \item[2.2.] Otherwise. For $i\in\singleton{1,\ldots,n}$ let $(\be^i_1,\ldots,\be^i_{m_i})$ be $\trs^{\ali}(\be)$
            provided $\be>\ali$, and set $m_i:=1$, $\be^i_1:=\ali\cdot\be$ if $\be\le\ali$.
            We first set the insertion index
            \[i_0:=\max\left(\singleton{1}\cup\set{j\in\singleton{2,\ldots,n}}{\be^j_1\le\mu^\tau_{\al_{j-1}}}\right),\]
            then define $\trst(\al):=(\ale,\ldots,\al_{i_0-1},\be^{i_0}_1,\ldots,\be^{i_0}_{m_{i_0}})$.
      \end{enumerate}
\end{enumerate}            
For technical convenience we set $\trs^0:=\trs^1$, and instead of $\trs^1$ we also simply write $\trs$.  \index{$\trst$!$\trs$}
\end{defi}

The following lemma shows that the above definition is sound, that the image of $\trst$ is contained in $\TSt$, 
and states basic properties of $\trst$. Its proof proceeds by straightforward induction along the definition of $\trst$, i.e.\ the
length of $\tau$-localization of multiplicative principal numbers and the number of factors in the multiplicative normal form of
additive principal numbers.

\begin{lem}[3.14 of \cite{CWc}]\label{trsbasicpropslem}
Let $\tau\in\Ezone$ and $\al\in[\tau,\tauinf)\cap\Hz$. Let further $(\ale,\ldots,\aln)$ be
$\trst(\al)$, the tracking sequence of $\al$ above $\tau$.
\begin{enumerate}
\item If $\al\in\Mz$ then $\aln=\al$ and $\trst(\ali)=(\ale,\ldots,\ali)$ for $i=1,\ldots,n$.
\item If $\al=_\mNF\eta\cdot\xi\not\in\Mz$ then $\aln\in\Hz\cap[\xi,\al]$ and $\aln=_\mNF\alnbar\cdot\xi$.
\item $(\ale,\ldots,\al_{n-1})$ is either empty or a strictly increasing sequence of epsilon numbers in the interval $(\tau,\al)$. 
\item For $1\le i\le n-1$ we have $\alie\le\mutali$, and if $\ali<\alie$ then $(\ale,\ldots,\alie)$ is a subsequence of the $\tau$-localization of $\alie$.
\end{enumerate}
We therefore have $\trst(\al)\in\TSt$.
\end{lem}

The following lemma is part of showing (in Theorems 3.19 and 3.20 of \cite{W17}) that $\trs$ is a $(\le,\kglex)$-isomorphism 
between $\ups_1\cap\Hz$ and $\TS$. We will generalize this isomorphism to all of $\Hz$. However, we will need to define the evaluation
function $\ov$ for additive principal numbers first, which will be shown to be the inverse of $\trs$.

\begin{lem}[3.15 of \cite{CWc}]\label{citedinjtrslem}
Let $\tau\in\Ezone$ and $\al,\ga\in[\tau,\tauinf)\cap\Hz$, $\al<\ga$. Then we have 
\[\trst(\al)\klex\trst(\ga).\] 
\end{lem}

Note that the proof of the above lemma given in \cite{CWc} is in fact an induction along the inductive definition of $\trst(\ga)$ with 
a subsidiary induction along the inductive definition of $\trst(\al)$.

The function $\ups$ gives rise to a segmentation of the ordinals into intervals $[\ups_\iota,\ups_{\iota+1})$, which we will use
for a generalization of the notion of tracking sequence. 

\begin{defi}[\boldmath$\upsseg$\unboldmath]\label{upssegdefi}
For $\al\in\On$ let $(\la,m)\in\Limnod\times\om$ be $\klex$-minimal such that $\al<\ups_{\la+m+1}$ and define $\upsseg(\al):=(\la,m)$,
the {\bf \boldmath$\ups$-segment\unboldmath} of $\al$. 
\end{defi}

Note that for $\tau\in\Ezone$ we have $\tauinf=\ups_{\la+m+1}$ where $(\la,m):=\upsseg(\tau)$.
Recall that by $(1/\ga)\cdot\al$ we denote the least ordinal $\de$ such that $\al=\ga\cdot\de$, whenever such an ordinal exists.
Note that in the extension of the function $\trs$ to all of $\Hz$ below, the original definition of $\trs$ is referred to as $\trs^0$.

The technical appearance of the following definition is owed to the fact that ordinal notations cover in each relativized instance 
at most an interval of the form $[\ups_\al,\ups_{\al+1})$. Note that in the case where $\ga<\al<\ga^\om$ we need to explicitly cancel 
the lead factor $\ga$, which is the one but last element of the tracking sequence, resulting in the final $\nu$-index $(1/\ga)\cdot\al$.

\begin{defi}[Extended domain \boldmath$\trs$\unboldmath]\label{latrsdefi}
For $\al\in\Hz$ let $(\la,m):=\upsseg(\al)$ be the $\ups$-segment of $\al$. 
For better readability we set $\gavec:=(\ups_{\la+1},\ldots,\ups_{\la+m})$ and $\ga:=\ups_{\la+m}$. We define
\[\trs(\al):=\left\{\begin{array}{ll}
                 \trs^{\ups_\la}(\al)&\mbox{if }m=0\mbox{, otherwise:}\\[2mm]
                 \gavec&\mbox{if }\al=\ga,\\[2mm]
                 \gavec^\frown(1/\ga)\cdot\al&\mbox{if }\ga<\al<\ga^\om,\\[2mm]
                 \gavec^\frown\trs^{\ga}(\al)&\mbox{if }\al\ge\ga^\om.
               \end{array}\right.
\]
\end{defi}

Having defined the assignment function $\trs$ of tracking sequences on $\Hz$, we need to calculate the evaluation function $\ov$
on $\TS$. The evaluation function $\ov$ was first given in \cite{CWc},
but redefined in \cite{W17}, which led to a substantial disentanglement that allowed to directly see elementary recursiveness.
For the verification proof that these functions are inverses of each other, we need to establish a definition of $\ov$ 
that provides expressions in multiplicative normal forms accompanied with the tracking sequences of the initial products of such
normal forms. Along the way we will characterize the fixed points of $\ov$, i.e.\ all $\alvec^\frown\be\in\TS$ such that 
$\ov(\alvec^\frown\be)=\be$. Clearly, fixed points are pivotal in the proof of order-isomorphism between $\Hz$ and $\TS$, and
decomposition into multiplicative normal form carries the proof on. More specifically,
given a sequence in $\alvec^\frown\be\in\TS$ with evaluation $\ga=_\NF\eta\cdot\xi\not\in\Mz$, we want to determine the tracking sequence 
of $\eta$. It turns out that $\xi=\lf(\be)$, i.e.\ either $\be$ itself or its last factor. The latter can occur if 
$\be=_\MNF\be_1\cdot\ldots\cdot\be_k$ with $k>1$, $\be_1\in\Ez^{>\aln}$, and $\be_2\le\mu_{\be_1}$, where $\alvec=(\ale,\ldots,\aln)$. 
The tracking sequence of $\eta$ is determined by an auxiliary function $\hop$, which in turn uses two functions $\sk$ 
(for \emph{skimmed sequence}) and $\mts$ (for \emph{minimal tracking sequence}). The former sequence, say $\skbe(\ga)$, extends a given 
point on a local main line (i.e.\ starting from a $\nu$-index $\ga$ and extending along nested $\nu$-indices through the $\letwo$-component
rooted in the start point indexed by $\ga$) as far as possible without obtaining evaluating last factors below 
$\be$. The latter, say $\mtsal(\be)$, produces a minimal $\mu$-covering from $\al$ to $\be$, obtaining large evaluating factors as
quickly as possible. For the reader's convenience we include all relevant (slightly modified) definitions, starting with $\mu$-coverings.
Correcting definitions in \cite{W17} technically, we need to allow $(\al)$ to be a $\mu$-covering of $\al$ itself and let a $\mu$-covering 
from $\al$ to $\be$ start with $\al$.

\begin{defi}[cf.\ 3.2 of \cite{W17}]\label{mucovering}
Let $\tau\in\Ezone$, $\al\in\Ez\cap(\tau,\tauinf)$, and $\be\in\Hz\cap[\al,\alinf)$.
A sequence $(\al_0,\dots,\aln)$ where $\al_0=\al$, $\aln=\be$, such that
$(\al_1,\ldots,\aln)\in\TSal$ and $\al<\al_1\le\mutal$ if $n>0$,
is called a {\bf \boldmath$\mu$-covering\unboldmath} from $\al$ to $\be$. 
\end{defi}

\begin{lem}[3.3 of \cite{W17}]\label{mucovloc}
Any $\mu$-covering from $\al$ to $\be$ is a subsequence of the $\al$-localization of $\be$.
\end{lem}

\begin{defi}[3.4 of \cite{W17}]\label{maxminmucov}
Let $\tau\in\Ezone$.
\begin{enumerate}
\item For $\al\in\Hz\cap(\tau,\tauinf)$ we define $\maxmucovtau(\al)$ to be the longest subsequence $(\ale,\ldots,\alne)$ of the $\tau$-localization of $\al$ 
which satisfies $\tau<\ale$, $\alne=\al$, and which is {\it $\mu$-covered}, i.e.\ which satisfies $\alie\le\mutali$ for $i=1,\ldots,n$.
\item For $\al\in\Ez\cap(\tau,\tauinf)$ and $\be\in\Hz\cap[\al,\alinf)$ we denote the shortest subsequence $(\be_0,\be_1,\ldots,\be_n)$
of the $\al$-localization of $\be$ which is a $\mu$-covering from $\al$ to $\be$ by $\minmucoval(\be)$, if such sequence exists. 
\end{enumerate}
\end{defi}

\begin{defi}[cf.\ 3.6 of \cite{W17}]\label{mts}
Let $\tau\in\Ezone$, $\al\in\Ez\cap(\tau,\tauinf)$, $\be\in\Hz\cap[\al,\alinf)$, and let $\al=\al_0,\ldots,\aln=\be$ be the $\al$-localization 
of $\be$.
If there exists the least index $i\in\{0,\ldots,n-1\}$ such that $\ali<\be\le\mutali$, then 
\[\mtsal(\be):=\mtsal(\ali)^\frown(\be),\]
otherwise $\mtsal(\be):=(\al)$.
\end{defi}

Note that $\mtsal(\be)$ reaches $\be$ if and only if it is a $\mu$-covering from $\al$ to $\be$. We will see a criterion for this to hold
in Lemma \ref{mtshatlem}.

\begin{lem}[3.7 of \cite{W17}]\label{covcharlem} Fix $\tau\in\Ezone$.
\begin{enumerate}
\item For $\al\in\Hz\cap(\tau,\tauinf)$ let $\maxmucovtau(\al)=(\ale,\ldots,\alne)=\alvec$. If $\ale<\al$ then $\alvec$ is a $\mu$-covering from $\ale$ to $\al$
and $\mtsale(\al)\subseteq\alvec$.
\item If $\al\in\Mz\cap(\tau,\tauinf)$ then $\maxmucovtau(\al)=\trst(\al)$.
\item Let $\al\in\Ez\cap(\tau,\tauinf)$ and $\be\in\Hz\cap[\al,\alinf)$. Then $\minmucoval(\al)$ exists if and only if $\mtsal(\be)$ is a $\mu$-covering from $\al$ to $\be$,
in which case these sequences are equal, characterizing the lexicographically maximal $\mu$-covering from $\al$ to $\be$.
\end{enumerate}
\end{lem}

The following ordinal operator provides a useful upper bound when calculating the reach of connectivity components. The subsequent 
Lemma \ref{mtshatlem} justifies the definition further.

\begin{defi}[cf.\ 3.16 of \cite{CWc}]\label{alhatdefi}\index{$\alhat$}
Let $\tau\in\Ezone$ and $\al\in(\tau,\tauinf)\cap\Ez$. We define
\[\alhat:=\min\set{\ga\in\Mz^{>\al}}{\trsal(\ga)=(\ga)\andsp\mutal<\ga}.\] 
For $\al=\ups_\xi$, $\xi>0$, we set \[\alhat:=\ups_{\xi+1}.\]
\end{defi}

Note that in the above context we have $\alhat\le\alplus$ if $\al\not\in\Image(\ups)$. 
As is the case with $\al^+$ we suppress the base $\tau$ in 
the notation $\widehat{\al}$ assuming that it will always be well understood from the respective context.

\begin{lem}[3.17 of \cite{CWc}]\label{widehatestimlem} Let $\tau,\al$ be as in the above definition, $\al\not\in\Image(\ups)$. Then  
\[\widehat{\be}\le\widehat{\al}\quad\mbox{ for any }\quad\be\in\Tal\cap\Ez\cap(\al,\mutal].\] 
We further have $\latal<\widehat{\al}$.
\end{lem}

\begin{lem}[cf.\ 3.8 of \cite{W17}]\label{mtshatlem}
Let $\tau\in\Ezone$, $\al\in\Ez\cap(\tau,\tauinf)$, and $\be\in\Mz\cap[\al,\alinf)$. Then $\mtsal(\be)$ is a $\mu$-covering from $\al$ to 
$\be$ if and only if  $\be<\alhat$. This holds if and only if either $\be=\al$ or for $\trsal(\be)=(\be_1,\ldots,\be_m)$ we have 
$\be_1\le\mutal$.
\end{lem}

\begin{defi}[cf.\ 4.9 of \cite{CWc}] For $\be\in\Mz^{>1}$ and $\ga\in\Ez\setminus\Image(\ups)$ 
let $\sk_\be(\ga)$\index{$\sk$} 
be the maximal sequence $\de_1,\ldots,\de_l$ such that (setting $\de_0:=1$)
\begin{itemize}
\item $\de_1=\ga$ and 
\item if $i\in\singleton{1,\ldots,l-1}\andsp\de_i\in\Ez^{>\de_{i-1}}\andsp\be\le\mu_{\de_i}$, then $\de_{i+1}=\overline{\mu_{\de_i}\cdot\be}$.
\end{itemize}
\end{defi}

Note that equivalently, one obtains the ordinal $\overline{\mu_{\de_i}\cdot\be}$ in the above definition by \emph{skipping} the factors
strictly below $\be$ from the multiplicative normal form of $\mu_{\de_i}$.
Lemma \ref{muestimlem} guarantees that the above definition terminates.
We have $(\de_1,\ldots,\de_{l-1})\in\RS$ and $(\de_1,\ldots,\de_l)\in\TS$. Notice that $\be\le\de_i$ for $i=2,\ldots,l$.
Note that if $\be=\om$, the sequence $\ga=\de_1,\ldots,\de_l$ is maximal such that $\de_{i+1}=\mu_{\de_i}$, $i=1,\ldots,l-1$.

After the above preparations we can now assemble the auxiliary sequence $\hga(\alvec^\frown\be)$ that plays a key role in the definition of 
the evaluation function $\ov$. 

\begin{defi}[cf.\ 3.11 of \cite{W17}]\label{hgamalbe}
Let $\alvec^\frown\ga\in\laRS$, where $\ga\not\in\Image(\ups)$, and $\be\in\Mz\cap(1,\gahat)$.
If $\be>\ga$ let $\mtsga(\be)=\etavec^\frown(\eps,\be)$.
\[\hbe(\alvecga):=\left\{\begin{array}{ll}
\alvec^\frown\etavec^\frown\skbe(\eps)&\mbox{ if }\ga<\be<\gahat\\[2mm]
\alvec^\frown\skbe(\ga)&\mbox{ if }\be\le\ga\mbox{ and }\be\le\muga\\[2mm]
\alvec^\frown\ga&\mbox{ if }\be\le\ga\mbox{ and }\be>\muga.
\end{array}\right.\]
\end{defi}

We are now going to extend Definition 3.14 of \cite{W17} (2.3 of \cite{W18}) to the extended class $\TS$ as the union over all $\laTS$. 
If $\alvec^\frown\be\in\laTS$ is of the form $(\ups_{\la+1},\ldots,\ups_{\la+m})$ where $m<\om$ or $\alvec^\frown\be=(1+\ups_\la)$, 
we will have $\ov(\alvec^\frown\be)=\be$, extending the class of fixed points of $\ov$. As mentioned earlier, the evaluation function
$\ov$ enables a smooth definition of the component enumerating functions $\ka$ and $\nu$ in the next subsection.

\begin{defi}[cf.\ 3.14 of \cite{W17}]\label{odef}
Let $\alvecbe\in\laTS$, where $\alvec=(\ups_{\la+1},\ldots,\ups_{\la+m})^\frown(\ale,\ldots,\aln)\in\laRS$ for the maximal such $m<\om$ and
$\be=_\MNF\be_1\cdot\ldots\cdot\be_k$. Let $\bepr:=1$ if $k=1$ and $\bepr:=\be_2\cdot\ldots\cdot\be_k$ otherwise. 
We set $\alnod:=1+\ups_{\la+m}$, $\alne:=\be$, $h:=\htarg{\alnod}(\ale)+1$, and
$\gavec_i:=\trs^{\alimin}(\ali)$ for $i=1,\ldots,n$, while
\[\gavec_{n+1}:=\left\{\begin{array}{ll}
            (\be)&\mbox{if }\be\le\aln\\[2mm]
            \trsaln(\be_1)^\frown\be_2&\mbox{if } k>1,\be_1\in\Ez^{>\aln}\andsp\be_2\le\mu_{\be_1}\\[2mm]
            \trsaln(\be_1)&\mbox{otherwise,}
            \end{array}\right.\]
and write $\gavec_i=(\ga_{i,1},\ldots,\ga_{i,m_i})$, $i=1,\ldots,n+1$.
We then define
\[\lSeq(\alvecbe):=(m_1,\ldots,m_{n+1})\in[h]^{\le h}\]
where $[h]^{\le h}$ is the set of sequences of natural numbers $\le h$ of length at most $h$, ordered lexicographically. 
We may now define $\ov(\alvecbe)$ recursively in $\lSeq(\alvecbe)$, as well as auxiliary parameters 
$n_0(\alvecbe)$ and $\ga(\alvecbe)$, which are set to $0$ where not defined explicitly.
\begin{enumerate}
\item If $\alvec=()$ and $\be=1+\ups_\la$, then $\ov((\be)):=\be$.
\item If $\alvec\not=()$\footnote{The corresponding condition of part 2 of Definition 2.3 in \cite{W18} should read $n\ge 1$, not $n>1$.} 
      and $\be_1\le\aln$, then $\ov(\alvecbe):=_\NF\ov(\alvec)\cdot\be$.
\item If $\be_1\in\Ez^{>\aln}$, $k>1$, and $\be_2\le\mu_{\be_1}$, then
set $n_0(\alvecbe):=n+1$, $\ga(\alvecbe):=\be_1$, and define
\[\ov(\alvecbe):=_\NF\ov(\hop_{\be_2}(\alvec^\frown\be_1))\cdot\bepr.\] 
\item Otherwise. Then setting
\[n_0:=n_0(\alvecbe):=\max\left(\{i\in\{1,\ldots,n+1\}\mid m_i>1\}\cup\{0\}\right),\]
define 
\[\ov(\alvecbe):=_\NF\left\{\begin{array}{ll}
            \be&\mbox{if } n_0=0\\[2mm]
            \ov(\hop_{\be_1}({\alvec_{\restriction_{n_0-1}}}^\frown\ga))\cdot\be&\mbox{if } n_0>0,
            \end{array}\right.\] 
where $\ga:=\ga(\alvecbe):=\ga_{n_0,m_{n_0}-1}$.
\end{enumerate}
\end{defi}
The just defined extension of the evaluation function $\ov$ to $\TS$ is easily seen to have the following desired properties.

\begin{theo}[cf.\ 3.19, 3.20, and 3.21 of \cite{W17}]\label{trsvallem} 
The (class) function $\trs:\Hz\to\TS$ is a $(<,\klex)$-order isomorphism with inverse $\ov$:
\begin{enumerate}
\item For $\al\in\Hz$ we have \[\ov(\trs(\al))=\al.\]
\item For $\alvec^\frown\be\in\TS$ we have \[\trs(\ov(\alvec^\frown\be))=\alvec^\frown\be.\]
\item $\ov$ is strictly increasing with respect to the lexicographic ordering on $\TS$ and 
continuous in the last vector component. 
\end{enumerate}
\end{theo}

For proofs and further details see Section 3 of \cite{W17}. Note that part 1 of the above theorem is proved by induction along the
definition of $\trs(\al)$, while part 2 is proved by induction on $\lSeq(\alvecbe)$ along the ordering $(\lSeq,\klex)$. 

\subsection{Connectivity components of \boldmath$\Rtwo$\unboldmath}
We can now define the complete system of enumeration functions of (relativized) connectivity components, using the evaluation function $\ov$.
To this end we define $\ka=\kvnod$ on all of $\On$, which will be conceptually justified in Remark \ref{globalkapparmk}, 
and define functions $\kval$ and $\nuval$ for $()\not=\alvec\in\laRS$ where $\la\in\Limnod$, 
extending Definition 4.1 of \cite{W17} (Definition 2.4 of \cite{W18}) to $\laRS$.
We first define these functions on the additive principal numbers.

\begin{defi}[cf.\ 2.4 of \cite{W18}]\label{kappanuprincipals}
Let $\alvec\in\laRS$ where $\la\in\Limnod$, $\alvec=(\ups_{\la+1},\ldots,\ups_{\la+m})^\frown(\ale,\ldots,\aln)$ 
for the maximal such $m<\om$,
and set $\alnod:=1+\ups_{\la+m}$.  
If $\alvec\not=()$, we define for $\be$ such that $\alvec^\frown\be\in\laTS$ \[\nu^\alvec_\be:=\ov(\alvecbe).\]
Now let $\be\in\Hz$ such that $\be\le\laaln$ if that exists (i.e.\ $\la+m>0$ if $n=0$) and $\be\le\ups_1$ otherwise.
Let $\be=_\MNF\be_1\cdot\ldots\cdot\be_k$ and set $\bepr:=(1/\be_1)\cdot\be$. 
We first define an auxiliary sequence $\gavec\in\xiRS$ for some $\xi\in\Limnod$, $\xi\le\la$.\\[2mm]
{\bf Case 1:} $\be\le\aln$. Here we consider two subcases:\\[2mm]
{\bf Subcase 1.1:} $n>0$ and there exists the maximal $i\in[0,\ldots,n-1]$ such that $\ali<\be$. Then we set 
\[\gavec:=(\ups_{\la+1},\ldots,\ups_{\la+m})^\frown(\ale,\ldots,\ali).\]
{\bf Subcase 1.2:} Otherwise. Then we have $\be\le\alnod$. 
Let $(\xi,l)\in\Limnod\times\om$ be $\kglex$-minimal such that $\be\le\ups_{\xi+l+1}$, so that
$(\xi,l)\kglex(\la,m)$, and set 
\[\gavec:=(\ups_{\xi+1},\ldots,\ups_{\xi+l}).\]
{\bf Case 2:} $\be>\aln$. Then $\gavec:=\alvec$.\\[2mm]
Writing $\gavec=(\ga_1,\ldots,\ga_r)$ we now define
\[\kvalbe:=\left\{\begin{array}{ll}
            \ov(\gavec)\cdot\bepr&\mbox{if } r>0,\; \be_1=\ga_r,\mbox{ and }k>1\\[2mm]
            \ov(\gavec^\frown\be)&\mbox{otherwise.}
            \end{array}\right.\]
For arbitrary $\be\in\Hz$ let $(\la,m)\in\Limnod\times\om$ be $\kglex$-minimal such that $\be\le\ups_{\la+m+1}$ and set 
$\alvec:=(\ups_{\la+1},\dots,\ups_{\la+m})$.
Writing $\kappa_\be$ instead of $\kappa^{()}_\be$, we define \[\ka_\be:=\kvalbe,\]
and call $\ka$ the {\bf global \boldmath$\ka$-function of $\Rtwo$\unboldmath}.
\end{defi}

\begin{rmk}\label{globalkapparmk} We still need to define the $\ka$- and $\nu$-functions for arguments that are additively decomposable. For orientation
and clarification we make some statements about the global $\ka$-function that follow from our results in Section \ref{structuresec}.
Recall that an ordinal $\be\ge\al$ is called $\al$-$\le_i$-minimal if $\ga\le\al$ for any $\ga$ such that $\ga<_i\be$. 

\begin{itemize}
\item The restriction of the global $\ka$-function of $\Rtwo$ to the initial segment $\oneinf=\ups_1$ enumerates all 
$\leo$-minimal ordinals, the largest one being $\ups_1$. 
\item For any $\la\in\Lim$ the $\ups_\la$-$\leo$-minimal ordinals are enumerated by $\be\mapsto\ups_\la+\ka_\be$ for $\be\le\ups_{\la+1}$,
and the greatest $\lo$-predecessor of $\ups_{\la+1}$ is $\ups_\la$.
\item And for any $\xi\in\Limnod$, $m<\om$, the $\ups_{\xi+m+1}$-$\leo$-minimal ordinals are enumerated by 
$\be\mapsto\ups_{\xi+m+1}+\ka_\be$ for $\be\le\ups_{\xi+m+1}$. We have $\ups_{\xi+m+1}\cdot 2\lo\infty$. 
\end{itemize}

Note that $\ka$ as defined here acts as the identity on $\Image(\ups)$. The global $\ka$-function can therefore equally be defined using
Definition \ref{upssegdefi}, since for $\be\in\Hz$ and $(\xi,m):=\upsseg(\be)$, we have $\ka_\be=\ka^\alvec_\be$, 
where $\alvec:=(\ups_{\xi+1},\ldots,\ups_{\xi+m})$, as $\ka^\alvec_\be=\ka^\alvecpr_\be$
for $\be=\ups_{\xi+m+1}$ and $\alvecpr=\alvec^\frown\be$.
 
Defining $\ka$ on all of $\On$ despite the fact that we have $\ups_1\lo\infty$,
wherefore the largest $\leo$-connectivity component in $\Rtwo$ is $[\ups_1,\infty)$,
is justified by the uniformity of $\ka$ through all $\ups$-segments in that it simply end-extends continously.
\end{rmk}

In order to extend the above 
definition to non-principal indices $\be$, we need some preparation.
We introduce a term measure on terms that use finitely many parameters from $\Image(\ups)$. 
Suppose $\la\in\Limnod$, $m<\om$, and $\al\in\T^{\ups_{\la+m}}$. This term representation uses finitely many parameters below $\ups_{\la+m}$,
each of which, in turn, can be represented in a system $\T^{\ups_\iota}$ for some $\iota<\la+m$ with parameters below $\ups_\iota$.
Resolving hereditarily (using transfinite recursion) all parameters results in a term representation of $\al$ that makes use of finitely 
many relativized $\thtnod$-functions $\thtt$ where $\tau\in(1+\ups_{\iota_1},\ldots,\ups_{\iota_l})$ for an increasing sequence of indices $0\le\iota_1,\ldots,\iota_l=\la+m$.
A term measure can therefore be defined elementary recursively relative to such resolved term representation for any ordinal $\al$. 
We adapt this motivation to the setting of $\laRS$ as follows. For $\la=m=0$ we obtain a parameter-free
representation for notations below $\ups_1$ as in \cite{W18}. We begin with a useful auxiliary notion. Note that we can in general assume that 
parameters of terms in $\Tt$ can itself be represented as terms of some suitable $\Ts$ where $\si,\tau\in\Ezone$, $\si<\tau$.

\begin{defi}
Let $\tauvec=(\tau_1,\ldots,\tau_{n+1})$ be a strictly increasing sequence of ordinals in $\Ezone$ and $\tau:=\tau_{n+1}$.
We say that $\al\in\Tt$ is given in $\Ttvec$-representation with parameter set $\Par^\tauvec(\al)$ if either 
$n=0$ and $\Par^\tauvec(\al)=\Par^\tau(\al)$ 
or $n>0$ and each parameter term $\be\in\Par^\tau(\al)$ is given in $\Tsivec$-representation where $\sivec=(\tau_1,\ldots,\tau_n)$
and \[\Par^\tauvec(\al)=\bigcup_{\be\in\Par^\tau(\al)}\Par^\sivec(\be).\]
\end{defi}

For clarification, in the case $n=0$, where $\tau=\tau_1$, we trivially have $\T^{(\tau)}=\Tt$ and $\Par^{(\tau)}(\al)=\Par^\tau(\al)$.
For $n>0$ we have $\Par^\tauvec(\al)\subseteq\tau_1$, which contains but in general can be much larger than $\Par^\tau(\al)\cap\tau_1$. 
Note that terms $\al$ in $\Ttvec$-representation as above can have $\tht^{\tau_i}$-subterms for $i=1,\ldots,n+1$, but for 
$i,j\in\{1,\ldots,n+1\}$, $i<j$, while a $\tht^{\tau_j}$-subterms can itself have $\tht^{\tau_i}$-subterms, this cannot happen the other 
way around.

\begin{defi}\label{resolvingseqdefi}
Setting $\alnod:=1+\ups_\la$ where $\la\in\Limnod$, let $\alvec\in\laRS$ be of a form 
$\alvec=(\ale,\ldots,\almn)$ where $m,n<\om$, $\ale=\ups_{\la+1},\ldots,\alm=\ups_{\la+m}$,
and $\alme\in(\ups_{\la+m},\ups_{\la+m+1})$ if $n>0$.
The term system $\laTalvec$ is obtained from $\Talmn$ by successive substitution of parameters from $(\ali,\alie)$ by their $\Tali$-representations,
for $i=m+n-1,\ldots,0$. The parameters $\ali$ are represented by the terms $\thtali(0)$, $0\le i\le m+n$. Remaining unresolved parameters
are below $\alnod$. More formally, we proceed as follows.
\begin{itemize}
\item $\al\in\T^{\almn}$ is a $\T^{(\almn)}$-term. $\Par^{(\almn)}(\al):=\Par^{\almn}(\al)$.
\item Let $\al$ be a $\T^{(\al_{m+n-i},\ldots,\al_{m+n})}$-term, $0\le i\le m+n-1$. Replace those parameters of $\al$ that are in the set 
$\Par^{(\al_{m+n-i},\ldots,\al_{m+n})}(\al)\cap[\al_{m+n-i-1},\al_{m+n-i})$ by $\T^{\al_{m+n-i-1}}$-terms. The resulting representation
of $\al$ is a $\T^{(\al_{m+n-i-1},\ldots,\al_{m+n})}$-term, and the set of parameters $\Par^{(\al_{m+n-i-1},\ldots,\al_{m+n})}(\al)$
is the set of remaining (and new) parameters in the new representation of $\al$ after replacement. All parameters are now
below $\al_{m+n-i-1}$.
\item We also call the $\T^{\alnod^\frown\alvec}$-representation of $\al$ the $\laTalvec$-term representation of $\al$ and the 
corresponding parameter set $\la$-$\Par^\alvec(\al)$. We have $\la$-$\Par^\alvec(\al)\subseteq\alnod$.
\end{itemize}

For $\al\in\Talmn$ in $\laTalvec$-representation, let $\iovec:=(1+\ups_{\iota_1},\ldots,\ups_{\iota_l})$ be the uniquely defined,
(possibly empty) finite increasing sequence below $\ups_\la$ needed to resolve all parameters of $\al$. 
This results in a $\Ttvec$-representation of $\al$, where 
$\tauvec=\iovec^\frown\alnod^\frown\alvec$, that uses relativized $\thtnod$-functions $\thtti$ for $i=1,\ldots,l+m+n+1$.
More formally, we proceed as follows.
\begin{itemize}
\item Suppose that $\al$ is given in $\Tsivec$-representation with parameter set $\Par^\sivec(\al)$ still containing nonzero elements.
Let $(\xi,k)$, where $\xi\in\Limnod$ and $k<\om$, be $\klex$-minimal such that $\Par^\sivec(\al)\subseteq\ups_{\xi+k+1}$, and set 
$\si:=1+\ups_{\xi+k}$.
Replace all parameters in $\al$ that are in the interval $[\si,\siinf)$ by $\Ts$-terms 
(using $\thts$-terms) to obtain $\al$ in $\T^{\si^\frown\sivec}$-representation with parameter set
$\Par^{\si^\frown\sivec}(\al)$ consisting of the remaining (and new) parameters, which are now below $\si$.
\item Starting from the $\laTalvec$-representation of $\al$, i.e.\ starting with $\sivec=\alnod^\frown\alvec$, after finitely many steps 
all nonzero parameters in $\al$ are resolved, and we call the resulting sequence $\tauvec$ the {\bf resolving sequence} for $\al$.
\end{itemize}
\end{defi}

\begin{defi}[cf.\ 4.3 of \cite{W17} (2.6 of \cite{W18})]\label{Ttauvec}
For $\al\in\Talmn$ in $\laTalvec$-representation, where $\la\in\Limnod$, $\alvec=(\ale,\ldots,\almn)\in\laRS$, and $\alnod:=1+\ups_\la$, 
let $\tauvec=(\tau_1,\ldots,\tau_r)$ be a sequence of strictly increasing ordinals in $\Ezone$ containing the resolving sequence 
for $\al$ as defined above.
The {\bf length} $\ltvec(\al)$ of the representation of $\al$ as $\Ttvec$-term $\al$ is defined by induction on the build-up of $\al$ 
as follows.
\begin{enumerate}
\item $\ltvec(0):=0$,
\item $\ltvec(\be):=\ltvec(\ga)+\ltvec(\de)$ if $\be=_\NF\ga+\de$, and 
\item $\ltvec(\tht(\eta)):=\left\{\begin{array}{l@{\quad}l}
1&\mbox{ if }\quad\eta=0,\\
\ltvec(\eta)+4&\mbox{ if }\quad\eta>0,
\end{array}\right.$\\[2mm] 
where $\tht\in\{\tht^{\tau_i}\mid 1\le i\le r\}\cup\{\tht_{i+1}\mid i\in\N\}$.
\end{enumerate}
\end{defi}

Note that in the above defintion the term $\al\in\T^{\al_{m+n}}$ in $\laTalvec$-representation determines uniquely, 
which relativized $\thtt$-terms occur in its resolved term notation. 
If $\tauvec$ is or contains the resolving sequence for $\al$, then for each such $\thtt$-term, $\tau$ is an element of $\tauvec$.
For $\la=m=0$ the above definition is compatible with Definition 4.3 of \cite{W17} and Definition 2.6 of \cite{W18}. In this case the 
sequence $1^\frown\alvec$ is resolving for $\al$, so the $0$-$\T^\alvec$-term representation is completely resolved and directly 
corresponds to the notion of $\Ttvec$-representation in \cite{W17,W18}.

\begin{lem}[cf.\ Subsection 2.1 of \cite{W18}]\label{ltveclem} 
For $\al\in\Talmn$ in $\laTalvec$-representation, where $\la\in\Limnod$, $\alvec=(\ale,\ldots,\almn)\in\laRS$, and 
$\alnod:=1+\ups_\la$, let $\tauvec=(\tau_1,\ldots,\tau_r)$ be a sequence of strictly increasing ordinals in $\Ezone$ containing 
the resolving sequence for $\al$ as defined above.
\begin{enumerate}
\item Setting $\tau:=\almn$, if $\al=\thtt(\De+\eta)\in\Ez$ such that $\al\le\mu_\tau$, then we have
\begin{equation*}\label{iotalen}
\ltvec(\De)=\ltvecal(\iota_{\tau,\al}(\De))<\ltvec(\al).
\end{equation*} 
\item If $\al\in\Ttvec\cap\Hz^{>1}\cap\Om_1$ let $\tau\in\{\tau_0,\ldots,\tau_r\}$ (where $\tau_0:=1$) be maximal such that $\tau<\al$. 
we have
\begin{equation*}\label{barlen} \ltvec(\albar) < \ltvec(\al),
\end{equation*}
and
\begin{equation*}\label{zelen}
\ltvec(\zetal) < \ltvec(\al).
\end{equation*}
In case of $\al\not\in\Ez$ we have
\begin{equation*}\label{loglen} \ltvec(\log(\al)), \ltvec(\log((1/\tau)\cdot\al)) < \ltvec(\al),
\end{equation*}
and for $\al\in\Ez$ we have
\begin{equation*}\label{lalen} 
\ltvecal(\latal) < \ltvec(\al).
\end{equation*}
\end{enumerate}
\end{lem}
{\bf Proof.} The inequalities follow directly from the term representations for the ordinal operations applied, see 
Lemma \ref{logendcharlem} and Definitions \ref{barop}, \ref{zetaldefi}, and \ref{lataldefi}.
\qed

Preparations are complete now to extend Definition \ref{kappanuprincipals} by the following two definitions. In order to do so, we need
to define the function $\dpval$, simultaneously with defining $\kval$. The \emph{depth} function $\dpf$ was first introduced in \cite{C99} 
for the analysis of $\Rone$. $\dpf$ satisfies the equation 
$\ka_{\al+1}=\ka_\al+\dpf(\al)+1$ for $\al<\ups_1$, 
i.e.\ $\lh(\ka_\al)=\ka_\al+\dpf(\al)$, where we let $\dpf$ operate on indices.
Clause 5 of the following definition applies when recovering and generalizing the recursion formula for $\Rone$ in the context of $\Rtwo$.
However, our extension of $\dpf$ to $\dpval$ for all of $\Rtwo$ does not any longer characterize the $\leo$-reach of an ordinal.
In case that the relativized $\leo$-component that $\dpval$ is applied to falls back onto a surrounding main line (recall the informal outline
given before Definition \ref{indicatorchi}), $\dpval$ leads us only to the point where that happens, as required for understanding the 
internal structure of local components, cf.\ clause 2 of Definition \ref{nuvaldefi}, rather than leading us to the 
$\leo$-reach of the component it is applied to. For an explanation in different words and some greater detail, see \cite{CWc} before 
Definition 4.4.
 
\begin{defi}[cf.\ 4.4 of \cite{CWc}]\label{kappadpdefi}\index{$\kval$}\index{$\dpval$}
Let $\alvec\in\laRS$ where $\alvec=(\ups_{\la+1},\ldots,\ups_{\la+m})^\frown(\al_1,\ldots,\al_n)$ and $\alnod:=1+\ups_{\la+m}$,
where $\ale\in(\ups_{\la+m},\ups_{\la+m+1})$ if $n>0$. 
We define global functions $\ka,\dpf:\On\to\On$, omitting superscripts $()$ for ease of notation, as well as, 
for $\alvec\not=()$, local functions $\kval,\dpval$ where 
$\domkval=[0,\laaln]$ and $\dom(\dpval)=\domkval$ if $n>0$ while $\dom(\dpval)=\ups_{\la+m+1}$ if $n=0$, 
simultaneously by recursion on $\ltvec(\be)$, where $\tauvec$ is the resolving sequence for $\be$ in $\laTalvec$-representation,
extending Definition \ref{kappanuprincipals}.
The clauses extending the definition of $\kval$ are as follows.  
\begin{enumerate}
\item $\kval_0:=0$, 
\item\label{kappapl} $\kvalbe:=\kvalga+\dpval(\ga)+\kvalde$ for $\be=_\NF\ga+\de$.
\end{enumerate}

\noindent $\dpval$ is defined as follows, using $\nu$ as already defined on $\laTS$.
\begin{enumerate}
\item $\dpf(\ups_\xi):= 0$ for all $\xi\in\On$.
\item $\dpval(0):=0$, $\dpval(1):=0$, and $\dpval(\aln):=0$ in case of $\alvec\not=()$,
\item $\dpval(\be):=\dpval(\de)$ if $\be=_\NF\ga+\de$,
\item\label{dpred} $\dpval(\be):=\dpf_{\alvecpr}(\be)$ if $\alvec\not=()$ for $\be\in\Hz\cap(1,\aln)$ where $\alvec={\alvecpr}^\frown\aln$,
\item for $\be\in\Hz^{>\aln}\setminus\Ez$ let $\ga:=(1/\aln)\cdot\be$ and $\log(\ga)=_\ANF\ga_1+\ldots+\ga_m$ and set 
      \[\dpval(\be):=\kval_{\ga_1}+\dpval(\ga_1)+\ldots+\kval_{\ga_m}+\dpval(\ga_m),\]
\item\label{dpeps} for $\be\in\Ez^{>\aln}$ let $\gavec:=(\ups_{\la+1},\ldots,\ups_{\la+m})^\frown(\ale,\ldots,\aln,\be)$, 
and define        
      \[\dpval(\be):=\nuvga_{\mu^\aln_\be}+\kvga_{\laalnbe}+\dpvga(\laalnbe).\]
\end{enumerate}
\end{defi}

\begin{defi}[cf.\ 4.4 of \cite{CWc}]\index{$\nuval$}\label{nuvaldefi}
Let $\alvec\in\laRS$ be of the form $\alvec=(\ups_{\la+1},\ldots,\ups_{\la+m})^\frown(\al_1,\ldots,\al_n)\not=()$ and set 
$\alnod:=1+\ups_{\la+m}$. 
We define the local function $\nuval$ on $[0,\mualn]$,
extending Definition \ref{kappanuprincipals} and setting $\al:=\ov(\alvec)$, by
\begin{enumerate}
\item $\nuval_0:=\al$, 
\item\label{nuple} $\nuval_{\be}:=\nuval_\ga+\kval_{\rhoalnga}+\dpval(\rhoalnga)+\chialncheck(\ga)\cdot\al$ if $\be=\ga+1$,
\item $\nuval_\be:=\nuval_\ga+\kval_{\rhoalnga}+\dpval(\rhoalnga)+\nuval_{\de}$ if $\be=_\NF\ga+\de\in\Lim$.
\end{enumerate}
\end{defi}

\begin{rmk} In Corollary \ref{nuimagecor} we will see that the image of $\nu^\alvec$ indeed consists of multiples of $\al$ and that infinite 
additive principal numbers in its domain are mapped to additive principal numbers greater than $\al$. It is not obvious, but a crucial
property of the indicator function $\chi$, that clause 2 of the above definition also yields a multiple of $\al$ if $\chi^\aln(\ga)=1$.
\end{rmk}

It is easy to see that the properties of $\ka$-, $\dpf$-, and $\nu$-functions established in Section 4 of both \cite{CWc} and \cite{W17}
extend as expected to the functions defined above. We essentially only need estimates of components in the interior of $\ups$-segments.
Of particular importance are monotonicity and continuity of $\ka$- and $\nu$-functions (Corollary \ref{kappanuhzcor}) and the verification
of agreement on the common domain with the definitions given in \cite{CWc} (Theorem \ref{agreementthm}).

\begin{lem}[cf.\ 4.6 of \cite{W17}]\label{kdpmainlem}
Let $\alvec=(\ups_{\la+1},\ldots,\ups_{\la+m})^\frown(\ale,\ldots,\aln)\in\laRS$ for the maximal such $m<\om$, and set $\alnod:=1+\ups_{\la+m}$.
\begin{enumerate}
\item Let $\ga\in\domkval\cap\Hz$, $\ga\not\in\Image(\ups)$. If $\ga=_\MNF\ga_1\cdot\ldots\cdot\ga_k\ge\aln$, setting
$\gapr:=(1/\ga_1)\cdot\ga$, we have 
\[(\kvalga+\dpval(\ga))\cdot\om=\left\{\begin{array}{ll}
            \ov(\alvec)\cdot\gapr\cdot\om&\mbox{if } \ga_1=\aln\\[2mm]
            \ov(\alvec^\frown\ga\cdot\om)&\mbox{otherwise.} 
            \end{array}\right.\]
If $\ga<\aln$ we have $(\kvalga+\dpval(\ga))\cdot\om<\ov(\alvec)$. 
\item For $\ga\in\domkval-(\Ez\cup\{0\})$ we have \[\dpval(\ga)<\kvalga.\] 
\item For $\ga\in\dom(\dpval)\cap\Ez^{>\aln}$ such that $\muga<\ga$ we have \[\dpval(\ga)<\ov(\alvec^\frown\ga)\cdot\muga\cdot\om.\]
\item For $\ga\in\domkval\cap\Ez^{>\aln}$, $\ga\not\in\Image(\ups)$, we have 
\[\kvalga\cdot\om\le\dpval(\ga) \quad\mbox{ and }\quad \dpval(\ga)\cdot\om=_\NF\ov(\homega(\alvecga))\cdot\om.\]
\item Let $\ga\in\domnuval\cap\Hz$,  $\ga=_\MNF\ga_1\cdot\ldots\cdot\ga_k\not\in\Image(\ups)$. We have
\[(\nuvalga+\ka^\alvec_{\varrho^\aln_\ga}+\dpval(\varrho^\aln_\ga))\cdot\om=\left\{\begin{array}{ll}
            \ov(\alvec)\cdot\ga\cdot\om&\mbox{if } \ga_1\le\aln\\[2mm]
            \ov(\homega(\alvec^\frown\ga))\cdot\om&\mbox{if } \ga\in\Ez^{>\aln}\\[2mm]
            \ov(\alvec^\frown\ga)\cdot\om&\mbox{otherwise.} 
            \end{array}\right.\]
\end{enumerate}
\end{lem}
{\bf Proof.}
The lemma is shown by simultaneous induction on $\ltvec(\ga)$, where $\tauvec$ is the resolving sequence for $\ga$, over all parts.
For a detailed proof see \cite{W17}.
\qed

\begin{cor}[cf.\ 4.7 of \cite{W17}]\label{kappanuhzcor} Let $\alvec\in\RS$. We have
\begin{enumerate}
\item $\kval_{\ga\cdot\om}=(\kvalga+\dpval(\ga))\cdot\om$ for $\ga\in\Hz$ such that $\ga\cdot\om\in\domkval$.
\item $\nuval_{\ga\cdot\om}=(\nuvalga+\ka^\alvec_{\varrho^\aln_\ga}+\dpval(\varrho^\aln_\ga))\cdot\om$ for $\ga\in\Hz$ such that $\ga\cdot\om\in\domnuval$.
\end{enumerate}
$\kval$ and for $\alvec\not=()$ also $\nuval$ are strictly monotonically increasing and continuous.
\end{cor}

\begin{theo}[cf.\ 4.8 of \cite{W17}]\label{agreementthm}
Let $\alvec=(\ups_{\la+1},\ldots,\ups_{\la+m})^\frown(\ale,\ldots,\aln)\in\laRS$ for the maximal such $m<\om$, and set $\alnod:=1+\ups_{\la+m}$.
For $\be\in\Hz$ let $\de:=(1/\bebar)\cdot\be$, so that $\be=_\NF\bebar\cdot\de$ if $\be\not\in\Mz$.
\begin{enumerate}
\item For all $\be\in\domkval\cap\Hz^{>\aln}$, $\be\not\in\Image(\ups)$, we have \[\kvalbe=\kval_{\bebar+1}\cdot\de.\]
\item If $\alvec\not=()$, then for all $\be\in\domnuval\cap\Hz^{>\aln}$, $\be\not\in\Image(\ups)$, we have 
\[\nuvalbe=\nuval_{\bebar+1}\cdot\de.\]
\end{enumerate}
The definitions of $\kappa,\nu$, and $\dpf$ extend the definitions given in \cite{CWc} and \cite{W17}.
\end{theo}

For the following estimates recall Definition \ref{alhatdefi}. Note that these estimates confirm that our relativized systems
$\Tt$ just suffice to contain relative connectivity components that occur between elements of $\Image(\ups)$. 

\begin{lem}[cf.\ 4.9 of \cite{W17}]\label{hatmainlem}
Let $\alvec=(\ale,\ldots,\aln)\in\RS$, $n>0$. 
\begin{enumerate}
\item For all $\be$ such that $\alvecbe\in\TS$ and $\be\not\in\Image(\ups)$ we have 
\[\ov(\alvecbe)<\ov(\alvec)\cdot\widehat{\aln}.\]
\item For all $\ga$ such that $\alvecga\in\RS$ and $\ga\not\in\Image(\ups)$ we have 
\[\ov(\homega(\alvecga))<\ov(\alvecga)\cdot\gahat.\]
\end{enumerate}
\end{lem}

\begin{cor}[cf.\ 4.10 of \cite{W17}]\label{kdpnuestimcor}
For all $\alvecga\in\RS$ such that $\ga\not\in\Image(\ups)$ the ordinal $\ov(\alvecga)\cdot\gahat$ is a strict upper bound of  
\[\Image(\kappa^{\alvecga}),\: \Image(\nu^{\alvecga}),\: \dpval(\ga),\: \mbox{and }
\nu^{\alvecga}_{\mu_\ga}+\kappa^{\alvecga}_{\la_\ga}+\dpf_{\alvecga}(\la_\ga).\]
\end{cor}
{\bf Proof.} This directly follows from Lemmas \ref{kdpmainlem} and \ref{hatmainlem}.
\qed

In order to formulate the assignment of tracking chains to ordinals in Subsection \ref{tcassignmentsubsec} we need to introduce a suitable 
notion of tracking sequence relative to a given context, as we did in \cite{CWc}. 
We first introduce an evaluation function for relativized tracking sequences.
\begin{defi}[cf.\ 4.13 of \cite{CWc}] Let $\alvec=(\ale,\ldots,\aln)\in\RS$ where $n>0$.
We define 
\[\TSalvec:=\set{\gavec\in\TSaln}{\ga_1\le\laalnminaln}\]\index{$\tst$@$\TSt$!$\TSalvec$}
and for $\gavec^\frown\be\in\TSalvec$
\[\ordvalalvec\left(\gavec^\frown\be\right):=\left\{\begin{array}{l@{\quad}l}
\kval_{\be} & \mbox{if } \gavec=()\\[2mm]
\nu^{\alvec^\frown\gavec}_\be & \mbox{otherwise.}
\end{array}\right.\]\index{$\ordvalue$!$\ordvalalvec$}
For convenience we identify $\TS^{()}$ with $\TS$ and
$\ordvalue^{()}$ with $\ordvalue$.
\end{defi}
{\bf Remark.} Note that this is well-defined thanks to part c) of Lemma \ref{rhomulamestimlem}. Notice also that $\TSalvec$ is a 
$\klex$-initial segment of $\TSaln$ and that in the case $\aln\in\Image(\ups)$ the sets $\TSalvec$ and $\TSaln$ coincide, 
since $\gavec\in\TSaln$ implies that $\ga_1<\al_n^\infty$. 
If $\aln\in\Image(\ups)$, the evaluation functions $\ordvalalvec$ and $\ordvalue$ agree.
We have the following 

\begin{lem}[cf.\ 4.14 of \cite{CWc}]\label{lamveclem}
Let $\alvec=(\ale,\ldots,\aln)\in\RS$ where $n>0$ and $\aln\not\in\Image(\ups)$.
Let $\la_1:=\mc(\la_\aln)$, and whenever $\la_i$ is defined and $\la_i\in\Ez^{>\la_{i-1}}$ (setting $\la_0:=\aln$), 
let $\la_{i+1}:=\mu_{\la_i}$. If we denote the resulting vector by $(\la_1,\ldots,\la_k)=:\vec{\la}$
then $\TSalvec$ is the initial segment of $\TSaln$ with $\klex$-maximum $\vec{\la}$.
We have 
\[\ordvalalvec(\vec{\la})=\mc(\ka^\alvec_{\la_\aln}+\dpf_\alvec(\la_\aln)).\]
\end{lem}
{\bf Proof.} The proof is by evaluation of $\mc(\ka^\alvec_{\la_\aln}+\dpf_\alvec(\la_\aln))$ using Lemmas \ref{kdpmainlem}, \ref{hatmainlem},
and Corollary \ref{kdpnuestimcor}.
\qed

The analogue to Lemma \ref{trsvallem} is as follows. Notice that we have to be careful regarding multiples of indices versus their evaluations.
\begin{lem}[cf.\ 4.15 of \cite{CWc}]\label{reltrsvallem}
Let $\alvec$ and $\gavec^\frown\be\in\TSalvec$ be as in the above definition and set $\al:=\ordvalue(\alvec)$.
\begin{enumerate}
\item We have 
\[\trsaln(\aln\cdot((1/\al)\cdot\ordvalalvec(\gavec^\frown\be)))=\gavec^\frown\be.\] 
\item For $\de\in\Hz\cap[\aln,\al_n^\infty)$ such that $\trsaln(\de)\in\TSalvec$ we have 
\[\ordvalalvec(\trsaln(\de))=\al\cdot((1/\aln)\cdot\de).\]
\item If $\aln\not\in\Image(\ups)$, setting $\la:=\aln\cdot((1/\al)\cdot\mc(\ka^\alvec_{\la_\aln}+\dpf_\alvec(\la_\aln)))$ we have 
\[\trsaln(\la)=\vec{\la}\in\TSalvec\]
for $\vec{\la}$ as defined in Lemma \ref{lamveclem}, and
the mapping $\trsaln$ is a $<$-$\klex$-order isomorphism of 
\[\set{\de\in\Hz\cap[\aln,\al_n^\infty)}{\trsaln(\de)\in\TSalvec}=[\aln,\la]\cap\Hz\] with $\TSalvec$.
\item If $\aln\in\Image(\ups)$, the mapping $\trsaln$ is a $<$-$\klex$-order isomorphism of 
$\Hz\cap[\aln,\al_n^\infty)$ with $\TSalvec$.
\end{enumerate}
\end{lem}
{\bf Proof.} Once the first claim of the lemma is shown by induction along $\klex$ on $\TSalvec$, the remaining claims follow using Lemmas \ref{citedinjtrslem} and \ref{lamveclem}. 
In proving the first claim for $\gavec^\frown\be\in\TSalvec$, say $\gavec=(\ga_1,\ldots,\ga_m)$ where $m\ge 0$, we proceed in analogy with the course of proof of Lemma \ref{trsvallem} (Theorem 3.20 of \cite{W17}),
replacing $\alvec$ with $\alvec^\frown\gavec$, $\aln$ with $\ga_m$ (setting $\ga_0:=\aln$), and $\al$ with $\ga:=\aln\cdot((1/\al)\cdot\ordvalue(\gavec))$ in the case $m>0$, where
by the i.h.\ we have $\trsaln(\ga)=\gavec$. 
\qed

\subsection{Extending the concept of tracking chains}
Every additive principal number $\al$ in $\Rtwo$ can be described uniquely by a tracking sequence
in the way we explained at the beginning of this Section, cf.\ (\ref{principalchain}).
Tracking chains, introduced first in Section 5 of \cite{CWc}, are sequeces of tracking sequences in the sense that each of these sequences
starts with a $\ka$-index and then possibly continues with $\nu$-indices. However, these $\ka$- and $\nu$-indices do not need to be
(infinite) additive principal numbers as is the case with tracking sequences. 
From the second sequence on the $\ka$-indices may have to be relativized to the local connectivity component specified by the preceding 
sequences. $\ka$-indices specify the (relativized local) $\leo$-component, while $\nu$-indices of suitably relativized $\nu$-functions
describe the selection of the nested $\letwo$-components.
Extending the \emph{address space} from the set of tracking sequences to the set of tracking chains allows us to describe
all ordinals uniquely in terms of indices of nested $\le_i$-components. Characterizing this general address space for ordinals 
in terms of conditions on the indices of nested $\le_i$-components is unfortunately more complicated than expected. 

Before providing the formal definitions determining the class $\TC$ of tracking chains for $\Rtwo$, we mention two simple examples for 
tracking chains. First, as mentioned above, the tracking chain of an additive principal number $\al$ is $(\alvec)$, where 
$\alvec:=\trs(\al)$ is the tracking sequence of $\al$.
Second, for an ordinal $\be=\ka_{\be_1}+\ldots+\ka_{\be_n}$ for suitable indices $\be_1>\ldots>\be_n$ has the tracking chain 
$((\be_1),\ldots,(\be_n))$. For $\be$ below $\epsn$ the conditions on the $\be_i$ are simply $\be_1<\epsn$ and $0<\be_{i+1}\le\logend(\be_i)$
for $i=1,\ldots,n-1$, as shown in \cite{C99}.

Strictly speaking, ordinals greater than $\ups_1$
drop out of the usual $2$-dimensional format of tracking chains for ordinals below $\oneinf=\ups_1$ as introduced in \cite{CWc}.
While a suitable address space for $\Rom$ would require a bookkeeping formalism handling $\om\times\om$-matrices of ordinals all but finitely
many entries of which are left blank and in which addresses for ordinals below $\ups_1$ would be represented in an equivalent but different way,
as long as we stay in $\Rtwo$, we can extend the address space $\TC$ of tracking chains without too many complications.

As we will see, ordinals $\ups_\la$ where $\la\in\Lim$ have $\ups_{\la+1}$-many $\ups_\la$-$\leo$-minimal successors. We will denote
the ordinal $\ups_\la$ by the chain $((\ups_\la))$, its least $<_2$-successor $\ups_\la\cdot 2$ by $((\ups_\la\cdot 2))$, and the
largest $\ups_\la$-$\leo$-minimal ordinal $\ups_{\la+1}$ by $((\ups_{\la+1}))$. The least $<_2$-successor of $\ups_2$ above $\ups_\la$
is denoted by $((\ups_\la+\ups_2))$. We will denote the ordinal $\ups_{\om^2+2}$, which we considered in the introduction in the context of 
$\Rthree$, simply by $((\ups_{\om^2+1},\ups_{\om^2+2}))$. 

Note that we lose a nice property of the original tracking chains: Namely that tracking chains of
all $<_2$-predecessors of an ordinal occur as initial chains of its tracking chain. However, in the presence of infinite $<_2$-chains this property cannot be kept anyway.

We are now going to extend Definition 5.1 of \cite{CWc} to a system of tracking chains for all of $\Rtwo$. The first sequence of a 
generalized tracking chain $\alvec$ will determine the $\ups$-segment in which the ordinal with address $\alvec$ is located. 
This will be called $\upsseg(\alvec)$, the $\ups$-segment of $\alvec$. For better accessiblity we are going to split up the 
definition of the class $\TC$ of tracking chains for all of $\Rtwo$ into several steps, as compared to the orginal definition in \cite{CWc}.
We begin with templates for tracking chains, which we call \emph{index chains}, and some useful general terminology.

\begin{defi}[Index chains, their domains, associated and initial chains]\label{indexchaindefi} \mbox{ }

\begin{enumerate}
\item An {\bf index chain} is a sequence $\alvec=(\alevec,\ldots,\alnvec)$, $n\in(0,\om)$, of ordinal vectors 
$\alivec=(\alcp{i,1},\ldots,\alcp{i,m_i})$ with $m_i\in(0,\om)$ for $1\le i\le n$.

\item We define $\dom(\alvec)$\index{$\dom(\alvec)$} to be the set of 
all index pairs of $\alvec$, that is
\[\dom(\alvec):=\set{(i,j)}{1\le i\le n\andsp1\le j \le m_i}.\]  

\item The vector $\tauvec=\left(\tauevec,\ldots,\taunvec\right)$ defined by $\taucp{i,j}:=\sumend(\alcp{i,j})$ for $(i,j)\in\dom(\alvec)$
(that is, $\taucp{i,j}$ is the least additive component of $\alcp{i,j}$) is called the 
{\bf\boldmath  chain associated with $\alvec$\unboldmath}.\index{chain associated with an index chain}

\item The {\bf initial chains $\alvecrestrarg{i,j}$ of $\alvec$}\index{initial chain}\index{$\alvecrestrarg{i,j}$}, where
$(i,j)\in\dom(\alvec)$, are 
\[\alvecrestrarg{i,j}:=((\alcp{1,1},\ldots,\alcp{1,m_1}),\ldots,(\alcp{i-1,1},\ldots,\alcp{i-1,m_{i-1}}),(\alcp{i,1},\ldots,\alcp{i,j})).\]
By $\alvecrestrarg{i}$\index{$\alvecrestri$!$\alvecrestrarg{i}$ for tracking chains} we abbreviate $\alvecrestrarg{i,m_i}$.
For convenience we set $\alvecrestrarg{i,0}:=()$ for $i=0,1$ and $\alvecrestrarg{i+1,0}:=\alvecrestrarg{i,m_i-1}$ for $1\le i<n$.
Initial chains of $\tauvec$ are defined in the same way.
\end{enumerate}
\end{defi}

Next we impose conditions familiar from tracking sequences, cf.\ definitions \ref{TSdefi}, \ref{RSdefi}, and \ref{laRSlaTSdefi}, in order to
introduce a notion of regularity with respect to $\nu$-indices.

\begin{defi}\label{nuregdefi}
Let $\alvec=(\alevec,\ldots,\alnvec)$, where $\alivec=(\alcp{i,1},\ldots,\alcp{i,m_i})$ for $i=1,\ldots,n$, 
be an index chain with associated chain $\tauvec$. 
We call $\alvec$ a {\bf \boldmath$\nu$-regular index chain\unboldmath} if for every $i\in\{1,\ldots,n\}$ such that $m_i>1$
\[\taucp{i,1}<\ldots<\taucp{i,m_i-1}\]
where $\taucp{i,j}\in\Ez$ and 
\[\alcp{i,j+1}\le\mu_\taucp{i,j}\]
for $j=1,\ldots,m_i-1$.
\end{defi}

For $\nu$-regular index chains the following definition becomes meaningful, which allows us to introduce a notion of \emph{$\ups$-segmentation} 
of such index chains. 

\begin{defi}[cf.\ Definition \ref{upssegdefi}]
For a $\nu$-regular index chain $\alvec=(\alevec,\ldots,\alnvec)$
let $\la\in\Limnod$ be maximal such that $\ups_\la\le\alcp{1,1}$ and $t<\om$ be maximal such that $\alvec_1$ is of the form 
\[\alvec_1=(\ups_{\la+1},\ldots,\ups_{\la+t},\ga_1,\ldots,\ga_l),\]
hence $\ga_1<\ups_{\la+t+1}$ if $l>0$.
Then $(\la,t)$ indicates the {\bf \boldmath$\ups$-segment\unboldmath} of $\alvec$, $\upsseg(\alvec):=(\la,t)$.
\end{defi}

\begin{defi}
Let $\alvec=(\alevec,\ldots,\alnvec)$, where $\alivec=(\alcp{i,1},\ldots,\alcp{i,m_i})$ for $i=1,\ldots,n$, be a $\nu$-regular index chain 
with associated chain $\tauvec$ and set $(\la,t):=\upsseg(\alvec)$. 
We say that $\alvec$ is {\bf \boldmath$\upsilon$-segmented\unboldmath} if the following two conditions hold.
\begin{enumerate}
\item $\alcp{1,1}\in[\ups_\la,\ups_{\la+1}]$, where
\[\alcp{1,1}=\ups_\la\mbox{ implies }n=m_n=1,\]
and unless $\alvec=((0))$, all indices $\alcp{i,j}$ are nonzero.
\item The sequences \[(\upsseg(\alcp{i,1}))_{2\le i\le n}\mbox{ and }(\upsseg(\taucp{i,1}))_{2\le i\le n}\]
are $\kglex$-weakly decreasing with upper bound $(\la,t)$.
\end{enumerate}
\end{defi}

\begin{defi}\label{segmentationdefi}
Let $\alvec=(\alevec,\ldots,\alnvec)$, where $\alivec=(\alcp{i,1},\ldots,\alcp{i,m_i})$ for $i=1,\ldots,n$, be an $\upsilon$-segmented 
$\nu$-regular index chain with associated chain $\tauvec$ and set $(\la,t):=\upsseg(\alvec)$. 
Let $s_1\in\{1,\ldots,n\}$ be minimal such that \[\taucp{s_1,1}<\ups_\la\] if that exists, in which case we 
further let $s_1<\ldots<s_p\le n$ indicate all indices $i\ge s_1$ where $\upsseg(\taucp{i,1})$ strictly decreases,
hence \[\upsseg(\taucp{s_j-1,1})\glex\upsseg(\taucp{s_j,1})\] for $j=2,\ldots,p$, otherwise we set $p:=0$ and $s_0:=1$. Let 
\[(\la_j,t_j):=\upsseg(\taucp{s_j,1})\] for $j=1,\ldots,p$ indicate the corresponding $\ups$-segments and set for convenience
$(\la_0,t_0):=(\la,t)$. 
We call 
\begin{enumerate}
\item $p$ the {\bf segmentation depth}, 
\item $(s_1,\ldots,s_p)$ the {\bf segmentation signature}, and
\item $(\la_i,t_i)_{0\le i\le p}$ the {\bf sequence of \boldmath$\ups$-segments\unboldmath}
\end{enumerate} 
of $\alvec$.
\end{defi}

Note that it is very well possible to have $p>0$ and $s_1=1$, a situation that occurs when $\alcp{1,1}>\ups_\la$ while $\taucp{1,1}<\ups_\la$
where $\la\in\Lim$.

The notion of \emph{unit} defined next allows us to trace back the greatest $\ktwo$-predecessor of an ordinal with tracking chain an
initial chain $\alvecrestrarg{i,1}$, $1\le i\le n$, of $\alvec$. 
It determines the setting of relativization (reference sequence) of 
the relativized $\ka$-function applicable in the $\letwo$-component rooted in that greatest $\ktwo$-predecessor of the ordinal with 
tracking chain $\alvecrestrarg{i,1}$.  

\begin{defi}[Units]\label{unitsdefi}
Let $\alvec=(\alevec,\ldots,\alnvec)$, where $\alivec=(\alcp{i,1},\ldots,\alcp{i,m_i})$ for $i=1,\ldots,n$, be an $\upsilon$-segmented 
$\nu$-regular index chain with associated chain $\tauvec$, segmentation depth $p$, segmentation signature $(s_1,\ldots,s_p)$,
and sequence of $\ups$-segments $(\la_i,t_i)_{0\le i\le p}$.
The {\bf\boldmath  $i$-th unit $\tauistar$\unboldmath}\index{unit}\index{$\tauistar$}
of $\alvec$ and its {\bf  index pair \boldmath$i^\star$\unboldmath}\index{index pair}\index{$*$@$i^\star$} for $1\le i\le n$ 
is defined as follows.
\[i^\star\mbox{, }\tauistar:=\left\{\begin{array}{ll}
   (l,j)\mbox{, }\;\;\taucp{l,j} & \mbox{ if $\klex$-max.\ } (l,j)\in\dom(\alvec) 
                             \mbox{ exists s.t.\ }(l,j)\klex(i,1)\mbox{, }j<m_l\mbox{ and }\taucp{l,j}\le\taucp{i,1}\\[2mm]
   (s_j,0)\mbox{, }1+\ups_{\la_j+t_j} & \mbox{ if otherwise max.\ }j\in\{1,\ldots,p\}\mbox{ exists s.t.\  }s_j\le i\\[2mm]
   (1,0)\mbox{, }\;1+\ups_\la & \mbox{ otherwise.}
\end{array}\right.\]
For technical convenience we set $\taucp{i,0}:=\tauistar$ for $i=1,\ldots,n$.
\end{defi}

Note that the definition of $\taucp{i,0}$ deviates from Definition 5.1 of \cite{CWc} but is conceptually more appropriate. 
Clearly, for $j=1,\ldots,p$ we have $(s_j)^\star=(s_j,0)$ by definition.
Note further that, unless $\alcp{1,1}=0$, for $i=1,\ldots,n$, by definition we have 
\[i^\star\klex(i,1),\quad \tauistar\le\taucp{i,1}, \quad\mbox{ and }\quad\upsseg(\taucp{i,1})=\upsseg(\tauistar).\]
For $\alvec=((1))$, i.e.\ $\alcp{1,1}=1$, we have $\tau_1^\star=1$, and for $\alvec=((\ups_\la))$ where $\la\in\Lim$ we have 
$\tau_1^\star=\ups_\la$. 

The more general notion of \emph{base} also considers greatest $\ktwo$-predecessors of ordinals with tracking chain $\alvecrestrarg{i,j}$
where $j\ge 1$ of a tracking chain $\alvec$. In the case $j>1$ the base is given by $\taucp{i,j-1}$, so the next definition simply
introduces a more convenient terminology to name bases.

\begin{defi}[Bases]\label{basesdefi}
Let $\alvec=(\alevec,\ldots,\alnvec)$, where $\alivec=(\alcp{i,1},\ldots,\alcp{i,m_i})$ for $i=1,\ldots,n$, be an $\upsilon$-segmented 
$\nu$-regular index chain with associated chain $\tauvec$.
For $(i,j)\in\dom(\alvec)$ we define the index pair $(i,j)^\prime$\index{$'$@$(i,j)^\prime$} by 
$i^\star$ if $j=1$ and by $(i,j-1)$ otherwise. 
The {\bf  base $\taucppr{i,j}$ of $\taucp{i,j}$ in $\alvec$}\index{base} is defined by
\index{$\taucppr{i,j}$}
\[\taucppr{i,j}:=\left\{\begin{array}{l@{\quad}l}
  \tauistar & \mbox{ if } j=1\\[2mm] 
   \taucp{i,j-1} & \mbox{ if } j>1.\end{array}\right.\]
For $1\le i\le n$ we define the {\bf\boldmath  $i$-th maximal base $\tauipr$\unboldmath} of 
$\alvec$\index{base!maximal base}\index{$\tauipr$} by 
\[\tauipr:=\taucppr{i,m_i}.\]
\end{defi}

With the above preparations we can determine how many $\al$-$\leo$-minimal $\lo$-successors an ordinal $\al$ represented by tracking 
chain $\alvec=(\alevec,\ldots,\alnvec)$ has. This type of condition also applies to the initial chains $\alvecrestrarg{i,1}$, $1\le i\le n$.
We call these upper bounds on the number of immediate $\lo$-successors the $i$-th critical indices. The formal definition is as follows. 

\begin{defi}[Critical \boldmath$\ka$-indices\unboldmath]\label{critkapdefi}
Let $\alvec=(\alevec,\ldots,\alnvec)$, where $\alivec=(\alcp{i,1},\ldots,\alcp{i,m_i})$ for $i=1,\ldots,n$, be an $\upsilon$-segmented 
$\nu$-regular index chain with associated chain $\tauvec$.
We define the {\bf\boldmath  $i$-th critical index of $\alvec$\unboldmath}\index{critical index of a tracking chain},
written as $\rhoi(\alvec)$, in short as $\rhoi$ if no confusion is likely, by
\[\rhoi:=\left\{\begin{array}{ll}
  \log\left((1/\tauistar)\cdot\taucp{i,1}\right)+1&\mbox{if }\;m_i=1\\[2mm]
  \rhoargs{\tauipr}{\taucp{i,m_i}}+\tauipr&\mbox{if }\;
     m_i>1\andsp\taucp{i,m_i}<\mu_{\tauipr}\andsp\chi^\tauipr(\taucp{i,m_i})=0\\[2mm]
  \rhoargs{\tauipr}{\taucp{i,m_i}}+1&\mbox{if }\;
     m_i>1\andsp\taucp{i,m_i}<\mu_{\tauipr}\andsp\chi^\tauipr(\taucp{i,m_i})=1\\[2mm]
  \la_{\tauipr}+1&\mbox{otherwise.}\end{array}\right.\]
\end{defi}

Note that $\al<\rhoi$ implies that $\al\le\taucp{i,1}$ if $m_i=1$, and $\al\le\la_{\taucp{i,m_i-1}}$ if $m_i>1$, since
according to part b) of Lemma \ref{rhomulamestimlem} we have $\rhoi\le\la_{\tauipr}+1$ where $\tauipr=\taucp{i,m_i-1}$. 
Note that the second clause is not a successor ordinal since we are approaching an ordinal that should be indexed by a $\nu$-index 
and therefore need to avoid ambiguity. This preference of $\nu$-indices over $\ka$-indices whenever possible also motivates the second
part of part 2 in the following definition, whereas part 1 is a necessity: If $\tauistar=\taucp{i,1}$ then we must have $m_i=1$ since
there do not exist $\lo$-sucessors of successor-$\ktwo$-successors (and hence we also must have $i=n$), cf.\ Lemma \ref{loalpllem}.
 
\begin{defi}[\boldmath$\ka$-regularity\unboldmath]\label{kapparegdefi}
Let $\alvec=(\alevec,\ldots,\alnvec)$, where $\alivec=(\alcp{i,1},\ldots,\alcp{i,m_i})$ for $i=1,\ldots,n$, be an $\upsilon$-segmented 
$\nu$-regular index chain with associated chain $\tauvec$.
$\alvec$ is called {\bf \boldmath $\ka$-regular\unboldmath} if the following conditions hold.
\begin{enumerate}
\item For all $i\in\{1,\ldots,n\}$ such that $m_i>1$ 
\[\tauistar<\taucp{i,1}.\]
\item For any  $i\in\singleton{1,\ldots,n-1}$
\[\alcp{i+1,1}<\rhoi,\]
and in case of $\tauipr<\taucp{i,m_i}\in\Ez$
\[\alcp{i+1,1}\not=\taucp{i,m_i}.\] 
\end{enumerate}
In short, such $\alvec$ is called a {\bf \boldmath$\ka\nu\ups$-regular\unboldmath} index chain. 
\end{defi}
 
Note that the class of $\ka\nu\ups$-regular index chains is closed under non-empty initial chains. The following lemma is technical,
but an important property of $\ka\nu\ups$-regular index chains. Essentially, this was shown in part b) of Lemma 5.8 (using Lemma 5.7) 
of \cite{CWc}, but we provide a new proof to enhance clarity.

\begin{lem}[cf.\ 5.7 and 5.8 of \cite{CWc}]\label{rserslem}
Let $\alvec=(\alevec,\ldots,\alnvec)$, where $\alivec=(\alcp{i,1},\ldots,\alcp{i,m_i})$ for $i=1,\ldots,n$, be a $\ka\nu\ups$-regular 
index chain with associated chain $\tauvec$.
For every $i$, $1\le i \le n$, we have \[\taucp{i,1}\le\la_\tauistar,\] provided that $\tauistar>1$.
\end{lem}
{\bf Proof.} 
We may assume that $i=n$ and $m_n=1$ since otherwise we can simply truncate $\alvec$ to $\alvecrestrarg{i,1}$. 
Suppose that $\taunstar>1$. We consider cases according to Definition \ref{unitsdefi}.
Let $p$, $s_1,\ldots,s_p$, and $(\la_l+t_l)_{0\le l\le p}$ be according to Definition \ref{segmentationdefi}.\\[2mm]
{\bf Case 1:} $n^\star=(j,k)\in\dom(\alvec)$. Then we have $j<n$ and $k<m_j$.\\[1mm]
{\bf Subcase 1.1:} $k+1<m_j$. Then $\taucp{n,1}<\taucp{j,k+1}\le\mu_\taucp{j,k}$ and $\taucp{j,k+1}\in\Ez^{>\taucp{j,k}}$.
    By part c) of Lemma \ref{rhomulamestimlem} we have \[\taucp{j,k+1}\le\la_\taucp{j,k}.\]
{\bf Subcase 1.2:} $k+1=m_j$. We have $\taucp{j+1,1}<\rho_j$, and since $\rho_j\le\la_\taucp{j,k}+1$ by part b) of Lemma        \
    \ref{rhomulamestimlem}, we have \[\taucp{j+1,1}\le\la_\taucp{j,k}.\]
{\bf Subcase 1.2.1:} $\min\{l\in(j,n)\mid m_l>1\}$ exists. Then we have \[\taucp{n,1}<\taucp{l,1}\le\ldots\le\taucp{j+1,1}\le\la_\taucp{j,k}.\]
{\bf Subcase 1.2.2:} Otherwise. Then we have \[\taucp{n,1}\le\ldots\le\taucp{j+1,1}\le\la_\taucp{j,k}.\]
{\bf Case 2:} $\taunstar=\ups_{\la_p+t_p}$ where $p>0$ and $\la_p+t_p>0$ and $n^\star=(s_p,0)$. Then according to the definition we have
\[\taucp{n,1}\in[\ups_{\la_p+t_p},\ups_{\la_p+t_p+1}).\] 
{\bf Case 3:} $\taunstar=\ups_\la$ where $n^\star=(1,0)$ and $\la\in\Lim$. Then according to the definition we have 
\[\taucp{n,1}\in[\ups_\la,\ups_{\la+1}).\]
Regarding cases 2 and 3, note that $\la_{\ups_\iota}=\ups_{\iota+1}$ for all $\iota>0$.
\qed

Next we consider stepwise (maximal) extension
of $\ka\nu\ups$-regular index chains and show that the class of such index chains is closed under maximal extension.
Index chains of the form $((\ups_{\la+1},\ldots,\ups_{\la+m}))$, which in principle could be maximally extended infinitely many
times, are excluded from this extension procedure as they do not play any role in the upcoming definition of tracking chains. 
For clarification, this notion of maximal extension does not change the local $\le_i$-component it originates from. In particular, it does
not jump from a local $\letwo$-component indexed by a submaximal $\nu$-index to another $\letwo$-component indexed by a larger $\nu$-index
in the same domain.

\begin{defi}[cf.\ 5.2 of \cite{CWc}]\label{maxextdefi}
Let $\alvec$ be a $\ka\nu\ups$-regular index chain with components $\alivec=(\alcp{i,1},\ldots,\alcp{i,m_i})$ for $1\le i\le n$ and associated
chain $\tauvec$.
The {\bf extension index}\index{extension index} for $\alvec$ is defined via the following cases, 
setting $\tau:=\taucp{n,m_n}$ and $\taupr:=\taunpr$:
\begin{enumerate}
\item[0.] $\taupr<\tau\in\Image(\ups)$: Then an extension index for $\alvec$ is not defined.
\item $m_n=1$: We consider three subcases:
    \begin{enumerate}
    \item[1.1.] $\taupr=\tau$: Then $\alvec$ is already maximal. An extension index for $\alvec$ does not exist.
    \item[1.2.] $\taupr<\tau\in\Ez$: Then the extension index for $\alvec$ is $\alcp{n,2}:=\mu_\tau$.
    \item[1.3.] Otherwise: Then the extension index for $\alvec$ is $\alcp{n+1,1}:=\log\left((1/\taupr)\cdot\tau\right)$.
    \end{enumerate}
\item $m_n>1$: We consider three subcases.
    \begin{enumerate}
    \item[2.1.] $\tau=1$: Then $\alvec$ is already maximal. An extension index for $\alvec$ does not exist.
    \item[2.2.] $\taupr<\tau\in\Ez$: Here we consider another two subcases.
        \begin{enumerate}
        \item[2.2.1.] $\tau=\mu_\taupr<\la_\taupr$: 
            Then the extension index for $\alvec$ is $\alcp{n+1,1}:=\la_\taupr$.
        \item[2.2.2.] Otherwise: Then the extension index for $\alvec$ is $\alcp{n,m_n+1}:=\mu_\tau$.
        \end{enumerate}
    \item[2.3.] Otherwise: We consider again two subcases.
        \begin{enumerate}
        \item[2.3.1.] $\tau<\mu_\taupr$: Then the extension index for $\alvec$ is $\alcp{n+1,1}:=\rhoargs{\taupr}{\tau}$.
        \item[2.3.2.] Otherwise: Then the extension index for $\alvec$ is $\alcp{n+1,1}:=\la_\taupr$.
        \end{enumerate}
    \end{enumerate}
\end{enumerate}
If the extension index for $\alvec$ is defined, we denote the extension of $\alvec$ by this index by $\ec(\alvec)$\index{$\ec$} 
and call this extended chain the {\bf\boldmath maximal $1$-step extension\unboldmath} of $\alvec$.\index{extension!maximal $1$-step extension}
If the iteration of maximal $1$-step extensions terminates after finitely many steps, we call the resulting index chain 
the {\bf maxplus extension} of $\alvec$ and denote it by $\mepl(\alvec)$.\index{$\me$!$\mepl$}
\end{defi}

\begin{rmk} In \cite{CWc} we called the extension of $\alvec$ by the extension index the \emph{extension candidate} for $\alvec$ since 
it might fail to be a tracking chain. Here we keep definitions essentially compatible with earlier work but avoid the notion of tracking 
chain since it has not been defined yet. In earlier work we called $\ec(\alvec)$ the maximal $1$-step extension only if it was a tracking chain.
The maxplus extension of the index chain $\alvec$, $\mepl(\alvec)$, might fail to be a tracking chain, cf.\ condition 2 
of Definition \ref{trackingchaindefi}. 

Regarding the formulation of the new clause 0, note that the extension index of the chain 
$\alvec=((\ups_1,\ups_1))$ is $\alcp{2,1}=\ups_1^2$ in the same way as $\Ga_0^2$ is the extension index of the chain $((\Ga_0,\Ga_0))$, by 
an application of clause 2.3.1. 
\end{rmk}

\begin{lem} The class of $\ka\nu\ups$-regular index chains is closed under maximal $1$-step extensions.
\end{lem}
{\bf Proof.} This immediately follows from the definitions. \qed

The following lemma settles that maximal extension is a finite process.
For an alternative, more general proof of termination regarding arbitrary $1$-step extensions, see Definition 5.3 and Lemma 5.4 of \cite{CWc}.

\begin{lem}\label{meterminationlem} 
The process of iterated maximal $1$-step extensions always terminates, hence $\mepl(\alvec)$ always exists.
\end{lem}
{\bf Proof.} The $\ltvec$-measure applied from the second extension index on strictly decreases if we omit those applications
of the $\mu$-operator (cases 1.2 and 2.2.2) that are directly followed by an application of the $\la$-operator (cases 2.2.1 and 2.3.2). 
Note that clause 2.3.1 can only be applied at
the beginning of the process of maximal extension since all $\nu$-indices occurring as maximally extending indices are maximal, i.e.\ obtained
by application of the $\mu$-operator. If the $\mu$-operator is applied twice in a row during maximal extension, say first extending
by $\mu_\tau$ and next by $\mu_{\mu_\tau}$, then according to the definition we have $\mu_\tau=\la_\tau$. Otherwise the next maximally extending
index after $\mu_\tau$ is the $\ka$-index $\la_\tau$. 
\qed

\begin{defi}
Let $\alvec$ be a $\ka\nu\ups$-regular index chain with components $\alivec=(\alcp{i,1},\ldots,\alcp{i,m_i})$ for $1\le i\le n$ and associated
chain $\tauvec$. 
The $\klex$-greatest index pair $(i,j)$ of $\alvec$ after which the elements of $\alvec$ fall onto the main line starting at $\alcp{i,j}$ 
is called the {\bf\boldmath  critical main line index pair of $\alvec$\unboldmath}.\index{index pair!critical main line index pair} 
The formal definition is as follows:

If there exists the $\klex$-maximal $(i,j)\in\dom(\alvec)$ such that $j<m_i$ and $\alcp{i,j+1}<\mu_\taucp{i,j}$,
and if $(i,j)$ satisfies the following conditions:
\begin{enumerate} 
\item $\chi^\taucp{i,j}(\taucp{i,j+1})=1$ and
\item $\alvec$
is reached by maximal $1$-step extensions starting from $\alvecrestrarg{i,j+1}$,
\end{enumerate}
then $(i,j)$ is called the critical main line index pair of $\alvec$, written as $\cml(\alvec)$. Otherwise $\alvec$ does not possess a critical main line index pair.\index{$\cml$}
\end{defi}

\begin{defi}[cf.\ 5.1 of \cite{CWc}]\label{trackingchaindefi}
Let $\alvec$ be a $\ka\nu\ups$-regular index chain with components $\alivec=(\alcp{i,1},\ldots,\alcp{i,m_i})$ for $1\le i\le n$ and associated
chain $\tauvec$. $\alvec$ is called a {\bf tracking chain}\index{tracking chain} if the following conditions hold.
\begin{enumerate}
\item All non-empty proper initial chains $\alvecrestrarg{i,j}$ of $\alvec$ are tracking chains.
\item If $m_n=1$ and if $\alvec$ possesses a critical main line index pair $\cml(\alvec)=(i,j)$, then 
\[\taucp{n,1}\not=\taucp{i,j}.\]
\end{enumerate}

By $\TC$\index{$\tc$@$\TC$} we denote the class of all tracking chains. 
For a tracking chain $\alvec$ with $\upsseg(\alvec)=(\la,t)$ we also write $\alvec\in\latTC$ or $\alvec\in\laTC$.

\medskip
\noindent {\bf Useful notation:}
\begin{enumerate}
\item By $(i,j)^+$\index{$+$@$(i,j)^+$} we denote the immediate $\klex$-successor of $(i,j)$ in $\dom(\alvec)$ if that exists and $(n+1,1)$ otherwise. For convenience we set $(0,0)^+:=(1,1)$ and $(i,0)^+:=(i,1)$ for $i=1,\ldots,n$.

\item Due to frequent future occurrences, we introduce the following notation for the modification of a tracking chain's last ordinal.
\[\alvec[\xi]:=\left\{\begin{array}{l@{\quad}l}
\alvecrestrarg{n-1}^\frown(\alcp{n,1},\ldots,\alcp{n,m_n-1},\xi) & \mbox{ if } \xi>0\mbox{ or } (n,m_n)=(1,1)\\
\alvecrestrarg{n-1}^\frown(\alcp{n,1},\ldots,\alcp{n,m_n-1}) & \mbox{ if } \xi=0\mbox{ and } m_n>1\\
\alvecrestrarg{n-1} & \mbox{ if } \xi=0, n>1,\mbox{ and } m_n=1.
\end{array}\right.\]\index{$\alvec[\xi]$}
\end{enumerate}
\end{defi}

\begin{rmk} Note that $\alvec[\xi]$ might not be a tracking chain. This has to be verified when this notation is used.
In the case $\xi\in(0,\alcp{n,m_n})$ the second part of condition 2 of Definition \ref{kapparegdefi} has to be checked.
\end{rmk}

\begin{defi}[Extension of tracking chains, \boldmath$\me(\alvec)$\unboldmath]\label{tcextensiondefi}
An {\bf extension}\index{extension} of a tracking chain $\alvec$ is a tracking chain of which $\alvec$ is an initial chain.
A {\bf\boldmath $1$-step extension\unboldmath}\index{extension!$1$-step extension} is an extension by exactly one additional index.
The {\bf maximal extension} of $\alvec$ is denoted by $\me(\alvec)$ and is the tracking chain obtained from $\alvec$ after the maximum 
possible number of maximal $1$-step extensions according to Definition \ref{maxextdefi}.
\end{defi}

\begin{rmk} Note that $\mepl(\alvec)$ and $\me(\alvec)$, where $\alvec\in\TC$, either coincide or differ by one index extending 
$\me(\alvec)$ to $\mepl(\alvec)$ that does not satisfy condition 2 of Definition \ref{trackingchaindefi}.
\end{rmk}

The following \emph{key} lemma and proof cover Lemma 5.5 and Corollary 5.6 of \cite{CWc}. For details see the alternative formulation and
proof of Lemma 5.5 of \cite{CWc}, however, the proof given below should suffice.

\begin{lem}[5.5 and 5.6 of \cite{CWc}]\label{cmlmaxextcor}
Let $\alvec\in\TC$ be maximal, i.e.\ $\alvec=\me(\alvec)$, with maximal index pair $(n,m_n)$ and associated chain $\tauvec$.
Suppose that $\cml(\alvec)=:(i,j)$ exists.
\begin{enumerate} 
\item We have $(i,j+1)\klex(n,m_n)$, $\ec(\alvec)=\mepl(\alvec)$ exists, 
and the extending index with index pair $(n+1,1)$ is a successor multiple of $\taucp{i,j}$ with 
$(n+1)^\star=(i,j)$.
\item For $(k,l)\in\dom(\alvec)$ such that $(i,j+1)\klex(k,l)$ and $l<m_k$ if $m_k>1$ we have $\taucp{i,j}<\taucp{k,l}$ and
$\chi^\taucp{i,j}(\taucp{k,l})=1$.
\end{enumerate}
\end{lem}
{\bf Proof.}
We explain first, why $\ec(\alvecrestrarg{i,j+1})$ always exists and always is a tracking chain (see the beginning of the proof 
of Lemma 5.5 of \cite{CWc}). 
If $\taucp{i,j}<\taucp{i,j+1}\in\Ez$, then case 2.2.2 of Definition \ref{maxextdefi} applies, and the extension is clearly a tracking chain, otherwise case 2.3.1 applies. In this latter case
$\alvec$ is extended by $\alcp{i+1,1}=\rhoargs{\taucp{i,j}}{\taucp{i,j+1}}=\taucp{i,j}\cdot\la$ where $\la:=\log(\taucp{i,j+1})$ is a limit ordinal, since $\taucp{i,j+1}\in\Lz^{\ge\taucp{i,j}}$
due to the assumption $\chi^\taucp{i,j}(\taucp{i,j+1})=1$. Hence $\taucp{i+1,1}>\taucp{i,j}=\taucppr{i+1,1}$, implying that $\ec(\alvec)\in\TC$. According to Lemma \ref{chiinvlem} we have 
$\chi^\taucp{i,j}(\taucp{i,j+1})=\chi^\taucp{i,j}(\la)=\chi^\taucp{i,j}(\sumend(\la))=1$ and hence also $\chi^\taucp{i,j}(\taucp{i+1,1})=1$.

Now consider the subset $J_0$ of $\dom(\alvec)\cup\{(n+1,1)\}$ of (index pairs of) maximally extending indices starting with the 
extending index of $\alvecrestrarg{i,j+1}$, which is obtained from an application of either case 2.2.2 or case 2.3.1 as mentioned before. 
Observe that the process of maximal extension, shown to always terminate in Lemma \ref{meterminationlem}, can only end with a $\ka$-index.
We reduce $J_0$ to the set $J$ cancelling those (index pairs of) indices obtained from applications of cases 1.2 and 2.2.2 that 
are immediately followed by applications of cases 2.2.1 or 2.3.2, i.e.\ applications of the $\mu$-operator immediately before application 
of the $\la$-operator, cf.\ the proof of Lemma \ref{meterminationlem}. The index pairs specified in part 2 then comprise the set
$J\setminus\{(n+1,1)\}$.

According to Lemma \ref{chiinvlem}, part a), we have 
\[\chi^{\taucp{i,j}}(\taucp{k,l})=\chi^{\taucp{i,j}}(\taucp{i,j+1})=1\]
for every $(k,l)\in J$, cf.\ the proof of Lemma 5.5 of \cite{CWc} and an alternative argumentation 
\footnote{Corrections to be made on p.\ 60 of \cite{W18}: $\alvecpl:=\mepl(\alvecpr)=\mepl(\alvec)$ in line 8, and in line -11 after defining
$\bspr_{l+1}$, add: $\bs_{k+1}$ is then defined to be ${\bspr_{k+1}}^\frown\si_{k+1}$ if $\si_{k+1}\in\Ez^{>\sipr_{k+1}}$, otherwise
$\bs_{k+1}:=\bspr_{k+1}$.}
in Subsection 2.2 of \cite{W18}.
Note that in particular $\taucp{i,j}\le\taucp{k,l}$ for all $(k,l)\in J$, where equality holds if and only if $(k,l)=(n+1,1)$, which
settles part 2. 
Due to Lemma \ref{chiinvlem}, part b), early termination, i.e.\ $\taucppr{k,l}=\taucp{k,l}$ for $(k,l)\in J\setminus\{(n+1,1)\}$,
is not possible.
For the index $\alcp{n+1,1}$ maximally extending $\alvec$ to 
$\ec(\alvec)=\mepl(\alvec)$ we have $\sumend(\alcp{n+1,1})=\taucp{n+1,1}=\taucp{i,j}$, completing the proof of the lemma.
\qed

In the following definition we will provide a notion of 
\emph{reference sequence} that will replace the notion of \emph{characteristic sequence} in Definition 5.3 of \cite{CWc} and allow us to 
see the analogues of Lemmas 5.7, 5.8, and 5.10 of \cite{CWc} without reiterating similar arguments given in the respective proofs.

\begin{defi}[cf.\ 5.1, 5.3 of \cite{CWc}]\label{charseqdefi}
Let $\alvec=\left(\alevec,\ldots,\alnvec\right)\in\TC$, where $\alivec=(\alcp{i,1},\ldots,\alcp{i,m_i})$, with associated chain 
$\tauvec$, segmentation depth $p$, segmentation signature $(s_1,\ldots,s_p)$, and sequence of $\ups$-segments $(\la_i,t_i)_{0\le i\le p}$
as in Definition \ref{segmentationdefi}.
\begin{enumerate}
\item For $i\in\{1,\ldots,n\}$ and $j\in\{0,\ldots,m_i\}$ the {\bf reference sequence $\rsij(\alvec)$ of $\alvec$ at $(i,j)$}
\index{characteristic sequence}\index{$\rs$} is defined by
\[\rsij(\alvec):=\left\{\begin{array}{rl}
   \rsistar(\alvec)^\frown\sivec&\mbox{ if }i^\star\in\dom(\alvec)\\[2mm] 
   (\ups_{\la_l+1},\ldots,\ups_{\la_l+t_l})^\frown\sivec &\mbox{ if }i^\star=(s_l,0)\mbox{ for some }l\in\{1,\ldots,p\}\\[2mm]
   \sivec&\mbox{ otherwise,}
  \end{array}\right.\]
where $\sivec:=(\taucp{i,1},\ldots,\taucp{i,j})$.
\item For $(i,j)\in\dom(\alvec)$ the {\bf reference index pair $\refcp{i,j}(\alvec)$ of $\alvec$ at 
$(i,j)$}\index{index pair!reference index pair}\index{$\refcp{}$} is  
\[\refcp{i,j}(\alvec):=\left\{\begin{array}{ll}
   (i,j-1) &\mbox{ if } (i,j-1)\in\dom(\alvec)\cup\{(1,0)\}\\[2mm]
   \refcp{i-1,m_{i-1}}(\alvec) &\mbox{ otherwise.} 
  \end{array}\right.\] 
\item Setting $(i_0,j_0):=\refcp{i,j}(\alvec)$, the {\bf evaluation reference sequence $\ersij(\alvec)$ of $\alvec$ at $(i,j)\in\dom(\alvec)$}
\index{ers} is 
\[\ersij(\alvec):=\left\{\begin{array}{ll}
() & \mbox{ if }(i_0,j_0)=(1,0)\\[2mm]
\rs_{i_0,j_0}(\alvec) & \mbox{ otherwise.}
\end{array}\right.\] 
For convenience we also define $\ers_{1,0}(\alvec):=()$.
\end{enumerate}
\end{defi}

\begin{lem}[cf.\ 5.8 of \cite{CWc}]\label{rslem}
In the situation of the above definition the following statements hold.
\begin{enumerate}
\item For all $(i,j)\in\dom(\alvec)$ we have \[\rsij(\alvec)=\rsistar(\alvec)^\frown(\taucp{i,1},\ldots,\taucp{i,j}),\]
and if $j>1$ we have $\rsij(\alvec)=\rs_{i,j-1}(\alvec)^\frown\taucp{i,j}$. 
\item For all $i\in\{1,\ldots,n\}$ and $j\in\{0,\ldots,m_i-1\}$ we have \[\rsij(\alvec)\in\RS.\]
\item For $(i,j)\in\dom(\alvec)$ let $\gavec:=\rs_{i,j-1}(\alvec)$. Then we have 
\[\taucp{i,j}\in\left\{\begin{array}{l@{\quad}l}\domkvga & \mbox{ if }j=1\\ \domnuvga & \mbox{ if }j>1.\end{array}\right.\]
\item For $(i,j)\in\dom(\alvec)$ we have \[\ersij(\alvec)\in\RS.\]
\end{enumerate}
\end{lem}
{\bf Proof.} 
Part 1 follows directly from the definition. Note that if $i^\star=(j,0)$ for some $j$, then $j^\star=(j,0)$.

According to Lemma \ref{rserslem} we have $\taucp{i,1}\le\la_\tauistar$ whenever $\tauistar>1$.
If $\tauistar>1$ and $\taucp{i,1}\in\Ez^{>\tauistar}$, part c) of Lemma \ref{rhomulamestimlem} yields $\taucp{i,1}\le\mu_\tauistar$ as well.
If $1\le j<m_i$ we have $\taucp{i,1}\in\Ez^{>\tauistar}$, so proceeding from $i=1$ up to $i=n$ we obtain part 2.

Parts 3 and 4 then readily follow.
\qed

The evaluation reference sequence $\gavec:=\ersij(\alvec)$ for $(i,j)\in\dom(\alvec)$ is defined to be the sequence in $\RS$ that matches the
correct setting of relativization needed to evaluate the index $\alcp{i,j}$ in terms of $\ka^\gavec$ (if $j=1$) or $\nu^\gavec$ (if $j>1$),
as intended in the definition of tracking chains. In order to have the index $\alcp{i,j}$ in the domain of the corresponding $\ka$- or 
$\nu$-function, we have defined the notion of reference index pair $\refcp{i,j}(\alvec)$. Note that we will make use of the global
$\ka$-function where possible, cf.\ Definition \ref{kappanuprincipals}.
The well-definedness of $\tauticp{i,j}$ and $\alticp{i,j}$ in part 1 of the following definition is warranted by $\ka$- and $\nu$-regularity,
see Definitions \ref{critkapdefi}, \ref{kapparegdefi}, and \ref{nuregdefi}. 

\begin{defi}[cf.\ 5.9 of \cite{CWc}]\label{trchevaldefi}
Let $\alvec=\left(\alevec,\ldots,\alnvec\right)\in\TC$, where $\alivec=(\alcp{i,1},\ldots,\alcp{i,m_i})$, with associated chain 
$\tauvec$.
\begin{enumerate}
\item The {\bf evaluations $\tauticp{i,j}$ and $\alticp{i,j}$}\index{evaluation function!evaluation of tracking chains}
\index{$\tauticp{i,j}$,$\alticp{i,j}$} for $(i,j)\in\dom(\alvec)$
are defined as follows, setting $\varsivec:=\ersij(\alvec)$.
\[\tauticp{i,1}:=\ka^\varsivec_\taucp{i,1},\quad\alticp{i,1}:=\ka^\varsivec_\alcp{i,1}\quad\quad\mbox{ if }j=1,\]
\[\tauticp{i,j}:=\nu^\varsivec_\taucp{i,j},\quad\alticp{i,j}:=\nu^\varsivec_\alcp{i,j}\quad\quad\:\mbox{ if }j>1.\]
For $i\in\{1,\dots,n\}$, setting $i^\star=:(k,l)$ we define $\tauticp{k,l}:=\tauistar$ in case of $l=0$. For convenience we define
$\tauticp{i,0}:=\tauticp{i^\star}$ for $i=1,\ldots,n$.
\item The {\bf initial values $\set{\ordcp{i,j}(\alvec)}{(i,j)\in\dom(\alvec)}$ of $\alvec$}\index{initial value} are defined,
setting for convenience $m_0:=0$, $\ordcp{0,0}(\alvec):=0$, and $\ordcp{1,0}(\alvec):=0$, for $i=1,\ldots,n$ by
\[\ordcp{i,1}(\alvec):=\ordcp{i-1,m_{i-1}}(\alvec)+\alticp{i,1}\]
and\index{$\ordvalue$!$\ordcp{i,j}$}
\[\ordcp{i,j+1}(\alvec):=\ordcp{i,j}(\alvec)+(-\tauticp{i,j}+\alticp{i,j+1})\mbox{ for }1\le j<m_i.\]
We define the {\bf  value of $\alvec$} by $\ordvalue(\alvec):=\ordcp{n,m_n}(\alvec)$\index{$\ordvalue$!$\ordvalue$ for tracking chains}
which is the terminal initial value of $\alvec$. For $\al\in\On$ and $\alvec\in\TC$ such that $\al=\ordvalue(\alvec)$ we call $\alvec$
a {\bf tracking chain for $\alvec$}.
\index{value of a tracking chain}
\end{enumerate}
\end{defi}

\begin{rmk}[\cite{CWc}]
The correction $-\tauticp{i,j}$ in the above definition avoids double summation: Consider the easy example of the chain $((\epsn,1))$ which codes $\epsn\cdot2$. 
Notice that $-\tauticp{i,j}+\alticp{i,j+1}$ is always a non-zero multiple of $\tauticp{i,j}$. We clearly have $\ordcp{i,j}(\alvec)=\ordvalue(\alvecrestrarg{(i,j)})$.
\end{rmk}

The sequences given by $\rsij(\alvec)$ directly provide the (possibly pruned) tracking sequences required when 
evaluating indices of the associated chain $\tauvec$ in their proper setting of relativization, as specified in the following lemma.

\begin{lem}[cf.\ Lemma 5.10 of \cite{CWc}]\label{evallem} In the situation of Definition \ref{trchevaldefi},
\begin{enumerate}
\item $\tauticp{i,1}=\ka^{\varsivecstar}_\taucp{i,1}$ where $\varsivecstar:=\rsistar(\alvec)$, 
\item $\trs(\tauticp{i,j})=\rsij(\alvec)$ for $(i,j)\in\dom(\alvec)$ such that $\taucp{i,1}\in\Ez^{>\tauistar}$ (if $\tauistar>1$) and 
$\taucp{i,j}>1$ (if $j>1$), and 
\item $(i,j)=(k,l)$ for $(i,j),(k,l)\in\dom(\alvec)$ such that $\trs(\tauticp{i,j})=\trs(\tauticp{k,l})$, 
$\taucp{i,1}\in\Ez^{>\tauistar}$, $\taucp{k,1}\in\Ez^{>\taukstar}$ $\andsp$ $\taucp{i,j},\taucp{k,l}>1$.
\end{enumerate}
\end{lem}
{\bf Proof.} 
Part 1: 
Let us assume that $\alvec\not=((0))$, which is the only possibility for the trivial case $\taucp{i,1}=0$ (where $i=1$).
Lemma \ref{rserslem} warrants (in the nontrivial case where $\tauistar>1$) that $\taucp{i,1}\in\dom(\ka^{\varsivecstar})$.
Observe that in the case $i^\star=(j,0)$ for some $j$ we have $\ka^{\varsivecstar}_\taucp{i,1}=\ka_\taucp{i,1}$
according to the definition of the global $\ka$-function, Definition \ref{kappanuprincipals}, and the fact that 
\[\ka^{\bevec^\frown\be}_\be=\ka^\bevec_\be\mbox{ whenver } \bevec^\frown\be\in\RS.\]
{\bf Case 1:} $\ers_{i,1}(\alvec)=()$. This implies that $i^\star\in\dom(\alvec)$ is impossible. Hence $\tauticp{i,1}=\ka_\taucp{i,1}$,
as observed above.\\[2mm]
{\bf Case 2:} $\varsivec:=\ers_{i,1}(\alvec)=\rs_{k,l}(\alvec)$ where $(k,l)=\refcp{i,1}(\alvec)\in\dom(\alvec)$.\\[2mm]
{\bf Subcase 2.1:} $(r,s):=i^\star\in\dom(\alvec)$. Then we have $(r,s)\kglex(k,l)$ and are done if equality holds.
Otherwise we have $\taucp{r,s}\le\taucp{i,1}<\taucp{k,l}$ and $\rs_{r,s}(\alvec)$ is a proper initial sequence of $\rs_{k,l}(\alvec)$,
because $(r,s)$ is $\klex$-maximal candidate of a $\taucp{u,v}\le\taucp{i,1}$ where $(u,v)\in\dom(\alvec)$, $(u,v)\klex(i,1)$, and 
$v<m_u$. By definition we obtain $\tauticp{i,1}=\ka^\varsivec_\taucp{i,1}=\ka^\varsivecstar_\taucp{i,1}$.\\[2mm]
{\bf Subcase 2.2:} $i^\star=(j,0)$ for some $j$. Then as observed above $\ka^{\varsivecstar}_\taucp{i,1}=\ka_\taucp{i,1}$.\\[2mm]
{\bf 2.2.1:} $j\le k$. Then $\rsistar(\alvec)$ is a proper initial sequence of $\rs_{k,l}(\alvec)$, say
$\rs_{k,l}(\alvec)=\rsistar(\alvec)^\frown(\be_1,\ldots,\be_r)$ where $\be_r=\taucp{k,l}$ and $\taucp{i,1}<\be_1$. As in case 2.1 it 
follows that $\tauticp{i,1}=\ka^\varsivec_\taucp{i,1}=\ka^\varsivecstar_\taucp{i,1}$, which is equal to $\ka_\taucp{i,1}$.\\[2mm]
{\bf 2.2.2:} $k<j$. Then all components of $\varsivec$ are strictly greater than $\taucp{i,1}$, and we obtain 
$\tauticp{i,1}=\ka^\varsivec_\taucp{i,1}=\ka_\taucp{i,1}$. This concludes the proof of part 1.

Part 2 follows for $j=1$ from part 1 since in the case $\varsivecstar:=\rsistar(\alvec)\not=()$, due to the assumption 
$\taucp{i,1}\in\Ez^{>\tauistar}$, we have 
\[\ka^{\varsivecstar}_\taucp{i,1}=\nu^{\varsivecstar}_\taucp{i,1},\] so that according to Theorem \ref{trsvallem} 
$\trs(\tauticp{i,1})=\varsivecstar^\frown\taucp{i,1}=\rs_{i,1}(\alvec)$. In the case $j>1$, setting $\varsivec:=\ersij(\alvec)$,
which is equal to $\rs_{i,j-1}(\alvec)$, we have $\tauticp{i,j}=\nu^\varsivec_\taucp{i,j}$, so that again according to Theorem \ref{trsvallem}
$\trs(\tauticp{i,j})=\varsivec^\frown\taucp{i,j}=\rsij(\alvec)$.
Note that for $\taucp{i,j}=1$ we have $\tauticp{i,j}=\ordvalue(\varsivec)\cdot 2$, which is not in the domain of $\trs$.

Part 3 is seen by induction on the length of $\trs(\tauticp{i,j})$, using part 2, according to which the assumption implies that
$\rsij(\alvec)=\rs_{k,l}(\alvec)$. 
Let $(\la_q,t_q)_{0\le q\le p}$ be the sequence of $\ups$-segments of $\alvec$ according to Definition \ref{segmentationdefi},
$(\la_0,t_0)$ being the $\ups$-segment and $p$ the segmentation depth of $\alvec$.
We consider the following cases.\\[2mm]
{\bf Case 1:} $1<j\le m_i$ and $1<l\le m_k$. Then the i.h.\ yields $(i,j-1)=(k,l-1)$ and we are done.\\[2mm]
{\bf Case 2:} $1<j\le m_i$ and $l=1$. We show that this contradicts the assumption and is therefore impossible.\\[2mm]
{\bf Subcase 2.1:} $k^\star\in\dom(\alvec)$. Then we have $\rs_{i,j-1}(\alvec)=\rs_{k^\star}(\alvec)$, and the i.h.\ yields $(i,j-1)=k^\star$.
On the other hand, by assumption we have $\taucp{i,j}=\taucp{k,1}$, which implies that $(i,j)\kglex k^\star$. Contradiction. 
\\[2mm]
{\bf Subcase 2.2:} $k^\star\not\in\dom(\alvec)$. In this situation $k^\star=(1,0)$ is not possible since according to the assumption,
$\rsij(\alvec)$ and $\rs_{k,l}(\alvec)$ have the same lenth. Therefore, $\rs_{k,1}(\alvec)$ is of the form 
$(\ups_{\la_r+1},\ldots,\ups_{\la_r+t_r},\taucp{k,1})$ where $r\in\{1,\ldots,p\}$, $\la_r+t_r<\la_0$, $\taukstar=1+\ups_{\la_r+t_r}$, and 
$\taucp{k,1}\in\Ez^{>\taukstar}$, and hence $\rs_{i,j-1}(\alvec)=(\ups_{\la_r+1},\ldots,\ups_{\la_r+t_r})$, which is impossible since 
$(i,j-1)\in\dom(\alvec)$.
\\[2mm]
{\bf Case 3:} $j=1$ and $1<l\le m_k$. We argue as in Case 2 to show that this case does not occur.\\[2mm]
{\bf Case 4:} $j=1$ and $l=1$. Then we have $\taucp{i,1}=\taucp{k,1}$ and need to show that $i=k$.\\[2mm]
{\bf Subcase 4.1:} $i^\star,k^\star\in\dom(\alvec)$. Then we have $\trs(\tauticp{i^\star})=\trs(\tauticp{k^\star})$, and hence by the i.h.\
$i^\star=k^\star$. Let us assume to the contrary that $i<k$. This implies $(i,1)\kglex k^\star$, since due to
$\ka$-regularity of $\alvec$ there exists $h\in\{i,\ldots,k-1\}$ such that $m_h>1$, hence $i^\star\klex k^\star$ contradicting the i.h. 
In the same way we see that $k<i$ cannot hold either.\\[2mm]
{\bf Subcase 4.2:} $i^\star,k^\star\not\in\dom(\alvec)$.
In this case $\rs_{i,1}(\alvec)=\rs_{k,1}(\alvec)$ is either equal to $(\taucp{i,1})$ or of the form 
$(\ups_{\la_r+1},\ldots,\ups_{\la_r+t_r},\taucp{i,1})$ where $r\in\{1,\ldots,p\}$. Assuming that $i<k$ we would either have
$m_i=\ldots=m_{k-1}=1$, which implies that $\taucp{i,1}=\ldots=\taucp{k,1}\in\Ez^{>\taucp{i,1}}$, contradicting $\taucp{i+1,1}<\taucp{i,1}$
by $\ka$-regularity, or there would exist the least $h\in\{i,\ldots,k-1\}$ such that $m_h>1$, so that we would obtain
$\taucp{k,1}=\taucp{i,1}\ge\ldots\ge\taucp{h,1}$ and hence $(h,1)\kglex k^\star$, which implies $k^\star\in\dom(\alvec)$ contradicting
our assumption. 
\\[2mm]
{\bf Subcase 4.3:} $i^\star\in\dom(\alvec)$ and $k^\star\not\in\dom(\alvec)$. Then $k^\star=(1,0)$ is impossible (as in Subcase 2.2),
and $\rs_{i,1}(\alvec)$ is of a form $(\ups_{\la_r+1},\ldots,\ups_{\la_r+t_r},\taucp{i,1})$ where $r\in\{1,\ldots,p\}$. Hence
$\rs_{i^\star}(\alvec)=(\ups_{\la_r+1},\ldots,\ups_{\la_r+t_r})$ while $i^\star\in\dom(\alvec)$, which is impossible as in Subcase 2.2.
\\[2mm]
{\bf Subcase 4.4:} $i^\star\not\in\dom(\alvec)$ and $k^\star\in\dom(\alvec)$. Then we argue as in Subcase 4.3.
\qed

Our intermediate goal is to establish an order isomorphism between tracking chains and their evaluations.
In order to formulate estimates on the values of tracking chains that will allow us to reach that goal we introduce a notion of 
\emph{depth} of a tracking chain by application of the function $\dpf$ introduced in Definition \ref{kappadpdefi} in accordance with the
concept of tracking chains. 

\begin{defi}[cf.\ 5.11 of \cite{CWc}]\label{dpfdefi}
Let $\alvec=\left(\alevec,\ldots,\alnvec\right)\in\TC$ where $\alivec=(\alcp{i,1},\ldots,\alcp{i,m_i})$ for $i=1,\ldots,n$, with associated
chain $\tauvec$. 
$\dpf(\alvec)$ is defined as follows. Let $\tau:=\taucp{n,m_n}$ and $\taupr:=\taunpr$. 
Let further $\tautipr$ be the evaluation of $\taupr$.
We set $\varsivec:=\rs_{n,m_n-1}(\alvec)$
and define
\[\dpf(\alvec):=\left\{\begin{array}{l@{\quad}l}
\dpf_\varsivec(\tau) & \mbox{ if } m_n=1\\[2mm]
\ka^\varsivec_{\rhoargs{ }{\tau}}+\dpf_\varsivec(\rhoargs{ }{\tau})+\chicheck^\taupr(\tau)\cdot\tautipr & 
\mbox{ if } m_n>1\andsp\tau<\mu_\taupr\\[2mm]
\ka^\varsivec_{\la_\taupr}+\dpf_\varsivec(\la_\taupr) & \mbox{ if } m_n>1\andsp\tau=\mu_\taupr. 
\end{array}\right.\]
\end{defi}

Note that $\taupr$ is equal to $\taucp{n,m_n-1}$ if $m_n>1$ and $\taunstar$ if $m_n=1$, that 
$\tautipr$ is equal to $\tauticp{n,m_n-1}$ if $m_n>1$ and $\tauticp{n^\star}$ if $m_n=1$, and that
$\varsivec$ as defined above is equal to $\trs(\tautipr)$ unless $\taunstar=1+\ups_{\la_r}$ for some $r\le p$ (assuming 
$(\la_i,t_i)_{0\le i\le p}$ to be the sequence of $\ups$-segments of $\alvec$), in which case we have $\varsivec=()$.
$\varsivec$ is the setting of relativization for the $\ka$-function enumerating $\al$-$\leo$-minimal ordinals, where 
$\al:=\ordvalue(\alvec)$.

The technical Lemma 5.12 of \cite{CWc} carries over for tracking chains $\alvec$ such that $\ov(\alvec)\not\in\Image(\ups)$, which
is just what we need here. For reader's convenience we present the content of Lemma 5.12 of \cite{CWc} in a way that may be easier to follow.
Each part is formulated separately and basically contains an observation about the interplay of (maximal) extensions of tracking chains
and the auxiliary function $\dpf$ that is useful in the context of order-isomorphisms of tracking chains and ordinals, in particular in
the next subsection, Lemma \ref{tcassgnmntlem} and its corollary, where we provide an assignment of tracking chains to ordinals. 
Basically, the consequences stated in the following lemma rest on the 
strict monotonicity of $\ka$- and $\nu$-functions verified in Corollary \ref{kappanuhzcor}.

\begin{lem}[cf.\ Lemma 5.12 of \cite{CWc}]\label{tcevalestimlem} Assume the settings of Definition \ref{dpfdefi}.
\begin{enumerate}
\item In case of $m_n>1$ and $\tau<\mu_\taupr$ we have
\[\ordvalue(\alvec)+\dpf(\alvec)=\ordvalue(\alvec[\alcp{n,m_n}+1]).\]
\item $\dpf(\alvec)=0$ if and only if there does not exist any proper extension of $\alvec$.
\item If $\ec(\alvec)$ exists, but $\ec(\alvec)\not\in\TC$, then 
\[\taucp{n+1,1}=\taucp{\cml(\alvec)}\quad\mbox{ and }\quad\sumend(\dpf(\alvec))=\tauticp{\cml(\alvec)},\]
where $\taucp{n+1,1}=\sumend(\alcp{n+1,1})$ and $\ec(\alvec)=\alvec^\frown(\alcp{n+1,1})$.
\item If $\alvecpl:=\ec(\alvec)\in\TC$ exists, there are the following three cases to consider.
\begin{enumerate}
\item $m_n>1$, $\tau<\mu_\taupr$, and $\chi^\taupr(\tau)=0$. Then 
\[\ordvalue(\alvecpl)+\dpf(\alvecpl)<\ordvalue(\alvec^\frown(\rhoargs{ }{\tau}+1))<\al +\dpf(\alvec).\]
\item $m_n>1$, $\tau<\mu_\taupr$, and $\chi^\taupr(\tau)=1$. Then $\cml(\me(\alvec))=(n,m_n-1)$, $\mepl(\alvec)\not\in\TC$, and
\[\ordvalue(\alvecpl)+\dpf(\alvecpl)=\al +\dpf(\alvec).\]
\item Otherwise. Then we again have
\[\ordvalue(\alvecpl)+\dpf(\alvecpl)=\al +\dpf(\alvec).\]
\end{enumerate}
\item For a non-maximal $1$-step extension $\alvecpl\not=\ec(\alvec)$ of $\alvec$ according to Definition \ref{tcextensiondefi}, 
$\alvecpl$ is of a form either 
\[\alvecpl=\alvec^\frown(\alcp{n+1,1})\quad\mbox{ or }\quad \alvecpl=(\alevec,\ldots,\alnminvec,\alnvec^\frown\alcp{n,m_n+1}),\]
and we set $\alcp{n,m_n+1}:=0$ if $\alvecpl$ is of the former, and $\alcp{n+1,1}:=0$ if $\alvecpl$ is of the latter form. 
We define 
\[\alvecpr:=\left\{\begin{array}{l@{\quad}l}
\alvec^\frown(\alcp{n+1,1}+1) & \mbox{ if } \alcp{n,m_n+1}=0\\[2mm]
(\alevec,\ldots,\alnminvec,\alnvec^\frown(\alcp{n,m_n+1}+1)) & \mbox { if } \alcp{n+1,1}=0,
\end{array}\right.\]
and consider two cases.
\begin{enumerate}
\item $\alvecpr\not\in\TC$: In this case we have $m_n>1$, $\tau=\mu_\taupr\in\Ez\cap(\taupr,\la_\taupr)$, $\alcp{n,m_n+1}=\mu_\tau$, and
\[\ordvalue(\alvecpl)+\dpf(\alvecpl)<\ordvalue(\alvec^\frown(\tau+1))\le\al +\dpf(\alvec).\]
\item $\alvecpr\in\TC$: Then we have
\[\ordvalue(\alvecpl)+\dpf(\alvecpl)\le\ordvalue(\alvecpr)\le\al +\dpf(\alvec)\quad\mbox{and}\quad\ordvalue(\alvecpl)+\dpf(\alvecpl)<\al +\dpf(\alvec).\]
\end{enumerate}
\item For any extension $\bevec$ of $\alvec$ we have
\begin{enumerate}
\item $\ordvalue(\bevec)+\dpf(\bevec)\le\al +\dpf(\alvec)$ and
\item $\ordvalue(\bevec)<\al +\dpf(\alvec)$ if $m_n>1$ and $\tau<\mu_\taupr$.
\end{enumerate}
\item The ordinal $\ordvalue(\me(\alvec))+\dpf(\me(\alvec))$ is calculated in the following three scenarios.
\begin{enumerate}
\item If $m_n>1$, $\tau<\mu_\taupr$, and $\chi^\taupr(\tau)=1$ we have \[\ordvalue(\me(\alvec))+\dpf(\me(\alvec))=\al+\dpf(\alvec)=\ordvalue(\alvec[\alcp{n,m_n}+1]).\]
\item If $\alvec$ does not possess a critical main line index pair $\cml(\alvec)$ then $\dpf(\me(\alvec))=0$ and
\[\ordvalue(\me(\alvec))=\left\{\begin{array}{l@{\quad}l}
\al + \dpf_\varsivec(\tau)&\mbox{ if }m_n=1\\[2mm]
\al + \ka^\varsivec_{\rhoargs{\taupr}{\tau}}+\dpf_\varsivec(\rhoargs{\taupr}{\tau})&\mbox{ if }m_n>1\mbox{ and }\tau<\mu_\taupr\\[2mm]
\al + \ka^\varsivec_{\la_\taupr}+\dpf_\varsivec(\la_\taupr)&\mbox{ otherwise,}
\end{array}\right.\]
which only deviates from $\al+\dpf(\alvec)$ in the case $m_n>1\andsp\tau<\mu_\taupr$.
\item If $\cml(\alvec)=:(i,j)$ exists then
\begin{eqnarray*}
\ordvalue(\me(\alvec))+\dpf(\me(\alvec))&=&\al+\dpf(\alvec)\\[2mm]
&=&\left\{\begin{array}{l@{\quad}l}
\al + \dpf_\varsivec(\tau)&\mbox{ if }m_n=1\\[2mm]
\al + \ka^\varsivec_{\rhoargs{\taupr}{\tau}}+\dpf_\varsivec(\rhoargs{\taupr}{\tau})&\mbox{ if }(n,m_n)=(i,j+1)\\[2mm]
\al + \ka^\varsivec_{\la_\taupr}+\dpf_\varsivec(\la_\taupr)&\mbox{ otherwise}
\end{array}\right.\\[2mm]
&=&\ordvalue(\alvecrestrarg{(i,j+1)}[\alcp{i,j+1}+1]),
\end{eqnarray*}
and 
\[\dpf(\me(\alvec))=\ka^\varsivecpr_{\taucp{i,j}(\xi+1)}\] where, say, $(r,k_r)$ is the $\klex$-greatest index pair of $\me(\alvec)$, 
$\varsivecpr:=\rs_{r,k_r-1}(\me(\alvec))$, 
and $\taucp{i,j}(\xi+1)$ for suitable $\xi$ is the extending index of $\ec(\me(\alvec))$.
\end{enumerate}
\end{enumerate}
\end{lem}
{\bf Proof.}
Part 1 follows from the definition of the $\nu$-functions.
Part 2 is a direct consequence of the definitions of tracking chain and $\dpf$.
Part 3 is a consequence of Lemma \ref{cmlmaxextcor}, using Lemma \ref{evallem}, 
as carried out in detail in the beginning of the proof of Lemma 5.12 of \cite{CWc}. 
Part 4 follows from the definitions involved with the aid of the monotonicity of $\ka$- and $\nu$-functions, 
and Lemma \ref{cmlmaxextcor} for part 4(b).
Part 5 is shown as in the proof of Lemma 5.12 of \cite{CWc} using monotonicity of $\ka$- and $\nu$-functions.
Part 6 is shown using the former parts regarding $1$-step extensions by induction on the number of $1$-step extensions that $\bevec$ 
results from (rather than the induction along $\klex$ on $\cspr(\alvec)$ as was stated in the proof of Lemma 5.12 of \cite{CWc}).
The first equality of part 7(a) follows from parts 6(a) and 4(b) and (c) where each extension step is chosen maximally. 
In part 7(c) the last stated equality 
was already shown by part 1 of Lemma \ref{cmlmaxextcor}.
Details of the verification of part 7 are given in the proof of Lemma 5.12 of \cite{CWc}. 
\qed

\begin{cor}\label{nuimagecor} Let $\alvec\in\RS$, $\alvec=(\ale,\ldots,\aln)$ where $n>0$, and $\al:=\ordvalue(\alvec)$. 
The image of $\nu^\alvec$ consists of multiples of $\al$. For $\be\in\dom(\nu^\alvec)$ we have $\nu^\alvec_\be\in\Hz^{>\al}$ if and only if
$\be\in\Hz^{>1}$.
\end{cor}
{\bf Proof.} Part 7(c) of the above Lemma \ref{tcevalestimlem} shows that for the only case in question, that is, $\ga\in\dom(\nu^\alvec)$ 
such that $\chi^{\aln}(\ga)=1$ the last summand of $\nu^\alvec_{\ga+1}$ is $\al$, which is seen as follows. 
Let $\bevec:=\me((\alvec^\frown\ga))$ with $\klex$-greatest index pair $(r,k_r)$ and set $\varsivecpr:=\rs_{r,k_r-1}(\bevec)$. 
The last summand of $\nu^\alvec_{\ga+1}$ is then $\ka^\varsivecpr_\aln=\ka^\alvec_\aln=\al$, since $\alvec\subseteq\varsivecpr$ by
Lemma \ref{cmlmaxextcor}. The second claim now readily follows. Note that $\nu^\alvec_0=\al$ and $\nu^\alvec_1=\al\cdot 2$.
\qed

\begin{defi}[5.13 of \cite{CWc}] We define a linear ordering $\ktc$\index{$<_tc$@$\ktc$,$\letc$} on $\TC$ as follows.
Let $\alvec, \bevec\in\TC$ be given, say, of the form
\[\alvec=((\alcp{1,1},\ldots,\alcp{1,m_1}),\ldots,(\alcp{n,1},\ldots,\alcp{n,m_n}))\]
and
\[\bevec=((\becp{1,1},\ldots,\becp{1,k_1}),\ldots,(\becp{l,1},\ldots,\becp{l,k_l})).\]
Let $(i,j)$ where $1\le i\le\min\singleton{n,l}$ and $1\le j\le\min\singleton{m_i,k_i}$ be
$\klex$-maxi\-mal such that $\alvecrestrarg{i,j}=\bevecrestrarg{i,j}$, if that exists,
and $(i,j):=(1,0)$ otherwise. 
\begin{eqnarray*}
\alvec\ktc\bevec\quad:\aeq&\quad&(i,j)=(n,m_n)\not=(l,k_l)\\
&\;\;\vee&(j<\min\singleton{m_i,k_i}\andsp\alcp{i,j+1}<\becp{i,j+1})\\
&\;\;\vee&(j=m_i<k_i\andsp i<n\andsp\alcp{i+1,1}<\taucp{i,j})\\
&\;\;\vee&(j=k_i<m_i\andsp i<l\andsp\taucp{i,j}<\becp{i+1,1})\\
&\;\;\vee&(j=k_i=m_i\andsp i<\min\singleton{n,l}\andsp\alcp{i+1,1}<\becp{i+1,1})\\[1ex]
\alvec\letc\bevec\quad:\aeq&&\alvec\ktc\bevec\veesp\alvec=\bevec.
\end{eqnarray*}
\end{defi}

Note that in the above definition, the first disjunction term covers the case where $\alvec$ is a proper initial chain of $\bevec$,
the second covers the situation where at $(i,j+1)$ the tracking chain $\bevec$ branches into a component with larger $\nu$-index, and 
the third, fourth, and fifth disjunction term cover the situation where $\bevec$ branches into a component with larger $\ka$-index.
In this latter situation the third term applies to the scenario where the branching $\ka$-index of $\bevec$ is not given explicitly but rather
by continuation with a $\nu$-index, whereas in the fourth term the lower $\ka$-index of $\alvec$ is not given explicitly but by continuation
with a $\nu$-index.   

\begin{lem}[5.14 of \cite{CWc}]\label{tcorderisolem} For all $\alvec,\bevec\in\TC$ we have
\[\alvec\ktc\bevec\quad\aeq\quad\ordvalue(\alvec)<\ordvalue(\bevec).\]
\end{lem}
{\bf Proof.} We first observe that $(\TC,\ktc)$ is a linear ordering.
The lemma then follows from the definitions of $\ktc$ and $\ordvalue$ using the strict monotonicity of $\ka$- and $\nu$-functions,
Corollary \ref{kappanuhzcor}, and Lemma \ref{tcevalestimlem}, by verifying that $\ordvalue(\alvec)<\ordvalue(\bevec)$ whenever 
$\alvec,\bevec\in\TC$ such that $\alvec\ktc\bevec$.
\qed

\begin{cor}[cf.\ 5.15 of \cite{CWc}]
For any $\al\in\On$ there exists at most one tracking chain for $\al$.\qed
\end{cor}

\subsection{Assignment of tracking chains to ordinals}\label{tcassignmentsubsec}
We will obtain an order isomorphism between $(\On,<)$ and $(\TC,\ktc)$ once we extend the inverse function $\tc$ from Definition 6.1 of
\cite{CWc} to an assignment of tracking chains to all ordinals. With the following adaptation of the notion of relative tracking sequence 
the extension of the assignment of tracking chains to ordinals from $\ups_1$ to all ordinals can be formulated conveniently. 

\begin{defi}[cf.\ Definition 4.16 of \cite{CWc}]\label{relchaindefi}
Let $\tauvec=(\tau_1,\ldots,\tau_n)\in\RS$ and $\be\in\Hz$, $\be\le\ka^\tauvec_{\la_{\tau_n}}+\dpf_\tauvec(\la_{\tau_n})$ if $n>0$. 
Denote the $\ups$-segment of $\be$ by $(\xi,u):=\upsseg(\be)$ and set $\tau_0:=\ups_{\xi+u}=:\tauti_0$. 
Let $k\in\singleton{0,\ldots,n}$ be maximal such that $\ordvalue(\tauvecrestrk)\le\be$, and set $\tauti_k:=\ordvalue(\tauvecrestrk)$ if $k>0$.
The {\bf tracking sequence $\trs[\tauvec](\be)$ of $\be$ relative to the reference sequence $\tauvec$} is defined by
\index{tracking sequence!tracking sequence of $\be$ relative to $\al$}\index{$\trst$!$\trsrelal$}
\[\trs[\tauvec](\be):=\left\{\begin{array}{l@{\quad}l}
    \trs^{\tau_k}(\tau_k\cdot(1/\tauti_k)\cdot\be)&\mbox{ if } k>0\\[2mm]
    \trs^{\tau_0}(\be)&\mbox{ if } k=0. \end{array}\right.
\]
\end{defi}
{\bf Remark.} $\trs[\tauvec]$ aims at a tracking sequence with starting point $\tauti_k$ instead of $0$.
In the above situation for $\trs[\tauvec](\be)$ to make sense, i.e.\ to be related to $\ordvalue(\tauvec)$, we should have 
$\be_1\le\la^{\tau_{k-1}}_{\tau_k}$ in case of $k>0$, where $\be_1$ is the first element of $\trs[\tauvec](\be)$.
It is easy to see (using Lemmas \ref{rhomulamestimlem}, \ref{citedinjtrslem}, and \ref{reltrsvallem}) that this holds if $k\in(0,n)$. 
However, in case of $k=n>0$ this holds if and only if $\be\le\ka^\tauvec_{\la_{\tau_n}}+\dpf_\tauvec(\la_{\tau_n})$ as follows from 
Lemmas \ref{lamveclem} and \ref{reltrsvallem}. Note that $\trs[()]$ and $\trs$ are not equal. 
 
\begin{lem}[cf.\ 4.17 of \cite{CWc}]\label{reltrslem}
Let $\tauvec=(\tau_1,\ldots,\tau_n)\in\RS$ and $\be,\ga\in\Hz$ such that $\be,\ga\le\ka^\tauvec_{\la_{\tau_n}}+\dpf_\tauvec(\la_{\tau_n})$ 
if $n>0$. Set $\tau_{n+1}:=\la_{\tau_n}+1$ if $n>0$ and $\tau_1:=\ups_{\xi+u+1}$ where $(\xi,u)=\upsseg(\be)$ otherwise.
\begin{enumerate}
\item[a)] With $k$ as in the above definition we have 
\[\tau_k\le\be_1<\tau_{k+1}\]
where $\be_1$ is the first element of $\trs[\tauvec](\be)$.
\item[b)] If $\be<\ga$ then
\[\trs[\tauvec](\be)\klex\trs[\tauvec](\ga).\]
\end{enumerate}
\end{lem}
{\bf Proof.} The lemma is proved by application of Lemmas \ref{citedinjtrslem} and \ref{reltrsvallem}, 
using part a) to show part b).
\qed

Recall part 4 of Lemma \ref{rslem}, according to which $\ersij(\alvec)\in\RS$ for all $(i,j)\in\dom(\alvec)$ where $\alvec\in\TC$. 
The assignment of tracking chains to ordinals given below is based on tracking sequences relativized to such evaluation reference sequences.

\begin{defi}[cf.\ 6.1 of \cite{CWc}]\label{tcassignmentdefi}
For $\al\in\On$ we define the {\bf\boldmath  tracking chain assigned to\unboldmath}
\index{tracking chain!tracking chain assigned to $\al$} $\al$, $\tc(\al)$\index{$\tc$}, by recursion on the length of the additive 
decomposition of $\al$.
We define $\tc(0):=((0))$, and if $\al\in\Hz$ we set $\tc(\al):=(\trs(\al))$. Now suppose $\tc(\al)=\alvec=(\alevec,\ldots,\alnvec)$ to be the
tracking chain already assigned to some $\al>0$, where $\alivec=(\alcp{i,1},\ldots,\alcp{i,m_i})$ for
$1\le i\le n$, with associated chain $\tauvec$, $(\la,t):=\upsseg(\alvec)$, and segmentation parameters $p,s_l,(\la_l,t_l)$ for $l=1,\ldots,p$
as in Definition \ref{segmentationdefi}, and let $\be\in\Hz$, $\be\le\sumend(\al)$.
For technical reasons, we set $\alcp{n+1,1}:=0$ and $m_{n+1}:=1$.
The definition of $\tc(\al+\be)$, the tracking chain assigned to $\al+\be$,
requires the following preparations.
\begin{itemize}
\item For $(i,j)\in\dom(\alvec)$ let
    \[(\beucp{i,j}_1,\ldots,\beucp{i,j}_{r_{i,j}}):=\trs[\ersij(\alvec)](\be).\]
For the tracking chain of $\be$ (according to Definition \ref{latrsdefi}) let \[(\be_1,\ldots,\be_r):=\trs(\be).\]
\item Let $(i_0,j_0)$, where $1\le i_0\le n$ and $1\le j_0< m_{i_0}$, be $\klex$-maximal with \[\alcp{i_0,j_0+1}<\mu_\taucp{i_0,j_0}\] if that exists, otherwise set $(i_0,j_0):=(1,0)$.
\item Let $(k_0,l_0)$ be either $(1,0)$ or satisfy $1\le k_0\le n+1$ and $1\le l_0\le m_{k_0}$, so that $(k_0,l_0)$ is $\klex$-minimal with 
$(i_0,j_0)\kglex(k_0,l_0)$ and 
\begin{enumerate}
\item for all $k\in\singleton{k_0,\ldots,n}$ we have
\[\alcp{k+1,1}+\beucp{k,m_k}_1\ge\rho_k\]
\item for all $k\in\singleton{k_0,\ldots,n}$ and all $l\in\singleton{1,\ldots,m_k-2}$ such that $(k_0,l_0)\klex(k,l)$ we have
\[\taucp{k,l+1}+\beucp{k,l+1}_1>\la_\taucp{k,l}.\]
\end{enumerate}
\end{itemize}
{\bf Case 0:} $(k_0,l_0)=(1,0)$.
Then $\tc(\al+\be)$ is defined by
       \[((\alcp{1,1}+\beucp{1,1}_1,\beucp{1,1}_2,\ldots,\beucp{1,1}_\rcp{1,1})).\]
{\bf Case 1:} $(i_0,j_0)=(k_0,l_0)\in\dom(\alvec)$. Then we set $(i,j):=(i_0,j_0+1)$ and consider three subcases:
\begin{itemize}
\item[{\bf 1.1:}] $\be<\tauticp{i_0,j_0}$. Then $\tc(\al+\be)$ is defined to be
    \[\alvecrestrarg{i,j}^\frown\left(\rhoargs{\taucp{i,j-1}}{\taucp{i,j}}+\beucp{i,j}_1,\beucp{i,j}_2,\ldots,\beucp{i,j}_\rcp{i,j}\right).\]
\item[{\bf 1.2:}] $\be=\tauticp{i_0,j_0}$. Then $\tc(\al+\be)$ is defined by 
       \[\alvecrestrarg{i,j}[\alcp{i,j}+1].\]
\item[{\bf 1.3:}] $\be>\tauticp{i_0,j_0}$. 
Then there is an $r_0\in(0,r)$ such that $\be_{r_0}=\taucp{i_0,j_0}=\taucp{i,j-1}$, and $\tc(\al+\be)$ is defined by
       \[\alvecrestrarg{i-1}^\frown(\alcp{i,1},\ldots,\alcp{i,j-1},
                \alcp{i,j}+\be_{r_0+1},\be_{r_0+2},\ldots,\be_r).\]
\end{itemize}
{\bf Case 2:} $(i_0,j_0)\klex(k_0,l_0)$. Then there are the following subcases:

\begin{itemize}
\item[{\bf 2.1:}] $k_0=n+1$ and $\beucp{n,m_n}_1=\taucp{n,m_n}\in\Ez^{>\taunpr}$. Then $\be=\tauticp{n,m_n}$, and $\tc(\al+\be)$ is defined by
        \[\alvecrestrarg{n-1}^\frown(\alcp{n,1},\ldots,\alcp{n,m_n},1).\]
\item[{\bf 2.2:}] $k_0\le n$, $l_0\in\singleton{1,\ldots,m_{k_0}-2}$ and $\taucp{k,l}+\beucp{k,l}_1\le\la_\taucp{k,l-1}$ for 
$(k,l):=(k_0,l_0+1)$. Then we define $\tc(\al+\be)$ by
 \[\alvecrestrarg{k,l}^\frown\left(\taucp{k,l}+\beucp{k,l}_1,\beucp{k,l}_2,\ldots,\beucp{k,l}_\rcp{k,l}\right),\]
 provided this vector satisfies condition 2 of Definition \ref{trackingchaindefi}, otherwise we have $(i_0,j_0)\in\dom(\alvec)$,
 $\be=\tauticp{i_0,j_0}$, and $\tc(\al+\be)$ is defined as in case 1.2.
\item[{\bf 2.3:}] Otherwise. Then $k_0>i_0$, $l_0=1$, and $\alcp{k+1,1}+\beucp{k,m_k}_1<\rho_k$ for $k:=k_0-1$, and $\tc(\al+\be)$ is defined by
\[\alvecrestrarg{k}^\frown\left(\alcp{k+1,1}+\beucp{k,m_k}_1,\beucp{k,m_k}_2,\ldots,\beucp{k,m_k}_\rcp{k,m_k}\right),\]
 provided this vector satisfies condition 2 of Definition \ref{trackingchaindefi}, otherwise we have $(i_0,j_0)\in\dom(\alvec)$,
 $\be=\tauticp{i_0,j_0}$, and $\tc(\al+\be)$ is defined as in case 1.2.
\end{itemize}
\end{defi}
\begin{rmk} 
Note that compared to \cite{CWc} we have changed indexing of the $\beucp{i,j}$-sequences:  $\beucp{i,j}$ here corresponds to  $\beucp{i,j-1}$ in \cite{CWc}. Regarding part 1 of the definition of $(k_0,l_0)$, note that for $0<k<n$ by definition of $\ers$ we have 
$\ers_{k+1,1}(\alvec)=\ers_{k,m_k}(\alvec)$. Writing more intuitively $\ers_{k+1,1}(\alvec)$ instead of $\ers_{k,m_k}(\alvec)$ would require
to set $\ers_{n+1,1}(\alvec)$ to $\ers_{n,m_n}(\alvec)$.

Case 0 corresponds to the scenario in which adding $\be$ to $\al$ means to jump into a larger ($\ups_\la$-)$\leo$-connectivity
component and used to be included in Case 1.3 in \cite{CWc}. Since in our more general setting here we no longer have $\tauticp{1,0}=0$, 
as was the case in \cite{CWc}, we decided to cover this case separately. Note that in Case 0 we use the tracking chain 
$\beucp{1,1}=\trs^{\ups_{\xi+u}}(\be)$ where $(\xi,u)=\upsseg(\be)$.
Case 1.3 now solely covers the situation of jumping into a larger (non-trivial) $\letwo$-connectivity component on the surrounding main line.
In Case 1.3 we use $\trs(\be)$ for the assignment.  
Case 2.1 takes care of condition 2 of $\ka$-regularity, Definition \ref{kapparegdefi}. 
\end{rmk}

\begin{lem}[cf.\ 6.2 of \cite{CWc}]\label{tcassgnmntlem} 
Let $\al\in\On$.
\begin{enumerate}
\item[a)] $\tc(\al)\in\TC$, i.e.\ $\tc(\al)$ meets all conditions of Definition \ref{trackingchaindefi}.
\item[b)] There exists exactly one tracking chain for $\al$, namely $\tc(\al)$ satisfies $\ordvalue(\tc(\al))=\al$.
\end{enumerate}
\end{lem}
{\bf Proof.} By adaptation of the proof of Lemma 6.2 of \cite{CWc}, including several corrections. Note that the proof in \cite{CWc} actually proceeds
by induction on the length of the additive decomposition of $\al$, rather than by induction on $\al$, as was stated there. 
The proof extensively utilizes Lemma \ref{tcevalestimlem} and the monotonicity of $\ka$- and $\nu$-functions.
The case $\al=0$ is trivial, and using Lemma \ref{trsvallem} we see that the claims hold whenever $\al\in\Hz$.
Now suppose the claims have been shown for some $\al>0$ with assigned tracking chain $\tc(\al)=\alvec$ as in the definition and suppose 
$\be\le\sumend(\al)$. We adopt the terminology of the previous definition and commence proving the inductive step for $\al+\be$ by showing 
the following preparatory claims.

\begin{claim}\label{tcassgnclmone} If $(i_0,j_0)\not=(1,0)$ and $\be\le\tauticp{i_0,j_0}$ then $\beucp{i_0,j_0+1}_1\le\taucp{i_0,j_0}$.
Equality holds if and only if $\be=\tauticp{i_0,j_0}$.
\end{claim}
In order to show the claim let us assume that $(i_0,j_0)\not=(1,0)$ and $\be\le\tauticp{i_0,j_0}$.
In the case $\beucp{i_0,j_0+1}_1=\taucp{i_0,j_0}$ the assumption implies $r_{i_0,j_0+1}=1$ and  $\be=\tauticp{i_0,j_0}$.
On the other hand, in case of $\be=\tauticp{i_0,j_0}$ we clearly have 
$\trs[\ers_{i_0,j_0+1}(\alvec)](\be)=(\taucp{i_0,j_0})$.  

Now assume $\be<\tauticp{i_0,j_0}$. Write $\sivec=(\si_1,\ldots,\si_s)$ for the sequence $\trs(\tauticp{i_0,j_0})$, which is equal to the 
reference sequence $\rs_{i_0,j_0}(\alvec)$, so that $\tauticp{i_0,j_0}=\ordvalue(\sivec)$ and $\taucp{i_0,j_0}=\si_s$. 
Let $r<s$ be maximal such that $\siti_r:=\ordvalue(\sivec{\scriptscriptstyle{\restriction {r}}})\le\be$.
If $r>0$, setting $\bepr:=\si_r\cdot(1/\siti_r)\cdot\be$, we have $\bepr\le\be$ and, by Lemma \ref{citedinjtrslem}, 
\[\beucp{i_0,j_0+1}=\trs^{\si_r}(\bepr)\klex\trs^{\si_r}(\tauticp{i_0,j_0})=(\si_{r+1},\ldots,\si_s),\] 
whence $\beucp{i_0,j_0+1}_1<\taucp{i_0,j_0}$.
If $r=0$ and hence $\be<\siti_1$, let $(\xi,u):=\upsseg(\be)$, so that \[\beucp{i_0,j_0+1}=\trs^{\ups_{\xi+u}}(\be)\klex\sivec\]
and thus $\beucp{i_0,j_0+1}_1<\taucp{i_0,j_0}$, using again Lemma \ref{citedinjtrslem} if necessary.
This concludes the proof of Claim \ref{tcassgnclmone}.\qed

\begin{claim}\label{tcassgnclmtwo} If $(i_0,j_0)\not=(1,0)$ and $\chi^\taucp{i_0,j_0}(\taucp{i_0,j_0+1})=1$ then $\be\le\tauticp{i_0,j_0}$ implies $(i_0,j_0)\klex(k_0,l_0)$.
\end{claim}
For the proof of this claim assume $(i_0,j_0)\not=(1,0)$, $\chi^\taucp{i_0,j_0}(\taucp{i_0,j_0+1})=1$, 
$\be\le\tauticp{i_0,j_0}$, and let $(i,j)$ be the $\kglex$-maximal index pair such that 
$\alvecrestrarg{(i,j)}$ is a common initial chain of $\alvec$ and $\gavec:=\me(\alvecrestrarg{i_0,j_0+1})$, 
hence $(i_0,j_0+1)\kglex(i,j)$. By Lemma \ref{cmlmaxextcor} we know that $\ec(\gavec)$ and hence $\ec(\alvecrestrarg{(i,j)})$ exists. 
Note also that by part 2 of Lemma \ref{cmlmaxextcor} we have $\ers_{i_0,j_0+1}(\alvec)\subseteq\ers_{i,j}(\alvec)$, which due to the assumption 
$\be\le\tauticp{i_0,j_0}$ entails $\beucp{i,j}=\beucp{i_0,j_0+1}$, and hence \[\beucp{i,j}_1\le\taucp{i_0,j_0}\]
according to Claim \ref{tcassgnclmone}.
In order to derive a contradiction we now assume that $(i_0,j_0)=(k_0,l_0)$ and discuss the possible cases in the
definition of $\ec(\alvecrestrarg{(i,j)})$. For convenience of notation we set $\tau:=\taucp{i,j}$ and $\taupr:=\taucppr{i,j}$.
\\[2mm]
{\bf Case 1:} $j=1$. Then we have $m_i=1$ by the maximality of $(i_0,j_0)$ and $(i,j)$, hence $i>i_0$ and $(i_0,j_0+1)\klex(i,j)$.
\\[2mm]
{\bf Subcase 1.1:} $\taupr<\tau\in\Ez$. By part 2 of Lemma \ref{cmlmaxextcor} and the assumption $k_0=i_0<i$ we then have 
\[\taucp{i_0,j_0}<\tau=\log((1/\taupr)\cdot\tau)<\rho_i\le\alcp{i+1,1}+\beucp{i,1}_1.\] 
However, we have already seen that $\beucp{i,1}_1\le\taucp{i_0,j_0}$, and by part 2 of Definition \ref{kapparegdefi} 
regarding $\ka\nu\ups$-regularity of tracking chains we have $\alcp{i+1,1}<\tau$, whence $\alcp{i+1,1}+\beucp{i,1}_1<\rho_i$. Contradiction.
\\[2mm]
{\bf Subcase 1.2:} Otherwise. Then by the maximality of $(i,j)$ the index $\alcp{i+1,1}$ is strictly less than $\log((1/\taupr)\cdot\tau)$, 
which is the extending index of $\ec(\alvecrestrarg{(i,j)})$ and according to Lemma \ref{cmlmaxextcor} a proper (non-zero) multiple of 
$\taucp{i_0,j_0}$. We run into the same contradiction as in Subcase 1.1.
\\[2mm]
{\bf Case 2:} $j>1$. Then $\taupr=\taucp{i,j-1}$.
\\[2mm]
{\bf Subcase 2.1:} $\taupr<\tau\in\Ez$. 
\\[2mm]
{\bf 2.1.1:} $\tau=\mu_\taupr<\la_\taupr$. Then $(i_0,j_0+1)\klex(i,j)$, which implies $(k_0,l_0)\klex(i,j-1)$. The extending index of 
$\ec(\alvecrestrarg{(i,j)})$ is then $\la_\taupr$, a proper multiple of $\taucp{i_0,j_0}$.
If $j<m_i$ we obtain the contradiction $\tau+\beucp{i,j}_1\le\la_\taupr$, otherwise we obtain the contradiction
$\alcp{i+1,1}+\beucp{i,j}_1<\rho_i=\la_\taupr+1$ in a similar fashion as in Case 1.
\\[2mm]
{\bf 2.1.2:} Otherwise. The extending index of $\ec(\alvecrestrarg{(i,j)})$ is then $\mu_\tau$, and $m_i=j$.
By part 2 of Definition \ref{kapparegdefi}, i.e.\ $\ka\nu\ups$-regularity of tracking chains, $\alcp{i+1,1}\not=\tau$. 
By the assumptions of this case and using Lemma \ref{cmlmaxextcor} we have $\taucp{i_0,j_0}<\tau$.
We first consider the case $(i,j)=(i_0,j_0+1)$. Then $\alcp{i+1,1}<\tau$ and $\rho_i=\tau+1$. We obtain the contradiction 
$\alcp{i+1,1}+\beucp{i,j}_1<\rho_i$.
Now assume $(i_0,j_0+1)\klex(i,j)$. Again we have $\rho_i=\tau+1$, $\alcp{i+1,1}<\tau$, and we run into the same contradiction.
\\[2mm]
{\bf Subcase 2.2:} Otherwise. Then again $m_i=j$. 
\\[2mm]
{\bf 2.2.1:} $\tau<\mu_\taupr$. This can only occur if $(i,j)=(i_0,j_0+1)$, thus $\taupr=\taucp{i_0,j_0}$ and $\tau=\taucp{i_0,j_0+1}$. 
The extending index of $\ec(\alvecrestrarg{(i,j)})$ is (the $\ka$-index) $\rhoargs{\taupr}{\tau}$, a proper multiple of $\taucp{i_0,j_0}$,
and $\rho_i=\rhoargs{\taupr}{\tau}+1$. 
We are then confronted with the contradiction 
$\alcp{i+1,1}+\beucp{i,j}_1<\rho_i$.
\\[2mm]
{\bf 2.2.2:} Otherwise, that is, $\tau=\mu_\taupr$. This implies $(i_0,j_0+1)\klex(i,j)$, and the extending index of 
$\ec(\alvecrestrarg{(i,j)})$ is $\la_\taupr$ which again is a proper multiple of $\taucp{i_0,j_0}$.
Thus $\alcp{i+1,1}+\beucp{i,j}_1\le\la_\taupr<\la_\taupr+1=\rho_i$. Contradiction.\\[2mm]
Our assumption $(i_0,j_0)=(k_0,l_0)$ therefore cannot hold true, which concludes the proof of Claim \ref{tcassgnclmtwo}. \qed
\\[2mm]
We are now prepared to verify the lemma for each of the seven clauses of the assignment of $\tc(\al+\be)$ to $\al+\be$. 
We check that $\tc(\al+\be)\in\TC$ and make a first step to calculate its ordinal value $\ordvalue(\tc(\al+\be))$.
However, a uniform argument exploiting the careful choice of the index pair $(k_0,l_0)$ will be given in the last part of this proof to 
complete the treatment of the single clauses.
\\[2mm]
{\bf Case 0:} $(k_0,l_0)=(1,0)$. Then $(i_0,j_0)=(1,0)$, hence all $\nu$-indices of $\alvec$ are maximal, i.e.\ given by the $\mu$-operator.
Letting $(\xi,u):=\upsseg(\be)$ and applying either Lemma \ref{trsvallem} in the case $u=0$ or Lemma \ref{reltrsvallem} to 
the reference sequence $(\ups_{\xi+1},\ldots,\ups_{\xi+u})$ if $u>0$, 
according to part 7(b) of Lemma \ref{tcevalestimlem} we obtain
\begin{eqnarray*}
\ordvalue(\tc(\al+\be))&=&\ordvalue(\me(\alvecrestrarg{(1,1)}))+\be,
\end{eqnarray*}
and will show later that this is equal to $\al+\be$.
\\[2mm]
{\bf Case 1:} $(i_0,j_0)=(k_0,l_0)\in\dom(\alvec)$. Here we have $\alcp{i_0,j_0+1}<\mu_\taucp{i_0,j_0}$. 
Let $\varsivec:=\ers(\taucp{i_0,j_0+1})$, which according to Lemma \ref{evallem} is equal to $\trs(\tauticp{i_0,j_0})$.
\\[2mm]
{\bf Subcase 1.1:} $\be<\tauticp{i_0,j_0}$. Then by Claim \ref{tcassgnclmone} $\beucp{i_0,j_0+1}_1<\taucp{i_0,j_0}$, 
and Claim \ref{tcassgnclmtwo} yields $\chi^\taucp{i_0,j_0}(\taucp{i_0,j_0+1})=0$, hence
$\tc(\al+\be)\in\TC$. Since $\alcp{i_0+1,1}+\beucp{i_0,m_{i_0}}_1\ge\rho_{i_0}$, we must have 
$m_{i_0}>j_0+1$ and hence $\rhoargs{\taucp{i_0,j_0}}{\taucp{i_0,j_0+1}}=\taucp{i_0,j_0+1}$. 
According to part 7(b) of Lemma \ref{tcevalestimlem} and Lemmas \ref{trsvallem} and \ref{reltrsvallem} we have
\begin{eqnarray*}
\ordvalue(\tc(\al+\be))&=&\ordvalue(\alvecrestrarg{(i_0,j_0+1)})+\dpf_\varsivec(\taucp{i_0,j_0+1})+\be\\
&=&\ordvalue(\me(\alvecrestrarg{(i_0,j_0+1)}))+\be.
\end{eqnarray*}
It remains to be shown that this is equal to $\al+\be$.
\\[2mm]
{\bf Subcase 1.2:} $\be=\tauticp{i_0,j_0}$. Again, by Claim \ref{tcassgnclmone} we have $\beucp{i_0,j_0+1}=(\taucp{i_0,j_0})$, and by 
Claim \ref{tcassgnclmtwo} $\chi^\taucp{i_0,j_0}(\taucp{i_0,j_0+1})=0$.
$\tc(\al+\be)\in\TC$ is immediate. We compute similarly as above
\begin{eqnarray*}
\ordvalue(\tc(\al+\be))&=&\ordvalue(\alvecrestrarg{(i_0,j_0+1)})+\ka^\varsivec_{\rhoargs{}{\taucp{i_0,j_0+1}}}+\dpf_\varsivec(\rhoargs{}{\taucp{i_0,j_0+1}})+\be\\
&=&\ordvalue(\me(\alvecrestrarg{(i_0,j_0+1)}))+\be,
\end{eqnarray*}
and again it remains to be shown that this is equal to $\al+\be$.
\\[2mm]
{\bf Subcase 1.3:} $\be>\tauticp{i_0,j_0}$.
Making use of Lemma \ref{tcevalestimlem} we observe that 
\[\tauticp{i_0,j_0}<\be\le\sumend(\al)=\sumend(\alticp{n,m_n})<\nu^\varsivec_{\mu_\taucp{i_0,j_0}},\]
which, realizing that due to Lemma \ref{trsvallem} 
$\trs(\nu^\varsivec_{\mu_\taucp{i_0,j_0}})=\varsivec^\frown\mu_\taucp{i_0,j_0}$, according to Lemma \ref{citedinjtrslem} 
implies \[\varsivec\klex\trs(\be)\klex\varsivec^\frown\mu_\taucp{i_0,j_0},\]
whence $\varsivec$ is a proper initial segment of $\trs(\be)$. Thus there is an $r_0\in(0,r)$ such that $\taucp{i_0,j_0}=\be_{r_0}$.
We now see that $\tc(\al+\be)\in\TC$. We will apply Lemma \ref{trsvallem} for the evaluation of $\trs(\be)$, considering two cases.\\[2mm]
{\bf 1.3.1:} $\chi^\taucp{i_0,j_0}(\taucp{i_0,j_0+1})=0$. Then by part 7(b) of Lemma \ref{tcevalestimlem} we have
\begin{eqnarray*}
\ordvalue(\tc(\al+\be))&=&\ordvalue(\me(\alvecrestrarg{(i_0,j_0+1)}))+\be.
\end{eqnarray*}
{\bf 1.3.2:} $\chi^\taucp{i_0,j_0}(\taucp{i_0,j_0+1})=1$. Then by part 7(a) of Lemma \ref{tcevalestimlem} we have
\begin{eqnarray*}
\ordvalue(\tc(\al+\be))&=&\ordvalue(\alvecrestrarg{(i_0,j_0+1)}[\alcp{i_0,j_0+1}+1])+\be\\
                       &=&\ordvalue(\me(\alvecrestrarg{(i_0,j_0+1)}))+\dpf(\me(\alvecrestrarg{(i_0,j_0+1)}))+\be.
\end{eqnarray*}
We leave the task of showing that this is equal to $\al+\be$ for later.
\\[2mm]
{\bf Case 2:} $(i_0,j_0)\klex(k_0,l_0)$.
\\[2mm]
{\bf Subcase 2.1:} $k_0=n+1$ and $\beucp{n,m_n}_1=\taucp{n,m_n}\in\Ez^{>\taunpr}$.
Since $\be\le\tauticp{n,m_n}$ we then have $\be=\tauticp{n,m_n}$, and $\tc(\al+\be)\in\TC$ is clear. Since $k_0=n+1$ we have 
$\taucp{n,m_n}<\rho_n$,
and realizing that $-\tauticp{n,m_n}+\nu^{\rsarg{n,m_n}(\alvec)}_1=\tauticp{n,m_n}$ we obtain \[\ordvalue(\tc(\al+\be))=\al+\be.\]
{\bf Subcase 2.2:} $k_0\le n$, $l_0\in\singleton{1,\ldots,m_{k_0}-2}$ and $\taucp{k,l}+\beucp{k,l}_1\le\la_\taucp{k,l-1}$ 
for $(k,l):=(k_0,l_0+1)$.
\\[2mm]
{\bf 2.2.1:} $\alvecrestrarg{(k,l)}^\frown\left(\taucp{k,l}+\beucp{k,l}_1,\beucp{k,l}_2,\ldots,\beucp{k,l}_\rcp{k,l}\right)$ satisfies 
condition 2 of Definition \ref{trackingchaindefi} for tracking chains.
Then this vector defines $\tc(\al+\be)$ and is easily seen to be a tracking chain. 
Note that since $\taucp{k,l}=\mu_\taucp{k,l-1}\in\Ez^{>\taucp{k,l-1}}\cap\la_\taucp{k,l-1}$ and 
$\alcp{k,l+1}=\taucp{k,l+1}=\mu_{\taucp{k,l}}$ we have
$\alvecrestrarg{k,l+1}\ktc\ec(\alvecrestrarg{k,l})$, implying that $\alvecrestrarg{k,l+1}$ does not possess a critical main line index pair. 
Part 7(b) of Lemma \ref{tcevalestimlem} therefore yields, setting $\varsivec:=\ers_{k,l+1}(\alvec)$,
\[\ordvalue(\me(\alvecrestrarg{k,l+1}))=
  \ordvalue(\alvecrestrarg{k,l+1})+\ka^\varsivec_{\la_\taucp{k,l}}+\dpf_\varsivec(\la_\taucp{k,l}).\]
Setting $\varsivecpr:=\ers_{k,l}(\alvec)$, we now compute using similarly as in Case 0 either Lemma \ref{trsvallem} or Lemma \ref{reltrsvallem}
\begin{eqnarray*}
\ordvalue(\tc(\al+\be))&=&\ordvalue(\alvecrestrarg{k,l})+\dpf_\varsivecpr(\taucp{k,l})+\be\\
&=&\ordvalue(\alvecrestrarg{k,l})+\nu^\varsivec_\taucp{k,l+1}+\ka^\varsivec_{\la_\taucp{k,l}}+\dpf_\varsivec(\la_\taucp{k,l})+\be\\
&=&\ordvalue(\alvecrestrarg{k,l+1})+\ka^\varsivec_{\la_\taucp{k,l}}+\dpf_\varsivec(\la_\taucp{k,l})+\be\\
&=&\ordvalue(\me(\alvecrestrarg{k,l+1}))+\be,
\end{eqnarray*}
and leave the task of showing this to be equal to $\al+\be$ for later.
\\[2mm]
{\bf 2.2.2:} Otherwise. Then $\tc(\al+\be)=\alvecrestrarg{i,j}[\alcp{i,j}+1]\in\TC$ where $(i,j):=(i_0,j_0+1)$.
According to the assumptions defining this case we have $r_{k,l}=1$, $\beucp{k,l}_1=\taucp{i_0,j_0}$, 
$\taucp{k,l}+\beucp{k,l}_1=\la_\taucp{k,l-1}$, which is the extending index of $\ec(\alvecrestrarg{k,l})\not\in\TC$, and
thus $\me(\alvecrestrarg{i,j})=\alvecrestrarg{k,l}$. Defining $\varsivecpr$ and $\varsivec$ as in the previous subcase 2.2.1,
setting $\varsivecnod:=\ers_{i,j}(\alvec)$, and
noticing that $\dpf_\varsivecpr(\la_\taucp{k,l-1})=0$, Lemmas \ref{trsvallem}, \ref{reltrsvallem}, and part 7(c) of Lemma \ref{tcevalestimlem} assure the computation
\begin{eqnarray*}
\ordvalue(\tc(\al+\be))&=&\ordvalue(\alvecrestrarg{i,j})+\ka^\varsivecnod_{\rhoargs{}{\taucp{i,j}}}+
  \dpf_\varsivecnod(\rhoargs{}{\taucp{i,j}})\\
&=&\ordvalue(\alvecrestrarg{k,l})+\ka^\varsivecpr_{\la_\taucp{k,l-1}}\\
&=&\ordvalue(\alvecrestrarg{k,l})+\dpf_\varsivecpr(\taucp{k,l})+\be\\
&=&\ordvalue(\alvecrestrarg{k-1}^\frown(\alcp{k,1},\ldots,\alcp{k,l},\mu_\taucp{k,l}))+\ka^\varsivec_{\la_\taucp{k,l}}+\dpf_\varsivec(\la_\taucp{k,l})+\be\\
&=&\ordvalue(\me(\alvecrestrarg{k,l+1}))+\be,
\end{eqnarray*}
where the last equality holds, since the tracking chain 
$\alvecrestrarg{k,l+1}=\alvecrestrarg{k-1}^\frown(\alcp{k,1},\ldots,\alcp{k,l},\mu_\taucp{k,l})$ does not possess a
critical main line index pair, according to part 7(b) of Lemma \ref{tcevalestimlem}. That this is equal to $\al+\be$ will be shown later.
\\[2mm]
{\bf Subcase 2.3:} Otherwise. Then $k_0>i_0$, $l_0=1$, and $\alcp{k+1,1}+\beucp{k,m_k}_1<\rho_k$ for $k:=k_0-1$.
\\[2mm]
{\bf 2.3.1:} The vector $\alvecrestrarg{k}^\frown\left(\alcp{k+1,1}+\beucp{k,m_k}_1,\beucp{k,m_k}_2,\ldots,\beucp{k,m_k}_\rcp{k,m_k}\right)$ satisfies condition 2 of Definition \ref{trackingchaindefi} for tracking chains.
Then $\tc(\al+\be)$ is defined by this vector and is easily seen to be a tracking chain, since we have already handled Subcase 2.1. 
Let us first assume that $k=n$.
Using Lemma \ref{trsvallem} or \ref{reltrsvallem} as before, we then have 
\[\ordvalue(\tc(\al+\be))=\al+\be.\]
Now we suppose $k<n$. We observe that $\alvecrestrarg{k+1,1}$ does not possess a critical main line index pair since 
$\alcp{k+1,1}<\rho_k\minusp1$, and it is only possible to have $m_k>1$ and $\taucp{k,m_k}<\mu_\taukpr$ if $(k,m_k)=(i_0,j_0+1)$. 
Now Lemmas \ref{trsvallem}, \ref{reltrsvallem}, and part 7(b) of \ref{tcevalestimlem} yield, setting $\varsivec:=\ers_{k,m_k}(\alvec)$,
\begin{eqnarray*}
\ordvalue(\tc(\al+\be))&=&\ordvalue(\alvecrestrarg{k+1,1})+\dpf_\varsivec(\taucp{k+1,1})+\be\\
&=&\ordvalue(\me(\alvecrestrarg{k+1,1}))+\be,
\end{eqnarray*} 
which will be shown to be equal to $\al+\be$.
\\[2mm]
{\bf 2.3.2:} Otherwise. Then $\tc(\al+\be)=\alvecrestrarg{i,j}[\alcp{i,j}+1]\in\TC$ where $(i,j):=(i_0,j_0+1)$.
In this final case we have $r_{k,m_k}=1$, $\beucp{k,m_k}=(\taucp{i,j-1})$, $\me(\alvecrestrarg{i,j})=\alvecrestrarg{k}$, and
$\alcp{k+1,1}+\taucp{i,j-1}=\rho_k\minusp1$ is the extending index of $\ec(\alvecrestrarg{k})\not\in\TC$.
The assumption $m_k>1$ and $\taucp{k,m_k}<\mu_\taukpr$ would imply $(k,m_k)=(i,j)$, but since $\ec(\alvecrestrarg{k})\not\in\TC$,
this would contradict Lemma \ref{cmlmaxextcor}, according to which $\ec(\alvecrestrarg{i,j})\in\TC$.
We therefore either have $\rho_k\minusp1=\log((1/\taukstar)\cdot\taucp{k,1})$ with $k>i$ in the case $m_k=1$, 
or we have $\rho_k\minusp1=\la_\taukpr$ with $\taucp{k,m_k}=\mu_\taukpr$ in the case $m_k>1$. 
Let $\varsivec:=\ers_{k,m_k}(\alvec)$.

$\alcp{k+1,1}$ must be a (possibly zero in the case $k=n$) multiple of $\taucp{i,j-1}$, which is seen as follows. 
Assume we had $0<\taucp{k+1,1}<\taucp{i,j-1}$.
Then $k<n$, and by part 6 of Lemma \ref{tcevalestimlem}, since $\alvec$ is an extension of $\alvecrestrarg{k+1,1}$,
\[\ordvalue(\alvecrestrarg{k+1,1})\le\al\le\ordvalue(\alvecrestrarg{k+1,1})+\dpf_\varsivec(\taucp{k+1,1})\]
and therefore by monotonicity of $\ka^\varsivec$
\[\tauticp{k+1,1}=\ka^\varsivec_\taucp{k+1,1}\le\ka^\varsivec_\taucp{k+1,1}+\dpf_\varsivec(\taucp{k+1,1})<\ka^\varsivec_{\taucp{k+1,1}+1}
<\ka^\varsivec_\taucp{i,j-1}=\tauticp{i,j-1}=\be,\]
so that either $\ordvalue(\alvecrestrarg{k+1,1})=\al$ and $\sumend(\al)=\tauticp{k+1,1}<\be$ or 
$\ordvalue(\alvecrestrarg{k+1,1})<\al$ and $\sumend(\al)\le\dpf_\varsivec(\taucp{k+1,1})<\be$, contradicting the assumption $\be\le\sumend(\al)$.

We are now prepared for another twofold application of Lemma \ref{tcevalestimlem}, first part 7(c), then part 7(b). 
In the case $k=n$ we are finished with the second equation, while otherwise we continue the computation as shown, where again
$\varsivecnod:=\ers_{i,j}(\alvec)$.
\begin{eqnarray*}
\ordvalue(\tc(\al+\be))&=&\ordvalue(\alvecrestrarg{i,j})+\ka^\varsivecnod_{\rhoargs{}{\taucp{i,j}}}+\dpf_\varsivecnod(\rhoargs{}{\taucp{i,j}})\\
&=&\ordvalue(\alvecrestrarg{k})+\ka^\varsivec_\alcp{k+1,1}+\dpf_\varsivec(\alcp{k+1,1})+\be\\
&=&\ordvalue(\alvecrestrarg{k+1,1})+\dpf_\varsivec(\taucp{k+1,1})+\be\\
&=&\ordvalue(\me(\alvecrestrarg{k+1,1}))+\be
\end{eqnarray*}
which in the case $k<n$ will be shown below to be equal to $\al+\be$.
\\[2mm]
We are going to show the equalities left open in the single cases. Notice that all cases where $k_0=n+1$ are finished already. We therefore assume $k_0\le n$ from now on,
whence $\beucp{n,m_n}_1\ge\rho_n$. In the first step we show that 
\begin{equation}\label{mealeq} 
\ordvalue(\me(\alvec))+\dpf(\me(\alvec))+\be=\al+\be.
\end{equation}
We have to consider three cases in each of which we use Lemma \ref{tcevalestimlem}. Let $\varsivec:=\ers_{n,m_n}(\alvec)$.
\\[2mm]
{\bf Case A:} $m_n=1$. Then $\al\le\ordvalue(\me(\alvec))+\dpf(\me(\alvec))=\al+\dpf_\varsivecstar(\taucp{n,1})$ according to part 7 of 
Lemma \ref{tcevalestimlem}, where $\varsivecstar:=\rs_{n^\star}(\alvec)$ is equal to $\rs_{n,0}(\alvec)$ and hence agrees with 
Definition \ref{dpfdefi}.
Since $\rho_n=\log((1/\taunpr)\cdot\taucp{n,1})+1$, we have
\[\dpf_\varsivecstar(\taucp{n,1})=\ka^\varsivec_{\log((1/\taunpr)\cdot\taucp{n,1})}+\dpf_\varsivec(\log((1/\taunpr)\cdot\taucp{n,1})).\]
By part b) of Lemma \ref{reltrslem} the assumption $\be\le\dpf_\varsivecstar(\taucp{n,1})$ would imply 
$\beucp{n,1}_1<\log((1/\taunpr)\cdot\taucp{n,1})+1$, which is not the case. 
\\[2mm]
{\bf Case B:} $m_n>1$ and $\taucp{n,m_n}<\mu_\taunpr$. This is only possible if $(n,m_n)=(i_0,j_0+1)$.
We then have 
\[\al\le\ordvalue(\me(\alvec))+\dpf(\me(\alvec))=\al+\ka^\varsivec_{\rhoargs{}{\taucp{i_0,j_0+1}}}+
  \dpf_\varsivec(\rhoargs{}{\taucp{i_0,j_0+1}}).\]
Here the assumption $\be\le\ka^\varsivec_{\rhoargs{}{\taucp{i_0,j_0+1}}}+\dpf_\varsivec(\rhoargs{}{\taucp{i_0,j_0+1}})$ 
would entail the contradiction $\beucp{i_0,j_0+1}_1<\rho_n$.
\\[2mm]
{\bf Case C:} 
Otherwise, i.e.\ $m_n>1$ and $\taucp{n,m_n}=\mu_\taunpr$. Then we have 
\[\al\le\ordvalue(\me(\alvec))+\dpf(\me(\alvec))=\al+\ka^\varsivec_{\la_\taunpr}+\dpf_\varsivec(\la_\taunpr),\] 
and the assumption $\be\le\ka^\varsivec_{\la_\taunpr}+\dpf_\varsivec(\la_\taunpr)$ would lead to the contradiction 
$\beucp{n,m_n}_1<\la_\taunpr+1=\rho_n$.\\[2mm] 
This concludes the verification of (\ref{mealeq}).\\[2mm] 
We now have to show that for index pairs $(i,j)\in\dom(\alvec)-\singleton{(n,m_n)}$ that are lexicographically greater than or equal to the 
index pair occurring in the respective case above, we have
\begin{equation}\label{mejialeq}
\ordvalue(\me(\alvecrestrarg{(i,j)}))+\dpf(\me(\alvecrestrarg{(i,j)}))+\be=
\ordvalue(\me(\alvecrestrarg{(i,j)^+}))+\dpf(\me(\alvecrestrarg{(i,j)^+}))+\be.
\end{equation}
This means that regarding the equations to be proven in Cases 0 and 1 we assume $(i_0,j_0+1)\kglex(i,j)$, regarding those to be shown in 
Case 2.2 we assume $(k_0,l_0+2)\kglex(i,j)$, and regarding Case 2.3 we assume $(k_0,1)\kglex(i,j)$. 

Let such an index pair $(i,j)$ be given. We may assume that $\me(\alvecrestrarg{(i,j)^+})\ktc\me(\alvecrestrarg{(i,j)})$, 
since in the case of equality there would be nothing to show, while $\me(\alvecrestrarg{(i,j)})\ktc\me(\alvecrestrarg{(i,j)^+})$ is not possible, 
for if this were the case we would have $(i,j)=(i_0,j_0+1)$
where $j_0>0$, $\chi^\taucp{i_0,j_0}(\taucp{i_0,j_0+1})=0$, $\rho_{i_0}=\rhoargs{}{\taucp{i_0,j_0+1}}+\taucp{i_0,j_0}$,
$(i,j)^+=(i+1,1)$, and $\alcp{i+1,1}=\rhoargs{}{i_0,j_0+1}+\xi$ for some $\xi\in(0,\taucp{i_0,j_0})$, which by Lemma \ref{tcevalestimlem} 
would imply that $\be\le\sumend(\al)<\tauticp{i_0,j_0}$, hence $\beucp{i_0,j_0+1}=\beucp{i_0,m_{i_0}}$ and according to 
Claim \ref{tcassgnclmone} $\beucp{i_0,j_0+1}_1<\taucp{i_0,j_0}$. Thus $k_0>i_0$, since $\alcp{i_0+1}+\beucp{i_0,m_{i_0}}_1<\rho_{i_0}$,
which implies Case 2 and the condition $i\ge k_0$, wherefore $i=i_0$ is not permitted.  
We therefore have \[\alvecrestrarg{(i,j)^+}\ktc\ec(\alvecrestrarg{(i,j)}),\] 
i.e.\ $\alvecrestrarg{(i,j)^+}$ is not a maximal $1$-step extension of $\alvecrestrarg{(i,j)}$, 
and consider the two possibilities for $(i,j)^+$:\\[2mm]
{\bf Case \Romannumeral{1}:} $(i,j)^+=(i+1,1)$. We then have $j=m_i$, set $\varsivec:=\ers_{i,m_i}(\alvec)$, and consider three subcases.
\\[2mm]
{\bf Subcase \Romannumeral{1}.1:} $m_i=1$. Then $\alcp{i+1,1}<\log((1/\tauipr)\cdot\taucp{i,1})=\rho_i\minusp 1$, hence 
$\alvecrestrarg{i+1,1}$ does not possess a critical main line index pair. 
Setting $\varsivec_{i}:=\rs_{i,0}(\alvec)$ and $\varsivec_{i+1}:=\rs_{i+1,0}(\alvec)$ (cf.\ Definition \ref{dpfdefi}), 
by part 7(b) of Lemma \ref{tcevalestimlem} we have
\begin{eqnarray*}
\ordvalue(\me(\alvecrestrarg{(i+1,1)}))&=&\ordvalue(\alvecrestrarg{(i+1,1)})+\dpf_{\varsivec_{i+1}}(\taucp{i+1,1})\\
&=&\ordvalue(\alvecrestrarg{(i,1)})+\ka^\varsivec_\alcp{i+1,1}+\dpf_\varsivec(\taucp{i+1,1})\\
&<&\ordvalue(\me(\alvecrestrarg{(i,1)}))+\dpf(\me(\alvecrestrarg{(i,1)}))\\
&=&\ordvalue(\alvecrestrarg{(i,1)})+\dpf_{\varsivec_i}(\taucp{i,1})\\
&=&\ordvalue(\alvecrestrarg{(i,1)})+\ka^\varsivec_{\log((1/\tauipr)\cdot\taucp{i,1})}+\dpf_\varsivec(\log((1/\tauipr)\cdot\taucp{i,1})).
\end{eqnarray*}
By setting $\de:=-\alcp{i+1,1}+\log((1/\tauipr)\cdot\taucp{i,1})$ and assuming $\be\le\ka^\varsivec_\de+\dpf_\varsivec(\de)$ we would obtain 
$\alcp{i+1,1}+\beucp{i,1}_1<\rho_i$, which because of $i\ge k_0$ is not the case.
Thus equation (\ref{mejialeq}) holds in the case $m_i=1$.\\[2mm]
{\bf Subcase \Romannumeral{1}.2:} $(i,m_i)=(i_0,j_0+1)$ where $j_0>0$. 
Then only Case 1 is possible, and it follows that $\alcp{i+1,1}\le\rhoargs{}{\taucp{i,m_i}}$, as we ruled out the situation where 
$\rho_{i_0}=\rhoargs{}{\taucp{i_0,j_0+1}}+\taucp{i_0,j_0}$ and $\alcp{i+1,1}=\rhoargs{}{\taucp{i_0,j_0+1}}+\xi$ 
for some $\xi\in(0,\taucp{i_0,j_0})$.
Lemma \ref{tcevalestimlem} supplies us with
\begin{eqnarray*}
\ordvalue(\me(\alvecrestrarg{i+1,1}))&=&\ordvalue(\alvecrestrarg{(i,m_i)})+\ka^\varsivec_\alcp{i+1,1}+\dpf_\varsivec(\alcp{i+1,1})\\
&<&\ordvalue(\me(\alvecrestrarg{(i,m_i)}))
\end{eqnarray*}
and
\[\ordvalue(\me(\alvecrestrarg{(i,m_i)}))+\dpf(\me(\alvecrestrarg{(i,m_i)}))=
  \ordvalue(\alvecrestrarg{(i,m_i)})+\ka^\varsivec_{\rhoargs{}{\taucp{i,m_i}}}+
  \dpf_\varsivec(\rhoargs{}{\taucp{i,m_i}}).\]
We now see that the assumption $\be\le\ka^\varsivec_\de+\dpf_\varsivec(\de)$, where $\de:=-\alcp{i+1,1}+\rhoargs{}{\taucp{i,m_i}}$, 
would have the consequence $\alcp{i+1,1}+\beucp{i,m_i}_1<\rho_{i}$, which (again because of $i\ge k_0$) is not the case.
We therefore have (\ref{mejialeq}) in this special case.\\[2mm]
{\bf Subcase \Romannumeral{1}.3:} $m_i>1$ and $(i_0,j_0+1)\klex(i,m_i)$. Then $\alcp{i+1,1}<\la_\taucp{i,m_i-1}=\rho_i\minusp 1$. 
Lemma \ref{tcevalestimlem} yields
\begin{eqnarray*}
\ordvalue(\me(\alvecrestrarg{(i+1,1)}))&=&\ordvalue(\alvecrestrarg{(i,m_i)})+\ka^\varsivec_\alcp{i+1,1}+\dpf_\varsivec(\alcp{i+1,1})\\
&<&\ordvalue(\me(\alvecrestrarg{(i,m_i)}))+\dpf(\me(\alvecrestrarg{(i,m_i)}))\\
&=&\ordvalue(\alvecrestrarg{(i,m_i)})+\ka^\varsivec_{\la_\taucp{i,m_i-1}}+\dpf_\varsivec(\la_\taucp{i,m_i-1}),
\end{eqnarray*}
and setting $\de:=-\alcp{i+1,1}+\la_\taucp{i,m_i-1}$ the assumption $\be\le\ka^\varsivec_\de+\dpf_\varsivec(\de)$ would again imply 
$\alcp{i+1,1}+\beucp{i,m_i}_1<\rho_i$.
Consequently, equation (\ref{mejialeq}) follows also in this situation.\\[2mm]
{\bf Case \Romannumeral{2}:} $(i,j)^+=(i,j+1)$.
Then we have $\alcp{i,j+1}=\mu_\taucp{i,j}$, since $(i_0,j_0)\klex(i,j)$.
Due to the fact that $\alvecrestrarg{(i,j+1)}$ is not a maximal $1$-step extension of $\alvecrestrarg{(i,j)}$, we have $j>1$, 
$\taucp{i,j}=\mu_\taucp{i,j-1}\in\Ez\cap(\taucp{i,j-1},\la_\taucp{i,j-1})$,
$\ec(\alvecrestrarg{(i,j)})=\alvecrestrarg{(i,j)}^\frown(\la_\taucp{i,j-1})$, and $(i_0,j_0+1)\klex(i,j)$. 
In particular, $\alvecrestrarg{(i,j+1)}$ does not possess a critical main line index pair. 
Setting $\varsivec:=\ers_{i,j+1}(\alvec)$ and $\varsivecpr:=\ersij(\alvec)$, part 7(b) of Lemma \ref{tcevalestimlem} yields
\begin{eqnarray*}
\ordvalue(\me(\alvecrestrarg{(i,j+1)}))&=&\ordvalue(\alvecrestrarg{(i,j+1)})+\ka^\varsivec_{\la_\taucp{i,j}}+\dpf_\varsivec(\la_\taucp{i,j})\\
&=&\ordvalue(\alvecrestrarg{(i,j)})+\nu^\varsivec_{\mu_\taucp{i,j}}+\ka^\varsivec_{\la_\taucp{i,j}}+\dpf_\varsivec(\la_\taucp{i,j})\\
&=&\ordvalue(\alvecrestrarg{(i,j)})+\dpf_\varsivecpr(\taucp{i,j}).
\end{eqnarray*}
Another extensive application of Lemma \ref{tcevalestimlem} provides us with
\[\ordvalue(\alvecrestrarg{(i,j)})+\dpf_\varsivecpr(\taucp{i,j})<\ordvalue(\me(\alvecrestrarg{(i,j)}))+\dpf(\me(\alvecrestrarg{(i,j)}))=
  \ordvalue(\alvecrestrarg{(i,j)})+\ka^\varsivecpr_{\la_\taucp{i,j-1}}+\dpf_\varsivecpr(\la_\taucp{i,j-1}).\] 
Now setting $\de:=-\taucp{i,j}+\la_\taucp{i,j-1}$, the assumption $\be\le\ka^\varsivecpr_\de+\dpf_\varsivecpr(\de)$ would imply, 
by Lemma \ref{reltrslem}, that $\beucp{i,j}_1\le\de$ and hence 
$\taucp{i,j}+\beucp{i,j}_1\le\la_\taucp{i,j-1}$ which is not the case: In Cases 0, 1 and 2.2 we always have $(k_0,l_0)\klex(i,j-1)$, 
while Case 2.3 presupposes (w.r.t.\ Case 2.2) that $\taucp{k_0,l_0+1}+\beucp{k_0,l_0+1}_1>\la_\taucp{k_0,l_0}$, 
which covers the only possibility where $(k_0,l_0)=(i,j-1)$. 
This concludes the proof of (\ref{mejialeq}).\\[2mm]
From the equations (\ref{mealeq}) and (\ref{mejialeq}) all claimed equalities follow, noticing that only in Subcase 1.3.2 the $\dpf$-term
is non-zero.
This completes the proof of Lemma \ref{tcassgnmntlem}.
\qed

\begin{cor}[cf.\ 6.5 of \cite{CWc}]\label{tcorderisocor} $\tc$ is a $<$-$\ktc$-order isomorphism between $\On$ and
$\TC$ with inverse $\ordvalue$. We thus have 
\[\tc(\ordvalue(\alvec))=\alvec\] 
for any $\alvec\in\TC$ and 
\[\al<\be\quad\aeq\quad\tc(\al)\ktc\tc(\be)\]
for all $\al,\be\in\On$.\qed
\end{cor}

\begin{cor}[corrected 6.6 of \cite{CWc}]\label{alpldpcor} 
Let $\al\in\On$ and $\alvec:=\tc(\al)$ with associated chain $\tauvec$, where $\alvec=(\alevec,\ldots,\alnvec)$ and 
$\alivec=(\alcp{i,1},\ldots,\alcp{i,m_i})$ for $1\le i\le n$. Then we have 
\[\tc(\al+\dpf(\alvec))=\left\{\begin{array}{l@{\quad}l}
\alvec_{\restriction_{i,j+1}}[\alcp{i,j+1}+1] & \mbox{ if }(i,j)=\cml(\alvec)\mbox{ or }
(i,j+1)=(n,m_n)\andsp\taucp{i,j+1}<\mu_\taucp{i,j}\\[2mm]
\me(\alvec)&\mbox{otherwise.}
\end{array}\right.\]
Let $\be\in\On$. Then $\tc(\be)$ is a proper extension of $\tc(\al)$ if and only if 
\[\be\in\left\{\begin{array}{l@{\quad}l}
\left(\al,\al+\dpf(\alvec)\right)&\mbox{ if }\cml(\alvec)\mbox{ exists or }m_n>1\andsp\taucp{n,m_n}<\mu_\taucp{n,m_n-1}\\[2mm]
\left(\al,\al+\dpf(\alvec)\right]&\mbox{ otherwise.}
\end{array}\right.\]  
\end{cor}
{\bf Proof.} The corollary directly follows from the Definition \ref{dpfdefi} of $\dpf$, Lemma \ref{tcevalestimlem}, 
and Corollary \ref{tcorderisocor}.
\qed

\subsection{Closed sets of tracking chains}\label{closedsetssubsec}
Here we introduce the notion of \emph{closed set of tracking chains}. This was first done in \cite{W18} in
order not just to locate an ordinal $\al$ within the core of $\Rtwo$, but to collect the tracking chains of ordinals needed to specify a 
pattern of minimal cardinality, the isominimal realization of which contains and therefore denotes $\al$. 
The corresponding result for $\Ronepl$ was established in \cite{W07c} and completed in Sections 5 and 6 of \cite{CWa}.
The generalization of closedness introduced here will be useful when verifying $<_3$-connections elsewhere, e.g.\ in \cite{W},
where this topic would have been too technical for the intended audience. However, in this article this subsection on
closedness is not used and can be skipped at first reading.

We are now going to generalize the notion of closed sets of tracking chains that was introduced in Section 3 of \cite{W18}. 
Closed sets are easily seen to be closed under the operation of maximal extension ($\me$) introduced in Definition \ref{maxextdefi}. 
We will need closedness to find all parameters from $\Image(\ups)$ involved in (tracking chains of) elements of $\Rtwo$.
These play a key role in handling all (finitely many) ``global'' $<_2$-predecessors needed to locate ordinals in $\Rtwo$. 
As we will see in the next section, ordinals of the form $\ups_{\la+m}$, where $\la+m>0$, $\la\in\Limnod$, and $m\in\N\setminus\{1\}$,
have arbitrarily large $<_2$-successors, namely all ordinals of the form $\ups_{\la+m}+\ups_{\la+m\minusp 1}\cdot(1+\xi)$, $\xi\in\On$. 
For illustration, consider the easy example of the ordinal $\ups_\om\cdot\ups_{17}$, the greatest $<_2$-predecessor of which is 
$\ups_\om$, while its $\leo$-reach is $\lh(\ups_\om\cdot\ups_{17})=\ups_\om\cdot\ups_{17}+\ups_{17}$, which in turn has the 
greatest $<_2$-predecessor $\ups_{18}$ (note that in this example, instead of $17,18$ any pair of natural numbers $k,k+1$, $k>1$ would do).
Another instructive example would be to consider the ordinal $\eps_{\ups_\om+\ups_{17}+1}$, where closure under $\bardot$
(see Section 8 of \cite{W07a} and Section 5 of \cite{CWa}) becomes 
essential, which holds for closed sets of tracking chains, cf.\ Lemma 3.19 of \cite{W18}.
The term decomposition of components of tracking chains in a closed set $M$ of tracking chains via the operations of additive decomposition,
logarithm, $\la$-, and $\bardot$-operator expose all bases of greatest $<_2$-predecessors of elements in $\ov[M]$, cf.\ also Lemma 3.20 of 
\cite{W18} in a more ambitious context in order to enable base minimization, cf.\ Definition 3.26 of \cite{W18}.

\begin{defi}\label{upsseqdefi}
Let $\tauvec\in\TS$. We call $\tauvec$ an {\bf \boldmath$\ups$-sequence} if it is of the form either $\tauvec=(\ups_\la)$ where $\la\in\Lim$ 
or $\tauvec=(\ups_{\la+1},\ldots,\ups_{\la+m})$ where $\la\in\Limnod$ and $m\in(0,\om)$.
A tracking chain $\alvec\in\TC$ is called an {\bf \boldmath$\ups$-sequence} if it is of the form $(\tauvec)$ where $\tauvec\in\TS$ is 
an $\ups$-sequence.
\end{defi}

\begin{defi}[cf.\ 3.7, 3.16 of \cite{W18}]\label{convexprincipalchaindefi}
Let $\alvec=(\alvec_1,\ldots,\alvec_n)\in\TC$, $\alvec_i=(\alcp{i,1},\ldots,\alcp{i,m_i})$, $1\le i\le n$, with associated chain $\tauvec$.
\begin{enumerate}
\item\label{convexpart} $\alvec$ is called {\bf convex} if every $\nu$-index of $\alvec$ is maximal, i.e.\ given by the $\mu$-operator. 
\item\label{principalchainpart} If $\alvec$ satisfies $m_n>1$ and $\alcp{n,m_n}=\mu_\tau$,  
where $\tau:=\taucp{n,m_n-1}$, then $\alvec$ is called a {\bf principal chain to base $\tau$}, 
and $\tau$ is called the {\bf base of $\alvec$}.
If $\alvec\in M$, where $M$ is some set of tracking chains, then we say that $\alvec$ is a {\bf principal chain in $M$} and that $\tau$ 
is a {\bf base in $M$}.
\end{enumerate}
\end{defi}

\begin{defi}[cf.\ 3.1, 3.2, 3.3, and 3.21 of \cite{W18}]\label{closedsetsdefi} 
Let $M\finsub\TC$. $M$ is {\bf closed} if and only if $M$ 
\begin{enumerate}
\item[1.] is {\bf closed under initial chains:} if $\alvec\in M$ and $(i,j)\in\dom(\alvec)$ then $\alvec_{\restriction_{(i,j)}}\in M$,
\item[2.] is {\bf $\nu$-index closed:} if $\alvec\in M$, $m_n>1$, $\alcp{n,m_n}=_\ANF\xi_1+\ldots+\xi_k$ then
\begin{enumerate}
\item[2.1.] $\alvec[\xi_1+\ldots+\xi_l]\in M$ for $1\le l\le k$ and 
\item[2.2.] $\alvec[\mu_{\taupr}]$, unless this is a $\ups$-sequence,
\end{enumerate}
\item[3.] {\bf unfolds minor $\letwo$-components:} if $\alvec\in M$, $m_n>1$, and $\tau<\mu_\taupr$ then: 
\begin{enumerate}
\item[3.1.] ${\alvec_{\restriction_{n-1}}}^\frown(\alcp{n,1},\ldots,\alcp{n,m_n},\mu_\tau)\in M$ in the case $\tau\in\Ez^{>\taupr}$, and
\item[3.2.] otherwise $\alvec^\frown(\varrho^\taupr_\tau)\in M$, provided that $\varrho^\taupr_\tau>0$, 
\end{enumerate}
\item[4.] is {\bf $\ka$-index closed:} if $\alvec\in M$, $m_n=1$, and $\alcp{n,1}=_\ANF\xi_1+\ldots+\xi_k$, then:
\begin{enumerate}
\item[4.1.] if $m_{n-1}>1$ and $\xi_1=\taucp{n-1,m_{n-1}}\in\Ez^{>\taucp{n-1,m_{n-1}-1}}$ then 
${\alvec_{\restriction_{n-2}}}^\frown(\alcp{n-1,1},\ldots,\alcp{n-1,m_{n-1}},\mu_{\xi_1})\in M$,
else ${\alvec_{\restriction_{n-1}}}^\frown(\xi_1)\in M$, and
\item[4.2.] ${\alvec_{\restriction_{n-1}}}^\frown(\xi_1+\ldots+\xi_l)\in M$ for $l=2,\ldots,k$,
\end{enumerate}
\item[5.] {\bf maximizes $\me$-$\mu$-chains:} if $\alvec\in M$ and $\tau\in\Ez^{>\taupr}$, then:
\begin{enumerate}
\item[5.1.] if $m_n=1$ then ${\alvec_{\restriction_{n-1}}}^\frown(\alcp{n,1},\mu_\tau)\in M$, unless this is a $\ups$-sequence, and
\item[5.2.] if $m_n>1$ and $\tau=\mu_\taupr=\la_\taupr$ then 
${\alvec_{\restriction_{n-1}}}^\frown(\alcp{n,1}\ldots,\alcp{n,m_n},\mu_\tau)\in M$, unless this is a $\ups$-sequence, 
\end{enumerate}
\item[6.] {\bf unfolds $\leo$-components:} for $\alvec\in M$, if $m_n=1$ and $\tau\not\in\Ez^{\ge\taupr}\cup\{1\}$ 
(i.e.\ $\tau=\taucp{n,m_n}\not\in\Ezone$, $\taupr=\taunstar$), let
\[\log((1/\taupr)\cdot\tau)=_\ANF\xi_1+\ldots+\xi_k,\]
if otherwise $m_n>1$ and $\tau=\mu_\taupr$ such that $\tau<\la_\taupr$ in the case $\tau\in\Ez^{>\taupr}$, let
\[\la_\taupr=_\ANF\xi_1+\ldots+\xi_k.\]
Set $\xi:=\xi_1+\ldots+\xi_k$, unless $\xi>0$ and $\alvec^\frown(\xi_1+\ldots+\xi_k)\not\in\TC$ 
(due to condition 2 of Definition \ref{trackingchaindefi}),
in which case we set $\xi:=\xi_1+\ldots+\xi_{k-1}$.
Suppose that $\xi>0$. Let $\alvecpl$ denote the vector $\{\alvec^\frown(\xi)\}$ if this is a tracking chain 
(condition 2 of Definition \ref{kapparegdefi}), 
or otherwise the vector ${\alvec_{\restriction_{n-1}}}^\frown(\alcp{n,1},\ldots,\alcp{n,m_n},\mu_\tau)$.
Then the closure of $\{\alvecpl\}$ under clauses 4 and 5 is contained in $M$.
\item[7.] {\bf supports bases:} if $\bevec$ is a principal chain in $M$ to base $\tau$ such that $\taubar\in(\taupr,\tau)$
then $\bevec^\frown(\taubar)\in M$.
\end{enumerate}
\end{defi}

Note that in clause 4.1, if $m_{n-1}=1$ and $\taucp{n-1,1}\in\Ez^{>\tau^\star_{n-1}}$, we have $\rho_{n-1}=\taucp{n-1,1}+1$,
and hence $\alcp{n,1}<\taucp{n-1,1}$ by condition 2 of Definition \ref{kapparegdefi}, 
so that the situation $\xi_1=\taucp{n-1,1}$ does not occur. If the conditions stated in clause 4.1 hold, the preference of $\nu$-indices
over $\ka$-indices applies, and the chain ${\alvec_{\restriction_{n-2}}}^\frown({\alvec_{n-1}}^\frown{\mu_{\xi_1}})$, which can not 
be a $\ups$-sequence, is taken instead of either ${\alvec_{\restriction_{n-1}}}^\frown(\xi_1)$, which is not a tracking chain, or
${\alvec_{\restriction_{n-2}}}^\frown({\alvec_{n-1}}^\frown{1})$, which in this context would be redundant.  

\begin{rmk} Due to the exclusion of $\ups$-sequences from closure in clauses 2 and 5 above it is easy to see that closure of a 
set $M\finsub\TC$ under clauses 1 -- 7 results in a finite set of tracking chains. This is due to decreasing $\htarg{}$- and 
$\operatorname{l}$-measures of the terms involved, cf.\ Definition 3.26 of \cite{W07a} and Definition \ref{Ttauvec}.
Closedness under clauses 1 -- 6 only results in $M$ being a {\bf spanning} set of tracking chains, 
first introduced in Section 5 of \cite{W17}.
\end{rmk}

\section{The Structure \boldmath$\Rtwo$\unboldmath}\label{structuresec}
We are now prepared to generalize Theorem 7.9 and Corollary 7.13 of \cite{CWc} to all ordinal numbers. Theorem 7.9 of \cite{CWc} 
provides the $\le_i$-predecessors ($i=1,2$) of ordinals below $\oneinf=\ups_1$, in particular the greatest $<_i$-predecessor of an 
ordinal in case such exists, while Corollary 7.13 of \cite{CWc} characterizes the $\le_i$-successors of ordinals below $\oneinf$.
The generalization carried out in this article consists of descriptions of $<_i$-pre- and $\le_i$-successorship within all of $\Rtwo$.
For this reason we may say that we \emph{display} the entire structure $\Rtwo$, as claimed in the abstract.
For an in-depth discussion of the initial segment of $\Rtwo$ below $\oneinf$ in arithmetical terms, as secured by Theorem 7.9
and Corollary 7.13 of \cite{CWc}, the reader is referred to Subsection 2.3 of \cite{W18}. There we called this arithmetical characterization 
$\Ctwo$, and in \cite{W17} we showed that it is an elementary recursive structure.

In order to proceed toward generalization of the arithmetical analysis established in \cite{CWc}, recall the notion of relativized
$\le_i$-minimality for $i\in\{1,2\}$: $\al$ is $\be$-$\le_i$-minimal if and only if there does not exist any $\ga\in(\be,\al)$ 
such that $\ga<_i\al$. Hence, $0$-$\le_i$-minimality is equivalent to $\le_i$-minimality. 
As in Definition 7.7 of \cite{CWc} we denote the greatest $<_i$-predecessor of an ordinal $\al$ by $\predec_i(\al)$
if that exists and set $\predec_i(\al):=0$ otherwise. Note that the latter case can have two reasons: either $\al$ is $\le_i$-minimal
or the order type of its $<_i$-predecessors is a limit ordinal. $\predecs_i(\al)$ denotes the set of all $<_i$-predecessors of $\al$,
$\succs_i(\al)$ denotes the class of all $\be$ such that $\al\le_i\be$, and $\lh_i(\al)$ denotes the maximum of $\succs_i(\al)$ if that
exists and $\infty$ otherwise, $\lh:=\lh_1$, where $\lh$ stands for \emph{length}. Note that $\lh_i(\al)$ is not defined to be the 
maximum $\be$ such that $\al\le_i\al+\be$ but rather to be the maximum $\be$ such that $\al\le_i\be$ (if such ordinal exists).

Recall Proposition \ref{letwocriterion} and Lemmas \ref{loalpllem}, \ref{ktwoinflochainlem}, and \ref{letwoupwlem} for basic but central
properties of relations $\leo$ and $\letwo$ in $\Rtwo$.
For the reader's convenience we also cite the notion of \emph{covering}, which is the natural notion of embedding for (pure) patterns,
and which plays a crucial role in the proof of Theorem \ref{maintheo}.

\begin{defi}[7.8 of \cite{CWc}]\label{coveringdef}
Given substructures $X$ and $Y$ of $\Rtwo$, a mapping $h:X\hookrightarrow Y$ is a \emph{covering of $X$ into $Y$}, if
\begin{enumerate}
\item $h$ is an injection of $X$ into $Y$ that is strictly increasing with respect to $\le$, and 
\item $h$ maintains $\le_i$-connections for $i=1,2$, i.e.\ $\forall \al,\be\in X\, (\al\le_i\be\imp h(\al)\le_i h(\be))$.
\end{enumerate}
We call $h$ a \emph{covering of $X$} if it is a covering from $X$ into $\Rtwo$. We call $Y$ a \emph{cover} of $X$ if there is 
a covering of $X$ with image $Y$.\index{cover, covering}
\end{defi}

As mentioned before, the following main theorem describes the structure $\Rtwo$ completely in terms of $\le_i$-predecessorship, $i=1,2$. 
As compared to Theorem 7.9 of \cite{CWc}, which only describes the initial segment $\ups_1$ of the structure $\Rtwo$ in this way, 
new cases arise in relation to the ordinals in $\Image(\ups)$.

The proof of Theorem \ref{maintheo}, of which Theorem \ref{maxchaintheo} from the introduction is an immediate consequence, 
is a modifying and generalizing rewrite of the proof of Theorem 7.9 of \cite{CWc} with 
several corrections and notational adjustments. We keep the proof structure and case numbering comparable to structure and numbering 
chosen in the proof of Theorem 7.9 of \cite{CWc}, however, with a more explicit numbering of subcases and some preference to deal with
cases involving translational isomorphism before cases that in general require base transformation as introduced in Subsection \ref{arithmeticsubsec}. 
The Special Case in Subcase 1.2, as well as Subcases 1.1.3 and 1.2.3, are new, due to the extended claim,
while other parts of the proof smoothly generalize. A correction of a part of the proof of Theorem 7.9 of \cite{CWc} 
is indicated. For further details see the \emph{proof map} at the beginning of the proof of Theorem \ref{maintheo}.

Recall the definition of tracking chain and $[\cdot]$-notation, Definition \ref{trackingchaindefi}, maximal extension $\me$, 
Definition \ref{tcextensiondefi}, the assignment of tracking chains to ordinals $\tc$ in Definition \ref{tcassignmentdefi}, and
of evaluation of (initial segments of) tracking chains, $\ordcp{i,j}(\alvec)$, as well as the evaluation at indices
$\alcp{i,j}$ of $\alvec$ and at indices $\taucp{i,j}$ of the associated chain $\tauvec$, $\alticp{i,j}$ and $\tauticp{i,j}$, respectively, 
see Definition \ref{trchevaldefi}. Recall the notation for \emph{units} $\tauistar$ introduced in Definition \ref{unitsdefi}. 

\begin{theo}[cf.\ 7.9 of \cite{CWc}]\label{maintheo}
Let $\al\in\On$ and $\tc(\al)=:\alvec$, where $\alivec=(\alcp{i,1},\ldots,\alcp{i,m_i})$ for $1\le i\le n$, 
with associated chain $\tauvec$ and segmentation parameters $(\la,t):=\upsseg(\alvec)$ and $p,s_l,(\la_l,t_l)$ for $l=1,\ldots,p$ 
as in Definition \ref{segmentationdefi}.
\begin{enumerate}
\item[a)] We have
\[\al\mbox{ is $\ups_\la$-$\leo$-minimal} \quad\aeq\quad (n,m_n)=(1,1)\]
and the greatest $\lo$-predecessor of $\al$ is
\[\predec_1(\al)=\left\{\begin{array}{l@{\quad}l}
\ups_\la & \mbox{if }(n,m_n)=(1,1)\mbox{ and }\ups_\la\in(0,\alcp{1,1})\mbox{ (hence $\la\in\Lim$)}\\[2mm]
\ordcp{n-1,m_{n-1}}(\alvec) & \mbox{if }m_n=1\mbox{ and } n>1\\[2mm]
\ordvalue(\alvec[\xi])& \mbox{if } m_n>1\mbox{, }\alcp{n,m_n}=\xi+1,
\mbox{ and }\chi^\taucp{n,m_n-1}(\xi)=0\\[2mm]
\ordvalue\left(\me\left(\alvec[\xi]\right)\right) & \mbox{if } m_n>1\mbox{, } \alcp{n,m_n}=\xi+1,
\mbox{ and }\chi^\taucp{n,m_n-1}(\xi)=1\\[2mm]
0 & \mbox{otherwise (a greatest $\lo$-predecessor does not exist).}
\end{array}\right.\]
The situation where the order type of $\predecs_1(\al)$, the set of $\lo$-predecessors of $\al$, is a limit ordinal is characterized 
by the following two cases:
\[\predecs_1(\al)=\left\{\begin{array}{ll}
\bigcup_{\xi\in(0,\la)}\;\;\,\predecs_1\left(\ups_\xi\right) & \mbox{ if }\al=\ups_\la>0\\[2mm]
\bigcup_{\xi<\alcp{n,m_n}}\predecs_1\left(\ordvalue(\alvec[\xi])\right) & \mbox{ if }m_n>1\mbox{ and }\alcp{n,m_n}\in\Lim.
\end{array}\right.\]
\item[b)] We have
\[\al\mbox{ is $\letwo$-minimal} \quad\aeq\quad m_n\le 2\mbox{ and }\taunstar=1,\]
and in terms of $\predec_2$ to denote the greatest $\ktwo$-predecessor we have
\[\predec_2(\al)=\left\{\begin{array}{l@{\quad}l}
\ordcp{n,m_n-1}(\alvec) & \mbox{ if }m_n>2\mbox{, otherwise:}\\[2mm]
\ordcp{i_0,j_0+1}(\alvec) & \mbox{ if } n^\star=:(i_0,j_0)\in\dom(\alvec)\\[2mm]
\ups_\la & \mbox{ if } n^\star=(1,0)\mbox{ and }\taunstar=\ups_\la\in(0,\al)\mbox{ (hence $\la\in\Lim$)}\\[2mm]
\ups_{\la_j} & \mbox{ if } n^\star=(s_j,0)\mbox{ and }\taunstar=\ups_{\la_j}\mbox{ for some }j\in\{1,\ldots,p\}\mbox{ where }\la_j\in\Lim\\[2mm]
\ups_{\la_j+t_j+1} & \mbox{ if } n^\star=(s_j,0)\mbox{ and }\taunstar=\ups_{\la_j+t_j}\mbox{ for some }j\in\{1,\ldots,p\}
\mbox{ where }t_j>0\\[2mm]
0 & \mbox{ otherwise (a greatest $\ktwo$-predecessor does not exist).}
\end{array}\right.\]
The order type of the set of $<_2$-predecessors of $\al$ is a limit ordinal if and only if $\al=\ups_\la>0$:
\[\al={\sup}^+\{\be\mid\be<_2\al\} \quad\aeq\quad \al=\ups_\la>0,\]
and if this is the case, we have
\[\predecs_2(\al)=\left\{\ups_{\ze+k}\mid\ze+k\in(0,\la)\mbox{ where }\ze\in\Limnod\mbox{ and }k\in\N\setminus \{1\}\right\}.\]   
\end{enumerate} 
\end{theo}
{\bf Proof.} The proof is by induction on $\al$. This means that according to the i.h.\ the relations $\leo$ and $\letwo$ look exactly
as claimed by the theorem on the set of ordinals $\al=\{\be\mid\be<\al\}$ and is the reason why the theorem is formulated in terms 
of $\le_i$-predecessors only.

\medskip

\emph{\underline{Proof map.}} Before getting into the proof technically and in detail, we provide an overview to facilitate better orientation.
Subcase 1.1 establishes the successor step of the proof, as its condition, namely $m_n=1$ and $\alcp{n,1}=\xi+1$ for some $\xi$, 
characterizes the situation where $\al$ is a successor ordinal $\be+1$. 
It contains Claim \ref{relleominclaim}, which implies $\de$-$\leo$-minimality of $\al$ as claimed in part a) of the
theorem, where $\de$ is as specified in the proof of Case 1, the case which generally pertains to the situation $m_n=1$, and
turns out to be equal to $\predec_1(\al)$ as claimed. 
It is shown that part b) of the theorem, which deals with $\letwo$-predecessors of $\al$, is trivial for successor ordinals $\al$. 
In the case $\de>0$ it is shown using the criterion provided by Proposition \ref{letwocriterion} that $\de\lo\al$.  

In Subcase 1.2 ($m_n=1$ and $\alcp{n,1}\in\Lim$) on the other hand, it is easy to see that part a) of the theorem immediately follows from
the i.h.\ for continuity reasons, as a limit of $\Sigma_i$-superstructures is again a $\Sigma_i$-superstructure.
Subcase 1.2, however, covers the situation where (successor-) $\ktwo$-connections in $\Rtwo$ need to be verified as claimed, 
and Subcase 1.2.1.2 particularly pertains to the situation where genuinely new $\ktwo$-relations arise in $\Rtwo$ 
(along increasing $\nu$-indices in the appropriate setting of relativization). 
An application of Claim \ref{relleominclaim} enables us to show that $\al$ does not have any $\ktwo$-predecessor greater than $\ga$, where
$\ga$ is defined to be either the greatest $\ktwo$-predecessor of $\al$ as claimed in the theorem in case such is claimed to exist at all,  
or otherwise $\ga:=0$.
Provided that $\ga>0$, the relation $\ga\ktwo\al$ is verified using the criterion given by Proposition \ref{letwocriterion}, which establishes
that $\predec_2(\al)=\ga$ as claimed.

Case 2, which contains Claim \ref{relletwominclaim}, covers the situation $m_n>1$ and most importantly needs to verify that
in Subcase 2.1 ($m_n>1$ and $\alcp{n,m_n}=\xi+1$ for some $\xi$) we actually have \[\alpr=\ordvalue(\alvec[\xi])\not\ktwo\ordvalue(\alvec)=\al.\]
Subcase 2.1.1, where $\predec_1(\al)=\alpr$, discusses the situation where the $\letwo$-component arising at $\alpr$ does not fall 
(non-trivially) back onto the mainline, which means that, unless $\xi$ is a successor ordinal and hence $\alvec[\xi]=\me(\alvec[\xi])$, 
we have \[\ordvalue(\me(\alvec[\xi]))\not\lo\al,\] and the condition for this scenario is $\chit(\xi)=0$.
Subcase 2.1.2 is the more involved situation using Claim \ref{relletwominclaim} where 
\[\alpr\ktwo\ordvalue(\me(\alvec[\xi]))=\predec_1(\al)\lo\al\] 
(the condition for this situation is that $\chit(\xi)=1$).
Subcase 2.2, where $\alcp{n,m_n}$ is a limit ordinal finally follows immediately from the i.h.\ for continuity reasons.

Further case distinction in the proof identifies situations where argumentation relies on the application of base transformation 
(see Subsection \ref{arithmeticsubsec}),  in particular Subcases 1.1.2.2 and 1.2.2.2.2, where an important correction of the proof of the 
corresponding Theorem 7.9 of \cite{CWc} takes place.
The argumentation in Subcases 1.1.1 and 1.2.2.1 is similar, namely by translational isomorphism as familiar from $\Rone$, 
as well as in Subcases 1.1.2.1 and 1.2.2.2.1 (also exploiting translation of isomorphic ordinal intervals). 
The more involved Subcases 1.1.2.2 and 1.2.2.2.2 are handled similarly using base transformation.
Base transformation also occurs in the proofs of Claims \ref{relleominclaim} (see Case \Romannumeral{2} there) 
and \ref{relletwominclaim} (not mentioned
explicitly, as the proof of Claim \ref{relleominclaim} is quite similar to the proof of Claim \ref{relletwominclaim}, in which
we focus on a situation (Subcase \Romannumeral{1}.2) that does not occur in the proof of Claim \ref{relleominclaim}).  

Claims \ref{relleominclaim} and \ref{relletwominclaim} are essential to confirm that certain $\le_i$-relations do not hold, as claimed
in the theorem and mentioned above. Finite sets $X$ and $Z$ are specified so that there does not exist any copy $\tildez$ of $Z$ such
that $X\cup\tildez\cong X\cup Z$ and $\tildez\subseteq\max(Z)$. The role of the sets $X$ is to force $\max(X)<\min(\tildez)$, but $X$ also
contains all existing greatest $\ktwo$-predecessors below $\max(X)$ of elements of $Z$. One can think of the sets $Z$ as sets that are 
\emph{incompressible} in the context or under the constraints provided by $X$.

The additional Subcases 1.1.3 and 1.2.3 as well as the Special Case at the beginning of Subcase 1.2 cover new situations involving ordinals
from $\Image(\ups)$ due to the generalization of the theorem to all ordinals as compared to Theorem 7.9 of \cite{CWc}, which covered
the initial segment $\ups_1$ only. 

\medskip 

\emph{\underline{Beginning of the formal proof.}} In the case $\al=0$, equivalently $\alvec=((0))$, there is nothing to show, so let us assume that 
$\al>0$, whence $\alcp{n,m_n}>0$. 
Defining 
\[\varsivec:=\ers_{n,m_n}(\alvec)\] 
in order to access the setting of relativized connectivity components in which $\al$ is located,
see Definitions \ref{charseqdefi} and \ref{trchevaldefi}, we distinguish between cases concerning $m_n$ and whether $\alcp{n,m_n}$ 
is a limit or a successor ordinal.\\[2mm] 
{\bf Case 1:} $m_n=1$. We define 
\[\de:=\left\{\begin{array}{ll}
  \ups_\la & \mbox{ if }n=1\\
  \ordcp{n-1,m_{n-1}}(\alvec) & \mbox{ if }n>1
  \end{array}\right.\] 
and consider cases regarding $\alcp{n,1}$.
\\[2mm]
{\bf Subcase 1.1:} $\alcp{n,1}$ is a successor ordinal, say $\alcp{n,1}=\xi+1$.
Thus $\taucp{n,1}=1$, $\taunstar=1$, $\al$ is a successor ordinal, say $\al=\be+1$, and clearly $\letwo$-minimal 
since we have Lemmas \ref{ktwoinflochainlem} and \ref{letwoupwlem}, according to which any $\ktwo$-predecessor would be supremum of an 
infinite $\lo$-chain, and finite patterns, such as for instance $\lo$-chains, below such a $\ktwo$-predecessor would reoccur cofinally 
below $\al$. 
According to Definitions \ref{kappanuprincipals}, \ref{kappadpdefi}, and \ref{trchevaldefi} we have
\[\be=\left\{\begin{array}{ll}
  \ka_\xi+\dpf(\xi) & \mbox{ if }n=1\\[2mm]
  \de+\ka^{\varsivec}_\xi+\dpf_{\varsivec}(\xi) & \mbox{ if }n>1,
\end{array}\right.\]
since $\xi\ge\ups_\la$ if $n=1$.
Note that the tracking chain of any ordinal in the interval
$[\de,\be]$ has the initial chain $\alvecrestrarg{n-1,m_{n-1}}$, see Corollary \ref{alpldpcor}. 
In the case $n=1$ we have to show that $\al$ is $\ups_\la$-$\leo$-minimal. This will be the special case $\de=\ups_\la$.
Generally, for $n\ge 1$ we now show that $\al$ is $\de$-$\leo$-minimal, which follows from the following Claim \ref{relleominclaim},  
since the existence of any $\ga\in(\de,\al)$ such that $\ga\lo\al$ would allow us to take the sets $X$ and $Z$ from Claim \ref{relleominclaim} and to $\leo$-reflect the
set $Z_2:=Z\cap[\ga,\al)$ down to a set $\tildez_2\subseteq\ga\le\be$ such that, setting $Z_1:=Z\cap\ga$, we would have $X, Z_1<\tildez_2$ and $X\cup Z_1\cup\tildez_2\cong X\cup Z_1\cup Z_2$, 
so that $\tildez:=Z_1\cup\tildez_2$ would satisfy $X<\tildez\subseteq\be$ and $X\cup\tildez\cong X\cup Z$, contradicting Claim \ref{relleominclaim}.

\begin{claim}\label{relleominclaim} 
There exists a finite set $Z\subseteq(\de,\al)$ such that there does not exist any cover $X\cup\tildez$
of $X\cup Z$ with $X<\tildez$ and $X\cup\tildez\subseteq\be$, where $X$ is the finite set that consists of $\de$ and all existing
greatest $<_2$-predecessors $\ga$ of elements in $Z$ such that $\ga\le\de$.
\end{claim}
{\bf Proof.} In order to prove the claim, let us first consider the case $\xi=\ups_\la$ if $n=1$ or $\xi=0$ if $n>1$. 
Then $\de=\be$ and $\al$ is clearly $\de$-$\leo$-minimal. We trivially choose $X:=\{\de\}$ and $Z:=\emptyset$. 

Now let us assume that $\xi=_\ANF\xi_1+\ldots+\xi_r>0$ such that $\xi>\ups_\la$ if $n=1$. 
Since $\alvec\in\TC$, we then have $\alvec[\xi]\in\TC$ if and only if condition 2 of Definition \ref{kapparegdefi} holds,
and accordingly set  
\[\gavec:=\left\{\begin{array}{l@{\quad}l}
\alvecrestrarg{n-2}^\frown({\alvec_{n-1}}^\frown\mu_{\taucp{n-1,m_{n-1}}})
&\mbox{if }n>1\mbox{ and } \xi=\taucp{n-1,m_{n-1}}\in\Ez^{>\taucppr{n-1}}\\
\alvec[\xi] & \mbox{otherwise.} 
\end{array}\right.\]
Let  $\tc(\be)=:\bevec$, where $\bevec_i=(\becp{i,1},\ldots,\becp{i,k_i})$ for $i=1,\ldots,l$, which according to 
part 7(b) of Lemma \ref{tcevalestimlem} is equal to $\me(\gavec)$ since due to the fact 
that $\alvec\in\TC$, $\gavec$ (and hence also $\bevec$) does not possess a critical main line index pair. 
Let $\sivec$ be the chain associated with $\bevec$ and set $k_0:=0$.
The i.h.\ yields $\de<\ga:=\ordvalue(\gavec)\leo\be$, $\de\lo\ga$ if $\de>0$, and we clearly have $k_l=1$ by the choice of $\gavec$ and the definition of $\me$.
Hence there exists a $\kglex$-minimal index pair $(e,1)\in\dom(\bevec)$ such that both $n\le e\le l$ and 
$\becp{e,1}\not\in\Ez^{>\siestar}$. 
Let 
\[\eta:=\left\{\begin{array}{ll}
  \ups_\la & \mbox{ if }e=1\\
  \ordcp{e-1,k_{e-1}}(\bevec) & \mbox{ if }e>1
  \end{array}\right.\] 
Notice that due to the minimality of $e$ the case $k_{e-1}=1$ can only occur if $e=n>1$, $m_{n-1}=1$, $\gavec=\alvec[\xi]$, 
and hence $\de=\eta$. 
Setting $\bevecpr:=\bevecrestrarg{(e,1)}$, $\bepr:=\ordvalue(\bevecpr)$, and setting $\varsivecpr:=\ers_{e,1}(\bevec)$
in general we have $\de\le\eta$,
\[\ga\le\bepr=\left\{\begin{array}{ll}
    \ka^{\varsivecpr}_\becp{e,1} & \mbox{ if }e=1\\[2mm]
    \eta+\ka^{\varsivecpr}_\becp{e,1} & \mbox{ if }e>1,
  \end{array}\right.\] 
and $\bepr+\dpf_\varsivecpr(\becp{e,1})=\be$.
Note that according to Lemma \ref{evallem} we have 
$\dpf_\varsivecpr(\becp{e,1})=\dpf_{\rs_{e^\star}(\bevec)}(\sicp{e,1})$.
We now consider cases regarding $\becp{e,1}$ in order to define in each case a finite set $Z_\eta\subseteq(\eta,\al)$ such that there 
does not exist any cover $X_\eta\cup\tildez_\eta$ of $X_\eta\cup Z_\eta$ with $X_\eta<\tildez_\eta$ and $X_\eta\cup\tildez_\eta\subseteq\be$,
where $X_\eta$ is the finite set that consists of $\eta$ and all existing greatest $<_2$-predecessors less than or equal to $\eta$ 
of elements in $Z_\eta$. 

{\small Case A below specifies the situation where $\be$ captures a (next) successor-$\lo$-successor of $\eta$ (simple example: $\al=\ups_\la+2$). 
Case B handles the occurrence of a successor-$\ktwo$-successor the greatest $\letwo$-predecessor of which is determined by $\siestar$ 
via the i.h. (simple example: $\al=\epsn\cdot(\om+1)+1$, $\becp{e,1}=\epsn$ where $e=2$. Note that in this example we have 
$\tc(\al)=((\epsn+1))$ and $\tc(\be)=((\epsn,\om),(\epsn))$). 
Case C is where $\be$ captures the $\becp{e,1}=_\NF\ze+\sicp{e,1}$-th $\eta$-$\leo$-minimal component ($\ze>0$), that is, 
a branching of $\eta$-$\leo$-component occurs (simple example: $\al=\om\cdot2+2$).
Case D is the situation where the least $\eta$-$\leo$-component is reached that itself $\lo$-connects to $\becp{e+1,1}$-many components, 
namely the $\becp{e,1}$-th component, that is, nesting of $\eta$-$\leo$-components occurs (simple example: $\al=\om^\om+\om+2$).}
\\[2mm]
{\bf Case A:} $\sicp{e,1}=1$. Then $l=e$, and by the i.h.\ applied to $\bepr=\be$, which is of the form $\be=\be^\circ+1$, there are
$\Xpr\finsub\eta+1$ and $\Zpr\finsub(\eta,\be)$ according to the claim, with the property that  
there does not exist any cover $\Xpr\cup\tildezpr$ of $\Xpr\cup\Zpr$ such that $\Xpr<\tildezpr$ and $\Xpr\cup\tildezpr\subseteq\be^\circ$.
Let
\[X_\eta:=\Xpr\quad\mbox{ and }\quad Z_\eta:=\Zpr\cup\singleton{\be}.\]
Clearly, if there were a set  $\tildez_\eta\subseteq(\eta,\be)$ such that $X_\eta\cup\tildez_\eta$ is a cover of $X_\eta\cup Z_\eta$ then $\Xpr\cup(\tildez_\eta\cap\max(\tildez_\eta))$ would be a cover of $\Xpr\cup \Zpr$ which is contained in $\be^\circ$.\\[2mm] 
{\bf Case B:} $\becp{e,1}=\siestar\in\Ez$. Then $\bevecpr$ is maximal, implying that $l=e$ and $\bepr=\be$.
Note that by monotonicity and continuity, Corollary \ref{kappanuhzcor},
\[\be=\sup\set{\ordvalue(\bevecpr[\ze])}{0<\ze<\siestar}.\] 
By the i.h.\ we see that $\be$ is a successor-$\ktwo$-successor of its greatest $\ktwo$-predecessor $\predec_2(\be)$. 
Accordingly,
\[X_\eta:=\{\predec_2(\be),\eta\}\quad\mbox{ and }\quad Z_\eta:=\singleton{\be}\]
has the requested property.\\[2mm]
{\bf Case C:} $\becp{e,1}=_\NF\ze+\sicp{e,1}$ where $\ze,\sicp{e,1}>1$. Since $\ze+1, \sicp{e,1}+1<\becp{e,1}$ we can apply the i.h.\ to 
$\be_\ze:=\ordvalue(\bevecpr[\ze+1])$ and $\be_\si:=\ordvalue(\bevecpr[\sicp{e,1}+1])$,
obtaining sets $X_1,Z_1$ and $X_2,Z_2$ according to the claim, respectively. We then set
\[X_\eta:=X_1\cup X_2\quad\mbox{ and }\quad Z_\eta:=Z_1\cup\left(\be_\ze+\left(-(\eta+1)+Z_2\right)\right).\]
Then $X_\eta$ and $Z_\eta$ have the desired property due to the fact that 
\[\be_\si\quad\cong\quad\eta+1\cup[\be_\ze,\al)\] 
which in turn follows from the i.h. Clearly, we exploit the i.h.\ regarding $\be_\ze$ and $\be_\si$ in order to see that a 
hypothetical cover of $X_\eta\cup Z_\eta$ contradicting the claim would imply the existence of
a cover of either $X_1\cup Z_1$ or $X_2\cup Z_2$ contradicting the i.h.\\[2mm]
{\bf Case D:} Otherwise. Then $\siestar<\sicp{e,1}=\becp{e,1}\not\in\Ezone$, and we have $k_e=1$, $(e+1,1)\in\dom(\bevec)$, and 
\[0<\becp{e+1,1}=\log((1/\siestar)\cdot\sicp{e,1})<\sicp{e,1}.\]
Note that $\bepr$ is a supremum of $\eta$-$\leo$-minimal ordinals $\be_\nu$ (exchanging the index $\becp{e,1}$ with $\nu$ in the definition
of $\bepr$ above) where $\siestar\le\nu<\becp{e,1}$ and
$\log((1/\siestar)\cdot\nu)<\becp{e+1,1}$, that is, $\becp{e,1}$ is the least index of an $\eta$-$\leo$-relativized component that 
$\leo$-connects to $\becp{e+1,1}$-many components. The constraint $\siestar\le\nu$ guarantees the same greatest $\letwo$-predecessors 
connecting to the $\nu$-th and $\becp{e,1}$-th components.
By the i.h.\ applied to $\be^\circ:=\ov(\bevecpr[\becp{e+1,1}+1])$ we obtain sets $\Xpr$ and $\Zpr\subseteq(\eta,\be^\circ)$ 
according to the claim. 
Define
\[X_\eta:=\Xpr\cup\{\predec_2(\bepr)\}\setminus\{0\}\quad\mbox{ and }\quad
  Z_\eta:=\singleton{\bepr}\cup\left(\bepr+(-\eta+\Zpr)\right).\]
Arguing toward contradiction, let us assume there were a set $\tildez_\eta\subseteq(\eta,\be)$ with  $X_\eta<\tildez_\eta$ such that $X_\eta\cup\tildez_\eta\subseteq\be$ is a cover of $X_\eta\cup Z_\eta$.
Since by the i.h.\ $\bepr\leo\be$, thus $\bepr\leo Z_\eta$ and hence $\mu:=\min(\tildez_\eta)\leo\tildez_\eta$, we find 
cofinally many copies of $\tildez_\eta$ below $\bepr$. We may therefore assume that $\tildez_\eta\subseteq(\eta,\bepr)$ and moreover
for some $\nu\in(0,\becp{e,1})$ such that $\nu\ge\siestar$ and $\log((1/\siestar)\cdot\nu)<\becp{e+1,1}$ (clearly satisfying $\bevecpr[\nu]\in\TC$)
\[\tildez^-_\eta:=\tildez_\eta\setminus\singleton{\mu}\subseteq\left(\be_\nu,\be_{\nu+1}\right)\]
where $\be_\nu:=\ordvalue(\bevecpr[\nu])$ and $\be_{\nu+1}:=\ordvalue(\bevecpr[\nu+1])$.
Setting
\[\tildezpr:=\eta+(-\be_\nu+\tildez^-_\eta)\]
and using that due to the i.h.\ we have
\[\eta+1\cup\left(\be_\nu,\be_{\nu+1}\right)\cong\eta+(-\be_\nu+\be_{\nu+1})\]
we obtain a cover $\Xpr\cup\tildezpr$ of $\Xpr\cup\Zpr$ with $\Xpr<\tildezpr$ and 
$\Xpr\cup\tildezpr\subseteq\ordvalue(\bevecpr[\becp{e+1,1}])$, which contradicts the i.h.\\[2mm]

Now, in the case $\de=\eta$ we are done, choosing $X:=X_\eta$ and $Z:=Z_\eta$. Let us therefore assume that $\de<\eta$. 
We claim that for every index pair $(i,j)\in\{(0,0)\}\cup\dom(\bevec)$ with 
$(n-1,m_{n-1})\kglex(i,j)\klex(e,1)$, where $m_0:=0$, setting $\eta_{0,0}:=\ups_\la$ and $\etaij:=\ordcp{i,j}(\bevec)$ for 
$(i,j)\in\dom(\bevec)$, 
there is $\Zij\finsub(\etaij,\al)$ such that there does not exist any cover $\Xij\cup\tildezij$ of $\Xij\cup\Zij$ 
with $\Xij<\tildezij$ and $\Xij\cup\tildezij\subseteq\be$, where $\Xij$ consists of $\etaij$ and all existing greatest $<_2$-predecessors 
less than or equal to $\etaij$ of elements of $\Zij$.
This is shown by induction on the finite number of $1$-step extensions from $\bevecrestrarg{(i,j)}$ to $\bevecpr$.
The initial step where $(i,j)=(e-1,k_{e-1})$ and $\etaij =\eta$ has been shown above. Now assume $(i,j)\klex(e-1,k_{e-1})$ and let 
$(u,v):=(i,j)^+$.
Let $\Xuv$ and $\Zuv\subseteq(\etauv ,\al)$ be according to the i.h. 
The i.h.\ provides us with knowledge of the $<_i$-predecessors of $\etauv $ ($i=1,2$), which in turn is in $\leo$-relation with every element in 
$\Zuv$. We consider cases regarding $(u,v)$. 

{\small Returning to our simple example for Case B above, where $\al=\epsn\cdot(\om+1)+2$ and $e=2$,
first Case \Romannumeral2 applies with $(i,j)=(1,1)$ and $\becp{u,v}=\becp{1,2}=\om$, and then Case \Romannumeral1 applies with $(i,j)=(0,0)$
and $\becp{u,v}=\becp{1,1}=\epsn$. The resulting sets are $X=X_{0,0}=\{0\}$ and $Z=Z_{0,0}=\{\epsn\cdot\om,\epsn\cdot(\om+1)\}$, 
where $0$ is contained in $X$ only for technical reasons, as $X=\emptyset$ would of course suffice.
Instructive variants of this example regarding Case \Romannumeral2 are $\al=\eps_1\cdot(\om+1)+2$, where $\overline{\eps_1}=\epsn$, 
or $\al=\eps_{\Ga_0+1}\cdot(\om+1)+2$, where $\overline{\eps_{\Ga_0+1}}=\Ga_0$.}
\\[2mm]
{\bf Case \Romannumeral{1}:} $(u,v)=(i+1,1)$ where $i\ge 0$. Then we have $\becp{u,v}=\sicp{u,v}\in\Ez^{>\si^\star_u}$ 
(by the minimality of $e$) and $(u,v)^+=(i+1,2)$ with $\becp{i+1,2}=\mu_{\sicp{u,v}}$.  We define
\[\Zij:=\Zuv\]
and observe that, setting $\Xpr:=\Xuv\setminus\{\etauv\}$, we have $\Xij=\Xpr\cup\{\etaij\}$ since for any greatest $<_2$-predecessor $\nu$ of an
element in $\Zij$ we have $\nu\le\etaij$.
Assume there were a cover $\Xij\cup\tildezij$ of $\Xij\cup\Zij$ with $\Xij<\tildezij$ and $\Xij\cup\tildezij\subseteq\be$.
We have $\etauv\lo\Zuv=\Zij$, and by the i.h.\ we may assume that $\tildezij\subseteq(\etaij,\etauv)$, since if necessary we may 
$\leo$-reflect the elements that are greater than or equal to $\etauv$ down below $\etauv$, into the interval $(\etaij,\etauv)$.
The i.h.\ shows that we have the following isomorphism
\[\etauv\quad\cong\quad\etaij+1\cup(\etauv,\ordvalue(\bevecrestrarg{i}^\frown(\becp{u,v},1))),\]
which shows that defining $\tildezijpr:=\etauv+(-\etaij+\tildezij)$ we obtain another cover $\Xij\cup\tildezijpr$ of $\Xij\cup\Zij$ with the assumed properties.
We now claim that $\Xuv\cup\tildezijpr$ is a cover of $\Xuv\cup\Zuv$ with $\Xuv<\tildezijpr$ and $\Xuv\cup\tildezijpr\subseteq\be$, contradicting the i.h.
Indeed, $\Xpr\cup\tildezijpr$ is a cover of $\Xpr\cup\Zuv$, and we have $\etauv\lo\Zuv,\tildezijpr$ and 
$\etauv\not\letwo\nu$ for any $\nu\in\Zuv\cup\tildezijpr$.\\[2mm]
{\bf Case \Romannumeral{2}:} $(u,v)=(i,j+1)$. Then $(i,j)\in\dom(\bevec)$, and letting $\si:=\sicp{i,j}$ and $\sipr:=\sicppr{i,j}$ we have, 
recalling Definition \ref{barop}, $\sipr\le\sibar$, since $\sipr<\si$ and, according to part 4 of Lemma \ref{trsbasicpropslem}, 
tracking sequences are subsequences of localizations, and $\becp{i,j+1}=\mu_\si$.
The i.h.\ applied to
$\be_\sibar:=\ordvalue(\bevecrestrarg{(i,j)}^\frown(\sibar+1))$ yields sets $X_\sibar$ and 
$Z_\sibar\subseteq(\etaij,\be_\sibar)$ according to the claim.
Setting $\be_\si:=\ordvalue(\bevecrestrarg{(u,v)}^\frown(\si))$ we now define
\[\Zij:=\singleton{\etauv }\cup(\etauv +(-\etaij +Z_\sibar))\cup\singleton{\be_\si}\cup\Zuv\]
and assume that there were a cover $\Xij\cup\tildezij$ of $\Xij\cup\Zij$ with $\Xij<\tildezij$ and $\Xij\cup\tildezij\subseteq\be$. 
Notice that the possibly redundant element $\be_\si$ is the least $<_2$-successor of $\etauv$ and provided explicitly for a practical reason, 
while $\Zuv$ must contain at least one $<_2$-successor of $\etauv$ which, however, we do not keep track of here. 
Thus, the image $\mu :=\min(\tildezij)$ of $\etauv$ must have a 
$<_2$-successor and therefore by the i.h.\ have a tracking chain ending with a limit $\nu$-index. 
Setting $\varsivecpr:=\ersij(\bevec)$ and $\siti:=\ka^\varsivecpr_\si$ we have 
$\tc(\etaij+\siti)=\bevecrestrarg{u,v}[1]\lo\be$, and since $\mu\leo\tildezij$, the set $\tildezij$ is contained in one component enumerated by 
$\ka^\varsivecpr$ starting from $\etaij$, so that, as we can see via $\lo$-downward reflection, the assumption can be fortified to assuming 
\[\tildezij\subseteq[\etaij+\ka^\varsivecpr_\ze,\etaij +\ka^\varsivecpr_{\ze+1})=:I\] 
for the \emph{least} $\ze$, which, using the i.h.\ and recalling that we incorporated a translation of the set $Z_\sibar$ into $\Zij$, 
can easily be seen to satisfy $\ze\in\Ez\cap(\sibar,\si)$ and (with the aid of Lemma \ref{tcevalestimlem})
\[\mu\lo\ordvalue(\me(\bevecrestrarg{(i,j)}^\frown(\ze)))=\etaij +\ka^\varsivecpr_{\ze}+\dpf_\varsivecpr(\ze)=\max(\tildezij).\]
The minimality of $\ze$ moreover allows us to assume that $\ordvalue(\bevecrestrarg{(i,j)}^\frown(\ze,\nu))\letwo\mu $ for some  index $\nu\le\mu_\ze$
for the following reasons: In case of $\mu <\ordvalue(\bevecrestrarg{(i,j)}^\frown(\ze,\mu_\ze))$ there is a least $\nu>0$ such that $\mu \lo\ordvalue(\bevecrestrarg{(i,j)}^\frown(\ze,\nu+1))$, and 
by the i.h.\ we have (making use of Lemma \ref{cmlmaxextcor}) 
$\ordvalue(\bevecrestrarg{(i,j)}^\frown(\ze,\nu))\letwo\predec_1(\ordvalue(\bevecrestrarg{(i,j)}^\frown(\ze,\nu+1)))$. 
If on the other hand $\mu \ge\ordvalue(\bevecrestrarg{(i,j)}^\frown(\ze,\mu_\ze))$ the assumption $\ordvalue(\bevecrestrarg{(i,j)}^\frown(\ze,\mu_\ze))\not\letwo\mu $ 
would imply, using the i.h.\ regarding $\letwo$-predecessors of $\mu$, that there is a least $q>i$ such that $\ordcp{q,1}(\me(\bevecrestrarg{(i,j)}^\frown(\ze)))\lo\mu $ with a corresponding $\ka$-index $\rho$ at $(q,1)$ 
such that $\sumend(\rho)<\ze$ - contradicting the minimality of $\ze$.
We may furthermore strengthen the assumption $\ordvalue(\bevecrestrarg{(i,j)}^\frown(\ze,\nu))\letwo\mu $ for some index $\nu\le\mu_\ze$ 
to actual equality
\[\mu=\ordvalue(\bevecrestrarg{(i,j)}^\frown(\ze,\nu)),\] 
since it is easy to check that this still results in a cover of $\Xij\cup\Zij$ with the assumed properties.

Since $\ze\in(\sibar,\si)$, setting $\phi:=\pizesiinv$ we have $\phi(\la_\ze)<\la_\si$ (cf.\ part 1 of Lemma \ref{finelocbasicpropslem}) 
and $\phi(\mu_\ze)\le\mu_\si$ by Lemma \ref{munoncofinlem}. The vectors in the $\ktc$-segment $\tc[I]$ of $\TC$ have a form
\[\iovec=\bevecrestrarg{(i,j)}^\frown(\zevec,\xivec_1,\ldots,\xivec_g)\]
where $\zevec=(\ze,\ze_1,\ldots,\ze_h)$ with $g, h\ge 0$, cf.\ again Corollary \ref{alpldpcor}.
Let  
\[\zevecpr:= \left\{\begin{array}{l@{\quad}l}
(\becp{i,1},\ldots,\becp{i,j},1) & \mbox{ if } h=0\\
(\becp{i,1},\ldots,\becp{i,j},1+\phi(\ze_1),\phi(\ze_2),\dots,\phi(\ze_h)) & \mbox{ otherwise.}
\end{array}\right.\] 
Let $g_0\in\singleton{1,\ldots,g}$ be minimal such that $\sumend(\xi_{g_0,1})<\ze$ if that exists, and $g_0=g+1$ otherwise.
If $g_0\le g$ let $\xivec_{g_0}=(\xi_{g_0,1},\ldots,\xi_{g_0,k})$ and define $\xivecpr:=(\phi(\xi_{g_0,1}),\xi_{g_0,2}\ldots,\xi_{g_0,k})$.
We can now define the base transformation of $\iovec$ by
\[t(\iovec):=\bevecrestrarg{i-1}^\frown\left(\zevecpr,\phi(\xivec_1),\ldots,
             \phi(\xivec_{g_0-1}),\xivecpr,\xivec_{g_0+1},\ldots,\xivec_g\right).\]
In order to clarify the definition, note that $t(\iovec)=\bevecrestrarg{i-1}^\frown\left(\zevecpr\right)$ in case of $g=0$.
The part $\left(\xivecpr,\xivec_{g_0+1},\ldots,\xivec_g\right)$, which is empty in case of $g_0=g+1$ and is equal to $(\xivecpr)$ if $g_0=g$,
refers to the addition of a parameter below $\ordvalue(\bevecrestrarg{(i,j)}^\frown(\ze))$
which is the reason why relevant indices are not subject to base transformation.
It is easy to see that $t(\iovec)\in\TC$ and therefore \[t: \tc[I]\to\TC,\mbox{ with }\ordvalue[\Image(t)]\subseteq[\etaij +\siti,\be).\]
Using $t$ and applying the i.h.\ in combination with the commutativity of $\phi$ with all operators acting on the indices, as shown in 
Lemma \ref{latallem} and Subsection \ref{ordopsec}, we obtain
\[\etaij+1\cup I\quad\cong\quad\etaij+1\cup\ordvalue[\Image(t)]\]
since thanks to $\sipr<\ze<\si$ it is easy to see that $\etaij +\ka^\varsivecpr_{\ze}$ and $\etaij +\siti$ have the same greatest 
$\ktwo$-predecessor (which then is less than or equal to $\etaij $) 
unless both are $\letwo$-minimal. The set $\dbltildezij:=\ordvalue\circ t\circ\tc[\tildezij]$ therefore gives rise to another cover of 
$\Xij\cup\Zij$ with the assumed properties. We have 
\[\muti:=\min(\dbltildezij)=\ordvalue(\bevecrestrarg{(u,v)}[\phi(\nu)]),\] 
corresponding to $\mu$. 
In the case $\phi(\nu)<\mu_\si=\becp{u,v}$ we may first assume (using iterated $\leo$-downward reflection if necessary) that the set
$\dbltildezij$ is contained in the interval 
\[[\muti,\ordvalue(\bevecrestrarg{(u,v)}[\phi(\nu)+1]))=:J,\] 
but then we may even assume that $\dbltildezij\subseteq[\etauv,\be)$ and $\muti=\etauv$, since otherwise, as seen directly from the i.h., 
we could exploit the translation isomorphism
\[\etaij+1\cup J\quad\cong\quad\etaij+1\cup(\etauv +(-\muti+J))\]
which shifts $J$ into the interval $[\etauv,\be)$. 

We have now transformed the originally assumed cover $\Xij\cup\tildezij$ to a cover $\Xij\cup\dbltildezij$ of $\Xij\cup\Zij$ which fixes 
$\etauv=\min(\dbltildezij)$ and still has the assumed property $\Xij\cup\dbltildezij\subseteq\be$.
Defining $\tildezuv$ to be the subset corresponding to $\Zuv$ in $\dbltildezij$
we obtain a cover $\Xuv\cup\tildezuv$ of $\Xuv\cup\Zuv$ that satisfies $\Xuv<\tildezuv$ and $\Xuv\cup\tildezuv\subseteq\be$,
contradiction.
This concludes the proof of Claim \ref{relleominclaim}.\qed

Clearly, in the case $\de=0$ (that is, $n=1$ and $\la=0$) we have verified that $\al$ is $\leo$-minimal and are done.
Assuming now that in the case $n=1$ we have $\la>0$ and hence $\de=\ups_\la>0$, we show next that $\de\lo\al$ as claimed in part a), 
using the criterion provided in Proposition \ref{letwocriterion}. 
Let finite sets $X\subseteq\de$ and $Y\subseteq[\de,\al)$ be given. Without loss of generality, we may assume that $\de\in Y$. 
We are going to define a set $\tildey$ such that $X<\tildey<\de$ and $X\cup\tildey\cong X\cup Y$,
distinguishing between two subcases, the second of which will require base transformation in its second part, 1.1.2.2.

{\small A simple example for Subcase 1.1.1 below is $\al=\om+1$, $\tc(\al)=((\om),(1))$, another easy but instructive example is 
$\al=\epsn\cdot\om^2+\epsn\cdot(\om+1)+1$, $\tc(\al)=((\epsn\cdot\om),(\epsn+1))$. An instructive example for Subcase 1.1.2.2 is
$\al=\eps_\om\cdot(\om+1)+1$, $\tc(\al)=((\eps_\om,\om),(\eps_\om+1))$, and for the easier (non-critical) Subcase 1.1.2.1 the example
$\al=\eps_\om\cdot\om+\epsn\cdot(\om+1)+1$, $\tc(\al)=((\eps_\om,\om),(\epsn+1))$ illustrates the difference between these two scenarios
in Subcase 1.1.2.
}
\\[2mm]
{\bf Subcase 1.1.1:} $n>1$ and $m_{n-1}=1$. 
Since \[\alcp{n,1}<\rho_{n-1}=\log((1/\taunminstar)\cdot\taucp{n-1,1})+1\]
we see that $\alcp{n-1,1}$ is a limit of ordinals $\eta<\alcp{n-1,1}$ such that $\log((1/\taunminstar)\cdot\sumend(\eta))\ge\xi$.
Now choose such an index $\eta$ large enough so that $\eta>\alcp{n-1,1}\minusp\taucp{n-1,1}$,  $\taunminstar\le\sumend(\eta)<\taucp{n-1,1}$, and $X<\ordvalue(\alvecrestrarg{n-1,1}[\eta])=:\ga$.
Notice that by the i.h.\  $\ga$ and $\de$ have the same $<_i$-predecessors, $i=1,2$.
We will define a translation mapping $t$ in terms of tracking chains that results in an isomorphic copy of the interval $[\de,\be]$ starting from $\ga$.
The tracking chain of an ordinal $\ze\in[\de,\be]$ has a form $\iovec:=\alvecrestrarg{(n-1,1)}^\frown\zevec$ where 
$\zevec=(\zevec_1,\ldots,\zevec_g)$, $g\ge0$, and $\zevec_i=(\zecp{i,1},\ldots,\zecp{i,w_i})$ for $1\le i\le g$. Let
\[\zevecpr:=\left\{\begin{array}{l@{\quad}l}
((\eta,1),\zevec_2,\ldots,\zevec_g)& \mbox{ if } g>0\andsp\zecp{1,1}=\sumend(\eta)\in\Ez^{>\taunminstar}\andsp w_1=1\\
((\eta,1+\zecp{1,2},\zecp{1,3},\ldots,\zecp{1,w_1}),\zevec_2,\ldots,\zevec_g)& \mbox{ if }g>0\andsp\zecp{1,1}=\sumend(\eta)\in\Ez^{>\taunminstar}\andsp w_1>1\\
((\eta),\zevec_1,\ldots,\zevec_g)& \mbox{ otherwise,}
\end{array}\right.\]
where the first two cases take care of condition 2 of Definition \ref{kapparegdefi} since the situation 
$\sumend(\eta)\in\Ez^{>\taunminstar}$ can actually occur, and define
\[t(\iovec):=\alvecrestrarg{(n-2,m_{n-2})}^\frown\zevecpr.\]
The mapping $t$ gives rise to the translation mapping 
\[\ordvalue\circ t\circ\tc:[\de,\be]\to[\ga,\ga+\ka^\varsivec_\xi+\dpf_\varsivec(\xi)],\] 
and by the i.h.\ we have
\[[0,\ga+\ka^\varsivec_\xi+\dpf_\varsivec(\xi)]\quad\cong\quad[0,\ga)\cup[\de,\be].\]
This shows that in order to obtain $X\cup\tildey\cong X\cup Y$ we may choose
\[\tildey:=\ga+(-\de+Y).\]
{\bf Subcase 1.1.2:} $n>1$ and $m_{n-1}>1$. Let $\tau:=\taucp{n-1,m_{n-1}}$, $\si:=\taupr=\taucp{n-1,m_{n-1}-1}$, 
$\sipr:=\taucppr{n-1,m_{n-1}-1}$, and let $\alcp{n-1,m_{n-1}}=_\NF\eta+\tau$ with $\eta=0$ in case of an additive principal number.
If $\alcp{n-1,m_{n-1}}\in\Lim$ let $\alpr$ be a successor ordinal in $(\eta,\eta+\tau)$, 
large enough to satisfy $\alstar:=\ordvalue(\alstarvec)>X$ where $\alstarvec:=\alvecrestrarg{n-1}[\alpr]\in\TC$, 
otherwise let $\alpr:=\alcp{n-1,m_{n-1}}\minusp1$ and set $\alstarvec:=\alvecrestrarg{n-1}[\alpr]$, which is equal to $\alvecrestrarg{n-1,m_{n-1}-1}$ if $\alpr=0$, according to the 
$[\cdot]$-notation introduced after Definition \ref{trackingchaindefi}. 
Notice that we have $\rho_{n-1}\ge\si$ and $\xi<\la_\si$. We consider the following subcases:\\[2mm]
{\bf 1.1.2.1:} $\xi<\si$. Here we can argue comfortably as in the treatment of Subcase 1.1.1: however, in the special case where 
$\chi^\si(\alpr)=1$ consider $\gavec:=\me(\alstarvec)$.
Using Lemma \ref{cmlmaxextcor} and part 7(c) of Lemma \ref{tcevalestimlem} we know that $\ec(\gavec)$ exists and has an extending index
of a form $\si\cdot(\ze+1)$ for some $\ze$ as well as that according to part 2 of Lemma \ref{cmlmaxextcor}
the maximal extension of $\alstarvec$ to $\gavec$ does not add epsilon bases (in the sense of Definition \ref{basesdefi}) 
between $\sipr$ and $\si$. In the cases where $\chi^\si(\alpr)=0$ we set $\gavec:=\alstarvec$.

Clearly, as under the current assumption we have $\xi<\si\in\Ez^{>\sipr}$, the ordinal $\si$ is a limit of ordinals $\eta$ such that 
$\log((1/\sipr)\cdot\sumend(\eta))=\xi+1$, which guarantees that $\sumend(\eta)>\sipr$, and $\eta$ can be chosen large enough so that setting
\[\nu:=\left\{\begin{array}{l@{\quad}l}
\si\cdot\ze+\eta&\mbox{ if }\chi^\si(\alpr)=1\\
\rhoargs{\si}{\alpr}+\eta&\mbox{ if }\alpr\in\Lim\andsp\chi^\si(\alpr)=0\\
\eta&\mbox{ otherwise}
\end{array}\right.\]
we obtain $X<\ordvalue(\gavec^\frown(\nu))=:\deti$. Observe that by the i.h.\ $\deti$ and $\de$ then have the same $\ktwo$-predecessors and the same $\lo$-predecessors
below $\deti$.
The i.h.\ shows that
\[\deti+\ka^{\varsivec}_\xi+\dpf_{\varsivec}(\xi)+1\quad\cong\quad\deti\cup[\de,\be]\]
whence choosing
\[\tildey:=\deti+(-\de+Y)\]
satisfies our needs.\\[2mm]
{\bf 1.1.2.2:} $\xi\ge\si$. Then, as $\xi<\lasi$, we have $\lasi>\si$, which implies $\si\in\Lim(\Ez)$, and $\rho_{n-1}>\xi+1>\si$, 
which entails $\alcp{n-1,m_{n-1}}\in\Lim$, hence $\tau>1$ and $\alpr$ is a successor ordinal.
According to Lemma \ref{cofinlem} 
$\si$ is a limit of $\rho\in\Ez$ with $\phi(\la^{\sipr}_\rho)\ge\xi$ where $\phi:=\pirhosiinv$.
The additional requirement $\rho>\sibar$ yields the bounds $\phi(\la_\rho)<\la_\si$ (cf.\ part 1 of Lemma \ref{finelocbasicpropslem}) and 
$\phi(\mu_\rho)\le\mu_\si$ (by Lemma \ref{munoncofinlem}). 

Note that for any $y\in Y$ the tracking chain $\tc(y)$ is an extension of $\tc(\de)$, see Lemma \ref{tcorderisolem} and 
Corollary \ref{alpldpcor}, and is of a form 
\[\tc(y)=\alvecrestrarg{n-2}^\frown(\alcp{n-1,1},\ldots,\alcp{n-1,m_{n-1}},\zecpy{0,1},\ldots,\zecpy{0,k_0(y)})^\frown\zevec^y\] 
where $k_0(y)\ge 0$, $\zevecy=(\zevecy_1,\ldots,\zevecy_{r(y)})$, $r(y)\ge 0$,  and $\zevecy_u=(\zecpy{u,1},\ldots,\zecpy{u,k_u(y)})$ with $k_u(y)\ge 1$  for $u=1,\ldots,r(y)$.
Notice that $k_0(y)>0$ implies that $\tau\in\Ez^{>\si}$ and $\xi\ge\tau$.
We now define $r_0(y)\in\singleton{1,\ldots,r(y)}$ to be minimal such that $\sumend(\zecpy{r_0(y),1})<\si$ if that exists, 
and $r_0(y):=r(y)+1$ otherwise.
For convenience let $\zecpy{r(y)+1,1}:=0$.
Using Lemma \ref{cofinlem} we may choose an epsilon number $\rho\in(\sibar,\si)$ satisfying $\tau,\xi\in\Tsr$ and 
$\la_\rho\ge\pi(\xi)$, where $\pi:=\pirhosi$, large enough so that 
\[\zecpy{r_0(y),1}, \zecpy{u,v}\in\Tsr\]
for every $y\in Y$, every $u\in[0,r_0(y))$, and every $v\in\singleton{1,\ldots,k_u(y)}$. 
We set 
\[\tildedevec:=\alstarvec^\frown(\rho,\pi(\tau)), \quad \deti:=\ordvalue(\tildedevec),\]
and easily verify using the i.h.\ that $\de$ and $\deti$ have the same $\ktwo$-predecessors in $\On$ and the same $\lo$-predecessors in $X$.
By commutativity with $\pi$ the ordinal $\deti$ has at least $\pi(\xi)$-many immediate $\le_1$-successors, 
i.e.\ $\rho_n(\tildedevec)\ge\pi(\xi)$.
Setting $\varsivecti:=\ers_{n,2}(\tildedevec)$ let
\[\beti:=\deti+\ka^\varsivecti_{\pi(\xi)}+\dpf_\varsivecti(\pi(\xi)).\]
Using $\phi$ we define the embedding
\[t:\beti+1\hookrightarrow\be+1\]
that fixes ordinals $\le\alstar$ and performs base transformation from $\rho$ to $\si$, thereby mapping $\deti$ to $\de$ and $\beti$ to $\be$,
as follows. For any tracking chain \[\zevec=\alstarvec^\frown(\zevec_1,\ldots,\zevec_r)\] of an ordinal in 
the interval $(\alstar,\beti]$, write $\zevec_u=(\zecp{u,1},\ldots,\zecp{u,k_u})$ and let $r_0\in\{1,\ldots,r\}$ be minimal such that
$\sumend(\zecp{r_0,1})<\rho$ if that exists and $r_0:=r+1$ otherwise. Then apply $\phi$ to every $\zecp{u,v}$ such that $u<r_0$
and $v\le k_u$ as well as to $\zecp{r_0,1}$ unless $r_0=r+1$. 
By the i.h.\ $t$ establishes an isomorphism between $\beti+1$ and $\Image(t)$, and
by our choice of $\rho$ we have $Y\subseteq\Image(t)$, so that defining
\[\tildey:=t^{-1}[Y]\]
we obtain the desired copy of $Y$.\\[2mm]
{\bf Subcase 1.1.3:} $n=1$ and $\la>0$. We then have $0<\de=\ups_\la<\xi+1=\alcp{1,1}$, and it is easy to see that we can choose an ordinal
$\nu<\la$ large enough so that $X\subseteq\ups_\nu$, all parameters occurring in (the tracking chains of) the elements of $Y$
are contained in $\ups_\nu$, and all existing greatest $<_2$-predecessors of elements in $Y$ are less than $\ups_\nu$. 
We may then apply straightforward base transformation $\pi_{\ups_\nu,\ups_\la}$ to produce the desired copy $\tildey$.\\[2mm]
\noindent
{\bf Subcase 1.2:} $\alcp{n,1}\in\Lim$.\\[2mm]
{\bf Special case:} $\al=\ups_\la>0$. According to the i.h.\ $\al$ is the supremum of the
infinite $<_1$-chain of ordinals $\ups_\xi$ where $\xi\in(0,\la)$ and of the infinite $<_2$-chain of ordinals $\ups_\xi$ where
$\xi\in(1,\la)$ is not the successor of any limit ordinal. This shows the claims for $\al$ in parts a) and b).\\[2mm]
{\bf Remaining cases:} $\al>\ups_\la$. 
By monotonicity and continuity in conjunction with the i.h.\ it follows that $\al$ is the supremum of ordinals 
either $\leo$-minimal as claimed for $\al$ or with the same greatest $<_1$-predecessor as claimed for $\al$, whence 
unless $\de=0$, $\de$ is the greatest $\lo$-predecessor of $\al$. This shows part a).

We now turn to the proof of part b).
If $\de=0$, i.e.\ $\al$ is $\leo$-minimal, we are done. We therefore assume that $\de>0$ from now on.
Let $\ga$ be the ordinal claimed to be equal to $\predec_2(\al)$. We first show that $\al$ is $\ga$-$\letwo$-minimal, meaning that we
show $\letwo$-minimality in case of $\ga=0$.
Arguing towards contradiction let us assume that there exists $\gastar$ such that $\ga<\gastar\ktwo\al$. Then clearly $\gastar\letwo\de$, 
and due to Lemma \ref{ktwoinflochainlem} we know that $m_{n-1}>1$ in case of $\gastar=\de$ and $n>1$.
Applying the i.h.\ to $\de$ we see that $\letwo$-predecessors of $\de$ either 
\begin{enumerate}
\item have a tracking chain of the form
$\alvec_{\restriction_{i,j+1}}$ where $(i,j)\in\dom(\alvec)$, $j<m_i$, and $i<n$, or 
\item are of the form $\ups_{\lapr}$ where
$\lapr\in\Lim$, $\lapr\le\la$, or
\item are of the form $\ups_{\lapr+\tpr}$ where
$\lapr\in\Limnod$ and $\tpr\in(1,\om)$ such that $\lapr+\tpr<\la$.
\end{enumerate}   
It follows that $\gastar$ must be of either of the forms just described. We define
\[\taustar:=\left\{\begin{array}{rl}
  \taucp{i,j}                & \mbox{ if }\gastar \mbox{ is of the form 1}\\[2mm]
  \ups_\lapr                 & \mbox{ if }\gastar \mbox{ is of the form 2}\\[2mm]
  \ups_{\lapr+\tpr\minusp 1} & \mbox{ if }\gastar \mbox{ is of the form 3,}
\end{array}\right.\] 
then $\taustar$ is the basis of $\gastar$, and by the assumption $\ga<\gastar$ we must have $\taucp{n,1}\in[\taunstar,\taustar)$. 
Let $\ze$ be defined by $\ordcp{i,j}(\alvec)$, if $\gastar$ is of the form 1,
by $\ups_\nu$ for the least $\nu\in(1,\lapr)$ such that $\taucp{n,1}<\ups_\nu$, if $\gastar$ is of the form 2,
and defined by $\ups_{\lapr+\tpr\minusp 1}$, if $\gastar$ is of the form 3.

In case of $\taucp{n,1}<\alcp{n,1}$ let $\eta$ be such that $\alcp{n,1}=_\NF\eta+\taucp{n,1}$, otherwise set $\eta:=0$, and define 
\[\be:=\left\{\begin{array}{rl}
\ka^\varsivec_\eta+\dpf_\varsivec(\eta) & \mbox{ if }n=1\\[2mm]
\de+\ka^\varsivec_\eta+\dpf_\varsivec(\eta) & \mbox{ otherwise,}
\end{array}\right.\] 
so that $\be+\tauticp{n,1}=\al$.
Set 
\[\varsivecstar:=\rs_{n^\star}(\alvec),\]
and note that $\ka^\varsivecstar_\taucp{n,1}=\tauticp{n,1}$ according to part 1 of Lemma \ref{evallem}.

Applying Claim \ref{relleominclaim} of the i.h.\ to the ordinal $\ze+\ka^\varsivecstar_{\taucp{n,1}+1}$ we obtain finite sets $X\subseteq\ze+1$
and $Z\subseteq(\ze,\ze+\ka^\varsivecstar_{\taucp{n,1}+1})$
such that there is no cover $X\cup\tildez$ of $X\cup Z$ with $X<\tildez$ and $X\cup\tildez\subseteq\ze+\tauticp{n,1}+\dpf_\varsivecstar(\taucp{n,1})$.
According to the i.h., there are copies $\tildez_\gastar$ of $Z$ cofinally below $\gastar$ such that $X\cup Z\cong X\cup\tildez_\gastar$. 
By Lemma \ref{letwoupwlem} and our assumption $\gastar\ktwo\al$ we now obtain copies $\tildez_\al$ of $Z$ cofinally below $\al$ (and hence above $\be$) such that 
$X\cup Z\cong X\cup\tildez_\al$.
The i.h.\ reassures us of the isomorphism
\[\ze+1+\tauticp{n,1}\quad\cong\quad\ze+1\cup(\be,\al).\]
This provides us, however, with a copy 
$\tildez\subseteq(\ze,\ze+\tauticp{n,1})$ of $Z$ such that $X\cup Z\cong X\cup\tildez$, contradicting our choice of $X$ and $Z$, whence $\gastar\ktwo\al$ is impossible.
The ordinal $\al$ is therefore $\ga$-$\letwo$-minimal.

We have to show that $\predec_2(\al)=\ga$ where $\ga$ is the ordinal according to the claim in part b).
From now on let us assume that $\taunstar>1$ and set $(i,j):=n^\star$, since $\taunstar=1$ immediately entails $\ga=0=\predec_2(\al)$.  
After having shown that $\al$ is $\ga$-$\letwo$-minimal, the next step is to verify that $\ga\ktwo\al$.
In the situation $\taunstar<\taucp{n,1}$ the ordinal $\al$ is a limit of $\ktwo$-successors of $\ga$ (the greatest $\ktwo$-predecessor of which
is $\ga$). This follows from the i.h.\ noticing that $\alcp{n,1}$ is a limit of indices which are successor multiples of $\taunstar$. 
It therefore remains to consider the situation \[\taunstar=\taucp{n,1}.\]
Here we show $\ga\ktwo\al$ using Proposition \ref{letwocriterion}. To this end let $X\finsub\ga$ and $Y\finsub[\ga,\al)$ be given.
Without loss of generality we may assume that $\ga\in Y$. 
Set $\tau:=\taunstar$ and 
\[(i_0,j_0):=\left\{\begin{array}{ll}
(i,j+1) & \mbox{ if }(i,j):=n^\star\in\dom(\alvec)\\[2mm]
(1,0) & \mbox{ otherwise.}
\end{array}\right.\]
We now check whether there is a $\klex$-maximal index pair $(k,l)\glex(i_0,j_0)$, after which $\alvec$ continues with a sub-maximal index:
let $(k,l)$ be the $\klex$-maximum index pair in $\dom(\alvec)$ such that $(i_0,j_0)\klex(k,l)\klex(n,1)$ and $\alcp{k+1,1}<\rho_k\minusp1$ in case of
$(k,l)^+=(k+1,1)$, i.e.\ $l=m_k$, whereas $\taucp{k,l}<\rho_k(\alvecrestrarg{(k,l)})\minusp1$ in case of $(k,l)^+=(k,l+1)$, i.e.\ $l<m_k$,
if that exists, and $(k,l):=(i_0,j_0)$ otherwise. 
We observe that 
\begin{enumerate}
\item $\alcp{u,v+1}=\mu_\taucp{u,v}$ whenever $(u,v), (u,v+1)\in\dom(\alvec)$ and $(k,l)\klex(u,v+1)$, and
\item $\alvec=\me(\alvecrestrarg{(k,l)^+}),$
\end{enumerate}
which is seen as follows. Assuming the existence of a lexicographically maximal $(u,v+1)$ violating property 1 of $(k,l)$, we can neither 
have $\chi^\taucp{u,v}(\taucp{u,v+1})=0$, as this would be in conflict with the maximality of $(k,l)$, nor can we have  
$\chi^\taucp{u,v}(\taucp{u,v+1})=1$, since $\alvec=\me(\alvecrestrarg{u,v+1})$ according to the definitions of $(k,l)$ and $\me$, while
$\taunstar=\taucp{n,1}<\taucp{u,v}$ by assumption and definition of $\taunstar$, which is in conflict with Lemma \ref{cmlmaxextcor}. 
Thus properties 1 and 2 follow hand in hand.
 
In case of $\taucp{k,l}<\alcp{k,l}$ let $\eta$ be such that $\alcp{k,l}=_\NF\eta+\taucp{k,l}$, otherwise $\eta:=0$.
We set 
\[\be:=\ordcp{k,l}(\alvec)\quad\mbox{ and }\quad\varsivecpr:=\ers_{k,l}(\alvec).\] 
For the reader's convenience we are going to discuss the following cases in full detail.
Subcase 1.2.1.2 below treats the situation in which a genuinely larger $\letwo$-connectivity component arises.
Subcase 1.2.2.2.2 is a correction of the corresponding subcase in \cite{CWc}.
Subcase 1.2.3 is new due to the extended claim of the theorem.

\medskip
{\small Simple examples involve the least $\ktwo$-pair $\epsn\cdot\om\ktwo\epsn\cdot(\om+1)$ in the example 
   $\tc(\epsn\cdot(\om+1))=((\epsn,\om),(\epsn))$ for Subcase 1.2.1.2, 
   $\tc(\varphi_{2,0}\cdot(\om^2+1))=((\varphi_{2,0},\om^2),(\varphi_{2,0}))$ for Subcase 1.2.1.1.,
   further instructive examples for these subcases are $\tc(\Ga_0^2+\Ga_0)=((\Ga_0,\Ga_0),(\Ga_0))$ (Subcase 1.2.1.1) and
   $\tc(\Ga_0^2\cdot(\om+1)+\Ga_0)=((\Ga_0,\Ga_0\cdot\om),(\Ga_0^2),(\Ga_0))$ (Subcase 1.2.1.2) for the pattern characterizing $\Ga_0$,
   namely $\Ga_0^2\cdot\om\ktwo\Ga_0^2\cdot(\om+1)+\Ga_0$ with the inner chain 
   $\Ga_0^2\cdot\om\ktwo\Ga_0^2\cdot(\om+1)\lo\Ga_0^2\cdot(\om+1)+\Ga_0$.
   Note that $\tc(\Ga_0^2\cdot2+\Ga_0)=((\Ga_0,\Ga_0+1))$ and not $((\Ga_0,\Ga_0),(\Ga_0^2),(\Ga_0))$, which is not a tracking chain,
   since $\Ga_0^2\cdot(2k+1)\ktwo\Ga_0^2\cdot(2k+2)$, but $\Ga_0^2\cdot(2k+1)\not\ktwo\Ga_0^2\cdot(2k+2)+\Ga_0$ for all $k<\om$.
   These latter examples also illustrate Case 2.1, in particular Subcase 2.1.2 with the application of Claim \ref{relletwominclaim} for
   $\tc(\Ga_0^2\cdot(2k+2)+\Ga_0)=((\Ga_0,\Ga_0\cdot(k+1)+1))$. 
   
   An example for Subcase 1.2.2.1 is $\al=\tau^\om+\tau^2\cdot\om+\tau$, $\tc(\al)=((\tau,\tau^\om),(\tau^2\cdot\om),(\tau))$,
   where $\tau:=\thtnod(\Om^2\cdot\om)=\thtnod(\thte(\thte(0)+1))$ and $\be=\tau^\om+\tau^2\cdot\om$. Note that $\lh(\tau)=\lh(\be)=\be+\tau+1$.
   An example for Subcase 1.2.2.2.1 is $\al=\be+\tau$ where $\tau=\BH=\thtnod(\thte(\tht_2(0)))$, $\si=\eps_{\BH+1}$, $\be=\si\cdot\om$,
   and $\tc(\al)=((\tau,\si,\om),(\tau))$. Note that $\lh(\tau)=\be+\si=\si\cdot(\om+1)$, completing the least $\ktwo$-chain of three ordinals.
   An easy example for Subcase 1.2.2.2.2 is $\al=\be+\si\cdot\tau+\tau$ where 
   $\tau:=|\mathrm{ID}_2|=\thtnod(\thte(\tht_2(\tht_3(0))))$,
   $\si=\BH_{\tau+1}=\thtt(\thte(\tht_2(0)))$, $\be=\eps_{\si+1}$, so that $\tc(\al)=((\tau,\si,\be),(\si\cdot\tau),(\tau))$. Note
   that $\lh(\tau)=\be\cdot(\om+1)$, completing the least $\ktwo$-chain of four ordinals.
   
   The simplest example for Subcase 1.2.3.1 is $\ups_2\ktwo\al=\ups_\om+\ups_1$, $\tc(\al)=((\ups_\om+\ups_1))$, 
   and for Subcase 1.2.3.2 $\ups_\om\ktwo\al=\ups_\om\cdot2$, $\tc(\al)=((\ups_\om\cdot2))$. 
}
\\[2mm]
{\bf Subcase 1.2.1:} $(k,l)=(i_0,j_0)=(i,j+1)$ where $(i,j)=n^\star\in\dom(\alvec)$. In this case we have $\be=\ga$.
Let $\varrho:=\alcp{i+1,1}$ if $(i,j+1)^+=(i+1,1)$ 
and $\varrho:=\taucp{i,j+1}$ (which then is an epsilon number greater than $\tau$) otherwise. Lemma \ref{chiinvlem}, 
allows us to conclude $\chit(\varrho)=1$ since maximal extension takes us to a successor multiple of
$\tau$ at $(n,1)$, and as verified by Lemma \ref{tcevalestimlem} we have \[\al=\ga+\ka^{\varsivecpr}_\varrho+\dpf_{\varsivecpr}(\varrho).\]
Let $\lapr\in\Limnod$ and $q<\om$ be such that $\logend(\alcp{i,j+1})=\lapr+q$, whence
by definition \[\alcp{i,j+1}=\eta+\om^{\lapr+q}\quad\mbox{ and }\quad\rhoargs{\tau}{\alcp{i,j+1}}=\tau\cdot(\lapr+q\minusp\chit(\lapr)).\]  
It follows from 
$\chit(\varrho)=1$ and $\varrho\le\rhoargs{\tau}{\alcp{i,j+1}}$ that 
$\varrho$ must have the form \[\varrho=\tau\cdot\xi\mbox{ for some }\xi\in(0,\lapr+q\minusp\chit(\lapr)]\]
where $\lapr+q>0$.\\[2mm]
{\bf 1.2.1.1:} $\varrho<\rhoargs{\tau}{\alcp{i,j+1}}$. In this case we are going to check that $\alcp{i,j+1}$ is a supremum of indices 
$\eta+\nu$ such that $\varrho\le\rhoargs{\tau}{\eta+\nu}$ and $\chit(\nu)=0$. Indeed, inspecting all cases we have
\[\alcp{i,j+1}=\sup\{\eta+\nu\mid\nu\in E\}\] 
where
\[E:=\left\{\begin{array}{l@{\quad}l}
  \{\om^{\ze+k}\mid k\in(0,\om),\:\ze\in\Limnod\cap\lapr,\mbox{ and }\ze+k\minusp\chit(\ze)\ge\xi\}&\mbox{ if }q=0\\[2mm]
  \{\om^\lapr\cdot r+\om^{\ze+k}\mid k,r\in(0,\om),\:\ze\in\Limnod\cap\lapr,\mbox{ and }\ze+k\minusp\chit(\ze)\ge\xi\}&\mbox{ if }
    \chit(\lapr)=1\andsp q=1\\[2mm]
  \{\om^{\lapr+q-1}\cdot r\mid r\in(0,\om)\}&\mbox{ if }\chit(\lapr)=0\andsp q>0\mbox{ or }\\ &\quad\;\chit(\lapr)=1\andsp q>1.
  \end{array}\right.\]
According to the definition we have
\[\rhoargs{\tau}{\alcp{i,j+1}}=\left\{\begin{array}{l@{\quad}l}
  \tau\cdot\lapr&\mbox{ if either }q=0\mbox{ or }\chit(\lapr)=1\andsp q=1\mbox{ (note: $\lapr>\xi$)}\\[2mm]
  \tau\cdot(\lapr+q)&\mbox{ if }\chit(\lapr)=0\andsp q>0\\[2mm]
  \tau\cdot(\lapr+q-1)&\mbox{ if }\chit(\lapr)=1\andsp q>1,
  \end{array}\right.\]
and obtain for $\nu\in E$ in the respective cases of the definition of $E$
\[\rhoargs{\tau}{\eta+\nu}=\left\{\begin{array}{l@{\quad}l}
  \tau\cdot(\ze+k\minusp\chit(\ze))&\mbox{ if either }q=0\mbox{ or }\chit(\lapr)=1\andsp q=1\mbox{ (note: $\lapr>\xi$)}\\[2mm]
  \tau\cdot(\lapr+q-1)&\mbox{ if }\chit(\lapr)=0\andsp q>0\\[2mm]
  \tau\cdot(\lapr+q-2)&\mbox{ if }\chit(\lapr)=1\andsp q>1.
  \end{array}\right.\]
Now it is easy to see that $\rhoargs{\tau}{\eta+\nu}\ge\varrho$, since $\varrho=\tau\cdot\xi<\rhoargs{\tau}{\alcp{i,j+1}}$ according to 
the assumption of this subcase.
By the i.h.\ we have 
\[\ga_\nu:=\ordvalue(\alvecrestrarg{(i,j+1)}[\eta+\nu])\ktwo\al_\nu:=\ga_\nu+\ka^{\varsivecpr}_\varrho+\dpf_{\varsivecpr}(\varrho)\] 
and
\begin{equation}\label{congalnu}
\al_\nu\quad\cong\quad\ga_\nu\cup[\ga,\al)
\end{equation}
for the $\nu$ specified above. 
Choose $\nu$ from $E$ large enough so that $X\subseteq\ga_\nu$ and let $Y_\nu\subseteq[\ga_\nu,\al_\nu)$ be the isomorphic copy of $Y$ according to 
isomorphism (\ref{congalnu}).
By the i.h.\ we obtain a copy $\tildey\subseteq\ga_\nu$ according to the criterion given by Proposition \ref{letwocriterion}.
Let $\tildey^+$ with $\tildey\subseteq \tildey^+\subseteq\ga$ be given, and set $U:=X\cup\tildey^+\cap\ga_\nu$, $V:=\tildey^+\setminus\ga_\nu$.
Since by the i.h.\ clearly $\ga_\nu\lo\ga$, we obtain a copy $\tildev$ such that $U<\tildev\subseteq\ga_\nu$ and $U\cup\tildev\cong U\cup V$.  
Setting $\tildey^+_\nu:=(\tildey^+\cap\ga_\nu)\cup\tildev$, hence $\tildey\subseteq\tildey^+_\nu\subseteq\ga_\nu$, the criterion yields
an appropriate extension $Y^+_\nu\subseteq\al_\nu$ such that $X\cup\tildey^+_\nu\cong X\cup Y^+_\nu$ extends $X\cup\tildey\cong X\cup Y_\nu$. 
Now let $Y^+$ be the isomorphic copy of $Y^+_\nu$ according to (\ref{congalnu}). 
This provides us with the extension of $Y$ according to $\tildey^+$ as required by Proposition \ref{letwocriterion}.
\\[2mm]
{\bf 1.2.1.2:} $\varrho=\rhoargs{\tau}{\alcp{i,j+1}}$. Recalling that we have $\chit(\varrho)=1$ this implies $(i,j+1)^+=(i+1,1)$ by  
condition 2 of Definition \ref{trackingchaindefi} and Lemma \ref{cmlmaxextcor}, which also shows that here the case $q=0$ does not occur,
since otherwise it would follow, invoking Lemma \ref{chiinvlem}, that $\chit(\alcp{i,j+1})=1$, whence $\cml(\al)=(i,j)$ and $\alvec\not\in\TC$. 
We now have
\[\alcp{i,j+1}=\sup\set{\eta+\om^{\lapr+q-1}\cdot r}{r\in(0,\om)},\]
and
\begin{eqnarray*}
\rhoargs{\tau}{\alcp{i,j+1}}
   &=&\left\{\begin{array}{l@{\quad}l}
      \tau\cdot(\lapr+q)  &\mbox{ if }\chit(\lapr)=0\\[2mm]
      \tau\cdot(\lapr+q-1)&\mbox{ if }\chit(\lapr)=1
      \end{array}\right.\\[2mm]
   &=&\left\{\begin{array}{l@{\quad}l}
      \rhoargs{\tau}{\eta+\om^\lapr\cdot r}&\mbox{ if }\chit(\lapr)=1\mbox{ and }q=1\\[2mm]
      \rhoargs{\tau}{\eta+\om^{\lapr+q-1}\cdot r}+\tau&\mbox{ otherwise}
      \end{array}\right.\\[2mm]
   &=&\varrho
\end{eqnarray*}
Let $r\in(0,\om)$ be large enough so that, setting $\nu:=\om^{\lapr+q-1}\cdot r$ and
$\ga_\nu:=\ordvalue(\alvecrestrarg{(i,j+1)}[\eta+\nu])$, we obtain $X\subseteq\ga_\nu$. 
Setting $\al_\nu:=\ordvalue(\alvecrestrarg{(i,j+1)}[\eta+\nu+1])$, by Lemma \ref{tcevalestimlem} we obtain 
\[\al_\nu=\left\{\begin{array}{l@{\quad}l}
  \ga_\nu+\ka^{\varsivecpr}_{\varrho}+\dpf_{\varsivecpr}(\varrho)&\mbox{ if }\chit(\lapr)=1\mbox{ and }q=1\\[2mm]
  \ga_\nu+\ka^{\varsivecpr}_{\rhoargs{\tau}{\eta+\nu}}+\dpf_{\varsivecpr}(\rhoargs{\tau}{\eta+\nu})+\tauti=
  \ga_\nu+\ka^{\varsivecpr}_{\varrho}&\mbox{ otherwise.}
  \end{array}\right.\]
Now the i.h.\ yields 
\begin{equation}\label{sndcongalnu}
\al_\nu\quad\cong\quad\ga_\nu\cup[\ga,\al),
\end{equation}
and we choose $\tildey$ to be the isomorphic copy of $Y$ under this isomorphism. Let $\tildey^+$ with $\tildey\subseteq\tildey^+\subseteq\ga$ be given.
Let $U:=X\cup\tildey^+\cap\al_\nu$ and $V:=\tildey^+\setminus\al_\nu$. Since by the i.h.\ we have $\al_\nu\lo\ga$ there exists $\tildev$ with $U<\tildev<\al_\nu$ and 
$U\cup\tildev\cong U\cup V$.
Now let $Y^+$ be the copy of $(\tildey^+\cap\al_\nu)\cup\tildev$ under (\ref{sndcongalnu}). 
This choice satisfies the requirements of Proposition \ref{letwocriterion}.
\\[2mm]
{\bf Subcase 1.2.2:} $(i_0,j_0)\klex(k,l)$. We argue similarly as in Subcases 1.1.1 and 1.1.2 above.
\\[2mm]
{\bf 1.2.2.1:} $l=1$. This subcase corresponds to Subcase 1.1.1.
Here we can only have $(k,l)^+=(k+1,1)$ and $\alcp{k+1,1}<\rho_k\minusp1=\log((1/\taukstar)\cdot\taucp{k,1})$, due to the maximality 
and property 1 of $(k,l)$.
We see that $\alcp{k,1}$ is a limit of ordinals $\eta+\nu<\alcp{k,1}$ such that $\taukstar<\sumend(\nu)<\taucp{k,1}$
and $\log((1/\taukstar)\cdot\sumend(\nu))\ge\alcp{k+1,1}$, and choosing $\nu$ large enough we may assume that 
$Y\cap\be\subseteq\ordvalue(\alvecrestrarg{(k,1)}[\eta+\nu])=:\be_\nu$. 
Using the i.h.\ and setting $\al_\nu:=\be_\nu+\ka^\varsivecpr_\alcp{k+1,1}+\dpf_\varsivecpr(\alcp{k+1,1})$ we now obtain the isomorphism
\begin{equation}\label{thrdcongalnu}
\al_\nu\quad\cong\quad\be_\nu\cup[\be,\al)
\end{equation}
via a mapping of the corresponding tracking chains defined similarly as in Subcase 1.1.1.
In fact, since $\ga\ktwo\al_\nu$ by the i.h., proving that $\ga\ktwo\al$ shows that this isomorphism extends to the suprema, that is, 
mapping $\al_\nu$ to $\al$. 
Exploiting (\ref{thrdcongalnu}) and using that the criterion holds for $\ga,\al_\nu$ we can now straightforwardly show that the criterion holds for $\ga,\al$.
\\[2mm]
{\bf 1.2.2.2:} $l>1$. Here we proceed in parallel with Subcase 1.1.2. Let $\xi:=\alcp{k+1,1}$ in case of $(k,l)^+=(k+1,1)$ and 
$\xi:=\taucp{k,l}$ otherwise, 
whence according to property 2 of $(k,l)$ and Lemma \ref{tcevalestimlem}
\[\al=\be+\ka^{\varsivecpr}_\xi+\dpf_{\varsivecpr}(\xi).\]
Let further $\si:=\taucp{k,l-1}$ and $\sipr:=\taucppr{k,l-1}$. In the case $\alcp{k,l}\in\Lim$ let $\alpr\in(\eta,\alcp{k,l})$ be a successor ordinal large enough so that, setting $\alstar:=\ordvalue(\alstarvec)$ where $\alstarvec:=\alvecrestrarg{(k,l)}[\alpr]$,
\[Y\cap[\alstar,\be)=\emptyset,\]
otherwise let $\alpr:=\alcp{k,l}\minusp1$ and define $\alstar$ and $\alstarvec$ as above.
Notice that we have $\rho_{k}(\alvecrestrarg{(k,l)})\ge\si$ and $\xi<\la_\si$. We consider two cases regarding $\xi$.\\[2mm]
{\bf 1.2.2.2.1:} $\xi<\si$. In the special case where $\chi^\si(\alpr)=1$ consider $\alvecpr:=\me(\alstarvec)$.
Using Lemma \ref{cmlmaxextcor} and part 7(c) of Lemma \ref{tcevalestimlem} we know that $\ec(\alvecpr)$ exists and is
of a form $\si\cdot(\ze+1)$ for some $\ze$ as well as that the maximal extension of $\alstarvec$ to $\alvecpr$ does not add epsilon bases between $\sipr$ and $\si$. In the cases where $\chi^\si(\alpr)=0$ we set $\alvecpr:=\alstarvec$.
Clearly, $\si$ is a limit of ordinals $\rho$ such that $\log((1/\sipr)\cdot\sumend(\rho))=\xi+1$, which guarantees that $\sumend(\rho)>\sipr$, and $\rho$ can be chosen large enough so that setting
\[\nu:=\left\{\begin{array}{l@{\quad}l}
\si\cdot\ze+\rho&\mbox{ if }\chi^\si(\alpr)=1\\[2mm]
\rhoargs{\si}{\alpr}+\rho&\mbox{ if }\alpr\in\Lim\andsp\chi^\si(\alpr)=0\\[2mm]
\rho&\mbox{ otherwise}
\end{array}\right.\]
we obtain, setting $\be_\nu:=\ordvalue(\alvecpr^\frown(\nu))$, $Y\cap[\be_\nu,\be)=\emptyset$. 
Observe that by the i.h.\ $\be_\nu$ and $\be$ then have the same $\ktwo$-predecessors and the same $\lo$-predecessors
below $\be_\nu$.
The i.h.\ shows that 
\[\al_\nu:=\be_\nu+\ka^{\varsivecpr}_\xi+\dpf_{\varsivecpr}(\xi)\quad\cong\quad\be_\nu\cup[\be,\al)\quad\mbox{ and }\ga\ktwo\al_\nu\]
which we can exploit to show that the criterion given by Proposition \ref{letwocriterion} holds for $\ga,\al$ from its validity for 
$\ga,\al_\nu$, implying that the above isomorphism extends to mapping $\al_\nu$ to $\al$.\\[2mm]
{\bf 1.2.2.2.2:} $\xi\ge\si$. Then we consequently have $\alcp{k,l}\in\Lim$, hence $\taucp{k,l}>1$, $\alpr$ is a successor ordinal,
$\si<\lasi$ and thus $\si\in\Lim(\Ez)$. We proceed as in Subcase 1.1.2.2 in order
to choose an epsilon number $\rho\in(\sibar,\si)$ suitable for base transformation. Clearly, $\taucp{k,l}$ takes the role of the ordinal $\tau$
in Subcase 1.1.2.2, and the role of $\de$ there is taken here by $\be$. Consequently, $\deti$ there will become $\be_\rho$ here, as defined later.

Parameters from $Y$ are treated as follows.
Note that for any $y\in Y^{>\alstar}$ the tracking chain $\tc(y)$ is an extension of $\tc(\be)$, and is of a form 
\[\tc(y)=\alvecrestrarg{k-1}^\frown(\alcp{k,1},\ldots,\alcp{k,l},\zecpy{0,1},\ldots,\zecpy{0,k_0(y)})^\frown\zevec^y\] 
where $k_0(y)\ge 0$, $\zevecy=(\zevecy_1,\ldots,\zevecy_{r(y)})$, $r(y)\ge 0$,  and $\zevecy_u=(\zecpy{u,1},\ldots,\zecpy{u,k_u(y)})$ with $k_u(y)\ge 1$  for $u=1,\ldots,r(y)$.
Notice that $k_0(y)>0$ implies that $\taucp{k,l}\in\Ez^{>\si}$ and $\xi\ge\taucp{k,l}$.
We now define $r_0(y)\in\singleton{1,\ldots,r(y)}$ to be minimal such that $\sumend(\zecpy{r_0(y),1})<\si$ if that exists, 
and $r_0(y):=r(y)+1$ otherwise.
For convenience let $\zecpy{r(y)+1,1}:=0$.
Using Lemma \ref{cofinlem} we may choose an epsilon number $\rho\in(\sibar,\si)$ satisfying $\taucp{k,l},\xi\in\Tsr$ and 
$\la_\rho\ge\pi(\xi)$, where $\pi:=\pirhosi$, large enough so that 
\[\zecpy{r_0(y),1}, \zecpy{u,v}\in\Tsr\]
for every $y\in Y^{>\alstar}$, every $u\in[0,r_0(y))$, and every $v\in\singleton{1,\ldots,k_u(y)}$. 

We may now map $\be$ to $\be_\rho:=\ordvalue(\bevec_\rho)$ where 
\[\bevec_\rho:=\alstarvec^\frown(\rho,\pi(\taucp{k,l})),\] 
easily verifying using the i.h.\ that $\be$ and $\be_\rho$ have the same
$\ktwo$-predecessors in $\On$ and the same $\lo$-predecessors in $X\cup(Y\cap\be_\rho)$.
Setting $\varsivecti:=\ers_{n,2}(\bevec_\rho)$, define
\[\al_\rho:=\be_\rho+\ka^\varsivecti_{\pi(\xi)}+\dpf_\varsivecti(\pi(\xi)).\]
In the same way as in Subcase 1.1.2.2 we can now define the embedding
\[t: \al_\rho\hookrightarrow\al\]
which fixes ordinals $\le\alstar$, so that by the i.h.\ 
\[\al_\rho\cong\Image(t)\quad\mbox{ and }\quad\ga\ktwo\al_\rho.\]
By our choice of $\alstar$ we have $X\cup(Y\cap\be)\subseteq\alstar$, and by our choice of $\rho$ we have 
$Y\subseteq\Image(t)$, hence $t^{-1}$ copies $Y^{>\alstar}$ into $[\be_\rho,\al_\rho)$, and applying the mapping $t$
we can now derive the validity of the criterion given by Proposition \ref{letwocriterion} for $\ga,\al$ from its validity for 
$\ga,\al_\rho$, which implies that the above isomorphism extends to mapping $\al_\rho$ to $\al$.
\\[2mm]
{\bf Subcase 1.2.3:} $(k,l)=(1,0)$. According to our assumptions, in particular $(i_0,j_0)=(1,0)$, $\taucp{n,1}=\taunstar>1$, 
and since the Special Case $\al=\ups_\la$ has been considered already, we have $(k,l)^+=(1,1)$, $\la\in\Lim$, 
$\alcp{1,1}\in(\ups_\la,\ups_{\la+1})$, $\alvec=\me(\alvecrestrarg{1,1})$ according to property 2 of $(k,l)$, and 
\[\ga=\left\{\begin{array}{ll}
\ups_\la & \mbox{ if }\taunstar=\ups_\la\\
\ups_{\la_j} & \mbox{ if }\taunstar=\ups_{\la_j}\mbox{ for some }j\in\{1,\ldots,p\}\mbox{ where }\la_j\in\Lim\\
\ups_{\la_j+t_j+1} & \mbox{ if }\taunstar=\ups_{\la_j+t_j}\mbox{ for some }j\in\{1,\ldots,p\}\mbox{ where }\la_j\in\Limnod\mbox{ and }t_j>0.
\end{array}\right.\]
{\bf 1.2.3.1:} $\ga<\ups_\la$. Here we argue as in Subcase 1.2.2, as $\alcp{1,1}$ is clearly a submaximal index. 
Cofinally in $\la$ we find $\rho$ such that $Y\cap[\ups_\rho,\ups_\la)=\emptyset$ and all parameters below $\ups_\la$ from components of
tracking chains of elements of $Y\setminus\ups_\la$ and of $\al=\ov(\me(\alvec_{\restriction_{1,1}}))$ are contained in $\ups_\rho$.
Applying base transformation $\pi:=\pi_{\ups_\rho,\ups_\la}$ to the elements of $Y\setminus\ups_\la$ then results in an isomorphic copy of $X\cup Y$
below $\al_\rho:=\pi(\al)$, which itself satisfies $\ga\ktwo\al_\rho$ by the i.h., so that again we can derive the validity of the criterion for 
$\ga\ktwo\al$ and $X\cup Y$ from its validity for $\ga\ktwo\al_\rho$ using the embedding from $\al_\rho$ into $\ups_\rho\cup[\ups_\la,\al)$ via
inverted base transformation $\pi^{-1}$.\\[2mm]   
{\bf 1.2.3.2:} $\ga=\ups_\la$. Setting $\xi:=\alcp{1,1}$ we have 
\[\al=\ka_\xi+\dpf(\xi)\mbox{, }\quad\quad\alvec=\me(((\xi)))\mbox{, }\quad\mbox{ and }\quad\chi^{\ups_\la}(\xi)=1,\]
invoking again Lemma \ref{chiinvlem}.
We again choose a sufficiently large $\rho<\la$, where now $\rho=\lapr+\tpr+2$ for suitable $\lapr\in\Limnod$, and $\tpr<\om$, 
such that $X<\ups_\rho$ and all parameters below $\ups_\la$ of $\al$ (equivalently, $\xi$) and of all components of tracking chains of the 
elements of $Y$ are contained in $\ups_\rho$. Setting $\pi:=\pi_{\ups_\rho,\ups_\la}$ and $\zevec:=(\ups_{\lapr+1},\ldots,\ups_{\lapr+\tpr+2})$, 
notice that by Lemma \ref{rhobasetrafolem} $\chi^{\ups_\rho}(\pi(\xi))=1$ and choose $\nu\in\Hz\cap\ups_{\rho+1}$ 
such that \[\rho_1((\zevec^\frown\nu))\minusp 1=\pi(\xi),\]
which is done as follows. Let $\xi=\ups_\la\cdot(\ze+l)$ where $\ze\in\Limnod$ and $l<\om$. 
Note that if $l=0$, we must have $\ze=(1/\ups_\la)\cdot\xi\in\Lim$ since $\xi>\ups_\la$, and 
$\chi^{\ups_\la}(\ze)=\chi^{\ups_\la}(\xi)=1$ according to Lemma \ref{chiinvlem}.
Also recall that base transformation commutes with $\om$-exponentiation (Lemma \ref{simplebarlem}) and the $\varrho$-operator 
(Lemma \ref{rhobasetrafolem}), which is useful to keep in mind during the following calculation.
Set 
\[k:=\left\{\begin{array}{ll}
    0&\mbox{ if }l=0\\
    l-1+\chi^{\ups_\la}(\ze)&\mbox{ if } l>0
  \end{array}\right.\quad\mbox{ and }\quad\nu:=\om^{\pi(\ze)+k}.\]
Now, if $\chi^{\ups_\rho}(\nu)=1$ it follows that $k=0$ and $\chi^{\ups_\rho}(\pi(\ze))=1$, hence $\chi^{\ups_\la}(\ze)=1$ and $l=0$,
whence $\rho_1(\zevec^\frown\nu)\minusp 1=\rhoargs{\ups_\rho}{\nu}=\ups_\rho\cdot\pi(\ze)=\pi(\xi)$. And if $\chi^{\ups_\rho}(\nu)=0$,
we have  $\rho_1(\zevec^\frown\nu)\minusp 1=\rhoargs{\ups_\rho}{\nu}+\ups_\rho$ and $k\minusp\chi^{\ups_\rho}(\pi(\ze))+1=l$, since 
if $k=0$ we must have $l>0$ as $\chi^{\ups_\la}(\xi)=1$. Thus $\rhoargs{\ups_\rho}{\nu}+\ups_\rho=\ups_\rho\cdot(\pi(\ze)+l)=\pi(\xi)$.

Setting $\mu:=\nu^\zevec_\nu$, we now obtain our master copy $\tildey$ of $Y$ by
\[\tildey:=\mu+(-\ups_\rho+\pi[Y]).\]
Now, let a finite set $\tildeypl$ such that $\tildey\subseteq\tildeypl\subseteq\ups_\la$ be given.
If necessary, let $\tildez$ be a copy of $\tildeypl\setminus\mu^+$, where $\mu^+:=\nu^\zevec_{\nu+1}$, below $\mu^+$ such that for 
$\tildeypr:=(\tildeypl\cap\mu^+)\cup\tildez$
\[X\cup\tildeypl\cong X\cup\tildeypr,\]
and set 
\[Y^+:=(\tildeypr\cap\mu)\cup\pi^{-1}[\ups_\rho+(-\mu+\tildeypr)].\]
It is now easy to see that the isomorphism of $X\cup\tildey$ and $X\cup Y$ extends to an isomorphism of $X\cup\tildeypl$ and $X\cup Y^+$.
\\[2mm]
\noindent
{\bf Case 2:} $m_n>1$. We now discuss the situation where the $\klex$-maximal index of the tracking chain of $\al$ is a $\nu$-index. 

{\small Simple examples for the subcases of Case 2 discussed below are $\tc(\epsn\cdot2)=((\epsn,1))$ (Subcase 2.1.1.1.1),
   $\tc(\eps_{\BH+1}\cdot2)=((\BH,\eps_{\BH+1},1))$ (Subcase 2.1.1.1.2), where $\BH$ stands of the Bachmann-Howard ordinal and 
   $\eps_{\BH+1}$ for the least epsilon number greater than the Bachmann-Howard ordinal, $\tc(\epsn\cdot3)=((\epsn,2))$ (Subcase 2.1.1.2),
   and $\tc(\varphi_{2,0}\cdot(\om+2))=((\varphi_{2,0},\om+1))$ (Subcase 2.1.1.3). 
   Examples for Subcase 2.1.2 were already given in the paragraph preceding 
   Subcase 1.2.1, and a simple example for the remaining Subcase 2.2 is $\tc(\epsn\cdot\om)=((\epsn,\om))$.
   Note that the characterizing pattern for $\BH$, $\BH\ktwo\eps_{\BH+1}\ktwo\eps_{\BH+1}\cdot(\om+1)$, is the least example for a
   $\ktwo$-chain of three ordinals.
}
\\[2mm]
{\bf Subcase 2.1:} $\alcp{n,m_n}$ is a successor ordinal, say $\alcp{n,m_n}=\xi+1$.
Let $\tau:=\taucp{n,m_n-1}$ and $\alpr:=\ordvalue(\alvec[\xi])$. We consider cases for $\chit(\xi)$:
\\[2mm]
{\bf Subcase 2.1.1:} $\chit(\xi)=0$. In order to verify part a) we have to show that $\predec_1(\al)=\alpr$. 
By monotonicity and continuity we have \[\al=\sup\set{\ordvalue(\alvec[\xi]^\frown(\rhoargs{\tau}{\xi}+\eta))}{\eta\in(0,\tau)},\]
which by the i.h.\ is a proper supremum over ordinals the greatest $\lo$-predecessor of which is $\alpr$.
\\[2mm]
We now proceed to prove part b) and consider cases regarding $\xi$.\\[2mm]
{\bf 2.1.1.1:} $\xi=0$. By part a) $\alpr=\ordcp{n,m_n-1}(\alvec)$ is the greatest $\lo$-predecessor of $\al$.\\[2mm]
{\bf 2.1.1.1.1:} $m_n=2$. By the i.h., $\alpr$ is either $\leo$-minimal or has a greatest $\lo$-predecessor, whence it does not have any
$<_2$-successor, and thus $\alpr\not\ktwo\al$ as claimed. 
Clearly, any $\ktwo$-predecessor of $\al$ then must be a $\ktwo$-predecessor of $\alpr$ as well. 
If $\predec_2(\alpr)>0$ then using the i.h.\ $\al$
is seen to be a proper supremum of $\ktwo$-successors of $\predec_2(\alpr)$ like $\alpr$ itself, hence $\predec_2(\al)=\predec_2(\alpr)$, 
as claimed.\\[2mm]
{\bf 2.1.1.1.2:} $m_n>2$. Then $\alpr\ktwo\al$, as according to the i.h.\ $\al$ then is the supremum of $\ktwo$-successors of $\alpr$, 
hence $\predec_2(\al)=\alpr$, as claimed.\\[2mm]
{\bf 2.1.1.2:} $\xi$ is a successor ordinal. Then by the i.h.\ $\alpr$ has a greatest $\lo$-predecessor, so $\alpr$ does not have any 
$<_2$-successor, in particular $\alpr\not\ktwo\al$. 
In the special case $m_n=2\andsp\taunstar=1$ the $\letwo$-minimality follows then
from the $\letwo$-minimality of $\alpr$, while in the remaining cases $\al$ is easily seen to be 
the supremum of ordinals with the same greatest $\ktwo$-predecessor as claimed for $\alpr$.\\[2mm]
{\bf 2.1.1.3:} $\xi\in\Lim$.
As $\alpr$ is its greatest $\lo$-predecessor, $\al$ is $\alpr$-$\letwo$-minimal, and showing that $\alpr\not\ktwo\al$ will imply the claim 
as above. Arguing toward contradiction, let us assume that $\alpr\ktwo\al$.
Let $X$ and $Z\subseteq(\alpr,\alpr+\ka^{\varsivec}_{\rhoargs{\tau}{\xi}+1})$ be sets 
according to Claim \ref{relleominclaim}, for which there does not exist any cover $X\cup\tildez$ such that $X<\tildez$ and
$X\cup\tildez\subseteq\alpr+\ka^{\varsivec}_{\rhoargs{\tau}{\xi}}+\dpf_{\varsivec}(\rhoargs{\tau}{\xi})$.
We set \[\Xpr:=X\setminus\singleton{\alpr}\quad\mbox{ and }\quad\Zpr:=\singleton{\alpr}\cup Z.\]
By part 2 of Lemma \ref{letwoupwlem} we obtain cofinally many copies $\tildezpr$ below $\alpr$ such that $\Xpr<\tildezpr$ and 
$\Xpr\cup\tildezpr\cong\Xpr\cup\Zpr$ with the property that
$\ordvalue(\alvec[1])\le\altipr:=\min\tildezpr\lo\alpr$. Let $\nu\in(0,\xi)$ be such that $\ordvalue(\alvec[\nu])\le\altipr<\ordvalue(\alvec[\nu+1])$. 
Choosing $\tildezpr$ accordingly we may assume that $\Xpr<\ordvalue(\alvec[\nu])$ and $\logend(\nu)<\logend(\xi)$, 
hence $\rhoargs{\tau}{\nu}\le\rhoargs{\tau}{\xi}$. Notice that if $\ordvalue(\alvec[\nu])<\altipr$
the i.h.\ yields $\chit(\nu)=1$ and $\altipr\le\predec_1(\ordvalue(\alvec[\nu+1]))=\me(\ordvalue(\alvec[\nu]))$, 
whence $\ordvalue(\alvec[\nu])\ktwo\altipr$.
We may therefore assume that $\altipr=\ordvalue(\alvec[\nu])$ since changing $\altipr$ to $\ordvalue(\alvec[\nu])$ would still result in a cover of 
$\Xpr\cup\Zpr$. Because $\ordvalue(\alvec[\nu+1])\lo\alpr$ by the i.h., we may further assume that $\tildezpr\subseteq\ordvalue(\alvec[\nu+1])$. 
Noticing that in the case $\rhoargs{\tau}{\nu}=\rhoargs{\tau}{\xi}$ we must have $\chit(\nu)=1$ and by the i.h.\  
$\ordvalue(\alvec[\nu])+\ka^{\varsivec}_{\rhoargs{\tau}{\nu}}\lo\alpr$,
we finally may assume that $\Xpr\cup\tildezpr\subseteq\ordvalue(\alvec[\nu]^\frown(\ze))$ for some $\ze<\rhoargs{\tau}{\xi}$ with $\min\tildezpr=\ordvalue(\alvec[\nu])$ so that 
 $\Xpr\cup\tildezpr$ is a cover of $\Xpr\cup\Zpr$. Since by i.h.
 \[\ordvalue(\alvec[\nu]^\frown(\ze))\quad\cong\quad\ordvalue(\alvec[\nu])\cup[\alpr,\ordvalue(\alvec[\xi]^\frown(\ze))),\]
setting \[\tildez:=(\alpr+(-\ordvalue(\alvec[\nu])+\tildezpr))-\singleton{\alpr}\] results in a cover $X\cup\tildez$ of $X\cup Z$ with $X<\tildez$ and 
$X\cup\tildez\subseteq\alpr+\ka^{\varsivec}_{\rhoargs{\tau}{\xi}}+\dpf_{\varsivec}(\rhoargs{\tau}{\xi})$. Contradiction.
\\[2mm]
{\bf Subcase 2.1.2:} $\chit(\xi)=1$. 
Part a) claims that, setting $\devec:=\me(\alvec[\xi])$,  $\predec_1(\al)=\ordvalue(\devec)=:\de$. 
It is easy to see that the extending index of $\ec(\devec)$ is of a form $\tau\cdot(\eta+1)$ for some $\eta$, shown explicitly in
part 7(c) of Lemma \ref{tcevalestimlem}. 
Notice that $\cml(\devec)=(n,m_n-1)$. By monotonicity and continuity we then have
\[\al=\sup\set{\ordvalue(\devec^\frown(\tau\cdot\eta+\ze))}{\ze\in(0,\tau)},\]
which by the i.h.\ is a proper supremum over ordinals the greatest $\lo$-predecessor of which is $\de$.
\\[2mm]
As to part b) we first show that $\al$ is $\alpr$-$\letwo$-minimal, arguing similarly as in the proof of (relativized) $\letwo$-minimality in Subcase 1.2., but providing the argument explicitly again for the reader's convenience.
We will then prove $\alpr\not\ktwo\al$ which as above implies the claim.

Let $\devec=(\devec_1,\ldots,\devec_r)$ where $\devec_i=(\decp{i,1},\ldots,\decp{i,k_i})$ for $1\le i\le r$ with associated chain $\sivec$. 
Then $r\ge n$, $k_n\ge m_n$, $\decp{n,m_n}=\xi$, and $\decp{i,j}=\alcp{i,j}$ for all $(i,j)\in\dom(\devec)$ such that $(i,j)\klex(n,m_n)$. 
Recall that we have $\predec_1(\al)=\ordvalue(\devec)=\de$, i.e.\ $\de$ is the greatest $\lo$-predecessor of $\al$, according to part a). 
For any $\ga\ktwo\al$ we therefore must have $\ga\le\de$.
According to the i.h.\ and with the aid of Lemma \ref{cmlmaxextcor} the maximal $\ktwo$-chain from $\alpr$ to $\de$ consists of ordinals
the tracking chains of which are initial chains of $\devec$ that extend $\tc(\alpr)=\alvec[\xi]$, in particular $(n,m_n)\klex(r,k_r)$ and
$\alpr\ktwo\de$, and as verified by part 7(c) of Lemma \ref{tcevalestimlem} we have, setting $\varsivecde:=\ers_{r,k_r}(\devec)$,
\[\al=\de+\ka^\varsivecde_{\tau\cdot(\eta+1)}=\be+\tauti,\mbox{ where }
  \be:=\de+\ka^\varsivecde_{\tau\cdot\eta}+\dpf_\varsivecde(\tau\cdot\eta)\mbox{ and }\tauti=\ka^\varsivecde_\tau.\] 

Now, arguing toward contradiction, let us assume that there exists a greatest $\ktwo$-predecessor $\ga$ of $\al$ such that 
$\alpr\ktwo\ga\letwo\de$, so that by the i.h.\ $\gavec:=\tc(\ga)$ is an initial chain of $\devec$, say $\ga=\ordcp{i,j+1}(\devec)$ 
for some $(i,j+1)\in\dom(\devec)$ with $j>0$ and $(n,m_n)\klex(i,j+1)$.
We then have $\sicp{i,j}>\tau$ by Lemma \ref{cmlmaxextcor} and set $\theta:=\ordcp{i,j}(\devec)$. 
According to Claim \ref{relleominclaim} of the i.h.\ for $\ov({\devec_{\restriction_{i,j}}}^\frown(\tau+1))=\theta+\tauti+1$ 
there exist finite sets $X\subseteq\theta+1$ and $Z\subseteq(\theta,\theta+\tauti+1)$
such that there does not exist any cover $X\cup\tildez$ of $X\cup Z$ with $X<\tildez$ and $X\cup\tildez\subseteq\theta+\tauti$.
By the i.h.\ we know that for every $\nu\in(0,\decp{i,j+1})$, setting $\de_\nu:=\ov(\devecrestrarg{(i,j+1)}[\nu])$ we have
\[\theta+1+\siticp{i,j}\quad\cong\quad\theta+1\cup(\de_\nu,\de_\nu+\siticp{i,j})\]
Since $\decp{i,j+1}=\mu_\sicp{i,j}\in\Hz$ and $\tauti<\siticp{i,j}$ (due to the monotonicity of an appropriately relativized $\ka$-function,
say, $\ka^\varsivecde$), we directly see that below $\ga$ there are cofinally many copies 
$\tildez_\ga$ of $Z$ such that $X\cup Z\cong X\cup\tildez_\ga$. 
By part 1 of Lemma \ref{letwoupwlem} and our assumption $\ga\ktwo\al$ we now obtain copies $\tildez_\al$ of $Z$ cofinally below $\al$ 
(and hence above $\be$) such that $X\cup Z\cong X\cup\tildez_\al$.
The i.h.\ reassures us of the isomorphism
\[\theta+1+\tauti\quad\cong\quad\theta+1\cup(\be,\al),\]
noting that (by the i.h.) the ordinals of the interval $(\be,\al)$ cannot have any $\ktwo$-predecessors in $(\theta,\be]$ 
and that the tracking chains of the ordinals
in $(\theta,\theta+\tauti)\cup(\be,\al)$ have the proper initial chain $\devecrestrarg{(i,j)}$. This provides us, however, with a copy 
$\tildez\subseteq(\theta,\theta+\tauti)$ of $Z$ such that $X\cup Z\cong X\cup\tildez$, contradicting our choice of $X$ and $Z$, whence $\ga\ktwo\al$ is impossible.
Therefore $\al$ is $\alpr$-$\letwo$-minimal.
\\[2mm]
We now show that $\alpr\not\ktwo\al$. In order to reach a contradiction let us assume to the contrary that $\alpr\ktwo\al$. Under this assumption we can prove the following variant of Claim \ref{relleominclaim}:
\begin{claim}\label{relletwominclaim} 
Assuming $\alpr\ktwo\al$, there exist finite sets $X$ and $Z\subseteq(\alpr,\al]$, where 
$X$ consists of $\alpr$ and all existing greatest $<_2$-predecessors $\ga$ of elements of $Z$ that satisfy $\ga\le\alpr$, such that 
there does not exist any cover $X\cup\tildez$ of $X\cup Z$ with $X<\tildez$ and $X\cup\tildez\subseteq\al$.
\end{claim}
{\bf Proof.} The proof of the above claim both builds upon Claim \ref{relleominclaim} and is similar to its proof, 
but for the reader's convenience we give it in detail, with an emphasis on a situation that is not particularly difficult but did not occur in the proof of 
Claim \ref{relleominclaim}, cf.\ Subcase \Romannumeral{1}.2.

We are going to show that for every index pair $(i,j)\in\dom(\devec)$ such that $(n,m_n)\kglex(i,j)\kglex(r,k_r)$, setting
$\etaij:=\ordcp{i,j}(\devec)$, there exists a finite set $\Zij\subseteq(\etaij,\al]$ such that 
for $\Xij$ consisting of $\etaij$ and all existing greatest $<_2$-predecessors below $\etaij$ of elements of $\Zij$ there does not exist 
any cover $\Xij\cup\tildezij$ of $\Xij\cup\Zij$ with $\Xij<\tildezij$ and $\Xij\cup\tildezij\subseteq\al$. 
We proceed by induction on the finite number of 1-step extensions from $\devecrestrarg{(i,j)}$ to $\devec$: 
The initial step is $(i,j)=(r,k_r)$, hence $\etaij=\de$. Recalling that $\al=\de+\ka^\varsivecde_{\tau\cdot(\eta+1)}$, we can apply Claim \ref{relleominclaim} of the i.h.\ to 
$\de+\ka^\varsivecde_{\tau\cdot\eta+1}$ to obtain sets $\Xpr$ and $\Zpr\subseteq(\de,\de+\ka^\varsivecde_{\tau\cdot\eta+1})$ such that there 
does not exist any cover $\Xpr\cup\tildezpr$ of 
$\Xpr\cup\Zpr$ with $\Xpr<\tildezpr$ and $\Xpr\cup\tildezpr\subseteq\de+\ka^\varsivecde_{\tau\cdot\eta}+\dpf_\varsivecde(\tau\cdot\eta)=:\depr$.
Defining
\[\Zij:=\Zpr\cup\singleton{\al}\]
and noticing that by our assumption we have $\alpr\ktwo\al$ and that by the i.h.\ there do not exist any $\ktwo$-successors of $\alpr$ 
in the interval $(\depr,\al)$, it is easy to check that $\Zij$ has the required property.
Let us now assume that $(i,j)\klex(r,k_r)$. We set $(u,v):=(i,j)^+$, $\varsivecpr:=\ers_{i,j}(\devec)$, and consider two cases.\\[2mm]
{\bf Case \Romannumeral{1}:} $(u,v)=(i+1,1)$. By Lemma \ref{cmlmaxextcor} we have $\sicp{u,v}>\tau\in\Ez$ and hence $\si^\star_u\ge\tau$. 
Notice that the case $\sicp{u,v}=\si^\star_u$ cannot occur since then $\ec(\devecrestrarg{(u,v)})$ would not exist.
We discuss the remaining possibilities for $\decp{u,v}$:\\[2mm]
{\bf Subcase \Romannumeral{1}.1:} $\decp{u,v}\in\Ez^{>\si^\star_u}$. 
We then argue as in the corresponding Case \Romannumeral{1} in the proof of Claim \ref{relleominclaim}. We therefore define
\[\Zij:=\Zuv.\]
That this choice is adequate is shown as in the proof of Claim \ref{relleominclaim}.\\[2mm]  
{\bf Subcase \Romannumeral{1}.2:} Otherwise. 
In case of $\decp{u,v}>\sicp{u,v}$ let $\ze$ be such that $\decp{u,v}=_\NF\ze+\sicp{u,v}$, otherwise set $\ze:=0$. 
If $\ze>0$ let $X_\ze$ and $Z_\ze\subseteq(\etaij,\etaij+\ka^{\varsivecpr}_{\ze+1})$ be the sets according to 
Claim \ref{relleominclaim} of the i.h.\ so that there does not exist any cover $X_\ze\cup\tildez_\ze$ of $X_\ze\cup Z_\ze$ with $X_\ze<\tildez_\ze$ and 
$X_\ze\cup\tildez_\ze\subseteq\etaij+\ka^{\varsivecpr}_\ze+\dpf_{\varsivecpr}(\ze)$, otherwise set $X_\ze:=\emptyset=:Z_\ze$. We now define
\[\Xij:=\{\etaij\}\cup X_\ze\cup(\Xuv\setminus\{\etauv\})\quad\mbox{ and }\quad\Zij:=Z_\ze\cup\singleton{\etauv}\cup\Zuv.\]
In order to show that this choice of $\Zij$ satisfies the claim let us assume to the contrary the existence of a set $\tildezij$ such
that $\Xij\cup\tildezij$ is a cover of $\Xij\cup\Zij$ with $\Xij<\tildezij$ and $\Xij\cup\tildezij\subseteq\al$.
Let $\Zpr:=\singleton{\etauv}\cup\Zuv$ and $\tildezpr$ be the subset of $\tildezij$ corresponding to $\Zpr$.
Due to the property of $Z_\ze$ in the case $\ze>0$ we have
\[\tildezpr\subseteq[\etaij+\ka^{\varsivecpr}_{\ze+1},\al),\]
and since $\etauv\lo\al$ there are cofinally many copies 
\[\tildezpr\subseteq[\etaij+\ka^{\varsivecpr}_{\ze+1},\etauv)\]
below $\etauv$ keeping the same $\ktwo$-predecessors. 
The ordinal $\mu:=\min\tildezpr$ corresponds to $\etauv$ in $\Zij$, and since $\mu\leo\tildezpr$ we see that there exists $\nu\in(\ze,\decp{u,v})$ 
such that setting $\eta_\nu:=\etaij+\ka^{\varsivecpr}_{\nu}$ and $\eta_{\nu+1}:=\etaij+\ka^{\varsivecpr}_{\nu+1}$ we have
\[\tildezpr\setminus\singleton{\mu}\subseteq(\eta_\nu,\eta_{\nu+1}),\]
which again we may assume to satisfy $\nu\ge\si^\star_u$ and
$\logend((1/\si^\star_u)\cdot\nu)<\log((1/\si^\star_u)\cdot\sicp{u,v})$. By the i.h.\ we have 
\[\eta_{\nu+1}\quad\cong\quad\eta_\nu\cup[\etauv,\etauv+(-\eta_\nu+\eta_{\nu+1}))\]
since $\etauv$ and $\eta_\nu$ have the same $\ktwo$-predecessors. Exploiting this isomorphism and noticing that 
$\Xuv\setminus\singleton{\etauv}\subseteq\Xij$ we obtain a copy 
$\tildezuv$ of $\tildezpr\setminus\singleton{\mu}$ such that
$\Xuv\cup\tildezuv$ is a cover of $\Xuv\cup\Zuv$ with $\Xuv<\tildezuv$ and $\Xuv\cup\tildezuv\subseteq\al$. Contradiction.\\[2mm]
{\bf Case \Romannumeral{2}:} $(u,v)=(i,j+1)$. Setting $\si:=\sicp{i,j}$ we then have $\decp{i,j+1}=\mu_\si$ 
and proceed as in the corresponding case in the proof of Claim \ref{relleominclaim}. 
Applying Claim \ref{relleominclaim} of the i.h.\ to 
$\ordvalue(\devecrestrarg{(i,j)}^\frown(\sibar+1))$ yields a set 
$Z_\sibar\subseteq(\etaij,\etaij+\ka^\varsivecpr_{\sibar+1})$ such that there does not exist a cover $\Xij\cup\tildez_\sibar$ of $\Xij\cup Z_\sibar$
with $\Xij<\tildez_\sibar$ and $\Xij\cup\tildez_\sibar\subseteq\etaij+\ka^\varsivecpr_\sibar+\dpf_\varsivecpr(\sibar)$.
We now define
\[\Zij:=\singleton{\etauv}\cup(\etauv +(-\etaij +Z_\sibar))\cup\singleton{\ordvalue(\devecrestrarg{(u,v)}^\frown(\si))}\cup\Zuv.\]
In order to show that $\Zij$ has the desired property we assume that there were a cover $\Xij\cup\tildezij$ of $\Xij\cup\Zij$ with 
$\Xij<\tildezij$ and 
$\Xij\cup\tildezij\subseteq\al$ and then argue as in the corresponding Case \Romannumeral{2} in the proof of 
Claim \ref{relleominclaim} in order to drive the assumption into a contradiction.\\[2mm] 
The final instance $(i,j)=(n,m_n)$ establishes Claim \ref{relletwominclaim}.\qed

We can now derive a contradiction similarly as in the previous subcase. Let $X,Z$ be as in the above claim. Without loss of generality me may assume that 
$\predec_1(\al)=\de\in Z$.
We set \[\Xpr:=X\setminus\singleton{\alpr}\quad\mbox{ and }\quad\Zpr:=\singleton{\alpr}\cup Z\setminus\singleton{\al}.\]
By part 2 of Lemma \ref{letwoupwlem} we obtain cofinally many copies $\tildezpr$ below $\alpr$ such that $\Xpr<\tildezpr$ and 
$\Xpr\cup\tildezpr\cong\Xpr\cup\Zpr$ with the property that
all $\leo$-connections to $\al$ are maintained. Let \[\altipr:=\min\tildezpr\] and notice that $\altipr\lo\al$ and $\altipr\letwo\gati$ 
for all $\gati\in\tildezprnod$, where $\tildezprnod$ is defined as the subset of $\tildezpr$ that consists of the copies of all 
$\ga\in\Zpr$ such that $\ga\lo\al$.
Let $\nu\in(0,\xi)$ be such that $\ordvalue(\alvec[\nu])\le\altipr<\ordvalue(\alvec[\nu+1])$. 
Choosing $\tildezpr$ accordingly we may
assume that both $\Xpr<\ordvalue(\alvec[\nu])$ and $\logend(\nu)<\logend(\xi)$. Notice that  
using the i.h.\ $\chit(\nu)=1$, hence $\rhoargs{\tau}{\nu}=\tau\cdot\logend(\nu)<\tau\cdot\logend(\xi)=\rhoargs{\tau}{\xi}$
(as $\nu=\nupr\cdot\om$ for some $\nupr$ would imply $\chit(\nu)=0$, simimarly for $\xi$), 
and $\altipr\leo\predec_1(\ordvalue(\alvec[\nu+1]))=\me(\ordvalue(\alvec[\nu]))$, whence 
$\ordvalue(\alvec[\nu])\letwo\altipr$.
We may therefore assume that $\altipr=\ordvalue(\alvec[\nu])$ since a replacement would still result in a cover of $\Xpr\cup\Zpr$. Because 
$\ordvalue(\alvec[\nu+1])\lo\alpr$ by the i.h., we may further assume that $\tildezpr\subseteq\ordvalue(\alvec[\nu+1])$ as elements of
$\tildezprnod$ are not affected. 
Noticing that since $\chit(\nu)=\chit(\xi)=1$ we have $\nu\cdot\om<\xi$, and setting 
\[\alstar:=\ordvalue(\alvec[\nu\cdot\om])+\ka^{\varsivec}_{\rhoargs{\tau}{\nu\cdot\om}}+\dpf_{\varsivec}(\rhoargs{\tau}{\nu\cdot\om})\]
we can use the isomorphism 
\[\ordvalue(\alvec[\nu+1])\quad\cong\quad\ordvalue(\alvec[\nu])\cup[\ordvalue(\alvec[\nu\cdot\om]),\alstar),\]
which is established by the i.h., in order to shift $\tildezpr$ by the translation 
\[\tildezstar:=\ordvalue(\alvec[\nu\cdot\om])+(-\ordvalue(\alvec[\nu])+\tildezpr).\] 
This results in the cover $\Xpr\cup\tildezstar$ of $\Xpr\cup\Zpr$. By the i.h.\ we know that 
\[\ordvalue(\alvec[\nu\cdot\om])\ktwo\alstar=\ordvalue(\alvec[\nu\cdot\om])+(-\ordvalue(\alvec[\nu])+\ordvalue(\alvec[\nu+1]))\]
and that for all $\gati\in\tildezprnod$ the corresponding element in $\tildezstar$ satisfies
\[\ordvalue(\alvec[\nu\cdot\om])+(-\ordvalue(\alvec[\nu])+\gati)\lo\alstar.\]
Since $\rhoargs{\tau}{\nu\cdot\om}=\tau\cdot\logend(\nu)<\rhoargs{\tau}{\xi}$, setting 
$\alti:= \alpr+\ka^{\varsivec}_{\rhoargs{\tau}{\nu\cdot\om}}+\dpf_{\varsivec}(\rhoargs{\tau}{\nu\cdot\om})$,  we may finally exploit the isomorphism
\[\alstar+1\quad\cong\quad\ordvalue(\alvec[\nu\cdot\om])\cup[\alpr,\alti]\]
so that setting
\[\tildez:=(\alpr+(-\ordvalue(\alvec[\nu\cdot\om])+(\tildezstar\cup\singleton{\alstar})))\setminus\singleton{\alpr}\]
we obtain the cover $X\cup\tildez$ of $X\cup Z$ which satisfies $X<\tildez$ and $X\cup\tildez\subseteq\al$. Contradiction.
\\[2mm]
{\bf Subcase 2.2:} $\alcp{n,m_n}\in\Lim$.
\\[2mm]
Part a) follows from the i.h.\ by monotonicity and continuity, according to which
\[\al=\sup\set{\ordvalue(\alvec[\xi])}{\xi\in(0,\alcp{n,m_n})}.\]
In order to see part b) we simply observe that according to part a) and the i.h.\ $(\ordvalue(\alvec[\xi]))_{\xi<\alcp{n,m_n}}$ 
is a $\lo$-chain of ordinals either $\letwo$-minimal as claimed for $\al$ or with the same greatest $\ktwo$-predecessor as claimed for $\al$.
\qed

The following Corollary \ref{maincor} applies Theorem \ref{maintheo} in order to completely describe the structure $\Rtwo$ in terms of
$\le_i$-successorship, $i=1,2$. As a preparation we need the following definition of \emph{greatest branching point} of a tracking chain,
denoted as $\gbo(\alvec)$, which is crucial in calculating $\lh(\al)$, i.e\ the maximum $\be$ such that $\al\leo\be$ if such ordinal exists.
Recall that we write $\lh_2(\al)$ for the maximum $\be$ such that $\al\letwo\be$ if such ordinal exists, and $\succs_i(\al)$ for the class
$\{\be\mid\al\le_i\be\}$, $i=1,2$. Recall also Definition \ref{charseqdefi} of \emph{reference sequence}, $\rsij(\alvec)$, and 
\emph{evaluation reference sequence}, $\ersij(\alvec)$.

\begin{defi}[7.12 of \cite{CWc}]\label{gbodefi}
Let $\alvec\in\TC$ where $\alivec=(\alcp{i,1},\ldots,\alcp{i,m_i})$ for $1\le i\le n$ and set \index{$\alvecst$}
\[\alvecst:=\left\{\begin{array}{l@{\quad}l}
\alvec & \mbox{if }m_n=1\\
\alvec[\mu_\taucp{n,m_n-1}] &\mbox{otherwise.}
\end{array}\right.\]
We define the (index pair of the) {\bf greatest branching point of}
\index{index pair!greatest branch-off index pair} $\alvec$, 
$\gbo(\alvec)$\index{$\gbo$}, by
\[\gbo(\alvec):=\left\{\begin{array}{l@{\quad}l}
\gbo(\alvecrestrarg{(i,j+1)}) & \mbox{if $(i,j):=\cml(\alvecst)$ exists}\\[2mm]
(n,m_n) & \mbox{otherwise.}
\end{array}\right.\]
\end{defi}

\begin{cor}[cf.\ 7.13 of \cite{CWc}]\label{maincor}
Let $\al\in\On$ with $\tc(\al)=\alvec$ where $\alivec=(\alcp{i,1},\ldots,\alcp{i,m_i})$ for $1\le i\le n$,
with associated chain $\tauvec$ and segmentation parameters $(\la,t):=\upsseg(\alvec)$ and $p,s_l,(\la_l,t_l)$ for $l=1,\ldots,p$ 
as in Definition \ref{segmentationdefi}.\\[2mm]
{\bf Case 1:} $\al\in\Image(\ups)$.
\begin{enumerate}
\item[{\bf 1.1:}] $\al=0$. Then $\lh_2(\al)=\lh(\al)=\al$.
\item[{\bf 1.2:}] $\al=\ups_{\la+1}$. Then $\lh_2(\al)=\al$ and $\lh(\al)=\infty$.
\item[{\bf 1.3:}] $\al=\ups_\la>0$. Then $\lh_2(\al)=\lh(\al)=\infty$ and \[\succs_2(\al)=\{\ups_\la\cdot(1+\xi)\mid\xi\in\On\}.\]
\item[{\bf 1.4:}] $\al=\ups_{\la+t^\prime+1}$, $t^\prime\in(0,\om)$. Then $\lh_2(\al)=\lh(\al)=\infty$ and 
            \[\succs_2(\al)=\{\al+\ups_{\la+t^\prime}\cdot\xi\mid\xi\in\On\}.\]
\end{enumerate}
{\bf Case 2:} $\al\not\in\Image(\ups)$. 
\begin{enumerate}
\item[a)] We first consider $\le_2$-successors of $\al$.\\[2mm] 
{\bf Subcase 2.1:} $m_n=1$. Then 
\[\succs_2(\al)=\singleton{\al}\quad\mbox{ and }\quad\lh_2(\al)=\al.\] 
{\bf Subcase 2.2:} $m_n>1$. Set $(i_0,j_0):=(n,m_n-1)$, $\tau:=\taucp{i_0,j_0}$, $\tauti:=\tauticp{i_0,j_0}$, 
$\varsivec:=\rs_{i_0,j_0}(\alvec)$, $\varrho:=\rhoargs{\tau}{\taucp{n,m_n}}$, and let
$\nu,\xi$ be such that 
\[\ka^{\varsivec}_\varrho+\dpf_{\varsivec}(\varrho)=\tauti\cdot\nu+\xi\mbox{ and }\xi<\tauti.\] 
Writing $\maxeta:=\nu\minusp\chit(\taucp{n,m_n})$ we then have
\[\succs_2(\al)=\set{\al+\tauti\cdot\eta}{\eta\le\maxeta}\quad\mbox{ and }\quad\lh_2(\al)=\al+\tauti\cdot\maxeta.\] 
\item[b)] Writing $(n_0,m_0):=\gbo(\alvec)$, $m:=m_0\minusp 2+1$, and setting
as in Definition \ref{trchevaldefi} $\varsivec:=\ers_{n_0,m}(\alvec)$
we have
\[\lh(\al)=\left\{\begin{array}{ll}
\ordcp{n_0,m}(\alvec)+\dpf_\varsivec(\taucp{n_0,m}) & \mbox{ if }\ordcp{n_0,m}(\alvec)\not\in\Image(\ups)\\[2mm]
\infty & \mbox{ otherwise.}
\end{array}\right.\]
\end{enumerate}
\end{cor}
{\bf Proof.} In Case 1 claims regarding $\lh$ follow directly from Theorem \ref{maintheo}. The claim regarding $\lh_2$ in Subcases 1.2 and 2.1 
follows from Lemma \ref{ktwoinflochainlem}, since according to Theorem \ref{maintheo} cofinal $\lo$-chains do not exist. 
We now consider the situations in Subcases 1.3, 1.4, and 2.2. 
If $\be$ is a $\ktwo$-successor of $\al$, according to the theorem either $\alvec\subseteq\bevec:=\tc(\be)$, where $m_n>1$, 
or $\al=\ups_{\la+t}$, where $\la\in\Limnod$, $t\in\N\setminus\{1\}$, and $\la+t>0$.
In the situation of Subcase 1.3 let $\tau:=\ups_\la=:\tauti$, in Subcase 1.4 let $\tau:=\ups_{\la+\tpr}=:\tauti$, and in Subcase 2.2 let 
$\tau$ and $\tauti$ be as defined there.
Note that $\le_i$-successorhip is closed under limits, $i=1,2$, so it is sufficient to consider 
successor-$\ktwo$-successors $\be$ of $\al$, which are immediate $\ktwo$-successors, as by the theorem non-immediate $\ktwo$-successors 
cannot be successor-$\ktwo$-successors.
Suppose therefore that $\bevec=(\bevec_1,\ldots,\bevec_r)$ where $\bevec_i=(\becp{i,1},\ldots,\becp{i,k_i})$ for $i=1,\ldots,r$ is a 
successor-$\ktwo$-successor of $\al$, whence by the theorem $\predec_2(\be)=\al$, $k_r=1$, and $\tau=\si^\star_r=\becp{r,1}>1$, 
where $\sivec$ is the chain associated with $\bevec$. Clearly, the converse of this latter implication holds as well.
Therefore, if $\be$ is a successor-$\ktwo$-successor of $\al$, then it must be of the form $\al+\tauti\cdot(\eta+1)$ for some $\eta$.
In Subcase 2.2 we see that such $\eta$ must be bounded as claimed, since according to the theorem $\bevec$ is then a proper extension of 
$\alvec$, whence with the aid of Corollary \ref{alpldpcor} $\be\le\al+\ka^{\varsivec}_\varrho+\dpf_{\varsivec}(\varrho)$, with strict 
inequality if $\chit(\taucp{n,m_n})=1$.
Note that any ordinal greater than $\al$ and bounded in this way has a tracking chain that properly extends $\alvec$, cf.\ Corollary 
\ref{alpldpcor}.

Having seen that all $\ktwo$-successors of $\al$ are of the claimed form, we now assume $\be$ to be of the form $\al+\tauti\cdot(\eta+1)$ 
for some $\eta$, bounded as stated in the situation of Subcase 2.2.
Let again $\tc(\be)=:\bevec=(\bevec_1,\ldots,\bevec_r)$ where $\bevec_i=(\becp{i,1},\ldots,\becp{i,k_i})$ for $i=1,\ldots,r$, 
with associated chain $\sivec$.
According to our assumption we have $\sumend(\siticp{r,k_r})=\tauti$. 

We first consider the situation of Subcase 2.2, where $\alvec\subseteq\bevec$ as we have seen above. 
Assuming that $\sicp{r,k_r}=1$, which implies $k_r>1$, we would have $\siticp{r,k_r-1}=\tauti$, hence by part 3 of Lemma \ref{evallem}
$(r,k_r-1)=(n,m_n-1)$ and thus $\alvec=\bevec$, which is not the case.
We therefore have $\sicp{r,k_r}>1$ and hence $\siticp{r,k_r}=\tauti$, which entails \[\trs(\siticp{r,k_r})=\trs(\tauti)\] 
due to Lemma \ref{trsvallem}, and by part 3 of Lemma \ref{evallem} it follows that
neither $k_r>1$ nor $k_r=1\andsp\sicp{r,1}\in\Ez^{>\sirstar}$, since otherwise $(r,k_r)=(n,m_n-1)$, which is not the case. 
By parts 1 and 2 of Lemma \ref{evallem} we have $\siticp{r,1}=\ka^\varsivecpr_\sicp{r,1}$ (where $\varsivecpr:=\rs_{r^\star}(\bevec)$)
and $\trs(\tauti)=\rs_{i_0,j_0}(\alvec)=\varsivec\in\RS$ (where by definition $(i_0,j_0)=(n,m_n-1)$), the membership relation according 
to part 2 of Lemma \ref{rslem}.
Thus $k_r=1$, and since the assumption $\sirstar<\sicp{r,1}$, where as we already know $\sicp{r,1}\not\in\Ez$, would imply 
$\trs(\siticp{r,1})\not\in\RS$ as this tracking sequence could not be a strictly increasing sequence of epsilon numbers, we obtain
$\sicp{r,1}=\sirstar$, so that $\trs(\siticp{r,1})=\varsivecpr$ and hence $\sirstar=\tau$, moreover $r^\star=(i_0,j_0)$ again
by part 3 of Lemma \ref{evallem}, whence by the theorem $\al\ktwo\be$.

Next, we consider the situation of Subcase 1.3, where $\al=\ups_\la>0$ and $\be=\ups_\la\cdot(1+\eta+1)$ for some $\eta$.
Then we have $\sumend(\siticp{r,k_r})=\ups_\la$, which implies $k_r=1$ as $\ups_\la$ can not be (the last summand of) an element of the 
image of any $\nu$-function, 
hence $\siticp{r,1}=\sicp{r,1}=\ups_\la$ as $\ka$-functions map additive principal numbers to addivitive principal numbers, 
and thus $\sirstar=\ups_\la$. The theorem now yields $\predec_2(\be)=\al$.

Finally, in Subcase 1.4, we have $\al=\ups_{\la+\tpr+1}$ and $\tau=\ups_{\la+\tpr}=\tauti$.
Then we have $\sumend(\siticp{r,k_r})=\tau<\al$, which again implies that $k_r=1$ as $\tau$ can only be 
(the last summand of an element) of the image of a $\nu$-function that occurs in tracking chains of ordinals below $\al$, 
namely if it is either equal to $\nu^\gavecpr_{\tau}=\nu^\gavec_0$, where $\gavecpr=(\ups_{\la+1},\ldots,\ups_{\la+\tpr-1})$ and 
$\gavec=\gavecpr^\frown\tau$, or (the last summand of) an element of the image of $\nu^\gavec$, cf.\ Corollaries \ref{kdpnuestimcor} and 
\ref{nuimagecor}. We again conclude that $\siticp{r,1}=\sicp{r,1}=\tau=\sirstar$, and hence by the theorem $\predec_2(\be)=\al$.

It remains to show part b) of the Corollary. 
We argue as in the corresponding proof of Corollary 7.13 of \cite{CWc}.
Let $\alvecpr:=(\alvecrestrarg{(n_0,m_0)})^\star$ using the ${}^\star$-notation from Definition \ref{gbodefi}, according to 
which the vector $\alvecpr$ does not possess a critical main line index pair. 
We set \[\alplus:=\ordvalue(\me(\alvecpr)).\]
If $\ordcp{n_0,m}(\alvec)\in\Image(\ups)$, we have $\alplus=\ov(\alvecpr)\in\Image(\ups)$, otherwise
inspecting definitions as was done in part 7(b) of Lemma \ref{tcevalestimlem} we obtain
\[\alplus=\ordcp{n_0,m}(\alvec)+\dpf_\varsivec(\taucp{n_0,m}).\]
We first show that 
\begin{equation}\label{leomeclaim}
\al\leo\alplus.
\end{equation} 
In the case $m_0=1$ we have $(n_0,m_0)=(n,m_n)$, and the claim follows directly from Theorem \ref{maintheo}. 
Now assume that $m_0>1$.
By Theorem \ref{maintheo} we have $\al\leo\ordvalue(\alvecst)\leo\ordvalue(\me(\alvecst))$.
If $\cml(\alvecst)$ does not exist, that is $\alvecpr=\alvecst$, we are done with showing (\ref{leomeclaim}). Otherwise
let $\cml(\alvecst)=:(i_1,j_1)$ and let $l_0$ be maximal so that for all $l\in(0,l_0)$ 
$\cml((\alvecrestrarg{(i_l,j_l+1)})^\star)=:(i_{l+1},j_{l+1})$ exists. Clearly, the sequence of index pairs we obtain in this way is $\klex$-decreasing, and by Definition \ref{gbodefi}
$(i_{l_0},j_{l_0}+1)=(n_0,m_0)$. Theorem \ref{maintheo} yields the chain of inequations
\[\al\leo\ordvalue(\alvecst)\leo\ordvalue((\alvecrestrarg{(i_1,j_1+1)})^\star)\leo\ldots\leo\ordvalue((\alvecrestrarg{(i_{l_0},j_{l_0}+1)})^\star)=\ordvalue(\alvecpr)\leo\alplus.\]
In the case $\ordcp{n_0,m}(\alvec)\in\Image(\ups)$ the claim follows directly from the theorem, so let us assume otherwise.
We claim that 
\begin{equation}\label{predoneclaim}
\predec_1(\alplus+1)<\ordvalue(\alvecpr).
\end{equation}
To this end note that $\tc(\alplus+1)$ must be of a form $\alvecrestrarg{i}^\frown(\alcp{i+1,1}+1)$ where $i\le n_0$. 
By Theorem \ref{maintheo} $\alplus+1$ is either $\ups_\la$-$\leo$-minimal or the greatest $\lo$-predecessor is $\ordvalue(\alcp{i-1,m_{i-1}})$.
Hence (\ref{predoneclaim}) follows, which implies that $\al\not\leo\alplus+1$. We thus have 
$\lh(\al)=\ordcp{n_0,m}(\alvec)+\dpf_\varsivec(\taucp{n_0,m})$.
\qed

\section*{Acknowledgements}
I express my gratitude to Professor Ulf Skoglund for encouragement and support of my research,
to Professor Wolfram Pohlers for helpful comments on an earlier draft,
and thank Dr.\ Steven D.\ Aird for editing the manuscript.

\end{document}